\documentclass[a4paper,12pt]{article}
\usepackage{latexsym,amssymb,amsmath,amsthm}

\newtheorem{Theorem}{Theorem}[section]
\newtheorem{Definition}[Theorem]{Definition}

\newtheorem{Proposition}[Theorem]{Proposition}
\newtheorem{Example}[Theorem]{Example}
\newtheorem{Remark}[Theorem]{Remark}

\newtheorem{Corollary}[Theorem]{Corollary}

\RequirePackage{amssymb}

\newcommand{\limpr}{\varprojlim}

\newcommand{\st}{\fontshape{n}\selectfont\rm}

\newcommand{\e}{{\st e}}

\renewcommand{\d}{{\st d}}
 \renewcommand{\c}{{\st c}}
\renewcommand{\r}{{\st r}}

\renewcommand{\t}{{\st t}}
\renewcommand{\l}{{\st l}}

\newcommand{\R}{\mathbb{R}}

\newcommand{\N}{\mathbb{N}}

\newcommand{\B}{\mathbb{B}}

\newcommand{\Z}{\mathbb{Z}}
\newcommand{\C}{\mathbb{C}}

\newcommand{\M}{\mathbb{M}}
\newcommand{\U}{\mathbb{U}}
\newcommand{\W}{\mathbb{W}}

\renewcommand{\L}{\mathbb{L}}

\newcommand{\bs}[1]{\boldsymbol{#1}}

\newcommand{\bsa}{{\bs{a}}}

\newcommand{\bsv}{{\bs{v}}}

\newcommand{\bsp}{{\bs{p}}}
\newcommand{\bsf}{{\bs{f}}}
\newcommand{\bsu}{{\bs{u}}}
\newcommand{\bsr}{{\bs{r}}}

\newcommand{\bsg}{\bs{g}}

\newcommand{\bst}{\bs{t}}

\newcommand{\bsP}{\bs{P}}
\newcommand{\bsI}{\bs{I}}

\newcommand{\bsG}{\bs{G}}
\newcommand{\bsA}{\bs{A}}

\newcommand{\bssK}{\bs{\s{K}}}
\newcommand{\bssS}{\bs{\s{S}}}

\newcommand{\wt}[1]{\widetilde{#1}}

\newcommand{\s}[1]{\mathcal{#1}}
\newcommand{\got}[1]{\mathfrak{#1}}

\newcommand{\stref}[1]{{\st \ref{#1}}}
\newcommand{\stcite}[1]{{\st \cite{#1}}}

\newcommand{\Hsum}[2]{\begin{matrix} #1 \\ \oplus \\ #2 \end{matrix}}

\newcommand{\Hsump}[2]{\begin{pmatrix} #1 \\ \oplus \\ #2 \end{pmatrix}} % umklammerte Hilbertsumme
\def\re{\mathop{\rm Re}\nolimits}
\def\Im{\mathop{\rm Im}\nolimits}

\newcommand{\dbar}{{d\hspace{-0,05cm}\bar{}\hspace{0,05cm}}}
\newcommand{\ol}{\overline}

\newcommand{\loc}{{\st loc}}
\newcommand{\sing}{{\st sing}}

\newcommand{\reg}{{\st reg}}
\newcommand{\comp}{{\st comp}}
\newcommand{\cone}{{\st cone}}

\newcommand{\elle}{{\st ell}}

\newcommand{\const}{{\st const}}
\newcommand{\grad}{{\st grad}}
\newcommand{\cl}{{\st cl}}

\renewcommand{\flat}{{\st flat}}

\DeclareMathOperator{\Diff}{Diff}

\DeclareMathOperator{\Vect}{Vect}
\DeclareMathOperator{\ord}{ord}
\DeclareMathOperator{\ind}{ind}
\DeclareMathOperator{\Int}{int}

\DeclareMathOperator{\op}{op}
\DeclareMathOperator{\Op}{Op}

\DeclareMathOperator{\id}{id}
\DeclareMathOperator{\im}{im}
\DeclareMathOperator{\coker}{coker}
\DeclareMathOperator{\corner}{corner}
\DeclareMathOperator{\codim}{codim}
\DeclareMathOperator{\supp}{supp}
\DeclareMathOperator{\proper}{proper}

\DeclareMathOperator{\symb}{symb}
\DeclareMathOperator{\dist}{dist}

\DeclareMathOperator{\luc}{luc}
\DeclareMathOperator{\ulc}{ulc}

\begin{document}

\title{The Structure of Operators on Manifolds with Polyhedral Singularities}
\author{ B.-W. Schulze}
\date{}
\maketitle

\begin{abstract}
We discuss intuitive ideas and historical
background of concepts  in the analysis
on configurations with singularities, here in connection
with our iterative approach for higher singularities.

\end{abstract}

\tableofcontents

\section*{Introduction}
\addcontentsline{toc}{section}{Introduction}
The analysis on configurations with singularities (e.g., conical ones, edges, corners, etc.)
is a classical area of mathematics, motivated by models of physics and the applied sciences and
also by structures  of geometry and topology.
The development goes back to (at least) the 19th century.

This paper is a survey on some aspects of the recent development and, at the same time, an introduction.
Moreover, it  is aimed at discussing new phenomena around the solvability of partial differential equations near singularities and the interaction of analytic, geometric and topological aspects.
The crucial point will be the concept of ellipticity of operators with respect to their symbolic structure.
After recalling the standard notion on a smooth manifold, based on the homogeneous principal symbol, we discuss the way of how a geometric singularity may contribute to other symbolic levels with associated notions of ellipticity.
Observations of that kind are basic for the construction of pseudo-differential algebras with symbolic hierarchies on manifolds with singularities.
The corresponding theories are voluminous; they may be found in several monographs and articles of the
author, see \cite{Schu2}, \cite{Schu20}, \cite{Schu27}.
The complexity of structures makes it desirable to explain the intuitive ideas in a separate exposition; this is just our motivation here.
Clearly this cannot be exhaustive.
On the one hand there is a vast variety of papers of different orientation, from concrete models of the applied sciences to operator algebra aspects and index theory.
On the other hand there are different believes on priorities and adequate approaches in the singular analysis.
We hope to illustrate the fascination of the structure insight connected with understanding
 and solving problems in this field.
Our conclusion will be that, despite of the  enormous experience through the work of the authors of the `singular community' (and also of the `regular' one), many important problems are unsolved and that the new challenges  open a bright future of the analysis on manifolds with singularities.

What concerns the literature we cannot give a complete review here.
We therefore content ourselves with a  list of references that have from different point of view  connections with 
this exposition.
In particular, we want to mention 
Gelfand \cite{Gelf2},
Agmon \cite{Agmo4},
Agranovich and Vishik \cite{Agra1}, 
Kohn and Nirenberg \cite{Kohn1},
Vishik and Eskin \cite{Vivs3}, 
Atiyah and Singer \cite{Atiy6},
Eskin \cite{Eski2},
Vishik and Grushin \cite{Vivs4},
Sternin \cite{Ster2}, \cite{Ster3},
Kondratyev \cite{Kond1}, 
Plamenevskij \cite{Plam1}, \cite{Plam4}, \cite{Plam2},
Rabinovich \cite{Rabi1},
Bolley and Camus \cite{Boll1},
Gramsch \cite{Gram6}, \cite{Gram3},
Gohberg and Sigal \cite{Gohb3},
Fedosov \cite{Fedo6},
Seeley \cite{Seel2},
Grushin \cite{Grus2},
Boutet de Monvel \cite{Bout1},
Atiyah, Patodi, and Singer \cite{Atiy8},
Maz'ja and Paneah \cite{Mazj3},
Shubin \cite{Shub1},
Parenti \cite{Pare1},
Cordes \cite{Cord1},  
Fichera \cite{Fich1},
Teleman \cite{Tele1},
Cheeger \cite{Chee1},
Melrose \cite{Melr4}, \cite{Melr2},
Melrose and Mendoza \cite{Melr1},
Rempel and Schulze \cite{Remp1}, \cite{Remp2},   
Grubb \cite{Grub1}, \cite{Grub4},
Kondratyev and Oleynik \cite{Kond2},
Grisvard \cite{Gris1},
Maz'ja and Rossmann \cite{Mazj1},
Dauge \cite{Daug1},
Chkadua and Duduchava \cite{Chka1}, \cite{Chka2},
Costabel and Dauge \cite{Cost1},
Shaw \cite{Shaw1},
Mazzeo \cite{Mazz1},
Roe \cite{Roe1},
Rozenblum \cite{Roze2},
Egorov and Schulze \cite{Egor1},
Booss-Bavnbek and Wojciechowski \cite{Boos1}, 
Mantlik \cite{Mant1},
Fedosov and Schulze \cite{Fedo2},
Lesch \cite{Lesc2},
Nistor \cite{Nist2},
Nistor, Weinstein, and Xu \cite{Nist1},
Witt \cite{Witt3},
Fedosov, Schulze, and Tarkhanov \cite{Fedo3}, \cite{Fedo5}, \cite{Fedo1},
Nazaikinskij and Sternin \cite{Naza10}, \cite{Naza11},
Savin and Sternin \cite{Savi2}, \cite{Savi3},
Grieser and Lesch \cite{Gries1},
Krainer \cite{Krai5}, \cite{Krai11},
Seiler \cite{Seil1}, \cite{Seil2},
Gil, Krainer, and Mendoza \cite{Gil6}, \cite{Gil7},
Ammann, Lauter, and Nistor \cite{Amma1},
Tarkhanov \cite{Tark4}.

  \section{Simple questions, unexpected answers}
  \label{s.10.1}

\markboth{Remarks and general comments}{Remarks and general comments}
 \begin{minipage}{\textwidth}
\setlength{\baselineskip}{0cm}
\begin{scriptsize}    
 We start from naive questions such as `what are the
 basic  questions or the right notions' around simple
 objects who everybody knows, e.g., on differential
 operators, their symbols, or the right function spaces. Other 
 questions concern classical objects from complex
 analysis who suddenly become obscure if we ask too
 much $\ldots$
 \end{scriptsize}
\end{minipage}

\subsection{What is ellipticity?}
\label{s.10.1.1}

%10.1.1 {Mi.,21.04.04}
The `standard' ellipticity of a differential operator
\begin{equation}
\label{1011.A2104.eq}
A = \sum_{|\alpha| \leq \mu} a_\alpha(x)D_x^\alpha
\end{equation}
in an open set $\Omega \subseteq \R^n$ with coefficients
$a_\alpha \in C^\infty(\Omega)$ refers to its
homogeneous principal symbol
\begin{equation}
\label{1011.sigmapsi2104.eq}
\sigma_\psi(A)(x, \xi) = \sum_{|\alpha|= \mu}
       a_\alpha(x)\xi^\alpha,
\end{equation}
$(x, \xi) \in \Omega \times (\R^n \setminus \{ 0 \})$.
More generally, if $M$ is a $C^\infty$ manifold, an
operator $A \in \Diff^\mu(M)$ has an invariantly defined
homogeneous principal symbol
\begin{equation}
\label{1011.glob2104.eq}
\sigma_\psi(A) \in C^\infty(T^* M \setminus 0).
\end{equation}
($\Diff^\mu(.)$ denotes the space of all differential operators of order $\mu$ with smooth coefficients on the manifold in parentheses.)
\begin{Definition}
\label{1011.Ellpsi2104.de}
The operator $A$ is called elliptic if $\sigma_\psi(A)
\not= 0$ on $T^*M \setminus 0$.
\end{Definition}

\begin{Remark}
\label{1011.hom2104.re}
Since $\sigma_\psi(A)$ is {\st{(}}positively{\st{)}}
homogeneous of order $\mu$, i.e.,
\begin{equation}
\label{hom2104.eq}
\sigma_\psi(A)(x, \lambda \xi) = \lambda^\mu
    \sigma_\psi(A)(x, \xi)
\end{equation}
for all $(x, \xi) \in T^* M \setminus 0$, $\lambda \in
\R_+$, we may equivalently require $\sigma_\psi(A)
\big|_{S^*M} \not= 0$, where $S^*M$ is the unit cosphere
bundle of $M$ {\st{(}}with respect to some fixed
Riemannian metric on $M${\st{)}}.
\end{Remark}

Clearly we can also talk about the complete symbol            
\begin{equation}
\label{1011.loc2204.eq}
\bs{\sigma}(A)(x, \xi) := a(x, \xi) = \sum_{|\alpha|\leq \mu}
     a_\alpha(x) \xi^\alpha
\end{equation}     
of an operator $A$, first on an open set $\Omega \subseteq \R^n$
and then on a $C^\infty$ manifold $M$. In the latter case by a
complete symbol we understand a system of local complete
symbols \eqref{1011.loc2204.eq} with respect to charts $\chi : U
\to \Omega$ when $U$ runs over an atlas on $M$.

The invariance of symbols refers to transition maps $\kappa
:= \tilde{\chi} \circ \chi^{-1}$ for different charts $\chi : 
U \to \Omega$, $\tilde{\chi} : U \to \widetilde{\Omega}$ which induce
isomorphisms
$\kappa_* : \Diff^\mu(\Omega) \to
      \Diff^\mu(\widetilde{\Omega})$     
(subscript
\textup{`$*$'}
denotes the push forward of an operator under the corresponding
diffeomorphism) and associated symbol push forwards
$a(x, \xi) \to \tilde{a}(\tilde{x}, \tilde{\xi})$
between the local complete symbols $a(x, \xi)$ and
$\tilde{a}(\tilde{x}, \tilde{\xi})$ of
$A$ and $\widetilde{A} = \chi_*A$, respectively. 
As is known we have
$\tilde{a}(\tilde{x}, \tilde{\xi}) \big|_{\tilde{x}=\kappa(x)} \sim
  \sum_{\alpha} \frac{1}{\alpha!} (\partial_\xi^\alpha a)
  (x, {}^\t d \kappa(x)\tilde{\xi})\Phi_\alpha(x, \tilde{\xi})$
for $\partial_\xi^\alpha = \partial_{\xi_1}^{\alpha_1} \cdot \ldots
\cdot \partial_{\xi_n}^{\alpha_n}$, $\Phi_\alpha(x, \tilde{\xi}) :=
D_z^\alpha e^{i \delta(x,z)\tilde{\xi}} \big|_{z=x}$ where
$\delta(x,z) := \kappa(z) - \kappa(x) - d \kappa(x)(z-x)$, with $d
\kappa(x)$  being the Jacobi matrix of $\kappa$, and the function 
$\Phi_\alpha(x, \tilde{\xi})$ is a polynomial in $\tilde{\xi}$ of degree $\leq | \alpha |/2$, $\alpha \in \N^n$, $\Phi_0 = 1$. 
In the asymptotic expression for $\tilde{a}(\tilde{x},
\tilde{\xi})$ we have equality for differential operators
(since the sum is finite) and an asymptotic sum of symbols in the
pseudo-differential case.
(Well known material on spaces $S_{(\cl)}^\mu(\Omega \times
\R^n)$ of pseudo-differential symbols of order $\mu \in \R$
(classical or non-classical) will be given in connection with
Definition \ref{1013.Sy2005.de} below; associated
pseudo-differential operators are discussed in Section 2.2).

In particular, for $\widetilde{A} := \kappa_*A$ it follows that
\[  
\sigma_\psi(\widetilde{A})(\tilde{x}, \tilde{\xi}) =
      \sigma_\psi(A)(x, \xi)\quad \text{for}\quad
      \tilde{x} = \kappa(x), \ \tilde{\xi} = ({}^\t d \kappa(x))^{-1} \xi, 
\]
which shows that $\sigma_\psi(A) \in C^\infty(T^*M \setminus 0)$.        

The ellipticity on the level of complete symbols
\eqref{1011.loc2204.eq} in local coordinates is the
condition that for every chart $\chi : U \to \Omega$ there is a
$p(x, \xi) \in S^{- \mu}(\Omega \times \R^n)$, $n = \dim M$,
such that
$p(x, \xi)a(x, \xi) - 1 \in S^{-1}(\Omega \times \R^n)$.

Recall that principal symbols and complete symbols have natural
properties with respect to various operations, for instance,
\[  \sigma_\psi(I) = 1, \; \sigma_\psi(AB) = 
    \sigma_\psi(A) \sigma_\psi(B)     \]
(with $I$ being the identity operator), and    
$\bs{\sigma}(I) = 1,  \bs{\sigma}(AB) =
    \bs{\sigma}(A) \# \bs{\sigma}(B)$,       
with the Leibniz product $\#$ between the 
local complete symbols $a(x, \xi)$ and $b(x, \xi)$ of 
the operators $A$ and $B$,
respectively,
$(a  \# b)(x, \xi) \sim \sum_{\alpha} \frac{1}{\alpha!} 
\bigl(\partial_\xi^\alpha
    a(x, \xi)\bigr) D_x^\alpha b(x, \xi)$
(the sum on the right hand side is finite in the case of a
differential operator $B$).      

A crucial (and entirely classical) observation is the following
result:

\begin{Theorem}
\label{1011.Ellpsi2204.th}
Let $M$ be a closed compact $C^\infty$ manifold and $A \in
\Diff^\mu(M)$. Then the following properties are equivalent:
\begin{enumerate}
\item The operator $A$ is elliptic with respect to 
      $\sigma_\psi$.
\item $A$ is Fredholm as an operator
      \begin{equation}
      \label{1011.neu1506.eq} 
      A : H^s(M) \to H^{s- \mu}(M)    
      \end{equation}
      for some fixed $s \in \R$.
\end{enumerate}
The property {\em(ii)} entails that {\em\eqref{1011.neu1506.eq}} is a Fredholm operator for every $s \in \R$. 

\end{Theorem}

Parametrices of elliptic differential operators are known to be
(classical) pseudo-differential operators of opposite order. Let
$L_{(\cl)}^\mu(M)$ denote the space of all pseudo-differential
operators on $M$ of order $\mu \in \R$; the manifold $M$ is 
not necessarily
compact (in this notation subscript \textup{`(\cl)'} indicates 
the classical
or the non-classical case). More generally, there are the spaces
$L_{(\cl)}^\mu(M; \R^l)$ of parameter-dependent pseudo-differential
operators on $M$ of order $\mu \in \R$ with the parameter $\lambda
\in \R^l$. In this case the local amplitude functions (in H\"ormander's classes)
contain $(\xi, \lambda) \in \R^n \times \R^l$ as covariables, the operator action (locally based on the Fourier
transform) refers to $(x, \xi)$, and the operators contain 
$\lambda$ as a parameter.

Every $A \in L_{\cl}^\mu(M)$ has a homogeneous principal symbol
$\sigma_\psi(A) \in C^\infty(T^*M \setminus 0)$ and 
a system of local complete symbols $\bs{\sigma}(A)(x,
\xi)$. More generally, for $A \in L_{\cl}^\mu(M; \R^l)$ there is a
corresponding principal symbol
\begin{equation}
\label{1011.1105Nr.eq}
\sigma_\psi(A)(x, \xi, \lambda) \in C^\infty(T^*M \times \R^l
      \setminus 0),
\end{equation}
homogeneous of order $\mu$ in $(\xi, \lambda) \not= 0$, and for $A
\in L^\mu(M; \R^l)$ we have complete parameter-dependent symbols.

The ellipticity of an operator $A(\lambda) \in L_{(\cl)}^\mu(M;
\R^l)$ is defined in an analogous manner as before (for $l > 0$ the
parameter is  treated as a component of the \textup{`covariable'}
$(\xi, \lambda)$).

Let $L^{- \infty}(M)$ denote the space of all operators $C$ on $M$
with kernel $c(x,x') \in C^\infty(M \times M)$, i.e., $C u(x) =
\int_{M} c(x,x') u(x') dx'$ ($dx'$ refers to a Riemannian 
metric on
$M$). Moreover, set $L^{- \infty}(M; \R^l) := {\cal S}(\R^l, L^{-
\infty}(M)$).

\begin{Theorem}
\label{1011.P1205.th}
Let $M$ be a closed compact $C^\infty$ manifold. An elliptic operator
$A \in L_{(\cl)}^\mu(M; \R^l)$, $\mu \in \R$, $l \in \N$, has a
parametrix $P \in L_{(\cl)}^{- \mu}(M; \R^l)$ in the sense
\begin{equation}
\label{1011.NrP1205.eq}
PA = I - C_\l, \quad AP = I - C_\r
\end{equation}
for operators $C_\l, C_\r \in L^{- \infty}(M; \R^l)$.

If $M$ is not compact we have an analogous result; in order to have
well defined compositions in \eqref{1011.NrP1205.eq} we may choose
$P$ in a suitable way, namely, 
{\st{}}\textup{`}properly supported\textup{'}, which is always
possible.
\end{Theorem}

\begin{Remark}
\label{1012.1205re.re}
Theorems {\st{\ref{1011.Ellpsi2204.th}}} and  
{\st{\ref{1011.P1205.th}}} are  true in analogous form for pseudo-differential operators acting between Sobolev spaces of distributional sections of {\st{(}}smooth
complex{\st{)}} vector bundles $E, F$ on $M$,
\begin{equation}
\label{1011.AH1205.eq}
A : H^s(M, E) \to H^{s- \mu}(M, F);
\end{equation}
the principal symbol is then a bundle morphism $\sigma_\psi(A) :
\pi^* E \to \pi^*F$ for $\pi : T^* M \setminus 0 \to M$, and
ellipticity means in this case that $\sigma_\psi(A)$ is an
isomorphism. 
There are then corresponding extensions of Theorems
{\st{\ref{1011.Ellpsi2204.th}}} and {\st{\ref{1011.P1205.th}}}.      
In addition the Fredholm property of \eqref{1011.AH1205.eq} for
a special $s = s_0 \in \R$ entails the Fredholm property for
all $s \in \R$.
\end{Remark}

\begin{Remark}
\label{1011.re1506.re}
We do not repeat all
elements of the classical calculus around pseudo-differential
operators and ellipticity on a smooth manifold. Let us only
mention that the index
$\ind A := \dim \ker A - \dim \coker A$
of the Fredholm operator \eqref{1011.neu1506.eq} is
independent of $s$. In fact, there are finite-dimensional
subspaces
\[  V \subset C^\infty(M,E), \quad W \subset C^\infty (M,F)   \]
such that $V = \ker A, W \cap \im A = \{ 0 \}$ and $W + \im A
= H^s(M,F)$ for all $s \in \R$.
\end{Remark}

It is a general idea to reduce interesting questions on the nature of an operator $A$ (as a map between spaces of distributions on $M$ or on the solvability of the equation $A u = f$) to the level of symbols which are much easier objects than operators. 
This is, of course, a general program, not only for elliptic operators, but also for other types of operators, e.g., parabolic or hyperbolic ones.

The aspect of connecting symbols with operators and vice versa plays
a role in wide areas of  mathematics. 
Key words in this connection are `index theory', `microlocal analysis', or `quantisation'. 
The symbolic structure of operators is
basic for many areas, e.g., in pseudo-differential and Fourier integral operators,
symplectic geometry, Hamiltonian mechanics, spectral theory,
operator algebras, or $K$-theory.

It is not the intention of our remarks to persuade the reader that
all this is relevant and useful. We want to focus here on the
analysis of operators on manifolds with singularities with
questions on the nature of symbols, ellipticity, homotopies, 
index,
and other natural objects. In the singular case those questions arise once
again from the very beginning, similarly as in the early days of
the microlocal analysis on smooth manifolds. Nevertheless, the
analysis on non-smooth and non-compact configurations has a long
history, and there is much experience of different generality with
the solvability of concrete elliptic (and also non-elliptic)
problems with singularities. The notions and inventions from
the smooth case might be a guideline, at least as a special 
case.
However, such an approach has a difficulty in principle: 
There is,
of course, no universal `true analysis' of (linear) partial
differential equations on a smooth manifold, and hence we cannot
expect anything like that in the singular case.

As noted at the beginning there exist different confessions 
 in the fields `ellipticity', or `index theory' on manifolds
with singularities. 
Our choice of aspects is motivated by
an iterative approach for manifolds with higher (regular)
singularities.

If we know nothing and want to see the smooth situation as a
special case we can start from an operator $A \in \Diff^\mu(\R^{n+1})$ (the dimension $n+1$ is taken here for convenience) and interpret the origin of $\R^{n+1}$ as a conical point. 
Introducing polar coordinates $(r, \phi) \in \R_+ \times S^{n}$ we
obtain $A \big|_{\R^{n+1} \setminus \{ 0 \}}$ (briefly denoted again
by $A$) as
\begin{equation}
\label{1011.A01905.eq}  
A = r^{- \mu} \sum_{j=0}^{\mu} a_j(r) \bigl(- r
    \frac{\partial}{\partial r}\bigr)^j     
\end{equation}
with coefficients $a_j \in C^\infty(\ol{\R}_+, \Diff^{\mu-j}(S^{n}))$.   
Note that the  operator
$A = \sum_{j=1}^{n+1} \wt{x}_j \frac{\partial}{\partial \wt{x}_j}$
in polar coordinates takes the form 
$r \displaystyle \frac{\partial}{\partial r}$.
Another example is the Laplace operator
$\Delta = \displaystyle \sum_{k=1}^{n+1} \frac{\partial^2}{\partial \wt{x}_k^2}$ in
polar coordinates:
$$
 \Delta = r^{-2}\Bigl( \bigl(r \frac{\partial}{\partial r}
    \bigr)^2 + (n-1) r \frac{\partial}{\partial r} + 
    \Delta_{S^{n}} \Bigr)     
$$
for the Laplace operator $\Delta_{S^{n}}$ on $S^{n}$. 
Setting
\[  \sigma_\c(A)(w) = \sum_{j=0}^\mu a_j(0)w^j,    \]
$w \in \C$,  we just obtain the so-called conormal symbol of $A$ of order
$\mu$ (with respect to the origin). In this way the operator $A$
suddenly has a second (operator-valued) principal symbolic
component, namely,      
$$
  \sigma_\c(A)(w) : H^s(S^{n}) \to H^{s- \mu}(S^{n}),  
$$  
$s \in \R$. This can be regarded as a component of a `principal
symbolic hierarchy'
\begin{equation}
\label{1011.sigma(A)2204.eq}
\sigma(A) := (\sigma_\psi(A), \sigma_\c(A))
\end{equation}
(with a natural compatibility property between 
$\sigma_\psi(A)$ and $\sigma_\c(A)$). 

For  the identity operator $I$ we obtain the constant family 
$\sigma_\c(I)(w)$
of
identity maps in Sobolev spaces, and the multiplicative rule including
conormal symbols has the form
\begin{equation}
\label{1011.2204eq.eq}
\sigma(AB) = (\sigma_\psi(A)\sigma_\psi(B),
            (T^\nu\sigma_\c(A))\sigma_\c(B))  
\end{equation}
if $A$ and $B$ are differential operators of order $\mu$ and
$\nu$, respectively; $(T^\nu f)(w) = f(w+ \nu)$. In order to
recognise $\sigma_\c(A)$ as a symbol of $A$ in a new context we
have to be aware of the following aspects:
\begin{enumerate}
\item the origin is singled out as a fictitious conical singularity
      {\st{(}}we could have taken any other point{\st{),}} and
      $\sigma_\c(A)$ also depends on the lower order terms of the
      operator $A$ {\st{(}}in any neighbourhood of $0${\st{)}};
\item $\sigma_\c(A)$ refers to a chosen conical structure in
      $\R^{n+1}$, i.e., to a splitting of variables $(r, \phi)
      \in \R_+ \times X$ for $X = S^{n}$ in which we express
      the operator $A$;
\item $\sigma_\c (A)$ is operator-valued, with values in operators on a smooth configuration which is of less singularity order than the conical case.      
\end{enumerate}    

We can pass to other splittings
$(\tilde{r}, \tilde{\phi}) \in \R_+ \times X$ of variables
when the transition diffeomorphism $\R_+ \times X \to \R_+ \times X$, 
$(r, \phi) \to (\tilde{x}, \tilde{\phi})$, is smooth up to $r = 0$. 
There is then a transformation rule of conormal symbols
which just expresses the invariance, cf. \cite{Kapa1}. 
Specific changes $\R_+ \times X \to \R_+ \times X$ which are smooth up to zero are generated by different diffeomorphisms $\R^{n+1} \to \R^{n+1}$ who preserve the origin. 
A Taylor expansion argument then shows that in such a case the  transformation of conormal symbols is induced by a linear isomorphism of $\R^{n+1}$.       
However, in the context of interpreting a point $v$ (here $v =0$) as a conical singularity, we admit arbitrary changes $(r, \phi) \to (\wt{r}, \wt{\phi})$ which are smooth up to $r=0$; then, in general, the transformed operator cannot be reduced to an operator with smooth coefficients across $v$ in the original Euclidean coordinates.

If $A \in \Diff^\mu(M)$ is a differential operator on a smooth
compact manifold $M$ we can fix any $v \in M$ as a fictitious
conical singularity and express $\sigma_\c(A)$ in local
coordinates under a chart $\chi : U \to \R^{n+1}$, $v \in U$, such that 
$\chi(v) = 0$. 
This gives us a conormal symbolic structure
$\sigma_\c(A)$ of operators $A \in \Diff^\mu(M)$. 
Together with the interior principal symbol we obtain a two-component symbolic hierarchy
\eqref{1011.sigma(A)2204.eq}. 
The same can be done for finitely many points 
$\{v_1, \ldots, v_N\} \subset M$; 
this gives us $N$ independent conormal symbols.
Let us restrict the discussion to the simplest case $N = 1$.

What is now ellipticity of $A$?  

We could refer to any classical exposition, and employ the technical background from there.
However, we want to develop the idea here from the point of view of a child who looks conciously for the first time to the sky and realises all the different stars, each of them representing another ellipticity and another index theory.

For the definition we have to foresee a kind of natural analogues
of Sobolev spaces in which the elliptic operators should act as
Fredholm operators. 
Considering $M$ as a manifold with conical
singularity $\{ v \}$ we have the associated stretched manifold $\M$.
By definition $\M$ is obtained from $M\setminus \{ v \}$ by (invariantly)  attaching a sphere $S^{n}$.
This produces a $C^\infty$ manifold with boundary $\partial \M \cong S^{n}$.
For instance, if $M$ is locally near $v$ identified with $\R^{n+1}$ (via a chart) and $v$ with the origin, then $\M$ is locally near $\partial \M$ identified with $\ol{\R}_+ \times S^{n}$ where $(r,\phi) \in \R_+ \times S^{n}$ correspond to polar coordinates in 
$\R^{n+1} \setminus \{ 0 \}$.
There is now a  scale of weighted Sobolev spaces 
${ \cal H}^{s, \gamma}(\M)$ of smoothness  $s \in \R$ and weight  $\gamma \in \R$, contained in
$H_{\loc}^s(M \setminus \{ v \})$. 
Locally near $v$ our operator
(in the chosen splitting of variables $(r, \phi)$) is a
polynomial in vector fields
\[ r \frac{\partial}{\partial r}, \frac{\partial}{\partial
    \phi_1}, \ldots, \frac{\partial}{\partial
    \phi_{n}}      \]
(for $n+1 = \dim M$), up to the weight factor $r^{- \mu}$.

By definition the  stretched manifold $\M$ is a $C^\infty$ manifold with boundary, 
and we can talk about all vector fields that are tangent to $\partial \M$. 
This is a motivation for a definition of the spaces ${\cal H}^{s,
\gamma}(\Int \M)$ for $s \in \N$ as
\begin{align*} 
{\cal H}^{s, \gamma}(\M) := & \{u \in {\cal H}^{0, \gamma}
   (\M) : D^\alpha u \in {\cal H}^{0, \gamma}(\M) \; 
   \text{for all $|\alpha| \leq s$}, \\ 
   & \text{for any tuple  $D = (D_1, \ldots, D_{n+1})$ of vector fields tangent to $\partial \M$ \}, }  
\end{align*}
where
$D^\alpha := (D^{\alpha_1}_1, \cdot \ldots \cdot D^{\alpha_{n+1}}_{n+1})$, and 
${\cal H}^{0, \gamma}(\Int \M)$ is a weighted $L^2$-space, locally 
near the boundary defined as $r^{\gamma-\frac{n}{2}}L^2$ 
$(\R_+ \times \partial \M)$. 
This definition immediately extends to an arbitrary
(compact) manifold $\M$ with boundary $\partial \M \cong X$ for any
closed compact $C^\infty$ manifold $X$, first for $s \in \N$ and
then, by duality and interpolation for all $s \in \R$. 

In particular, for $M = \R^{n+1}$ and $\M = \ol{\R}_+ \times S^{n}$ we
also set
\[  {\cal H}^{s, \gamma}(\R^{n+1} \setminus \{ 0 \}) = {\cal H}^{s,
       \gamma}(\M).      \]
The role of the weight $\gamma \in \R$ may appear somehow
mysterious at the first glance. To give a motivation we 
observe that the conormal symbol
\begin{equation}
\label{1011.sigmac2204.eq}
\sigma_\c(A)(w) : H^s(X) \to H^{s- \mu}(X)
\end{equation}
represents a holomorphic family of Fredholm operators, cf.
Theorem \ref{1011.Ellpsi2204.th}. The ellipticity of $A$ with
respect to $\sigma_\c(A)$ should have the meaning of some
invertibility of the maps \eqref{1011.sigmac2204.eq}, because a
parametrix in the pseudo-differential sense is expected to 
be associated with
the inverse symbol. However, in general, there exists a
non-trivial set $D_A \subset \C$ of points such that
\eqref{1011.sigmac2204.eq} is not invertible. What we know is
that
\begin{equation}
\label{1011.DA2204.eq}
D_A \cap \{ w \in \C : c \leq \re w \leq c' \} 
\end{equation}
is finite for every $c \leq c'$. 
This is a consequence of the parameter-dependent ellipticity of of 
$\sigma_\c(A) \big|_{\Gamma_\beta}$ as a family of operators on $X$ with
parameter $\Im w$, for  every $\beta \in \R$; here 
$\Gamma_\beta := \{ w \in \C : \re w = \beta\}$. 
In our definition of ellipticity we should exclude the set \eqref{1011.DA2204.eq} and
feed in an extra weight information:

\begin{Definition}
\label{1011.Ellc2204.de}
An operator $A \in \Diff^\mu(M)$ is called elliptic with respect
to the symbol
\[  \sigma_\gamma(A) := \Bigl( (\sigma_\psi(A), \sigma_\c(A)
     \big|_{\Gamma_{\frac{n+1}{2} - \gamma}} \Bigr)     \]
for some given weight $\gamma \in \R$, if $ A$ is elliptic with
respect to $\sigma_\psi(A)$, cf. Definition 
{\st{\ref{1011.Ellpsi2104.de}}},
and if
\begin{equation}
\label{1011.For0705.eq}  
\sigma_\c(A)(w) : H^s(X) \to H^{s- \mu}(X)    
\end{equation}
is a family of isomorphisms for all $w \in \Gamma_{\frac{n+1}{2} -
\gamma}$ and some $s \in \R$.     
\end{Definition}

The justification lies in the following result.

\begin{Theorem}
\label{1011.2204th.th}
For an operator $A \in \Diff^\mu(M)$ on $M$ {\st{(}}regarded as a
manifold with conical singularity $v \in M${\st{)}} the following
properties are  equivalent:
\begin{enumerate}
\item The operator $A$ is elliptic with respect to $(\sigma_\psi,
      \sigma_\c \big|_{\Gamma_{\frac{n+1}{2} - \gamma}})$.     
\item A is Fredholm as an operator
     \[  A : {\cal H}^{s, \gamma}(\M) \to {\cal H}^{s- \mu, 
           \gamma - \mu}(\M)          \]
      for some fixed $s \in \R$.
\end{enumerate}
The property {\em(ii)} for a specific $s$ entails the same for all $s \in \R$.
\end{Theorem}

Theorem {\st{\ref{1011.2204th.th}}} extends to the case of a general
stretched manifold $\M$ belonging to a manifold $M$ with conical singularity $\{v \}$. 
As noted before $\M$ is to be replaced in this case by an 
arbitrary compact $C^\infty$ manifold with boundary $\partial \M
\cong X$ (where $X$ is not necessarily a sphere). 
Then $M := \M / \partial \M$ (the quotient space in which $\partial \M$ is collapsed
to a point $v$) is a manifold with conical singularity, cf. Section 3 below. 
The operators $A$ in this case are assumed to belong to
$\Diff_{\deg}^\mu(\M)$ which is defined to be the subspace of all $A \in
\Diff^\mu(\M \setminus \partial \M)$ that are of the form 
\begin{equation}
\label{eq.1.16}
 A =r^{- \mu} \sum_{j=0}^{\mu} a_j(r)\bigl(- r 
       \displaystyle \frac{\partial}{\partial r} \bigr)^j
\end{equation}
in a collar neighbourhood $\cong [0,1) \times X$ of the boundary, 
with $a_j \in C^\infty([0,1), \Diff^{\mu
-j}(X))$. Operators of that kind will also be called of Fuchs type. 
For this situation there exists a pseudo-differential
algebra in analogy to the algebra of pseudo-differential operators on
a $C^\infty$ manifold, here with a principal symbolic hierarchy
\[  \sigma(A) = (\sigma_\psi(A), \ \sigma_\c(A)),    \]
ellipticity, parametrices, etc., cf. also Section 3 below.

Let us now return to differential operators in the Euclidean
space and ask whether there are other natural notions of
ellipticity. 
First, under suitable conditions on the coefficients
of an operator $A \in \Diff^\mu(\R^m)$, $m := n+1$, we have continuity 
$A : H^s(\R^m) \to H^{s- \mu}(\R^m)$ between 
Sobolev spaces globally in $\R^m$ for every $s \in \R$. 
For instance, if an operator
\[  A = \sum_{|\alpha|\leq \mu} a_\alpha(x) D_x^\alpha    \]
has coefficients $a_\alpha(x) \in S_{\cl}^0(\R_x^m)$ we are 
in the situation of the  calculus of operators globally in $\R^m$, cf. Parenti \cite{Pare1}, Cordes \cite{Cord1}, with the principal symbols
\[  \sigma(A) = (\sigma_\psi(A), \sigma_\e(A), 
    \sigma_{\psi,\e}(A)),    \]
see also the notation in Section 3.3 below, or \cite[Section 1.4]{Schu20}.
More generally, there is an analogous notion of ellipticity on
an arbitrary manifold with conical exits to infinity. 
We do not repeat once again the elements of that theory.
 Let us only recall that when we introduce the origin of $\R^m$ as a conical
singularity we have a combination of the principal symbolic
structure near $0$ from the cone calculus and of the exit
symbolic structure near $\infty$, with a principal symbolic
hierarchy
\begin{equation}
\label{1011.sigmaX2304.eq}
\sigma(A) = (\sigma_\psi(A), \sigma_\c(A), \sigma_\e(A),   
             \sigma_{\psi,\e}(A)).
\end{equation}
The adequate scale of weighted Sobolev spaces in this case is
${\cal K}^{s, \gamma}(\R^m \setminus \{ 0 \})$, $s, \gamma \in
\R$, defined by
\begin{equation}
\label{1011.calK1305.eq} 
{\cal K}^{s, \gamma}(\R^m \setminus \{ 0 \}) := \{ \omega u + (1-
   \omega)v : u \in {\cal H}^{s, \gamma}(\R^m \setminus \{ 0 \}), \
   v \in H^s(\R^m) \}
\end{equation}   
for any $\omega \in C_0^\infty(\R^m)$ such that $\omega \equiv 1$ in
a neighbourhood of zero.
In the behaviour with respect to ellipticity there is, of
course, no kind of symmetry under the transformation $(r,
\phi) \to (r^{-1}, \phi)$ when $(r, \phi)$ are
polar coordinates in $\R^m \setminus \{ 0 \}$.	
Similarly, if $X$ is a closed compact $C^\infty$ manifold, we
have a class of natural operators $A$ on the infinite stretched
cone $X^\land := \R_+ \times X \ni (r,x)$ and the principal 
symbolic structure \eqref{1011.sigmaX2304.eq}, see \cite[Section 1.4]{Schu20}.

\begin{Definition}
\label{1011.Elle2304.de}
Let $A \in \Diff^\mu(X^\land)$ be an operator of the form
\begin{equation}
\label{1011.Fuchs0705.eq}  
r^{- \mu} \sum_{j=0}^{\mu} a_j(r) \bigl(- r
    \frac{\partial}{\partial r} \big)^j,      
\end{equation}
with coefficients $a_j \in C^\infty(\ol{\R}_+, \Diff^{\mu-j}(X))$  such that the coefficients locally in $(r,x)$ in a conical subset of $\R^m$ for $r \to \infty$ $(n=m-1 = \dim X)$ admit the pair of exit symbols $\sigma_\e (A)$, $\sigma_{\psi, \e} (A)$ {\em(}this is the case, for instance, when the coefficients are independent of $r$ for large $r$, cf. also 
Section {\em 3.3} below, or 
{\em\cite[Section 1.4.5]{Schu20}).}
Then $A$ is called elliptic with respect to     
\begin{equation}
\label{1011.sigmagamma2304.eq}
\sigma_\gamma(A) := (\sigma_\psi(A), \sigma_\c(A)
\big|_{\Gamma_{\frac{n+1}{2} - \gamma}}, 
   \sigma_\e(A), \sigma_{\psi, \e}(A))
\end{equation}
if all components are elliptic. 
For $\sigma_\psi(A)$ that means $\sigma_\psi(A) \not= 0$ on $T^*X^\land \setminus 0$
and, in local coordinates $x$ on $X$ with covariables $\xi$, and
$\wt{\sigma}_\psi (A) (r,x,\rho,\xi):=
           r^\mu \sigma_\psi(A)(r,x,r^{- 1} \varrho, \xi) \not= 0$ 
for $(\varrho, \xi) \not= 0$, up to $r = 0$. 
For $\sigma_\c(A)$ the condition is that
\eqref{1011.sigmac2204.eq} is bijective for all $w \in
\Gamma_{\frac{n+1}{2} - \gamma}$ and any $s \in \R$. The
ellipticity condition for the exit symbols $\sigma_\e(A)$ and
$\sigma_{\psi,\e}(A)$ will  be explained in detail in
Section {\em 3.3}.
\end{Definition}

For an arbitrary closed compact $C^\infty$ manifold $X$ there is an
analogue of the spaces \eqref{1011.calK1305.eq}, namely,
\begin{equation}
\label{1011.1305calKs.eq}
{\cal K}^{s, \gamma}(X^\land) = \{ \omega u + (1- \omega) v : u \in
    {\cal H}^{s, \gamma}(X^\land), \ v \in H_{\cone}^s(X^\land) \}
\end{equation}
for an arbitrary cut-off function $\omega$. 
The space $H_{\cone}^s(X^\land)$ is locally modelled on the standard 
Sobolev spaces for $r \to \infty$.   

More generally, it may be reasonable to consider the spaces
\begin{equation}
\label{eq.new10}
\s{K}^{s,\gamma; g} (X^\wedge) :=
    \langle r \rangle^{-g} \s{K}^{s,\gamma} (X^\wedge)
\end{equation}
with an extra weight $g \in \R$ at infinity.

Recall that for $s \in \N$, $\gamma \in \R$, we set
\begin{equation}
\label{1011.calH2605.eq}
{\cal H}^{s, \gamma}(X^\land) := \{ u \in r^{\gamma -
   \frac{n}{2}} L^2(\R_+\times X) : (r \partial_r)^k D_x^\alpha
   u \in r^{\gamma - \frac{n}{2}} L^2(\R_+ \times X) \quad
   \textup{for all}\quad k + | \alpha | \leq s \},
\end{equation}
where $D_x^\alpha$ runs over the space of all differential
operators of order $| \alpha |$ on $X$; $n = \dim X$. 
The definition of the spaces ${\cal H}^{s,0}(X^\land)$ for
arbitrary real $s$ follows by duality with respect to the $r^{-
\frac{n}{2}} L^2(\R_+ \times X)$-scalar product and interpolation, and then we set 
${\cal H}^{s,\gamma}(X^\wedge) = r^\gamma {\cal H}^{s,0} (X^\wedge)$ for $\gamma \in \R$.   
   
\begin{Remark}
\label{1011.1305re.re}
The spaces ${\cal K}^{s, \gamma}(X^\land)$ are Hilbert spaces with
suitable scalar products; in particular, we have natural
identifications
${\cal K}^{0,0}(X^\land) = {\cal H}^{0,0}(X^\land) = r^{-
       \frac{n}{2}} L^2(\R_+ \times X)$
for $n = \dim X$. Setting
\[  \kappa_\lambda u(r,x) = \lambda^{\frac{n+1}{2}} 
      u(\lambda r, x), \ \lambda \in \R_+,    \]
for $u \in {\cal K}^{s, \gamma}(X^\land)$ we obtain a strongly
continuous group of isomorphisms
\[  \kappa_\lambda : {\cal K}^{s, \gamma}(X^\land) \to
       {\cal K}^{s- \mu, \gamma - \mu}(X^\land),    \]
$\lambda \in \R_+$.
More generally, on $\s{K}^{s,\gamma; g} (X^\wedge)$ we can consider
\begin{equation}
\label{eq.new11}
\kappa^{g}_\lambda u (r,x) :=
    \lambda^{g+\frac{n+1}{2}} u(\lambda r,x), \ \ \ 
         \lambda \in \R_+.
\end{equation}
\end{Remark}                 

\begin{Theorem}
\label{1011.Elle2304.th}
For an $A \in \Diff^\mu(X^\land)$ as in Definition
{\st{\ref{1011.Elle2304.de}}} the following properties are equivalent:
\begin{enumerate}
\item The operator $A$ is elliptic with respect to
      \eqref{1011.sigmagamma2304.eq}.
\item A is Fredholm as an operator
      $A : {\cal K}^{s, \gamma}(X^\land) \to {\cal K}^{s-\mu,
      \gamma - \mu}(X^\land)$
      for some fixed $s \in \R$.
\end{enumerate}
\end{Theorem}      

Similarly as in Theorem \ref{1011.2204th.th} the property (ii) for a specific $s$ entails the same for all $s\in \R$.
Let us now consider operators $A \in \Diff^\mu(\R_{\wt{x},y}^{m+q})$,
\ $q > 0$, from the point of view of polar coordinates 
$\R^m_{\wt{x}} \setminus \{ 0 \} \ni \wt{x} \to
       (r, \phi) \in \R_+ \times S^{n}$ in $A \big|_{(\R^m \setminus \{ 0 \}) \times \R^q}$ (briefly denoted again by $A$). 
Similarly as at the beginning we obtain $A$ in the form
\begin{equation}
\label{1011.A1905.eq}  
A = r^{- \mu} \sum_{j+|\alpha|\leq \mu} a_{j \alpha}
    (r,y) \bigl(- r\frac{\partial}{\partial r} \bigr)^j 
    (rD_y)^\alpha,       
\end{equation}
now with coefficients $a_{j \alpha} \in C^\infty(\ol{\R}_+ \times
\R^q, \Diff^{\mu-(j+|\alpha|)}(S^{n}))$ for $n=m-1$. 
For instance, for the Laplace operator 
$\Delta =\displaystyle \sum_{k=1}^{n+1}
\frac{\partial^2}{\partial \wt{x}_k^2} + 
\sum_{l=1}^{q} 
\displaystyle \frac{\partial^2}{\partial y^2_l}$ we obtain
$$
    \Delta = r^{-2} \Bigl( \bigl( r \frac{\partial}{\partial r}
    \bigr)^2 + (n-1)  r \frac{\partial}{\partial r} +
    \Delta_{S^{n}} + \sum_{l=1}^{q} r^2
    \frac{\partial^2}{\partial y_l^2} \Bigr).      
$$
This case generates a new operator-valued principal symbol,
namely,
\begin{equation}
\label{1011.sigmaland2304.eq}
\sigma_\land(A)(y, \eta) := r^{- \mu} \sum_{j+|\alpha|\leq \mu}       
   a_{j \alpha}(0,y) \bigl(- r \frac{\partial}{\partial r}
   \bigr)^j (r \eta)^\alpha,
\end{equation}   
$(y,\eta) \in \R^q \times (\R^q \setminus \{ 0 \})$, which is called the principal edge symbol of the operator $A$ with $\R^q$ being interpreted as an edge. 
\eqref{1011.sigmaland2304.eq} represents a family of continuous 
operators
\[  \sigma_\land(A)(y, \eta) : {\cal K}^{s, \gamma}(X^\land)
    \to {\cal K}^{s- \mu, \gamma - \mu}(X^\land),     \]
$X = S^{n}$, for every $s, \gamma \in \R$.

Our new principal symbolic hierarchy here has two components

\begin{equation}
\label{1011.sigmaW2304.eq}    
\sigma(A) = (\sigma_\psi(A), \sigma_\land(A)).
\end{equation}
The second component is homogeneous in the sense
\begin{equation}
\label{1011.kappa1105.eq}  
\sigma_\land(A)(y, \lambda \eta) = \lambda^\mu
       \kappa_\lambda \sigma_\land(A)(y, \eta) 
       \kappa_\lambda^{-1}     
\end{equation} 
for all  $\lambda \in \R_+, (y, \eta) \in \R^q \times (\R^q
\setminus \{ 0 \})$. For the identity operator $I$ we have
$\sigma_\land(I) = \id$ and for the composition 
\[  \sigma_\land(A B) = \sigma_\land(A) \sigma_\land(B)     \]
when $A$ and $B$ are differential operators of order $\mu$ and
$\nu$, respectively.    

Operators of the form \eqref{1011.A1905.eq} including their symbols
\eqref{1011.sigmaW2304.eq} are meaningful on $\R_+ \times X \times
\R^q$ for an arbitrary closed compact $C^\infty$ manifold $X$. In
this connection $\R^q$ is regarded as the edge of the (open
stretched) wedge $X^\land \times \R^q$ with the (open stretched)
model cone $X^\land = \R_+ \times X$. Such operators are called
edge-degenerate. This notation comes from the connection with
`geometric' wedges
\[  W = X^\Delta \times \R^q,  \]
with non-trivial cones
\[  X^\Delta := (\ol{\R}_+ \times X)/(\{ 0 \} \times X)    \]
(in the quotient space $\{ 0 \} \times X$ is identified with a point,
the tip of the cone).

\begin{Remark}
\label{1011.Delta11905.re}
The Laplace-Beltrami operator on $\R_+ \times X \times \R^q \ni (r,x,y)$ belonging
to a Riemannian metric of the form
\[  dr^2 + r^2g_X(r) + dy^2   \]
for a family of Riemannian metrics $g_X(r)$ on a $C^\infty$ manifold
$X$ {\st{(}}smooth in $r \in \ol{\R}_+$ up to $0${\st{)}} is 
edge-degenerate. In
particular, for $q = 0$ we obtain an operator of Fuchs type.
\end{Remark}

\begin{Remark}
\label{1011.re1905.re}
As we see from the preceding discussion, differential operators $A$
in $\R^{n+1}$ {\st{(}}with their standard principal symbolic structure
$\sigma_\psi${\st{)}} secretly belong to several distinguished
societies, namely,
\begin{enumerate}
\item the class of Fuchs type operators with respect to any
      {\st{(}}fictitious{\st{)}} conical singularty;
\item the class of edge-degenerate operators with respect to any
{\st{(}}fictitious{\st{)}} edge {\st{(}}when $n \geq 2${\st{)}}.  
\end{enumerate}   
\end{Remark}

The ellipticity with respect to $\sigma_\land$ in the edge-degenerate
case is a longer story, and we return later on to this point, cf.
Section 2.1.

The question is now whether our operators have other hidden qualities
that we did not notice so far.

The answer is `yes' (when the dimension is not too small).

The operator \eqref{1011.A2104.eq} (in the dimension $n+1$ rather than $n$) written in the form
\eqref{1011.A01905.eq} for $X = S^{n}$ or \eqref{1011.A1905.eq} 
(when the original dimension is equal to $m + q$) for $X = S^{m-1}$ allows us to repeat the game, namely,
to introduce once again fictitious conical points or edges on the
sphere and to represent the coefficients $a_j(r)$ or $a_{j
\alpha}(r,y)$ in Fuchs or edge-degenerate form. 
In order to make the effects more visible we slightly change the transformation of \eqref{1011.A2104.eq} to operators of the form \eqref{1011.A01905.eq} or \eqref{1011.A1905.eq} by
\begin{equation}
\label{1011.+1905.eq}
(x_1, \ldots, x_m) \to (x', r)
\end{equation}
for $x' := (x_1, \ldots, x_{m-1}), r := x_m$. The orthogonal
projection of $S_+^{m-1} := \{ x \in S^{m-1} : x_m > 0 \}$ to
$$
  B := \{ x' \in \R^{m-1} : |x'| < 1 \}    
$$  
is one-to-one; thus $x' \in B$ can be taken as local coordinates on
$S_+^{m-1}$. 
Clearly for local representations it suffices to
consider the hemisphere $S_+^{m-1}$ (up to a rotation). 
The substitution of \eqref{1011.+1905.eq} into \eqref{1011.A2104.eq} 
(when the original dimension is equal to $m+q$)
gives us the operator $A$ in the form \eqref{1011.A1905.eq} with coefficients 
\begin{equation}
\label{1011.aB1905.eq}
a_{j \alpha}(r, y) \in C^\infty(\ol{\R}_+ \times \Omega, \Diff^{\mu
     -(j+ | \alpha |)}(B))
\end{equation}     
(here we
assume $q > 0$; the case $q = 0$ is simpler and corresponds to the
Fuchs type case). 

Now, splitting up the variables in $B$ into $x'= (x'', t, z)$ for
$x'' := (x_1, \ldots, x_{k-1})$, $t = x_k$, $z = (x_{k+1}, \ldots,
x_{m-1})$, the substitution $x' \to (x'',t)$ turns
\eqref{1011.aB1905.eq} into
\begin{equation}
\label{1011.+z1905.eq}
a_{j \alpha}(r,y) = t^{- \mu + j + | \alpha |} \sum_{l+ |\beta|\leq
    \mu -(j+ |\alpha |)} b_{j \alpha; l \beta}(r,y,t,z) 
    \Bigl( - t \frac{\partial}{\partial t} \Bigr)^l (t D_z)^\beta
\end{equation}    
with $\Diff^{\mu-(j+| \alpha |) - (l+ | \beta |)}(C)$ -valued
coefficients $b_{j \alpha; l \beta}$, smooth in $(r, t,y,z)$ (up to
$r = 0$ and $t = 0$) for $C := \{ x'' \in \R^{k-1} : | x'' | <
\frac{1}{2} \}$ and $(t,y)$ varying in a neighbourhood of $(t,y) = 0$.
Inserting \eqref{1011.+z1905.eq} into \eqref{1011.A1905.eq} we
obtain a differential operator of the form
\begin{equation}
\label{1011.Adeg1905.eq}
A = r^{-\mu}t^{- \mu} \widetilde{A}
\end{equation}
where $\widetilde{A}$ is a polynomial of degree $\mu$ in the vector
fields
\begin{equation}
\label{1011.Vect1905.eq}
rt \partial_r, \partial_{x_1}, \ldots, \partial_{x_{k-1}},
   rt \partial_{y_1}, \ldots, rt \partial_{y_q},
   t \partial_t, t \partial_{z_1}, \ldots, t \partial_{z_{q_1}}
\end{equation}
for $q_1 := m - k - 1$ with smooth coefficients in $(r,t,x'', y,z)$
up to $r = 0$, $t = 0$. The operator \eqref{1011.Adeg1905.eq} is
degenerate in a specific way. There are two axial variables $t$ and
$r$, and the principal symbolic hierarchy of $A$ in this case
consists of 3 components    
\[  \sigma(A) = (\sigma_\psi(A)(x, \xi), \sigma_{\land_1}(A)(r,y,z,
      \varrho, \eta, \zeta), \sigma_{\land_2}(A)(y, \eta))    \]
with the standard principal symbol $\sigma_\psi(A)$, the edge symbol
$\sigma_{\land_1}(A)$ of first generation and the edge symbol
$\sigma_{\land_2}(A)$ of second generation.    

The ellipticity of $A$ with respect to $(\sigma_\psi,
\sigma_{\land_1}, \sigma_{\land_2})$ cannot be explained in a few
words. What we mainly need as a new ingredient is an analogue of the
${\cal K}^{s, \gamma}$-spaces on infinite cones whose base spaces are
manifolds with edges. 
This will be discussed later on in Section 6.2.
 Another important point are extra edge conditions which are
necessary both for $\sigma_{\land_1}$ and $\sigma_{\land_2}$; they
also require separate constructions, cf. Section 2.1 for a very
simple model situation.

\begin{Remark}
\label{1011.Delta21905.re}
The Laplace-Beltrami operator on $\R_+ \times \{ \R_+ \times X \times
\R^{q_1} \} \times \R^q \ni (r,t,x,z,y)$ belonging to a Riemannian
metric of the kind
\[  dr^2 + r^2 \{ dt^2 + t^2 g_X (r,t,z,y) + dz^2 \} + dy^2    \]
for a family of Riemannian metrics $g_X(r,t,z,y)$ on a $C^\infty$
manifold $X$ {\st{(}}smooth in the variables up to $r = 0$, $t =
0${\st{)}} has the form
\begin{equation}
\label{1011.LA1905.eq}
t^{- \mu} r^{- \mu} \sum_{j+| \alpha | + l+ | \beta |\leq \mu}
    b_{j \alpha; l \beta}(r,y,t,z) \Bigl(- t \frac{\partial}{\partial
    t} \Bigr)^l (t D_z)^\beta \Bigl(- r t \frac{\partial}{\partial r}
    \Bigr)^j(r t D_y)^\alpha
\end{equation}
{\st{(}}for $\mu = 2${\st{)}} with $\Diff^{\mu-(j+
|\alpha|+l+|\beta|)}(X)$ -valued coefficients which are smooth up to $r = 0$ 
and $t= 0$. Operators of the kind \eqref{1011.LA1905.eq} are 
called
corner-degenerate of second generation.   
\end{Remark}  

It is now clear that the constructions which lead from
\eqref{1011.A2104.eq} to \eqref{1011.A1905.eq} and then to
\eqref{1011.Adeg1905.eq} can be iterated as often as 
we want (only limited by the total dimension). Every time we produce
new types of degenerate operators with higher principal symbolic
structures. As the Remarks \ref{1011.Delta11905.re} and
\ref{1011.Delta21905.re} show, operators with such degeneracies are
connected with higher corner geometries, not merely with fictitious
edges and corners.

Other variants of degenerate operators appear when we introduce in
\eqref{1011.A1905.eq} polar coordinates in different hypersurfaces
not only with respect to the $x$-variables on $S^{m-1}$ but also with
respect to the $y$-variables in $\R^q$. This leads again to new
principal symbolic structures and new ellipticities (provided 
that the
concepts of ellipticity for such higher-degenerate operators 
are developed far enough).

Summing up we see that the process of iteratively blowing up
singularities produces a large variety of degenerate operators, the
ellipticity of which (including their Fredholm property, in which
Sobolev spaces?) was never studied before.

Operators with analogous degeneracies are natural on manifolds with
edge and corner geometries in  general. In the following sections we
develop step by step more ideas, motivation and technicalities
around operators on corner manifolds.

The surprising answer to the question `what is ellipticity' is that
there are many ellipticities, according to the chosen symbolic
structures, most of them being unknown in detail, including all the
consequences for the analysis of the corresponding operators and
their index theory.

In the above examples we saw that the additional principal symbolic components, apart from the standard homogeneous principal symbol on the `main stratum', are contributed by lower-dimensional (here fictitious)  strata.
Since the latter ones are special cases of `real' strata we see that the already derived minimal information has to be a part of the elliptic story also in cases with general polyhedral singularities.

\subsection{Meromorphic symbolic structures}
\label{s.10.1.2}
As we saw in the preceding section differential operators may have
many kinds of symbols, not only the standard homogeneous symbol.
Each of those symbols controls another kind of ellipticity, the
Fredholm property in different scales of (weighted) Sobolev spaces,
and parametrices. One of the most substantial novelties are the
conormal symbols who consist of parameter-dependent operators, in
simplest cases on a closed manifold $X$, the base of the local model
cone (a sphere when the conical point is fictitious).

As Definition \ref{1011.Elle2304.de} shows the conormal symbol
$\sigma_\c(A)(w)$ of an operator $A$ of the form
\eqref{1011.Fuchs0705.eq} refers to a chosen weight $\gamma \in \R$
which is admissible in the sense of the bijectivity of
\eqref{1011.For0705.eq} for all $w \in \Gamma_{\frac{n+1}{2} - \gamma}$.
Nevertheless, the conormal symbol may be of interest in the whole complex plane
as a (for differential operators) holomorphic operator family. The
inverse (in the elliptic case) exists as a meromorphic family of
Fredholm operators between the corresponding Sobolev spaces on $X$.

There are now several interesting questions.
\begin{enumerate}
\item Which is the role of the poles (including Laurent expansions)
      of $\sigma_{\c}(A)^{-1}$ for the operator $A$ or for the 
      nature of solutions $u$ of $A u = f$?
\item Can we control spaces of meromorphic operator functions as
      spaces of conormal symbols in analogy to the scalar symbol
      spaces?
\end{enumerate}

Concerning (i), as we shall illustrate below, there are many 
properties
of solvability that depend on poles and zeros (i.e.,
non-bijectivity points) of the conormal symbols. The main aspects
are asymptotics of solutions and the Fredholm index (especially, the
relative index when we change weights).      

For (ii) we have to specify the meaning of \textup{`control'}. 
The point is
that every operator $A$ generates a pattern of poles and zeros of
its conormal symbol $\sigma_\c(A)$ which is individually determined
by $A$. Spaces of such meromorphic symbols contain all possible
patterns of that kind, i.e., such symbol spaces encode 
the asymptotic behaviour of solutions, the relative index
behaviour  and other effects, influenced by the conormal
symbols for all possible operators $A$ at the same time.

This is far from being a purely
\textup{`administrative'} discussion on the structure of the calculus. 
In fact, if we pass to edge singularities and
edge-degenerate operators $A$ we have subordinate conormal symbols
\[  \sigma_{\c}\sigma_\land(A) (y,w) = \sum_{j=0}^{\mu} a_{j0}
       (0,y)w^j      \]
which are families varying with the edge variable $y$, and, of
course, all data connected with meromorphy (including the position
and multiplicities of poles and zeros) depend on the variable $y$.

Let us have a look at a very simple example which shows how the
operator determines individual asymptotics of solutions near $r =
0$.

Let
\begin{equation}
\label{1012.Eq0705.eq}
Au := \sum_{j=0}^{\mu} a_j \bigl( - r
     \frac{\partial}{\partial r} \bigr)^j u(r) = f
\end{equation}     
be an equation of Fuchs type on $\R_+$ with constant 
coefficients (any weight factor as in \eqref{eq.1.16} in front the operator
is not really essential in the conical case).
Then, for
\begin{equation}
\label{1012.h2505.eq}  
h(w) := \sum_{j=0}^{\mu} a_j w^j     
\end{equation}
the equation \eqref{1012.Eq0705.eq} takes the form
$\op_M(h) u = f$.     
Here $\op_M(h) u = M^{-1}h M$ is the pseudo-differential operator
based on the Mellin transform $M$ in $L^2(\R_+)$, $Mu(w) =
\int_{0}^{\infty} r^{w-1} u(r)dr$.
Under the ellipticity condition
\begin{equation}
\label{1012.ell0805.eq}
\sigma_{\c}(A)(w) = h(w) \not= 0\quad \text{on}\quad
                        \Gamma_{\frac{1}{2}}
\end{equation}
we  can realise $\op_M(h^{-1})$ as a continuous operator
$L^2(\R_+) \to L^2(\R_+)$, and we find the solution in the form
\begin{equation}
\label{1012.sol0805.eq}  
u(r) = \op_M(h^{-1})(f)(r) = M_{w \to r}^{-1} \bigl(h^{-1}(w)
           M(f)(w)\bigr).     
\end{equation}
The Mellin transform $M$ is operating not only on $L^2(\R_+)$ but on
subspaces $L_P^2(\R_+)$ of functions with asymptotics of type
\begin{equation}
\label{1012.P0805.eq}  
P = \{ (p_j, m_j) \}_{j \in \N}.      
\end{equation}
Here $p_j \in \C$, $m_j \in \N$, $\re p_j < \frac{1}{2}$, $\re p_j
\to - \infty$ as $j \to \infty$. The space $L_P^2(\R_+)$ is defined
to be the subspace of all $u \in L^2(\R_+)$ such that for every
$\beta \in \R$ there is an $N = N(\beta)$ with
\[  \omega(r) \bigl\{ u(r) - \sum_{j=0}^{N} \sum_{k=0}^{m_j} c_{jk}
     r^{- p_j}  \log^k r\} \in r^\beta L^2(\R_+)   \]
with coefficients $c_{jk} \in \C$ depending on $u$, for any
cut-off function $\omega$ (i.e., an element of $C_0^\infty ({\ol{\R}}_+)$
that is equal to $1$ near $r = 0$).

\begin{Theorem}
\label{1012.th0705.th}
Let $A$ satisfy the conditions $a_\mu \not= 0$ and \eqref{1012.ell0805.eq}, and let 
$f
\in L^2(\R_+)$. Then the equation $A u = f$ has a unique solution $u
\in L^2(\R_+)$. Moreover, $f \in L_Q^2(\R_+)$ for some asymptotic
type $Q$ entails $u \in L_P^2(\R_+)$ for some resulting asymptotic
type $P$.
\end{Theorem}     			

We have, of course, more
regularity of solutions than in $L^2$ (cf. Remark \ref{1042.re2110.re} below).

The meromorphic function $h^{-1}(w)$ belongs to a category of
spaces that are defined as follows.

Let
\begin{equation}
\label{1012.R0805.eq}
R = \{ (r_j, n_j) \}_{j \in \Z}
\end{equation}
be a sequence of pairs $\in \C \times \N$, such that $|\re r_j|
\to \infty$ as $|j| \to \infty$. Set $\pi_\C R := \bigcup_{j \in
\Z} \{ r_j \}$.
A $\pi_\C R -$excision function is any $\chi_R \in
C^\infty(\C)$ such that $\chi_R(w) = 0$ for $\dist(w, \pi_\C R) <
c_0$, $\chi_R(w) = 1$ for $\dist(w, \pi_\C R) > c_1$ for
certain $0 < c_0 < c_1$.

\begin{Definition}
\label{1012.neuM0507.de}
Let ${\cal M}_R^\nu$ denote the space of all meromorphic
functions $f$ in the complex plane with poles at $r_j$ of
multiplicity $n_{j} +1$ such that
$\chi_R(w)f(w) |_{\Gamma_\beta} \in S_{\cl}^\nu(\Gamma_\beta)$
for every $\beta \in \R$ uniformly in compact $\beta$-intervals; 
here $\chi_R$ is any $\pi_\C R$-excision function, and    
 $S_{\cl}^\nu(\Gamma_\beta)$ is the space of all classical
symbols of order $\nu$ in  the covariable 
$\Im w$ for $w \in \Gamma_\beta$ with constant coefficients, cf.
Definition {\st{\ref{1013.Sy2005.de}}} below. For $\pi_\C R = \emptyset$
the corresponding space will be denoted by ${\cal M}_{{\cal
O}}^\nu$.
\end{Definition}

In our example, if $a_\mu \not=0$, we have
\begin{equation}
\label{1012.R2505.eq}  
h^{-1}(w) \in {\cal M}_R^{-\mu}     
\end{equation}
for some $R$ of the kind \eqref{1012.R0805.eq} determined by the
zeros of $h(w)$ in the complex plane.

In order to obtain the regularity result of Theorem
\ref{1012.th0705.th} with asymptotics we consider the solutions
\eqref{1012.sol0805.eq} and observe that the space
\[  M_{r \to w} L_Q^2(\R_+)    \]
for an asymptotic type $Q = \{(q_j, l_j)\}_{j \in \N}$,
$\pi_\C Q \subset \{ w: \re w < \frac{1}{2} \}$, can be
characterised as the space ${\cal A}_Q^0$ of those meromorphic
functions $m(w)$ in the half-plane $\re w < \frac{1}{2}$ 
with poles at $q_j$ of multiplicity $j+1$, such  that for
every $\pi_\C Q$-excision function $\chi_Q$ we have
\begin{equation}
\label{1012.L20805.eq}
\chi_Q(w)m(w)|_{\Gamma_\beta} \in L^2(\Gamma_\beta)
\end{equation}
for all $\beta \leq \frac{1}{2}$ (the meaning for $\beta =
\frac{1}{2}$ is that $\chi_Q(\beta + i \varrho) m(\beta + i
\varrho)$ has an $L^2(\R_\varrho)$-limit for $\beta \nearrow
\frac{1}{2}$), and \eqref{1012.L20805.eq} holds uniformly in
compact $\beta$-intervals $\subset (- \infty, \frac{1}{2}]$.

In other words, $M : L^2(\R_+) \to L^2(\Gamma_{\frac{1}{2}})$
restricts to an isomorphism
\[  M : L_Q^2(\R_+) \to {\cal A}_Q^0 \quad \text{for every}\quad Q.
                                              \]
Now $M f \in {\cal A}_Q^0$ entails
$h^{-1}(w) M f(w) \in {\cal A}_P^0$
for some asymptotic type \eqref{1012.P0805.eq}. Then the
relation \eqref{1012.sol0805.eq} gives us immediately $u \in
L_P^2(\R_+)$.

This consideration shows by a very simple example how the
regularity of solutions near $r = 0$ is influenced by the
operator $A$. Namely, the resulting asymptotic type is
determined (apart from $Q$ on the right hand side) by
\[ R \bigl|_{\re w < \frac{1}{2}} \quad \text{for}\quad
      h^{-1} \in {\cal M}_R^{- \mu}.     \]
Here $R \big|_{\re w < \delta} := \{ (q,n) \in R : 
\re q < \delta \}$.
The same questions can be asked for $r \to \infty$, or both for $r \to
0$ and $r \to \infty$. Let
\begin{equation}
\label{1012.P02505.eq}  
P^0 = \bigl\{ (p_j^0, m_j^0) \bigr\}_{j \in \N}, \quad
    P^\infty = \{ (p_j^\infty, m_j^\infty) \}_{j \in \N}   
\end{equation}
be asymptotic types, $P_0$ responsible for $r \to 0$ as before
and $P_\infty$ for $r \to \infty$ (where $\re p_j^\infty >
\frac{1}{2}$, $\re p_j^\infty \to \infty$ as $j \to \infty$).
Let $L_{P^0, P^\infty}^2(\R_+)$ be the subspace of all $u \in
L_{P^0}^2(\R_+)$ such that for every $\beta \in \R$ there is an
$N = N(\beta)$ such that
\[ (1- \omega(r))
              \Big\{ u(r) - \sum_{j = 0}^{N} \sum_{k=0}^{m_j^\infty}
    d_{jk} r^{- p_j^\infty} \log^k r \Big\} \in r^{- \beta}
    L^2(\R_+)     \]
for some coefficients $d_{jk}$ depending on $u$, and a cut-off
function $\omega(r)$.            

Then a simple generalisation of the regularity result of
Theorem \ref{1012.th0705.th} with asymptotics is that
\[  A u = f \in L_{Q^0,Q^\infty}^2(\R_+) \Rightarrow
       u \in L_{P^0,P^\infty}^2(\R_+)     \]
for every pair $(Q^0, Q^\infty)$ of asymptotic types with some
resulting $(P^0, P^\infty)$.

The correspondence
\begin{equation}
\label{1012.As2505.eq}
Q^0 \to P^0  \quad \text{comes from} \quad R \big|_{\re w <
                          \frac{1}{2}}\quad \text{and}\quad
Q^\infty \to P^\infty \quad \text{from} \quad
                          R \big|_{\re w > \frac{1}{2}}   
\end{equation}
by a simple multiplication of meromorphic functions in the
complex Mellin $w$-plane.	

In other words, the asymptotic type $R$ of the Mellin symbol
$h^{-1}(w)$ is subdivided into parts in different half-planes,
responsible for the asymptotics of solutions for $r \to 0$ and
$r \to \infty$.		  			        

Let us now slightly change the point of view and ask solutions of the equation 
\eqref{1012.Eq0705.eq} for $f \in r^\gamma L^2(\R_+) := L^{2,
\gamma}(\R_+)$ rather than
$L^2(\R_+)$, for some  weight $\gamma \in \R$.

To this end we first recall that the Mellin transform 
$M u = \int_{0}^{\infty} r^{w-1}u(r)dr
       \big|_{\Gamma_{\frac{1}{2} - \gamma}}$,
$u \in C^\infty_0 (\R_+)$, extends to an isomorphism
\[  M_\gamma : L^{2,\gamma} (\R_+) \to L^2(\Gamma_{\frac{1}{2} -
    \gamma})      \]
for every $\gamma \in \R$ (which is equal to $M$ for $\gamma = 
0$). Then, having a Mellin symbol \eqref{1012.h2505.eq}, we can
form the associated operator
\begin{equation}
\label{1012.neu12206.eq}  
u \to M_\gamma u \to h \big|_{{\Gamma_{\frac{1}{2} - \gamma}}}
     M_\gamma u \to M_\gamma^{-1} 
     \Bigl(h \big|_{{\Gamma_{\frac{1}{2}-
     \gamma}}}\Bigr) M_\gamma u =: \op_M^\gamma(h)u.    
\end{equation}       
We also write
\begin{equation}
\label{1012.neu22206.eq}
\op_M(.) := \op_M^0(.).
\end{equation}

Observe that
\[  \op_M(h) u = r^\gamma \op_M
      (T^{- \gamma} h)r^{- \gamma}u      \]
for $(T^{- \gamma} h)(w) := h(w- \gamma)$, for arbitrary 
$\gamma \in \R$, and $u \in C_0^\infty(\R_+)$. Considering the 
equation \eqref{1012.Eq0705.eq} for
$f := r^\gamma f_0$, $u := r^\gamma u_0$, for a given element $f_0
\in L^2(\R_+)$ it follows that
\begin{equation}
\label{1012.Agamma255.eq}
\sum_{j=0}^{\mu} a_j \Bigl( - r \frac{\partial}{\partial r}
      \Bigr)^j (r^\gamma u_0) = r^\gamma f_0 = r^\gamma
      \sum_{j=0}^{\mu} a_j \Bigl( - r \frac{\partial}{\partial r} -
      \gamma \Bigr)^j u_0.
\end{equation}                 
Thus the equation 
\[  
\sum_{j=0}^{\mu} a_j 
      \Bigl( - r \frac{\partial}{\partial r} \Bigr)^j 
           (r^\gamma u_0)  = r^\gamma f_0,\quad f_0 \in L^2(\R_+)    \]
is equivalent to
\begin{equation}
\label{1012.A02505.eq}
\op_M(T^{- \gamma} h)u_0 = f_0,
\end{equation}
and solutions $u = r^\gamma u_0$ of \eqref{1012.Agamma255.eq}
follow from solutions $u_0$ of \eqref{1012.A02505.eq}.

Analogously as before we form spaces of the kind
\[  L_{P_\gamma}^{2, \gamma}(\R_+)\quad \text{or}\quad
    L_{P_\gamma^0, P_\gamma^\infty}^{2,\gamma}(\R_+)     \]
for asymptotic types $P_\gamma$ or $(P_\gamma^0, P_\gamma^\infty$) 
defined in a similar manner (and with a similar meaning) as
before. More precisely, we have
\[  P_\gamma^0 = \bigl\{ (p_j^0 - \gamma, m_j^0 \bigr\}_{j 
       \in \N},\quad
    P_\gamma^\infty = \bigl\{ (p_j^\infty - \gamma, m_j^\infty)
       \bigr\}_{j \in \N}    \]
for sequences as in \eqref{1012.P02505.eq}; then 
$u(r) \in
      L_{P_\gamma^0, P_\gamma^\infty}^2(\R_+)$ has asymptotics of type
$P_\gamma^0$ for $r \to 0$ and of type $P_\gamma^\infty$ for
$r \to \infty$.      

As a corollary of Theorem \ref{1012.th0705.th} we now  obtain the following result:

\begin{Theorem}
\label{1012.th2505.th}
Let $A$ satisfy the conditions $a_\mu \not= 0$ and 
\[  \sigma_\c(A)(w) = h(w) \not= 0\quad \text{on}\quad
       \Gamma_{\frac{1}{2} - \gamma}.    \]
Then the equation $A u = f \in L^{2, \gamma}(\R_+)$ has a unique
solution $u \in L^{2, \gamma}(\R_+)$. Moreover, $f \in
L_{Q_\gamma^0, Q_\gamma^\infty}^2(\R_+)$ for some asymptotic types
$(Q_\gamma^0, Q_\gamma^\infty)$ entails $u \in L_{P_\gamma^0,
P_\gamma^\infty}^{2, \gamma}(\R_+)$ for resulting asymptotic types
$(P_\gamma^0, P_\gamma^\infty)$. Analogously, we can ignore
asymptotics at $\infty$ and conclude from $f \in
L_{Q_\gamma}^{2, \gamma}(\R_+)$ solutions $u \in
L_{P_\gamma}^{2, \gamma}(\R_+)$ for every $Q_\gamma$ with some
resulting $P_\gamma$.
\end{Theorem}

This is immediate from the reformulation of
\eqref{1012.Agamma255.eq} as \eqref{1012.A02505.eq}.     

What we also see in analogy of \eqref{1012.As2505.eq}  in the
weighted case is that the transformation
\begin{equation}
\label{1012.Asgamma255.eq}
Q_\gamma^0 \to P_\gamma^0\quad \text{comes from}\quad 
   R\big|_{\re w < \frac{1}{2} - \gamma}\quad \text{and}\quad
   Q_\gamma^\infty \to P_\gamma^\infty\quad \text{from}\quad
   R \big|_{\re w > \frac{1}{2} - \gamma}.
\end{equation}

\begin{Remark}
\label{1012.re2505.re}
In principle, the generalisation of Theorem
{\st{\ref{1012.th0705.th}}} from $\gamma_0 = 0$ to arbitrary $\gamma
\in \R$ is completely trivial. Nevertheless, something very
strange happed during the change to the new weight. Comparing
\eqref{1012.As2505.eq} and \eqref{1012.Asgamma255.eq} we see that
some part of the \textup{`meromorphic information'} of the
inverted conormal symbol $\sigma_\c(A)^{-1}$, encoded by $R$,
which is responsible for the asymptotics of solutions in 
$L^{2,\gamma_0}(\R_+)$ for $r \to \infty$ may suddenly be responsible
for the asymptotics of solutions in $L^{2, \gamma}(\R_+)$ for $r
\to 0$, and vice versa, according to the specific position of $R$ 
relative to the weight lines $\Gamma_{\frac{1}{2} - \gamma_0}$ and
$\Gamma_{\frac{1}{2} - \gamma}$, respectively. In the extremal
case, since $\pi_\C R$ {\em(}in the case of a differential operator{\em)} is finite for our differential operator $A$,
we may have
\begin{equation}
\label{1012.gamma02505.eq}
\pi_\C R \subset \{ w: \re w < \frac{1}{2} - \gamma_0 \},
\end{equation}
or
\begin{equation} 
\label{1012.gamma12505.eq}
\pi_\C R  \subset \{ w : \re w > \frac{1}{2} - \gamma \}.
\end{equation}

In  the case \eqref{1012.gamma02505.eq} there is no influence of $R$ to
the asymptotics of solutions in $L^{2, \gamma_0}(\R_+)$ for $r \to
\infty$ but \textup{`very much'} for $r \to 0$, in the case 
\eqref{1012.gamma12505.eq} for solutions in $L^{2, \gamma}(\R_+)$
it is exactly the opposite.
\end{Remark}    

The situation becomes even more mysterious if we pass from the
operator $A$ to its formal adjoint $A^*$ with respect to the
scalar product of $L^2(\R_+)$. Writing $A$ in the form $A =
\op_M(h)$ it follows that
\[  A^* = \op_M(h^*)\quad \text{for}\quad
       h^*(w) = \ol{h}(1- \ol{w}).      \]
Similarly as \eqref{1012.R0805.eq} there is then an $R^* = \{
(r_j^*, n_j^*) \}_{j \in \Z}$ such that
$(h^*)^{-1}(w) \in {\cal M}_{R^*}^{- \mu}$,     
and it is obvious in this case that $w \in \pi_\C R$ is equivalent
to $1 - \ol{w} \in \pi_\C R^*$, cf. the relation
\eqref{1012.R2505.eq}.       

\begin{Remark}
\label{1012.2505re.re}
The influence of $R$ to the asymptotics of solutions of $Au = f
\in L^2(\R_+)$ for $r \to \infty$ $(r \to 0)$ is
translated to an influence of to the asymptotics of solutions of
$A^*v = g \in L^2(\R_+)$ for $r \to 0$ $(r \to \infty)$. In fact,
there is a natural bijection $R \to R^*$ induced by $w \to 1 -
\ol{w}$ under which $\pi_\C R \cap \{ \re w \lessgtr \frac{1}{2}
\}$ is transformed to $\pi_\C R^* \cap \{ \re w \gtrless
\frac{1}{2} \}$, cf. the relation \eqref{1012.As2505.eq}.
\end{Remark}

Let us now pass to operators of the form
\begin{equation}
\label{1012.A2605.eq}
A = \sum_{j=0}^{\mu} a_j \bigl( - r
\frac{\partial}{\partial r} \bigr)^j
\end{equation}
with coefficients $a_j \in \Diff^{\mu-j}(X)$ for an arbitrary
closed compact $C^\infty$ manifold $X$ of dimension $n$ (it also makes sense to admit $a_j \in
C^\infty(\ol{\R}_+, \Diff^{\mu-j}(X))$, together  with some control of
the $r$-dependence for $r \to \infty$, cf. also Remark \ref{1042.re2110.re} below).

As noted in Section 1.1 we have the pair of symbols
\eqref{1011.A2104.eq}, especially, the conormal symbol
$\sigma_\c(A)(w) = \sum^\mu_{j=0} a_j w^j$ which represents a family of continuous
operators
\begin{equation}
\label{eq.50}
  \sigma_\c(A)(w) : H^s(X) \to H^{s- \mu}(X)     
\end{equation}
for every $s \in \R$, holomorphic in $w \in \C$.
The generalisation of the discussion before to the case $n =
\dim X > 0$ gives rise to some substantial new aspects.

Assuming $\sigma_\psi$-ellipticity of $A$ in the sense that the standard homogeneous principal symbol
$\sigma_\psi (A) (r,x,\rho,\xi)$ does not vanish on $T^* X^\wedge \setminus 0$ and that
\[  \tilde{\sigma}_\psi(A)(x, \varrho, \xi) :=
       \sigma_\psi(A)(r,x, r^{-1}\varrho, \xi)     \]
satisfies the condition
\[  \tilde{\sigma}_\psi(A)(x, \varrho, \xi) \not= 0\quad
      \textup{for all}\quad (\varrho, \xi) \not= 0, \quad
      \textup{up to}\quad r= 0,      \]      
the operators \eqref{eq.50} are parameter-dependent elliptic on $X$ with the parameter $\Im w$ for
$w \in \Gamma_\beta = \{ w \in \C : \re w = \beta \}$ for every $\beta \in \R$.
The operator function \eqref{eq.50} belongs to a space ${\cal M}^\mu_{\cal O}(X)$ which is defined as follows.
First let $L^\mu_\cl (X; \R^l)$ denote the space of all parameter-dependent classical pseudo-differential operators on $X$ of order $\mu \in \R$ with the parameter $\lambda \in \R^l$, that is, the local amplitude functions $a(x,\xi,\lambda)$ are classical symbols in the covariables $(\xi, \lambda) \in \R^{n+l}$, while
$L^{-\infty} (X; \R^l) := {\cal S} (\R^l, L^{-\infty} (X))$ with $L^{-\infty}(X)$ being the space of smoothing operators on $X$, i.e., with kernels in $C^\infty (X \times X)$ (and identified with $L^{-\infty} (X)$, including the Fr\'echet space structure from $C^\infty (X \times X)$).
For $l=0$ we simply write $L^\mu_\cl (X)$.

Then ${\cal M}^\mu_{\cal O} (X)$ is the space of all 
$h(w) \in {\cal A} (\C_w, L^\mu_\cl (X))$ (i.e., holomorphic $L^\mu_\cl (X)$-valued functions) such that
$f(\beta +i \rho) \in 
       L^\mu_\cl (X; \R_\rho)$ for every $\beta \in \R$, uniformly in finite $\beta$-intervals.

For an operator \eqref{1012.A2605.eq}  we then have $\sigma_\c (A) \in {\cal M}^\mu_{\cal O} (X)$.
In addition, the $\sigma_\psi$-ellipticity of $A$ has the consequence that \eqref{eq.50} is invertible for all $w \in \C \setminus D$ for a certain discrete set $D$.
In order to describe $\sigma^{-1}_{\c} (A) (w)$ we define sequences
\begin{equation}
\label{eq.51}
R = \{ (r_j, n_j, N_j) \}_{j\in \Z},
\end{equation}
where $\pi_\C  R := \{ r_j \}_{j \in \N} \subset \C$
is a subset such that $|\re r_j| \to \infty$ as $|j| \to \infty$, 
$n_j \in \N$, and 
$N_j \subset L^{-\infty} (X)$ are  finite-dimensional subspaces of operators of finite rank.

If $E$ is a Fr\'echet space and $U \subseteq \C$ an open set, by ${\cal A} (U,E)$ we denote the space of all holomorphic functions in $U$ with values in $E$.
\begin{Definition}
\label{d.}
\begin{enumerate}
\item
Let ${\cal M}^{-\infty}_R (X)$ denote the space of all 
$f(w) \in {\cal A}(\C \setminus \pi_\C R, L^{-\infty}_\cl (X))$ that are meromorphic
with poles at $r_j$ of multiplicity $n_j +1$ and Laurent coefficients at
$(w -r_j)^{-(k+1)}$ belonging to $N_j$ for all $0 \leq k \leq n_j$, such that for any $\pi_\C R$-excision function $\chi_R(w)$ we have
$( \chi_{R} f) (\beta + i \rho) \in L^{-\infty} (X; \R_\rho)$
for every $\beta \in \R$, uniformly in compact $\beta$-intervals;
\item
for $\mu \in \R$ we define
\begin{equation}
\label{eq.new20}
{\cal M}^\mu_R (X) := 
     {\cal M}^\mu_{\cal O} (X) + {\cal M}^{-\infty}_R (X).
\end{equation}
\end{enumerate}
\end{Definition} 
 
In order to describe regularity and asymptotics of solutions to elliptic equations $Au =f$ we can introduce subspaces $\s{H}^{s,\gamma}_{P^0,P^\infty} (X^\wedge)$ of 
$\s{H}^{s,\gamma} (X^\wedge)$ of elements with asymptotics of types
\begin{equation}
\label{eq.52}
P^0 = \{ (p^0_j, m^0_j, L^0_j) \}_{j \in \N},   \ \ \
P^\infty = \{ (p^\infty_j, m^\infty_j, L^\infty_j) \}_{j \in \N},   
\end{equation}
where the meaning is quite similar as before for the case $\dim X=0$.
In \eqref{eq.52} we assume $p^0_j, p^\infty_j \in \C$,
$m^0_j, m^\infty_j \in \N$,
$\re p^0_j < \frac{n+1}{2} -\gamma$,
$\re p^\infty_j > \frac{n+1}{2} -\gamma$ for all $j$, and
$\re p^0_j \to -\infty$, $\re p^\infty_j \to +\infty$ for $j \to \infty$;
moreover $L^0_j, L^\infty_j \subset C^\infty (X)$ are subspaces of finite dimension.
Then $u(r,x) \in \s{H}^{s,\gamma}_{P^0, P^\infty} (X^\wedge)$
means that there are coefficients
$c^0_{j,k} \in L^0_j$ and $c^\infty_{j,k} \in L^\infty_j$ for all $0 \leq k \leq m^0_j$ and
$0 \leq k \leq m^\infty_j$, resprectively, $j \in \N$, such that for every $\beta \in \R$ there exists an $N = N(\beta)$ such that
\begin{equation}
\label{eq.p0}
\omega(r)
 \Big\{ u(r,x) -
   \sum^N_{j=0} \sum^{m^0_j}_{k=0}
       c^0_{jk} (x) r^{-p^0_j} \log^k r \Big\} \in
 \omega (r) \s{H}^{s,\gamma+\beta} (X^\wedge)
\end{equation}
and
\begin{equation}
\label{eq.pinfty}
(1 - \omega(r))
 \Big\{ u(r,x) -
   \sum^N_{j=0} \sum^{m^\infty_j}_{k=0}
       c^\infty_{jk} (x) r^{-p^\infty_j} \log^k r \Big\} \in
 (1- \omega (r)) \s{H}^{s,\gamma-\beta} (X^\wedge).
\end{equation}
Here $\omega (r)$ is any cut-off function.

We can also consider subspaces of elements $u \in \s{H}^{s,\gamma} (X^\wedge)$ of the kind
$\s{H}^{s,\gamma}_{P^0} (X^\wedge)$ (or $\s{H}^{s,\gamma}_{P^\infty} (X^\wedge)$) where we observe asymptotics of type $P^0$ (or $P^\infty$) only for $r \to 0$ by requiring \eqref{eq.p0} or 
($r \to \infty$ by \eqref{eq.pinfty}).
Now a general theorem which summarises several features on operators of the kind \eqref{1012.A2605.eq} is the following.
First, let us write
\begin{equation}
\label{eq.53}
\op^\delta_M (h) u := r^\delta \op_M (T^{-\delta} h) r^{-\delta}
\end{equation}
for any $h(w) \in \s{M}^\mu_R (X)$ and $\delta \in \R$ such that 
$\pi_\C R \cap \Gamma_{\frac{1}{2}-\delta} = \emptyset$
(observe that the notation \eqref{eq.53} is not a contradiction to \eqref{1012.neu22206.eq},
because for $h \in \s{M}^\mu_{\s{O}} (X)$ we have $\op_M (h)u = \op^\delta_M (h) u$ for 
$u \in C^\infty_0 (X^\wedge))$.
\begin{Theorem}
\label{t.1.21}
Let \eqref{1012.A2605.eq} be $\sigma_\psi$-elliptic and write $h(w) = \sum^\mu_{j=0} a_j w^j$.
Then we have 
$h^{-1} (w) \in \s{M}^{-\mu}_R (X)$ for some $R$ as in \eqref{eq.51}.
For every $\gamma \in \R$ such that $\pi_\C R \cap \Gamma_{\frac{n+1}{2}-\gamma}= \emptyset$
the operator $A$ induces an isomorphism
\begin{equation}
\label{eq.new1-21}
A : \s{H}^{s,\gamma} (X^\wedge) \to 
     \s{H}^{s-\mu, \gamma} (X^\wedge)
\end{equation}
for every $s\in \R$, and the inverse has the form
$A^{-1} = \op^{\gamma-\frac{n}{2}}_M (h^{-1})$.
Moreover, for every pair of asymptotic types $Q^0, Q^\infty$ as in \eqref{eq.52} there is an analogous pair $P^0,P^\infty$ such that
\begin{equation}
\label{eq.new2-21}
Au \in \s{H}^{s-\mu,\gamma}_{Q^0,Q^\infty} (X^\wedge) \Rightarrow
    u \in \s{H}^{s,\gamma}_{P^0,P^\infty} (X^\wedge)
\end{equation}
for  every $s \in \R$.
\end{Theorem}
\begin{Remark}
\label{r.1.22}
For the asymptotics of solution of the equation $Au=f$ for a $\sigma_\psi$-elliptic operator 
\eqref{1012.A2605.eq}
we have a simple analogue of the Remarks {\em\ref{1012.re2505.re}} and {\em \ref{1012.2505re.re}}, now referring to the spaces 
$\s{H}^{s,\gamma}(X^\wedge)$ and subspaces with asymptotics, with $R$ being as in Theorem {\em\ref{t.1.21}} {\em(}this is a more precise information also for $\dim X=0$ compared with the discussion in $L^2$ spaces on the half-axis{\em)}.
\end{Remark}
The latter results have natural analogues for the case of $\sigma_\psi$-elliptic operators \eqref{1012.A2605.eq} when the coefficients $a_j$ depend on $r$ in a controlled manner (e.g., smooth up to $r=0$).
Instead of unique solvability we then have a Fredholm operator 
\eqref{eq.new1-21} (under an analogous condition on the weight $\gamma$), and the relation \eqref{eq.new2-21} can be interpreted as a result on elliptic regularity in spaces with asymptotics.
Results of that kind exist in many variants, e.g., for finite asymptotic types,  or so-called continuous asymptotic types,  cf. \cite{Schu32}.
In the framework of pseudo-differential parametrices which exist in the cone algebra, acting as continuous operators in weighted Sobolev spaces and subspaces with asymptotics, it is important to stress the conormal symbolic structure, i.e., the spaces of meromorphic operator function in 
\eqref{eq.new20}.
The Mellin asymptotic types \eqref{eq.51} in those spaces of symbols vary over all possible configurations of that kind which is an enormous input of a priori information in the corresponding cone calculus with asymptotics.
Starting from a specific asymptotic type $R$, known by the inverse of the conormal symbol of the operator $A$, the correspondence
$$
(Q^0, Q^\infty) \to (P^0, P^\infty)
$$
in the sense of \eqref{eq.new2-21} is completely determined.
However, to really compute $R$ may be a difficult task in concrete cases.
For every individual operator $A$ we have to solve a corresponding non-linear eigenvalue problem, and the asymptotic information $(P^0, P^\infty)$ on the solution is not merely defined by the homogeneous principal symbol $\sigma_\psi (A)$ of the elliptic operator $A$ but by the global spectral behaviour of operators on the base $X$ of the cone which is also influenced by the lower  order terms.

Similar observations are true when we are only interested in the asymptotics for $r \to 0$ (or $r \to \infty$) alone.
In this connection later on we shall employ $\s{K}^{s,\gamma}$-spaces and weighted Schwartz spaces with asymptotics.
Set
$$
 \s{K}^{s,\gamma}_P (X^\wedge) :=
   \{ \omega u + (1-\omega)v :
      u \in \s{H}^{s,\gamma}_P (X^\wedge), 
       v \in \s{K}^{s,\gamma} (X^\wedge)   \}, 
$$
cf. the formula \eqref{1011.1305calKs.eq}.
Here $P = \{ (p_j, m_j, L_j )\}_{j \in \N}$
is an asymptotic type as the first one in the formula \eqref{eq.52}, i.e., responsible for $r \to 0$.
Moreover, we define 
\begin{equation}
\label{eq.n56}
 \s{S}^{\gamma}_P (X^\wedge) :=
   \{ \omega u + (1-\omega)v :
      u \in \s{H}^{\infty,\gamma}_P (X^\wedge), 
       v \in \s{S} (\ol{\R}_+, C^\infty (X))  \}.
\end{equation}
The spaces $\s{H}^{s,\gamma}_P (X^\wedge)$, $\s{K}^{s,\gamma}_P (X^\wedge)$, etc., are Fr\'echet in a natural way.

%%%%%%%%%%%%%%%%%%%%%%%%%%%%%%%%%%%%%%%%%%%%%%%%%%%%%%%%%%%%%%%%%%%%%%

\subsection{Naive and edge definitions of Sobolev spaces}
\label{s.10.1.3}

Sobolev spaces certainly belong to the prominent institutions in
the field of partial differential equations. The present modest
remarks do not reveal anything new as far as they concern the
\textup{`classical'} context. In fact, we content ourselves with
spaces based on $L^2$ norms and Fourier transforms. However,
things suddenly become much more uncertain if we ask the nature of
analogous spaces on manifolds with geometric singularities (cf.
also the considerations in Section 6.2 below). As is known the
\textup{`standard'} role of Sobolev spaces in elliptic PDE is to
encode the elliptic regularity of solutions. For instance, if $B$
is the unit ball in $\R^3$, solutions $u$ to the Dirichlet problem
$\Delta u = f \in H^{s-2}(B)$, $u \big|_{\partial B} = g \in H^{s-
\frac{1}{2}}(\partial B)$ for $s > \frac{3}{2}$ belong to
$H^s(B)$.

Now let $S \subset B$ be the hypersurface $S = \{ (x_1, x_2, x_3)
\in \R^3 : x_3 = 0$, $|x_1| + |x_2| \leq \frac{1}{2} \}$. What
can we say about the \textup{`Sobolev'} regularity of solutions of
$\Delta u = f$ in $B \setminus S$ with $u |_{\partial B} \in
H^{s- \frac{1}{2}}(\partial B)$, $u |_{\Int S_+} = g_+$,
$\displaystyle \frac{\partial u}{\partial x_3} \big|_{\Int S_-} = g_-$ 
(with
$\big|_{\Int S_\pm}$ denoting the limits at $\Int S$ from $x_3 > 0$
and $x_3 < 0$, respectively)? The question includes the choice of
\textup{`natural'} spaces for the boundary values on $\Int S_\pm$
as well as of the right notion of ellipticity in this case.
The critical zone is, of course, a neighbourhood of $\partial S$.
Problems of that kind occur, for instance, in crack theory.

Another question is the regularity of solutions to the
Zaremba problem
\[  \Delta u = f \; \text{in}\quad B, \quad
      u \big|_{S_+^2} = g_+, \ \frac{\partial u}{\partial \nu}
        \Big|_{S_-^2} = g_-,    \]
$S_\pm^2 := \partial B \cap \{ x_3 \gtrless 0 \}$, where
$\displaystyle \frac{\partial}{\partial \nu}$ denotes the derivative in direction
of the inner normal to the sphere. 
Also here the right notion of
ellipticity and the choice of analogues of Sobolev spaces (with
respect to their behaviour near	the interface $S^2 \cap 
\{ x_3 = 0 \}$) is far from being evident.

Problems with  singularities are meaningful also in
the pseudo-differential context.
Parametrices of elliptic problems for differential operators are
pseudo-differential, and questions then do not only concern the
spaces but also (hopefully manageable) quantisations, cf.
also Section 2.2 below.

For instance, we can ask the nature of solvability of the equation
\begin{equation}
\label{1013.Omega0805.eq}
\r^+ A \e^+ u = f
\end{equation}
in a (say, bounded) domain $\Omega \subset \R^n$, when $A =
\Op(a)$ is an elliptic pseudo-differential operator in $\R^n$ with
homogeneous principal symbol 
\[  \sigma_\psi(A)(\xi) = | \xi |^\mu     \]
for some $\mu \in \R$. In \eqref{1013.Omega0805.eq} by $\e^+$ we
mean the extension of distributions on $\Omega$ to zero outside
$\Omega$, and by $\r^+$ the operator of restriction to $\Omega$.
Even if $\partial \Omega$ is smooth and $\mu \not\in 2 \Z$ the
answer is not trivial. For $\mu \in 2 \Z$ we are in the frame of
operators with the transmission property at the boundary, cf. Section 2.1 below.

Another category of problems is the solvability of equations
$A u = f$,  
say, in $\R_+ \times \R^m$, when $A$ is a polynomial in vector
fields of some specific behaviour. In Section 10.1.1 we already
saw examples, such as vector fields of the form
\begin{equation}
\label{1013.F0805.eq}
r \partial_r, \partial_{x_1}, \ldots, \partial_{x_n}
\end{equation} 
for $m = n$ when $(r, x)$ are the coordinates in 
$\R_+ \times \R^n$
or, for the case $m = n + q$ with the 
coordinates $(r,x,y)$ in $\R_+ \times \R^n \times \R^q$
\begin{equation}
\label{1013.K0805.eq}
r \partial_r, \;\; \partial_{x_1}, \ldots, \partial_{x_n}, \;\;
    r \partial_{y_1}, \ldots, r \partial_{y_q}.
\end{equation}

Polynomials in vector fields \eqref{1013.F0805.eq} and
\eqref{1013.K0805.eq} just produce Fuchs type and edge-degenerate
operators, respectively (without weight factors in front of the
operators which we found natural in Section 10.1.1).

In the case \eqref{1013.F0805.eq} for $s \in \N, \gamma \in \R$,
we  can form the spaces
\begin{align}
\label{1013.FX0805.eq}
{\cal H}^{s, \gamma}(\R_+ \times \R^n) := & \; \{ u(r,x) \in 
     r^{-  \frac{n}{2} + \gamma} L^2(\R_+ \times \R^n) : 
      (r \partial_r)^k D_x^\alpha u(r,x) \in \\
      &r^{- \frac{n}{2} + \gamma} 
       L^2(\R_+ \times \R^n)\quad \text{for all} \quad
         k \in \N, \ \alpha \in \N^n, \ k + |\alpha| \leq s \}.    
                            \nonumber
\end{align}

In the case \eqref{1013.K0805.eq} for $s \in \N, \gamma \in \R$ we
might take	
\begin{align}
\label{1013.KY0805.eq}
{\cal H}^{s, \gamma}(\R_+ \times \R^{n+q}) := &\;
    \{ u(r,x,y) \in r^{- \frac{n}{2} + \gamma} L^2(\R_+ \times
          \R^{n+q}) :		    \\
     &\;(r \partial_r)^k D_x^\alpha (r D_y)^\beta u (r,x,y) \in
      r^{- \frac{n}{2} + \gamma} L^2(\R_+ \times \R^{n+q})
                                \nonumber\\
     & \; \text{for all}\quad k \in \N, \alpha \in \N^n, \beta \in \N^q,
       k+|\alpha|+|\beta| \leq s \}.
         \nonumber
\end{align}
Corresponding spaces for arbitrary real $s$ can be obtained by duality and
interpolation.
The spaces \eqref{1013.FX0805.eq} have natural invariance
properties and can be defined also on an open stretched cone
$\R_+ \times X =: X^\land$ for a (say, closed compact) $C^\infty$
manifold $X$. The resulting spaces are denoted by 
${\cal H}^{s,\gamma}(X^\land)$, cf. the formula
\eqref{1011.calH2605.eq}. 

Also the spaces \eqref{1013.KY0805.eq} have analogues in the
manifold case, namely, on $\Int \W$, where $\W$ is a (say, compact)
$C^\infty$ manifold with boundary $\partial \W$, such that
$\partial \W$ is an $X$-bundle over another $C^\infty$ manifold
$Y$. The corresponding spaces are denoted by 
${\cal H}^{s,\gamma}(\Int \W)$. In particular, for $\W = \ol{\R}_+
\times X \times \R^q$, we have the spaces
\begin{equation}
\label{1013.calH2005.eq}
{\cal H}^{s, \gamma}(X^\land \times \R^q).
\end{equation}
Note that ${\cal H}^{s, \gamma}(X^\land \times \R^q) = r^\gamma {\cal
H}^{s, 0}(X^\land \times \R^q)$ for all $s, \gamma \in \R$.
The spaces ${\cal H}^{s, \gamma}(X^\land)$ or their analogues
${\cal H}^{s, \gamma}(\Int \B)$ on a (compact) stretched manifold
$\B$ with conical singularities (that is, a compact $C^\infty$
manifold with boundary $\partial \B \cong X$) are common in the
investigation of elliptic operators on a manifold with conical
singularities. Also the spaces ${\cal H}^{s, \gamma}(\Int \W)$ for
$q > 0$ are  taken in many investigations in the literature 
when the operators are generated by the vector fields
\eqref{1013.K0805.eq}.	   

However, for nearly all purposes that we have in mind here, for
instance, the problems mentioned at the beginning of this
section, or also for geometric (edge-degenerate) operators with the
typical weight factor, we find the above mentioned definition of 
${\cal H}^{s, \gamma} (\Int \W)$-spaces for $ \dim Y > 0$ not really convenient (which says nothing on whether the spaces themselves are adequate).
That is why we talk in this connection about a \textup{`naive'} definition of
Sobolev spaces, in contrast to other ones which are more efficient for establishing calculi on manifolds with (regular)
geometric singularities. Also to express \textup{`canonical'}
singular functions of edge asymptotics another choice of Sobolev
spaces seems to be indispensable.

\textup{`Naive'} and \textup{`non-naive'} definitions of corner Sobolev spaces
are also possible for more than one axial variable. Higher
generations of Sobolev spaces in that sense will be considered in
Section 10.6.2 below.

In order to give some motivation for an alternative choice of Sobolev
scales we look at what happens when we
formulate a boundary value problem for an elliptic differential
operator with smooth coefficients in $\R^n$ in a smooth domain
$\Omega \subset \R^n$ with boundary.

A priori a Sobolev space, given in $\R^n$, has no relation to
possible boundaries of a domain. A boundary contributes some
anisotropy into the consideration, and tangential and normal
directions play a different role. More generally,
smooth (or non smooth) hypersurfaces of arbitrary codimension should
interact with isotropic Sobolev distributions in a specific manner.
We want to discuss this point in terms of a representation of the
Euclidean space as a \textup{`wedge'} $\R^n = \R^m \times \R^q \ni
(z,y)$ with edge $\R^q$ and model cone $\R^m$ (with the origin as a
fictitious conical singularity).

Recall that the standard $L^2$-spaces have the property
\begin{equation}
\label{1013.L21305.eq}
L^2(\R^m \times \R^q) = L^2(\R^q, L^2(\R^m)).
\end{equation}
More generally, we might try to employ Sobolev spaces taking values
in another Sobolev space.

Let $E$ be a Hilbert space, and let
$H^s(\R^q, E)$
denote the completion of ${\cal S}(\R^q, E)$ (the Schwartz space of
functions with value in $E$) with respect to the norm
 $\| u \|_{H^s(\R^q, E)} = \Bigl\{ \int \langle \eta \rangle^{2s}
    \| \hat{u}(\eta) \|_E^2 d \eta \Bigr\}^{\frac{1}{2}}$,    
$s \in \R$, $\langle \eta \rangle := (1+ | \eta |^2)^{1/2}$, with 
$\hat{u}(\eta) = (F_{y \to \eta} u)(\eta)$ being the Fourier transform
in $\R^q \ni y$.

Clearly we have
$$
  H^s(\R^m \times \R^q) \not= H^s(\R^q, H^s(\R^m))     
$$
unless $s = 0$. The question is how to find the \textup{`right'}
anisotropic reformulation of $H^s(\R^n)$.

The answer comes from the notion of \textup{`abstract'} edge Sobolev
spaces.

\begin{Definition}
\label{1013.1305de.de}
Let $E$ be  a Hilbert space, equipped with a strongly
continuous group of isomorphisms   
$\kappa_\lambda : E \to E, \quad \lambda \in \R_+$,
such that $\kappa_\lambda \kappa_\delta = \kappa_{\lambda \delta}$ for
all $\lambda, \delta \in \R_+$ {\em(}strongly
continuous means $\{ \kappa_\lambda e \}_{\lambda \in \R_+} \in
C(\R_+, E)$ for every $e \in E${\em)}.
In that case we will speak about a group action on $E$.
The abstract edge Sobolev space ${\cal W}^s(\R^q, E)$ of 
smoothness $s \in \R$, modelled on
a Hilbert space $E$ with group action $\{ \kappa_\lambda \}_{\lambda
\in \R_+}$, is defined to be the completion of ${\cal S}(\R^q, E)$
with respect to the norm
\begin{equation}
\label{1013.w1305.eq}
\| u \|_{{\cal W}^s(\R^q, E)} = \Bigl\{ \int \langle \eta \rangle^{2s}
     \| \kappa_{\langle \eta \rangle}^{-1} \hat{u}(\eta)
     \|_E^2 d \eta \Bigr\}^{\frac{1}{2}}.
\end{equation}
If $E$ is a Fr\'echet space with group action $\{ \kappa_\lambda \}_{\lambda \in \R_+}$, i.e., 
$E =\varprojlim_{j\in \N} E^j$ for a chain of Hilbert spaces with continuous embeddings 
$\ldots  \hookrightarrow E^{j+1} \hookrightarrow E^{j} \hookrightarrow \ldots  \hookrightarrow E^0$, where
$\{ \kappa_\lambda \}_{\lambda \in \R_+}$ on $E^0$ restricts to a group action on every $E^j$, $j \in \N$, then we set $\s{W}^s (\R^q,E) := \varprojlim_{j\in \N} \s{W}^s (\R^q, E^j)$.
\end{Definition}
For $E = \C$ with the trivial group action we recover the 
scalar Sobolev spaces $H^s(\R^q)$, i.e.,
$H^s(\R^q) = {\cal W}^s(\R^q, \C)$.
\begin{Remark}
\label{1013.re1305.re}
For $E = H^s(\R^m)$ and $\kappa_\lambda u(x) := \lambda^{\frac{m}{2}}
u(\lambda x)$, $\lambda \in \R_+$, we have a canonical isomorphism
\begin{equation}
\label{1013.Hs1305.eq}
H^s(\R^m \times \R^q) = {\cal W}^s(\R^q, H^s(\R^m))
\end{equation}
for every $s \in \R$. The group $\{ \kappa_\lambda \}_{\lambda \in
\R_+}$ is unitary on $L^2(\R^m) = H^s(\R^m)$; thus
\eqref{1013.Hs1305.eq} is compatible with the relation
\eqref{1013.L21305.eq}. 
\end{Remark}      

\begin{Remark}
\label{1013.re2005.re}
Writing
$\| u \|_{{\cal W}^s(\R^q, E)} = \Bigl\{ \int \langle \eta 
     \rangle^{2s}
    \| F\bigl(F^{-1} \kappa_{\langle \eta \rangle}^{-1} F \bigr) u
    \|_E^2 d \eta \Bigr\}^{\frac{1}{2}}$,
$F = F_{y \to \eta}$, we see that the
operator $T := F^{-1} \kappa_{\langle \eta \rangle}^{-1} F$ induces
an isomorphism
\[  F^{-1} \kappa_{\langle \eta \rangle}^{-1} F :
     {\cal W}^s(\R^q, E) \to H^s(\R^q, E)    \]
for every $s \in \R$.   
Given a closed subspace $L \subseteq E$, not necessarily invariant under the group action 
$\{ \kappa_\lambda \}_{\lambda \in \R_+}$, we can form 
$H^s (\R^q, L)$ and then
\begin{equation}
\label{eq.V}
\s{V}^s (\R^q, L) :=
    T^{-1} H^s (\R^q, L).
\end{equation}
In the case that $\{ \kappa_\lambda \}_{\lambda \in \R_+}$ induces a group action in $L$ {\em(}by restriction{\em)} we have, of course,
$\s{V}^s (\R^q, L) = \s{W}^s (\R^q, L)$; in any case $\s{V}^s (\R^q, L)$ is a closed subspace of $\s{W}^s (\R^q, E)$.
If we have a direct decomposition  $E = L \oplus M$ into closed subspaces, we get a direct decomposition
\begin{equation}
\label{eq.()}
\s{W}^s (\R^q, E) =
  \s{V}^s (\R^q, L) \oplus \s{V}^s (\R^q, M).
\end{equation}
\end{Remark}

Recall that the space ${\cal K}^{s, \gamma}(\R^m \setminus \{ 0 \})$
has the property 
\[  (1- \omega) {\cal K}^{s, \gamma}(\R^m \setminus 
       \{ 0 \}) = (1- \omega)H^s(\R^m)   \]
for every $\omega \in C_0^\infty(\R^m)$ which is equal to $1$ in a
neighbourhood of $0$. 
For $E := {\cal K}^{s, \gamma}(\R^m \setminus  \{0 \})$ with the group action $\kappa_\lambda u(z) =
        \lambda^{\frac{m}{2}} u(\lambda z), \ \lambda \in \R_+$, we can form
the corresponding edge Sobolev space and observe that
\[  (1- \omega){\cal W}^s(\R^q, {\cal K}^{s, \gamma}(\R^m \setminus \{ 0
      \})) = (1- \omega) H^s(\R^{m+q})     \]
for any such $\omega$. This implies
\begin{equation}
\label{eq.65.}
 H_{\comp}^s((\R^m \setminus \{0 \}) \times \R^q) \subset
     {\cal W}^s(\R^q, {\cal K}^{s, \gamma}(\R^m \setminus \{ 0 \}))
     \subset H_{\loc}^s((\R^m \setminus \{ 0 \}) \times \R^q).   
\end{equation}
 
\begin{Remark}
\label{r.g}
We can also form the spaces
\begin{equation}
\label{eq.w}
\s{W}^s (\R^q, \s{K}^{s,\gamma;g} (X^\wedge))
\end{equation}
based on the group action  \eqref{eq.new11}.
Those satisfy an analogue of the relation \eqref{eq.65.} for all $s,\gamma, g \in \R$.
For $g := s-\gamma$ the spaces have particularly natural properties.
\end{Remark}        

The relation \eqref{1013.Hs1305.eq} shows that classical Sobolev spaces are special examples of edge spaces in the sense of Definition \ref{1013.1305de.de}, where a hypersurface $\R^q$ of arbitrary codimension $\geq 1$ can be interpreted as an edge.
Such an anisotropic description of `isotropic' Sobolev spaces also makes sense with respect to any other (smooth) hypersurface of a $C^\infty$ manifold, cf. the articles \cite{Dine1}, \cite{Dine3}.
The anisotropic interpretation is particularly reasonable on a $C^\infty$ manifold with boundary; in this case the boundary is locally identified with $\R^q$, and $\R_+$ (the inner normal with respect to a chosen Riemannian metric in product form near the boundary) is the substitute of $\R^n$.
This gives us the possibility to encode various properties of regularity of distributions up to the boundary, not only $C^\infty$ in terms of 
$\s{W}^\infty (\R^q, \s{S} (\ol{\R}_+))$ but  other asymptotics, e.g.,
$\s{W}^\infty (\R^q, \s{S}^\gamma_P ({\R}_+))$ for an asymptotic type 
$P= \{ (p_j, m_j)\}_{j\in\N}$ as in Section 1.2, with 
$\pi_\C R \subset \{ \re w < \frac{1}{2} - \gamma \}$.
More generally, asymptotics of type $P= \{ (p_j, m_j, L_j) \}_{j\in \N}$, cf. the first sequence of \eqref{eq.52}, with
$\pi_\C P \subset \{ \re w < \frac{n+1}{2} - \gamma \}$ for $n=\dim X$, gives rise to edge asymptotics of distributions $u(r,x,y)$ on the stretched wedge $X^\wedge \times \R^q \ni (r,x,y)$, modelled on 
$$
\s{W}^s (\R^q, \s{K}^{s,\gamma}_P (X^\wedge)),
$$
when (which can be done) $\s{K}^{s,\gamma}_P (X^\wedge)$ is written as a projective limit of Hilbert subspaces of $\s{K}^{s,\gamma} (X^\wedge)$, endowed with the group action 
$(\kappa_\lambda u) (r,x) = \lambda^{\frac{n+1}{2}} u(\lambda r, x)$, $\lambda \in \R_+$.

Sobolev spaces $H^s(\R^q)$ described in terms of the Fourier
transform are perfectly adapted to pseudo-differential symbols in
H\"ormander's classes $S_{(\cl)}^\mu(U \times \R^q) \ni a(y, \eta)$,
cf. Definition \ref{1013.Sy2005.de} below. In particular, $\langle
\eta \rangle^s = (1+ | \eta |^2)^{s/2}$ is a classical symbol of
order $s$, and we have
\[  \| u \|_{H^s(\R^q)} = \| \langle \eta \rangle^s \hat{u}(\eta) 
    \|_{L^2(\R^q)}.     \]
This relation can be seen as a continuity result for the
pseudo-differential operator $A = \Op(\langle \eta \rangle^s)$,
\[  \Op(a)u(y) = \iint e^{i(y-y')\eta} a(y,y', \eta) u(y') dy' \dbar
        \eta,      \]
$\dbar \eta = (2 \pi)^{-q}d \eta$. A simplest version tells us that
for $a(y,y', \eta) \in S^\mu(\R_{y,y'}^{2q} \times \R_\eta^q)$ under
suitable conditions on the dependence on $(y,y')$ for $|(y,y')| \to
\infty$, we have continuity of the associated operator
\begin{equation}
\label{1011.cont2005.eq}
\Op(a) : H^s(\R^q) \to H^{s- \mu}(\R^q)
\end{equation}
for all $s \in \R$. (The conditions are satisfied, for instance,
when $a(y,y', \eta)$ is independent of $(y,y')$ for large
$|(y,y')|$, and, of course, in much more general cases.)

Analogously, the abstract edge Sobolev spaces of Definition
\ref{1013.1305de.de} have a counterpart in terms of 
operator-valued symbols.

\begin{Definition}
\label{1013.Sy2005.de}
\begin{enumerate}
\item The space $S^\mu(U \times \R^q; E, \widetilde{E})$ for open $U
      \subseteq \R^p$ and Hilbert spaces $E$ and $\widetilde{E}$, 
      endowed
      with group actions $\{ \kappa_\lambda \}_{\lambda \in \R_+}$ 
      and $\{\tilde{\kappa}_\lambda \}_{\lambda \in \R_+}$, 
      respectively, is defined to be the set of all
      $a(y,\eta) \in C^\infty(U \times \R^q, {\cal L}(E, 
                 \widetilde{E}))$
      such that
    \[ \sup_{\substack{y \in K \\ \eta \in \R^q}} \langle \eta 
       \rangle^{-
      \mu + | \beta |} \| \tilde{\kappa}_{\langle \mu \rangle}^{-1}
      \{ D_y^\alpha D_\eta^\beta a(y, \eta) \}
      \kappa_{\langle \eta \rangle} \|_{{\cal L}(E, \widetilde{E})}
      < \infty      \]
      for all $\alpha \in \N^p$, $\beta \in \N^q$, and $K \Subset U$; 
      $\langle \eta \rangle := (1+|\eta|^2)^{1/2}$;       		 	     
\item the subspace of classical symbols $S_{\cl}^\mu(U \times \R^q;
      E, \widetilde{E}) \ni a(y, \eta)$ is defined by the condition
      that there are elements $a_{(\mu -j)}(y, \eta)
      \in C^\infty(U\times (\R^q \setminus \{ 0 \}), {\cal L}(E,
      \widetilde{E}))$, $j \in \N$, satisfying the homogeneity 
      condition
      \begin{equation}
      \label{1013.tw2105.eq}  
      a_{(\mu -j)}(y, \lambda \eta) = \lambda^\mu
          \tilde{\kappa}_\lambda a_{(\mu -j)}(y, \eta)
          \kappa_\lambda^{-1}        
      \end{equation}
      for all  $\lambda \in \R_+$, such that
      \begin{equation}
      \label{1013.ho2005.eq}
      a(y, \eta) - \chi(\eta) \sum_{j=0}^{N} a_{(\mu -j)}(y, \eta)
         \in S^{\mu -(N+1)}(U \times\R^q; E, \widetilde{E})
      \end{equation}
      for all $N \in \N$ {\st{(}}here $\chi(\eta)$ is any excision
      function in $\eta${\st{)}}. The relation 
      \eqref{1013.tw2105.eq}
      is also referred to as twisted homogeneity {\st{(}}of order
      $\mu - j${\st{)}}.	 
\end{enumerate}
\end{Definition}

For $E = \widetilde{E} = \C$ and $\kappa_\lambda =
\tilde{\kappa}_\lambda = \id$ for all $\lambda \in \R_+$ we obtain
the scalar symbol spaces $S_{(\cl)}^\mu(U \times \R^q)$
(subscript `$(\cl)$' means that we are speaking about the
classical or the general case). By $S_{(\cl)}^\mu(\R^q; E,
\widetilde{E})$ we denote the subspaces of $S_{(\cl)}^\mu(U \times
\R^q; E, \widetilde{E})$ of $y$-independent elements (i.e., symbols
with constant coefficients).   

\begin{Remark}
\label{1013.hom2505.re}
Let $a(y, \eta) \in C^\infty(U \times \R^q, {\cal L}(E,
\widetilde{E}))$ satisfy the relation
\[  a(y, \lambda \eta) = \lambda^\mu \tilde{\kappa}_\lambda a(y,
      \eta) \kappa_\lambda^{- 1}     \] 
for all $\lambda \geq 1$, $(y, \eta) \in U \times \R^q$, $|\eta |
\geq C$ for a constant $C > 0$. Then we have $a(y, \eta) \in
S_{\cl}^\mu(U \times \R^q; E, \widetilde{E})$.     
\end{Remark}      

\begin{Remark}
\label{1013.re0206.re}
Definition {\st{\ref{1013.Sy2005.de}}} has a straightforward generalisation
to the case of Fr\'echet spaces $E$, {\st{(}}and/or{\st{)}}
$\widetilde{E}$, equipped with group actions in the
sense that the spaces are projective limits of Hilbert spaces with
corresponding group actions, cf. the corresponding notation in the second part of Definition 
{\em\ref{1013.1305de.de}}.
\end{Remark}

There are many beautiful and unexpected examples of
operator-valued symbols.

\begin{Example}
\label{1013.G0206.exa}
An important category of  examples are the Green, potential and trace
symbols in the calculus of boundary value problems with the
transmission property.

Consider functions
\begin{align*}
f_G(t,t'; y, \eta) 
& 
    \in {\cal S}(\ol{\R}_+ \times \ol{\R}_+, 
                     S_{\cl}^{\mu+1}(\Omega \times \R^q)),   \\
f_K(t;y,\eta) 
& 
    \in {\cal S}(\ol{\R}_+, S_{\cl}^{\mu + \frac{1}{2}}
                 (\Omega \times \R^q)), \\   
f_B(t'; y, \eta) 
& 
   \in {\cal S}(\ol{\R}_+, S_{\cl}^{\mu +
                   \frac{1}{2}}(\Omega \times \R^q)),
\end{align*}
and form the operator families		   		 		     
\begin{eqnarray*}
g(y,\eta)u(t) 
& := &
\int_{0}^{\infty} f_G(t[\eta],t'[\eta];y, \eta)  u(t')dt', 
\ \ \ \ \ \ u \in L^2(\R_+),   \\
k(y, \eta)c
&  := &
f_K(t[\eta];y, \eta)c,  
\ \ \ \ \ \  c \in \C,     \\
b(y, \eta)u 
& :=&
 \int_{0}^{\infty} f_B(t'[\eta]; y, \eta)u(t')dt',
\ \ \ \ \ \ u \in L^2(\R_+).
\end{eqnarray*}

Here $\eta \to [\eta]$ is any $C^\infty$ function in $\eta \in \R^q$, strictly positive, and $[\eta] = |\eta|$
for $|\eta| > C$ for some $C >0$.
We then obtain operator-valued symbols
\begin{align*}
g(y, \eta) & \in S_{\cl}^\mu(\Omega \times \R^q; L^2(\R_+), {\cal
                S}(\ol{\R}_+)), \\
k(y, \eta) & \in S_{\cl}^\mu(\Omega \times \R^q; \C, {\cal
                S}(\ol{\R}_+)),   \\
b(y, \eta) & \in S_{\cl}^\mu(\Omega \times \R^q; L^2(\R_+), \C),
\end{align*}
called Green, potential, and trace symbols, respectively, of order
$\mu \in \R$ and type $0$. Green and trace symbols of type $\d \in \N$
are defined as linear combinations
\[  g(y, \eta) = \sum_{j=0}^{\d} g_j(y, \eta)
      \frac{\partial^j}{\partial t^j}\quad \textup{and}\quad
      b(y, \eta) = \sum_{j=0}^{\d} b_j(y, \eta)
      \frac{\partial^j}{\partial t^j}     \]
with $g_j$ and $b_j$ being of order $\mu -j$ and type $0$ {\st{(}}with argument 
functions 
belonging to $H^s(\R_+)$ for $s > \d - \frac{1}{2}${\st{)}}.      
\end{Example}

The associated pseudo-differential operators
$\Op(g)$, $\Op(k)$ and $\Op(b)$
are called Green, potential, and trace operators (of the respective
types in the Green and trace case).

\begin{Example}
\label{1013.ed2105.exa}
Let $X$ be a closed compact $C^\infty$ manifold, $\Omega
\subseteq \R^q$ an open set, and consider an operator $A$ of the
form \eqref{1011.A1905.eq}  on $X^\land \times \Omega$, $X^\land
= \R_+ \times X$, with coefficients
\begin{equation}
\label{1013.co2105.eq}
a_{j \alpha}(r, y) \in C^\infty(\ol{\R}_+ \times \Omega, 
     \Diff^{\mu-(j+ | \alpha |)}(X))
\end{equation}
that are independent of $r$ for $r > R$ for some $R > 0$. Then we
have
\begin{equation}
\label{1013.a2105.eq}
a(y, \eta) := r^{- \mu} \sum_{j+| \alpha | \leq \mu} a_{j
     \alpha}(r,y) \Bigl(- r \frac{\partial}{\partial r}\Bigr)^j
     (r \eta)^\alpha
     \in S^\mu(\Omega \times \R^q; {\cal K}^{s, \gamma}(X^\land),
     {\cal K}^{s- \mu, \gamma - \mu}(X^\land))
\end{equation}         
for every $s, \gamma \in \R$ {\st{(}}the group action is as in
Remark {\st{\ref{1011.1305re.re}}}{\st{)}}. If the coefficients
\eqref{1013.co2105.eq} are independent of $r$, then $a(y, \eta)$
is classical, and we have $a_{(\mu)}(y, \eta) =
\sigma_\land(A)(y, \eta)$, cf. the expression
\eqref{1011.sigmaland2304.eq}.  
\end{Example}

In analogy to \eqref{1011.cont2005.eq} every $a(y,y', \eta) \in
S^\mu(\R_{y,y'}^{2q} \times \R_\eta^q; E, \widetilde{E})$ (again
under suitable conditions on the dependence on $(y,y')$ for
$|(y,y')| \to \infty$) induces continuous operators
\begin{equation}
\label{1013.contW2105.eq}
\Op(a) : {\cal W}^s(\R^q,E) \to {\cal W}^{s- \mu}(\R^q,
          \widetilde{E})
\end{equation}
for all $s \in \R$.	  

Applying that to the symbol \eqref{1013.a2105.eq} (for $\Omega =
\R^q$, under a corresponding condition on the coefficients
for $|y| \to \infty$, say, to be constant with respect to $y$
for large $|y|$) we see that the associated edge-degenerate
operator $A = \Op(a)$ (cf. the formula \eqref{1011.A1905.eq})
induces a continuous operator
\[  A : {\cal W}^s(\R^q, {\cal K}^{s, \gamma}(X^\land)) \to
           {\cal W}^{s- \mu}(\R^q, {\cal K}^{s- \mu, \gamma -
	   \mu}(X^\land))     \]
for every $s, \gamma \in \R$. Recall that we also have an
alternative scale of spaces, namely, ${\cal H}^{s, \gamma}(X^\land
\times \R^q)$, cf. the formula \eqref{1013.calH2005.eq}. 
Since the operator $A$ has the form $r^{- \mu} \widetilde{A}$,
where $\widetilde{A}$ is (locally with respect to $X$) a
polynomial of order $\mu$ in the  vector fields
\eqref{1013.K0805.eq}, the operator $A$ is also continuous in the sense
\[  A : {\cal H}^{s, \gamma}(X^\land \times \R^q) \to
         {\cal H}^{s- \mu, \gamma - \mu}(X^\land \times	\R^q) \]
for every $s, \gamma \in \R$. 
The question is now which is the more natural definition of spaces in connection with edge-degenerate
operators,
\begin{equation}
\label{1013.H2105.eq}
{\cal W}^{s, \gamma}(X^\land \times \R^q)\quad \text{or}\quad
     {\cal H}^{s, \gamma}(X^\land \times \R^q)?
\end{equation}

There are, of course, many other questions, for instance,
\begin{equation}
\label{1013.N2105.eq}
\textup{`what is natural'?}
\end{equation}
The question \eqref{1013.H2105.eq} seems
to have a natural answer in favour of the spaces ${\cal H}^{s,
\gamma}(X^\land \times \R^q)$, because, up to the weight factor,
the operators are polynomials in the typical vector fields.
Authors who employ this definition of weighted spaces in connection with 
configurations with edge singularities probably share this opinion. 

This is now an excellent opportunity to found a new
sect who believes the other truth. 

A wise outcome of the discussion would be that both parties
have their own right to exist; the various approaches might be
(to some extent) equivalent anyway, or point out different aspects
of the same phenomena . 
The trivial solution would  be that the different spaces are the same  at all (at least locally near $r = 0$). 
The latter, however, is not the case when
$\gamma \not= s$.

An answer is that the spaces 
$\s{W}^s (\R^q, \s{K}^{s,\gamma} (X^\wedge))$ belong to a continuously parametrised family of spaces
$\s{W}^s (\R^q, \s{K}^{s,\gamma; g} (X^\wedge))$ for $g \in \R$, cf. Remark \ref{r.g}.
All these spaces are possible choices for a consistent calculus with the same edge algebra.
However, for $g=s-\gamma$ the spaces
$\s{W}^s (\R^q, \s{K}^{s,\gamma; g} (X^\wedge))$  and 
$\s{H}^{s,\gamma} (X^\wedge \times \R^q)$ agree close to $r=0$.
This coincidence is a hidden effect and an aspect of what we  call a non-naive (or edge-) definition of weighted spaces, see also \cite{Tark4} and \cite{Schu53}.

Let us now have a look at another category of operator-valued symbols in
the sense of Definition \ref{1013.Sy2005.de}, which play a role
in the description of the internal properties of standard
Sobolev spaces.      	     

Let us write $\R^n = \R^m \times \R^q$.
Recall the well known fact that the operator of restriction
$\r' : {\cal S}(\R^n) \to {\cal S}(\R^q)$,    
$\r' u := u \big|_{\{0\} \times \R^q}$ extends to a continuous
operator 
\begin{equation}
\label{1013.co2505.eq}  
\r' : H^s(\R^n) \to H^{s- \frac{m}{2}}(\R^q)    
\end{equation}
for all $s > \frac{m}{2}$. This can easily be interpreted as a
continuity result of the kind \eqref{1013.contW2105.eq} for some
special operator-valued symbol.

In fact, writing $\r_0' u := u(0)$ for $u \in {\cal S}(\R^n)$, 
there is an extension to a continuous operator
\[  \r_0' : H^s(\R^m) \to \C    \]
for every $s > \frac{m}{2}$. If we endow $H^s(\R^m)$ with the group
action $\kappa_\lambda u(t) = \lambda^{\frac{m}{2}} u(t)$ for
$\lambda \in \R_+$ and $\C$ with the trivial group action,       
         from
\begin{equation}
\label{1013.rh2505.eq}
\r_0' \in C^\infty(\R^q, {\cal L}(H^s(\R^m), \C))\quad
             \text{and}\quad
      \r_0' = \lambda^{\frac{m}{2}} \r_0' \kappa_\lambda^{-1}\quad
      \text{for all}\quad \lambda \in \R_+
\end{equation}   
it follows that
\[ \r_0' \in S_{\cl}^{\frac{m}{2}}(\R^q; H^s(\R^m), \C)     \]
for every $q \in \N$, cf. Remark \ref{1013.hom2505.re} (i.e., $\r_0'$
is a classical symbol in the covariable $\eta \in \R^q$,
although it is independent of $\eta$). 
Then \eqref{1013.co2505.eq} is a
consequence of the continuity of
\[  \r' = \Op(\r_0') : {\cal W}^s(\R^q, H^s(\R^m)) \to
          {\cal W}^{s-\frac{m}{2}}(\R^q, \C) = H^{s-\frac{m}{2}} (\R^q),      \]
cf. the relations \eqref{1013.tw2105.eq} and 
\eqref{1013.KY0805.eq}.	     	     
More generally, for every $\alpha \in \N^m$ we can form the composition
$$
\r'_0 D^\alpha_x : H^s (\R^m) \to \C
$$
which is continuous for $s -|\alpha| > \frac{m}{2}$ and defines a symbol
$$
\r'_0 D^\alpha_x \in 
   S^{\frac{m}{2} +|\alpha|} 
       (\R^q; H^s (\R^m), \C)
$$
which is even homogeneous in the sense
$\r'_0 D^\alpha_x =
     \lambda^{\frac{m}{2} +|\alpha|}
        \r'_0 D^\alpha_x \kappa^{-1}_\lambda$ for all
$\lambda \in \R_+$.
From \eqref{1013.contW2105.eq} we then obtain a corresponding continuity of the associated pseudo-differential operator.
A certain counterpart of such symbols are potential symbols of the form
$k(\eta) : \C \to \s{S}(\R^m)$, defined by  
$c \to [\eta]^{\nu +\frac{m}{2}}
       (x [\eta])^\alpha \omega (x[\eta]) c$
for any $\nu \in \R$, $\alpha \in \N^m$, and a function 
$\omega \in C^\infty_0 (\R^m)$ that is equal to $1$ in a neighbourhood of $x=0$.
They define an element
$$
k(\eta) \in
    S^\nu_\cl (\R^q; \C, H^s (\R^m))
$$
for every $s \in \R$ and satisfy the homogeneity relation
$k (\lambda \eta) = 
      \lambda^{\nu}  \kappa_\lambda k(\eta)$ for 
$\lambda \geq 1$, $|\eta| \geq c$ for some constant $c>0$.

%%%%%%%%%%%%%%%%%%%%%%%%%%%%%%%%%%%%%%%%%%%%%%%%%%%%%%%%%%%%%%%%%%%%%%
\section{Are regular boundaries harmless?}
\label{s.10.2}
\begin{minipage}{\textwidth}
\setlength{\baselineskip}{0cm}
\begin{scriptsize}    
Classical boundary value problems (such as the Dirichlet or the
Neumann problem for the Laplace operator in a smooth bounded domain) are
well understood from the point of view of elliptic regularity up to
the boundary, the Fredholm index in Sobolev spaces, the nature of
pseudo-differential parametrices, etc. Regular boundaries in this
context are harmless in the sense that non-regular boundaries require
much more specific insight (even for the simplest case of
conical singularities). Of course, also for problems with smooth
boundaries there are interesting aspects, worth to be considered up
to the present, for instance, in connection with the index of elliptic
operators who do not admit Shapiro-Lopatinskij elliptic
conditions, or around the spectral behaviour.

However, this is not the idea of the discussion here. We want to see
how pseudo-differential operators behave near a smooth boundary and
show some connections to the edge calculus.
\end{scriptsize}
\end{minipage}

\subsection{What is a boundary value problem?}
\label{s.10.2.1}

In an exposition on operators on manifolds with higher singularities we should ask `what is an edge problem' or `what is a higher corner problem';
however, this  will be answered anyway in Chapter 5 below.
The structures and inventions for the higher corner calculus should derive their motivation from something very common, namely, boundary value problems.
Boundary value problems have something to do with the values of a solution at the boundary, i.e., with boundary conditions.
That leads to one of the basic ingredients also for the analysis on a polyhedral configuration near lower-dimensional strata, namely, to additional conditions along those strata, with a specific contribution to the symbolic structure and associated operators, in general, of trace and potential type.

In this section we are interested in the behaviour of pseudo-differential operators with smooth symbols in a smooth domain in $\R^n$ (or on a $C^\infty$
manifold with boundary). 
Moreover, we ask the nature of solvability
near the boundary when the operator is elliptic. 
For convenience, we
first consider a smooth bounded domain $\Omega \subset \R^n$ and a
classical pseudo-differential operator $A$ in a neighbourhood of $X =
\ol{\Omega}$.
In the simplest case $A$ is a differential operator,
\begin{equation}
\label{1021.Amu1105.eq}
A = \sum_{|\alpha|\leq \mu} a_\alpha(x) D_x
\end{equation}
with coefficients $a_\alpha \in C^\infty(\R^n)$. If $\ol{\Omega}$ is
locally modelled on the half-space
\[  \ol{\R}_+^n = \{ x = (x_1, \ldots, x_n) \in \R^n : x_n > 0 \}  \]
we also write $x = (y,t)$ for $y = (x_1, \ldots, x_{n-1}), \ t :=
x_n$, with corresponding covariables $\xi = (\eta, \tau)$. Then
$\sigma_\psi(A)(x, \xi)$, the homogeneous principal symbol of $A$ of
order $\mu$, cf. the formula \eqref{1011.sigmapsi2104.eq}, generates
a parameter-dependent family of differential operators on $\R_+ \ni
t$, namely,
\begin{equation}
\label{1021.sigmpartial1105.eq}
\sigma_\partial(A)(y, \eta) := \sigma_\psi(A)(y,0, \eta, D_t),
\end{equation}
$(y, \eta) \in T^*\R^{n-1} \setminus 0$, regarded as a 
family of maps
\begin{equation}
\label{1021.sigmaH1105.eq}
\sigma_\partial(A)(y, \eta) : H^s(\R_+) \to H^{s-\mu}(\R_+)
\end{equation}
between the Sobolev spaces $H^s(\R_+) := H^s(\R) \big|_{\R_+}$, $s \in
\R$. We call \eqref{1021.sigmpartial1105.eq} the homogeneous
principal boundary symbol of the operator $A$. This is another
example of a symbolic structure of operators in $\R^n$, not
explicitly mentioned in Section 1.1. 
Note that we have
\begin{equation}
\label{1021.comp2805.eq}  
\sigma_\partial(I) = \id, \quad \text{and}\quad
      \sigma_\partial(AB) = \sigma_\partial(A) 
      \sigma_\partial(B)  
\end{equation}
for differential operators $A$ and $B$ of order $\mu$ and $\nu$, 
respectively.

\begin{Remark}
\label{1021.hom1105.re}
The homogeneity of the principal boundary symbol refers to a strongly
continuous group $\{ \kappa_\lambda \}_{\lambda \in \R_+}$ of
isomorphisms on the $H^s(\R_+)$-spaces, defined by
\[  \kappa_\lambda u(t) := \lambda^{\frac{1}{2}} u(\lambda t) \quad
       \text{for} \quad \lambda \in \R_+.     \]
For an operator \eqref{1021.Amu1105.eq} we have
$\sigma_\partial(A)(y, \lambda \eta) = \lambda^\mu \kappa_\lambda
      \sigma_\land(A)(y, \eta) \kappa_\lambda^{-1}$
for all $\lambda \in \R_+$, $(y, \eta) \in T^* \R^{n-1} \setminus 0$
{\st{(}}cf., similarly, the formula \eqref{1011.kappa1105.eq}
{\st{)}}.      
\end{Remark}      

If $\Omega \subseteq \R^n$ is a smooth bounded domain and
\eqref{1021.Amu1105.eq} a differential operator in $\R^n$, we can
replace $\Omega$ locally near $\partial \Omega$ by the half-space
$\R_+^n \ni (y,t)$ and calculate the boundary symbol
$\sigma_\partial(A)(y, \eta)$. This is then invariantly defined as a
family of operators \eqref{1021.sigmaH1105.eq} for $(y, \eta) \in
T^*(\partial \Omega) \setminus 0$.

If we recognise the boundary symbol $\sigma_\partial(A)$ of an
operator $A$ as another principal symbolic level, i.e., interpret
the pair
\begin{equation}
\label{1013.sigmaneu2505.eq}  
\sigma(A) = (\sigma_\psi(A), \sigma_\partial(A))     
\end{equation}
as the \textup{`full'} principal symbol of $A$, then ellipticity
should be defined as the invertibility of both components. However,
the second component is not necessarily bijective, as we see by the
following theorem.
 
\begin{Theorem}
\label{1021.th1105.th}
Let \eqref{1021.Amu1105.eq} be elliptic with respect to
$\sigma_\psi$. Then \eqref{1021.sigmaH1105.eq} is a surjective family
of Fredholm operators for every $s > \mu - \frac{1}{2}, 
(y, \eta) \in T^*(\partial \Omega) \setminus 0$.  
\end{Theorem}

\begin{Remark}
\label{1021.re1105.re}
By virtue of Remark {\st{\ref{1021.hom1105.re}}} we have
\[  \dim \ker \sigma_\partial(A)(y, \eta) = 
    \dim \ker \sigma_\partial(A)(y, \frac{\eta}{|\eta|}).   \]
\end{Remark}

Simplest examples show what happens when we look at
$\sigma_\partial(A)$ for a $\sigma_\psi$-elliptic operator $A$: 
Let $A = \Delta$ be the Laplacian with its principal symbol 
$\sigma_\psi(\Delta) = - | \xi |^2$. Then
\begin{equation}
\label{1021.neu0206.eq}  
\sigma_\partial(\Delta) (\eta) = - | \eta |^2 +
         \frac{\partial^2}{\partial t^2} : H^s(\R_+) \to
	 H^{s- \mu}(\R_+)    
\end{equation}
has the kernel
$\ker \sigma_\partial(\Delta)(\eta) = \{ c e^{- | \eta |t} :
        c \in \C \}$
which is of dimension $1$ (the other solution $ce^{| \eta |t}$ of
$(- | \eta |^2 + \displaystyle \frac{\partial^2}{\partial  t^2}) 
u(t) = 0$ 
does not belong to $H^s(\R_+)$ on the positive half-axis).		  

In order to associate with $\sigma_\partial(A)$ a family of
isomorphisms we can try to enlarge the boundary symbol to a family of
isomorphisms
\begin{equation}
\label{1021.Delta2605.eq}  
\begin{pmatrix}
    \sigma_\partial(A) & \sigma_\partial(K) \\
    \sigma_\partial(T)& \sigma_\partial(Q)  
    \end{pmatrix} \; (y, \eta) :
       \Hsum{H^s(\R_+)}         {\C^{j_-}}     \to
       \Hsum{H^{s- \mu}(\R_+)}   {\C^{j_+}}      
\end{equation}
by entries $\sigma_\partial(T), \sigma_\partial(K),
\sigma_\partial(Q)$ of finite rank, cf. the discussion in this section below. 
In the case \eqref{1021.neu0206.eq} it suffices to set $j_- = 0$, $j_+ = 1$
and to take
\[  \sigma_\partial(T) := \r_0'      \]
with the restriction operator $\r_0' : H^s(\R_+) \to \C$ for $s >
\frac{1}{2}$, cf. analogously, Section 1.3, in particular, the
homogeneity relation \eqref{1013.rh2505.eq}. In other words, with
the Laplacian we can associate the family of isomorphisms
\begin{equation}
\label{1021.neuD1108.eq}  
\begin{pmatrix}
\sigma_\partial(\Delta) \\ \r_0'    \end{pmatrix}(\eta)\; :       
 H^s(\R_+) \to
          \Hsum{H^{s-2}(\R_+)}     {\C}       
\end{equation}
for $s > \frac{1}{2}$ which is just the boundary symbol of the
Dirichlet problem. Analogously, $\sigma_\partial(T) := \r_0' \circ
\displaystyle \frac{\partial}{\partial t}$ gives us the 
boundary symbol of the Neumann problem. 

Let us calculate the inverse of \eqref{1021.neuD1108.eq}.
Writing 
\[  l_\pm(\eta) := | \eta | \pm i \tau    \]
we have $- l_-(\tau)l_+(\tau) = - (|\eta|^2 + \tau^2)$ and
\[  \sigma_\partial(\Delta)(\eta) = - \op^+(l_-)(\eta)
       \op^+(l_+)(\eta) = - |\eta|^2 + 
       \frac{\partial^2}{\partial t^2}.    \]
The operator $\op^+(l_-)(\eta) : {\cal S}(\ol{\R}_+) \to
{\cal S}(\ol{\R}_+)$ is an isomorphism for every $\eta \not= 0$
where $(\op^+(l_-)(\eta))^{-1} = \op^+(l_-^{-1})(\eta)$, and
$\op^+(l_+)(\eta) : {\cal S}(\ol{\R}_+) \to
      {\cal S}(\ol{\R}_+)$
is surjective for every $\eta \not= 0$ with
\[  
\ker \op^+(l_+)(\eta) = \{ c e^{- | \eta |t} : c \in \C \}.
\]
Let us form the map $k(\eta) : \C \to {\cal S}(\ol{\R}_+)$ by
$k(\eta) c := c e^{- | \eta |t}$. 
Then we have
\[  \begin{pmatrix}
    \op^+(l_+)(\eta) \\ \r_0'\end{pmatrix}
          \begin{matrix}
          (\op^+(l_+^{-1})(\eta) \quad  % & \\ &             
          \quad k(\eta)) = \!\!\!\! &\\& \end{matrix}
     \begin{pmatrix}
     1 & 0  \\  0 & 1
     \end{pmatrix}              \]
and     
\[    \begin{matrix}
    (\op^+(l_+^{-1})(\eta) \quad k(\eta)) \!\!\!\! &\\&
    \end{matrix}
      \begin{pmatrix}   
      \op^+(l_+)(\eta) \\ \r_0'  
     \end{pmatrix}  = 1.       \] 
Thus, because of
\[  \begin{pmatrix}
    \sigma_\partial(\Delta)(\eta) \\  \r_0'
    \end{pmatrix} = \begin{pmatrix}
                    - \op^+(l_-)(\eta) & 0 \\  0 & 1
		    \end{pmatrix}
	 \begin{pmatrix}
	 \op^+(l^+)(\eta)  \\   \r_0'
	 \end{pmatrix}                 \]
it follows that
\begin{align*}  
    \begin{pmatrix}
    \sigma_\partial(\Delta)(\eta)  \\  \r_0'
    \end{pmatrix}^{-1}& = 
        \begin{matrix} 
	(\op^+(l_+^{-1})(\eta) \quad k(\eta))\!\!\!\!  &\\&	 	    
        \end{matrix}
     \begin{pmatrix}
     -   \op^+(l_-^{-1})(\eta)  &  0   \\   0  &  1
     \end{pmatrix}   \\
      & = \bigl(- \op^+(l_+^{-1})(\eta) \op^+(l_-^{-1})(\eta)
             \quad k(\eta)\bigr).    
\end{align*}

\begin{Remark}	     
\label{1021.re1208.re}
The potential part $k(\eta)$ gives rise to an operator-valued symbol in
the sense of Definition {\st{\ref{1013.Sy2005.de}}} {\st{(}}and
Remark {\st{\ref{1013.re0206.re}}}{\st{)}}, namely,
$\chi(\eta) k(\eta) \in S_{\cl}^{- \frac{1}{2}}(\R_\eta^{n-1}; \C,
       {\cal S}(\ol{\R}_+))$
for any excision function $\chi$.
Moreover, the operator function $\chi(\eta) g(\eta)$ for
\[  
 g(\eta) := - \op^+(l_+^{-1})(\eta) \op^+(l_-^{-1})(\eta) +
      \op^+(l_+^{-1} l_-^{-1})(\eta)   
 \]
is a Green symbol of order $- 2$ and type $0$, cf. the
terminology in Example {\st{\ref{1013.G0206.exa}}},
\[  
\chi(\eta) g(\eta) \in S_{\cl}^{-2}(\R_\eta^{n-1}; L^2(\R_+),
     {\cal S}(\ol{\R}_+)),     
\]
and the $\eta$-wise $L^2(\R_+)$-adjoint has the property
\[  
(\chi g)^*(\eta) \in S_{\cl}^{-2}(\R_\eta^{n-1}; L^2(\R_+), {\cal
         S}(\ol{\R}_+)).     \]
\end{Remark}     
           
In general, if $A$ is an elliptic differential operator, by virtue
of Theorem \ref{1021.th1105.th} we expect that $j_- = 0$ is
adequate and that we can complete $\sigma_\partial(A)$ by 
$j_+ := \dim \ker \sigma_\partial(A)$ trace conditions to a 
family of isomorphisms
\begin{equation}
\label{1021.neu2805.eq}  
\begin{pmatrix}
    \sigma_\partial(A) \\ \sigma_\partial(T)
\end{pmatrix} \; (y, \eta) : H^s(\R_+) \to
        \Hsum{H^{s- \mu}(\R_+)}       {\C^{j_+}},     
\end{equation}
$(y, \eta) \in T^* \R^{n-1} \setminus 0$.
In this case $\C^{j_+}$ is interpreted as the fibre of a
vector bundle over $T^*\R^{n-1} \setminus 0$; it may be regarded 
as the pull back of a vector bundle $J_{+,1}$ on the sphere 
bundle
$\R^{n-1} \times S^{n-2}$ under the projection $(y, \eta) \to (y,
\eta / | \eta |)$.	

In order to be able to interpret $\sigma_\partial(T)$ as the
boundary symbol of a trace operator
\[  T : H^s(\R_+^n) \to \oplus_{l=1}^{j_+}
        H^{s-m_j - \frac{1}{2}}(\R^{n-1})      \]
(with orders $m_j$, according to the $\kappa_\lambda$ 
homogeneity of the
components of $T = {}^\t(T_1, \ldots, T_{j_+}))$ 
the bundle $J_{+,1}$ has to be the pull back under the projection
$(y, \eta) \to y$ of a vector bundle
$J_+$ on the boundary $\R^{n-1}$ itself. This is an assumption
that we now
impose, although it may be too restrictive in some cases, 
cf. the
discussion of Section 5.3 below in connection with the Atiyah-Bott
obstruction.

However, for the Dirichlet or the Neumann problem for the Laplace
operator as well as for many other interesting problems this
obstruction vanishes; this is enough for the purposes of this
section	(it turns out that the insight from this situation is very
useful also for the general case, cf. \cite{Schu37}).
	 	  
From the classical analysis of boundary value problems for a
differential operator \eqref{1021.Amu1105.eq} of order $\mu = 2m$
it is known that additional trace operators $T = {}^\t(T_1,
\ldots, T_m)$ may have the form
\begin{equation}
\label{1021.Tj2805.eq}
(T_j u)(y) := \r'B_j u(y)
\end{equation}
for differential operators $B_j = \sum_{|\beta|\leq m_j} b_{j
\beta}(x) D_x^\beta$ of different orders $m_j$, with $(\r' v)(y) :=
v(y,0)$, such that, when we set $\sigma_\partial(T_j)(y, \eta) =
\r_0' \sigma_\psi(B_j)(y,0, \eta, D_t)$, $j = 1, \ldots, m$, 
with
$\sigma_\psi(B_j)$ being the homogeneous principal symbol of $B_j$
of order $m_j$, the operators
\[  \sigma_\partial(T)(y, \eta) := {}^\t(\sigma_\partial(T_j)(y,
       \eta))_{j=1, \ldots, m}      \]
complete $\sigma_\partial(A)(y, \eta)$ to a family of isomorphisms       		  
\eqref{1021.neu2805.eq} for all sufficiently large real $s$ (and
$j_+ = m$).

Globally, an elliptic boundary value problem for a (scalar)
differential operator $A$ on a (say, compact) $C^\infty$ manifold
$X$ with boundary $\partial X$ is represented by a column matrix 
\begin{equation}
\label{1021.calA2805.eq}
{\cal A} := \binom{A}{T}
\end{equation}
consisting of the elliptic operator $A$ itself, and a column
vector $T$ of trace operators, with entries $T_j$ of the form
\eqref{1021.Tj2805.eq}, with differential operators $B_j$ in a
neighbourhood of $\partial X$, $m_j = \ord B_j$.

It is often convenient to unify the orders by passing to
compositions
$\widetilde{T}_j := R^{\mu - m_j - \frac{1}{2}} T_j$ 
with order reducing isomorphisms $R^{\mu - m_j - \frac{1}{2}}$ 
(as mappings $H^s(\partial X) \to H^{s- \mu + m_j +
\frac{1}{2}}(\partial X)$, $s \in \R$) belonging to $L_{\cl}^{\mu
- m_j - \frac{1}{2}}(\partial X)$. We can find such operators with
homogeneous principal symbol $| \eta |^{\mu - m_j - \frac{1}{2}}$.
Then we have
$\sigma_\partial(\widetilde{T}_j)(y, \eta)  =
        | \eta |^{\mu - m_j - \frac{1}{2}} 
	\sigma_\partial(T_j)(y, \eta)$
and
$\sigma_\partial(\widetilde{T}_j)(y, \lambda \eta)  =
        \lambda^\mu \sigma_\partial(\widetilde{T}_j)
	(y, \eta) \kappa_\lambda^{-1}$
for all $\lambda \in \R_+$, $(y, \eta) \in T^* \R^{n-1} \setminus
0$.		
Of course, we can reach any other order of the trace operators by
composition from the left with a suitable order reducing
isomorphism.

If we now assume that the trace operators are defined from the
very beginning in combination with order reductions from the left
and denote the trace operators again by $T$ (rather than
$\widetilde{T}$) our boundary value problems
\eqref{1021.calA2805.eq} induces continuous operators
\begin{equation}
\label{1021.calAJ2805.eq}
{\cal A} = \binom{A}{T} : H^s(\Int X) \to
   \Hsum{H^{s- \mu}(\Int X)}
        {H^{s- \mu}(\partial X,J_+)}
\end{equation}
for sufficiently large $s \in \R$. In this notation $J_+$ is a
(smooth complex) vector bundle on $\partial X$, similar to the
above one in the half-space case, and $H^r(\partial X, J_+)$ is
the space of distributional sections in $J_+$ of Sobolev
smoothness $r \in \R$.	

The boundary symbol $\sigma_\partial({\cal A})$ of the operator
${\cal A}$ is  a bundle morphism
\begin{equation}
\label{1021.sigmapart2805.eq}
\sigma_\partial({\cal A}) : \pi_{\partial X}^* H^s(\R_+) \to 
      \pi_{\partial X}^*
    \Hsump{H^{s - \mu}(\R_+)}       {J_+},
\end{equation}
the global analogue of \eqref{1021.neu2805.eq}. Here 
$\pi_{\partial X} :
T^*(\partial X) \setminus 0 \to \partial X$ is the canonical
projection, and $\pi_{\partial X}^*$ denotes the  pull backs of vector
bundles, here with the corresponding infinite-dimensional fibres.

Homogeneity of $\sigma_\partial({\cal A})$ means
\begin{equation}
\label{1021hom2805.eq}
\sigma_\partial({\cal A})(y, \lambda \eta) = \lambda^\mu
  \begin{pmatrix}
  \kappa_\lambda & 0 \\ 0 & 1
  \end{pmatrix} \sigma_\partial({\cal A})(y, \eta)
  \kappa_\lambda^{-1}
\end{equation}
for all $\lambda \in \R_+$, $(y, \eta) \in T^*(\partial X)
\setminus 0$, where 1 indicates the identity operator;
$(\kappa_\lambda u)(t) = \lambda^{\frac{1}{2}} u(\lambda t)$,
$\lambda \in \R_+$.

The expectation that the composition of operators gives rise to
the composition of the associated principal symbols is not so easy
to satisfy in the case of boundary value problems, since there is
no reasonable composition between the corresponding column matrices
(although we have the relations \eqref{1021.comp2805.eq}).
However, as we shall see, such a composition property is true of 
block
matrices when the number of rows and columns in the
middle fits together. This is
a natural concept in an operator algebra in a suitably generalised
sense. An access to this construction is what we obtain from the
ellipticity. 

Before we give the definition we want to make a remark on the
nature of symbols of operators $A$ on a $C^\infty$ manifold $X$
with boundary $\partial X$. 
Without loss of generality we may assume that $A$ has the form
\begin{equation}
\label{1021.r+0906.eq}
A = \r^+ \widetilde{A} \e^+
\end{equation}
for a differential operator $\widetilde{A}$ in a neighbouring
$C^\infty$ manifold $\widetilde{X}$ (for instance, the double of
$X$), where $\e^+$ is the operator of extension by zero from
$\Int X$ to $\widetilde{X}$ and $\r^+$ the restriction to $\Int
X$. Operators of the from \eqref{1021.r+0906.eq} also make sense
for arbitrary $\widetilde{A} \in L_{\cl}^\mu(\Int X)$ (of
course, also for non-classical pseudo-differential operators) as
continuous operators
\begin{equation}
\label{1021.cont0906.eq}
A = \r^+ \widetilde{A} \e^+ : C_0^\infty(\Int X) \to
               C^\infty(\Int X).
\end{equation}

In this section we content ourselves with integer orders $\mu$.
Let $\tilde{a}_{(\mu-j)}(y,t, \eta, \tau)$, $j \in \N$, denote
the sequence of homogeneous components of order $\mu - j$
belonging to a representation of $\widetilde{A}$ in local
coordinates $(y, t) \in \Omega \times \R$ near the boundary,
$\Omega \subseteq \R^{n-1}$ open. Then $\widetilde{A}$ is said
to have the transmission property at the boundary if
\begin{equation}
\label{1021.neu1106.eq}  
\tilde{a}_{(\mu-j)}(y,t,- \eta, - \tau) - e^{i \pi(\mu-j)}
      \tilde{a}_{(\mu-j)}(y,t,\eta, \tau)     
\end{equation}
vanishes to the infinite order on the set of non zero normal
covectors to the boundary
\[  \{ (y,t,\eta,\tau) \in T^*(\Omega \times \R) : t = 0, \eta =
       0, \tau \not= 0 \}      \]
for all $j \in \N$. This is an invariant condition; so it makes
sense as a property of $\widetilde{A}$ globally on $X$ near the
boundary. Since the condition is satisfied if and only if all
$a_{(\mu-j)}(y,t,\eta,\tau) := \tilde{a}_{(\mu-j)}
(y,t,\eta,\tau) \big|_{\Omega \times \ol{\R}_+ \times (\R^n
\setminus \{ 0 \})}$ have this property, we also talk about the
transmission property of the operator $A$ itself.   

\begin{Remark}   
\label{1021.0906re.re}
A differential operator $A$ {\st{(}}with smooth coefficients up
to the boundary{\st{)}} has the transmission property at the
boundary. Writing $A$ in the form \eqref{1021.r+0906.eq} the
ellipticity of $A$ entails the ellipticity of $\widetilde{A}$ in
a neighbourhood of $\partial X$. Then, if we form a parametrix
$\widetilde{P}$ in $L_{\cl}^{- \mu}(\widetilde{X})$
{\st{(}}i.e., $I - \widetilde{A} \widetilde{P}, 
I - \widetilde{P}
\widetilde{A} \in L^{- \infty}(X)${\st{)}}, also
$\widetilde{P}$ has the transmission property at $\partial X$.
\end{Remark}

\begin{Remark}
\label{1021.con2206.re}
Let $S^*X$ denote the unit cosphere bundle on $X$ {\st{(}}with
respect to a fixed Riemannian metric on $X${\st{)}}, and let
$N^*$ denote the bundle of covectors normal to
the boundary that are of length $\leq 1$. Set
\begin{equation}
\label{1021.Psi2206.eq}
\Xi := S^*X |_{\partial X} \cup N^*
\end{equation}
which is a fibre bundle on $\partial X$ with fibres being
\textup\{unit spheres\} $\cup$ \{straight connection of 
south and north poles\}, where the south and north poles are 
locally
representend by $(y,0,0, - 1)$ and $(y, 0,0, + 1)$,
respectively. Then, if $\widetilde{A} \in
L_{\cl}^0(\widetilde{X})$ is an operator with the transmission
property at the boundary, the homogeneous principal symbol
$\sigma_\psi(A)$ of \eqref{1021.r+0906.eq} extends from $S^*X$
to a continuous function  $\bs{\sigma}_\psi(A)$ on $S^*X \cup
N^*$ {\st{(}}including the zero section of $N^*$ which is
represented by $\partial X${\st{)}}. The ellipticity of $A$
with respect to $\sigma_\psi(.)$, i.e.,
$\sigma_\psi(A) \not= 0 \quad \textup{on}\quad T^*X 
      \setminus  0$,  
entails the property
\begin{equation}
\label{1021.bssigma2206.eq}
\bs{\sigma}_\psi (A) \not= 0 \quad 
     \textup{on}\quad \Xi.
\end{equation}	       
\end{Remark}

\begin{Definition}
\label{1021.de2805.de}
The operator \eqref{1021.calA2805.eq} is called elliptic, if both
components of its principal symbol
\[  \sigma({\cal A}) = (\sigma_\psi({\cal A}), 
      \sigma_\partial({\cal A}))     \]
are bijective, i.e., for the principal interior symbol
$\sigma_\psi({\cal A}) := \sigma_\psi(A)$ we have
$\sigma_\psi({\cal A}) \not= 0$ on $T^* X \setminus 0$, and the
principal boundary symbol $\sigma_\partial({\cal A})$ defines an
isomorphism \eqref{1021.sigmapart2805.eq} for any 
{\st{(}}sufficiently large{\st{)}} $s$.      
\end{Definition}

\begin{Theorem}
\label{1021.Ell2805.th}
Let $X$ be a compact $C^\infty$ manifold with boundary, and
\eqref{1021.calA2805.eq} an operator of the described structure.
Then the following properties are equivalent:
\begin{enumerate}
\item The operator ${\cal A}$ is elliptic;
\item ${\cal A}$ induces a Fredholm operator
      \eqref{1021.calAJ2805.eq} for any fixed 
      {\st{(}}sufficiently large{\st{)}} $s$.
\end{enumerate}
The property {\em(ii)} entails the Fredholm property \eqref{1021.calAJ2805.eq} for all 
{\em(}sufficiently large{\em)} $s$.
\end{Theorem}

As a Fredholm operator \eqref{1021.calAJ2805.eq} the elliptic
operator ${\cal A}$ has a parametrix
\begin{equation}
\label{1021.calP0206.eq}
{\cal P} = (P \quad K)
\end{equation}
in the functional analytic sense, and it is interesting to
characterise the nature of ${\cal P}$. The operator $P$ should 
belong to
$L_{\cl}^{- \mu}(\Int X)$. As noted before, since the original elliptic differential
operator $A$ can be seen as the restriction of an elliptic differential operator
$\widetilde{A}$ in a neighbouring $C^\infty$ manifold
$\widetilde{X}$ to $X$, we can form a parametrix $\widetilde{P} \in
L_{\cl}^{- \mu}(\widetilde{X})$ of $\widetilde{A}$ and ask the
relationship between $\widetilde{P}$ and the operator $P$ in 
the formula
\eqref{1021.calP0206.eq}. An answer was given in Boutet de 
Monvel's
paper \cite{Bout1}, not only of this point, but about the
pseudo-differential structure of ${\cal P}$ itself. We do not 
repeat here all the
details; there are many expositions on Boutet de Monvel's 
theory
of pseudo-differential boundary value problems, see, for instance,
Rempel and Schulze \cite{Remp2}, Grubb \cite{Grub1}, Schulze
\cite{Schu50}.
We want to observe here some specific features and `strange'
points of the pseudo-differential calculus of boundary value
problems.
If we form $\r^+ \widetilde{P} \e^+$ we obtain a
continuous operator
\begin{equation}
\label{1021.P0206.eq}
\r^+ \widetilde{P} \e^+ : H^{s- \mu}(\Int X) \to H^s(\Int X)
\end{equation}
for every $s > \mu - \frac{1}{2}$. 

For our differential operator
we have $A = \r^+ \widetilde{A} \e^+$ and
$(\r^+ \widetilde{A} \e^+)(\r^+ \widetilde{P} \e^+) =
       \r^+ \widetilde{A} \widetilde{P} \e^+$     
which is the identity map modulo an operator with kernel in
$C^\infty(X \times X)$.  The composition $(\r^+ \widetilde{P}
\e^+)(\r^+ \widetilde{A} \e^+)$ has a more complicated structure;
it is equal to the identity modulo a smoothing operator $G$ in
$\Int X$, however, not with a kernel in $C^\infty(X \times X)$. The
operator $G$ is called a Green operator, and it is locally in a
collar neighbourhood of $\partial X$ of the form $\Op(g)$ for a
Green symbol $g(y, \eta)$ of some type $\d$, cf. Example
{\st{\ref{1013.G0206.exa}}}.

\begin{Theorem}
\label{1021.calP0206.th}
Let ${\cal A} := \displaystyle \binom{A}{T}$ be an elliptic boundary value
problem for the differential operator $A$. Then there is a
{\st{(}}two-sided{\st{)}} parametrix ${\cal P}$ of ${\cal A}$ of
the form \eqref{1021.calP0206.eq} for
\begin{equation}
\label{1021.neu0806.eq}  
P = \r^+ \widetilde{P} \e^+ + G   
\end{equation}
for some Green operator $G$, and a potential 
operator $K$, cf.
Example {\st{\ref{1013.G0206.exa}}}.      
\end{Theorem}

The result of Theorem \ref{1021.calP0206.th} is remarkable for
several reasons. First of all, if we accept the operator
\eqref{1021.calP0206.eq} as a `boundary value problem' for the
pseudo-differential operator $P$ (whatever its precise structure near the
boundary is) instead of boundary conditions we have potential 
conditions, represented by the
operator $K$. The symbol of $K$ is associated with the second
component of the inverse of \eqref{1021.sigmapart2805.eq} which is a
row matrix
\[  \sigma_\partial({\cal A})^{-1}(y, \eta) = (\sigma_\partial(P)(y,
      \eta)\quad \sigma_\partial(K)(y, \eta)).      \]
In this case
\begin{equation}
\label{1021.sigmaP0806.eq}
\sigma_\partial(P)(y, \eta) : H^{s- \mu}(\R_+) \to H^s(\R_+)
\end{equation}
is necessarily injective but not surjective, and the operators of
finite rank
\begin{equation}
\label{1021.eq0806.eq}
\sigma_\partial(K)(y, \eta) : J_{+, y} \to H^s(\R_+)
\end{equation}
fill up the family \eqref{1021.sigmaP0806.eq} to a family of
isomorphisms (here $J_{+,y}$ denotes the fibre of $J_+$ over the
point $y$). The local structure of $\sigma_\partial(K)(y, \eta)$ is
just as in Example \ref{1013.G0206.exa}; in fact
$\sigma_\partial(K)(y, \eta)$ is the vector of homogeneous
principal components of order $- \mu$ of potential symbols of the
kind $f_K(t[\eta]; y, \eta)$. Concerning the structure of
\eqref{1021.sigmaP0806.eq} we have, according to
\eqref{1021.neu0806.eq},
\begin{equation}
\label{1021.sigmapart(P)0806.eq}
\sigma_\partial(P)(y, \eta) = \sigma_\partial(\widetilde{P})(y,
                   \eta) + \sigma_\partial(G)(y, \eta),     
\end{equation}
where $\sigma_\partial(G)(y, \eta)$ is the homogeneous principal
component of order $- \mu$ of a Green symbol in the sense of
Example \ref{1013.G0206.exa}, while	
\[  
\sigma_\partial(\widetilde{P})(y, \eta) = \r^+ \sigma_\psi(\widetilde{P})(y,
           0, \eta, D_t) \e^+,	  
 \]
with $\sigma_\psi(\widetilde{P})(y,t,\eta, \tau)$ being the
homogeneous principal symbol of $\widetilde{P}$ 	   
near the boundary in the splitting of variables $x =(y,t)$, and $\e^+$ is the operator
of extension by zero from $\R_+$ to $\R$, and $\r^+$ the
restriction from $\R$ to $\R_+$.

What we see is the following. Given an elliptic pseudo-differential
operator $\widetilde{A}$ of order $\mu$ in a neighbouring manifold
$\widetilde{X}$ of a $C^\infty$ manifold $X$ with boundary (with the
transmission property at $\partial X$)  we can form
the operator
\[  A = \r^+ \widetilde{A} \e^+ : H^s(\Int X) \to H^{s- \mu}(\Int
                          X)       \]
(say, for $s > \max(\mu,0) - \frac{1}{2}$). Its boundary 
symbol
\begin{equation}
\label{1021.Fred0806.eq}
\sigma_\partial(A)(y, \eta) = \r^+ \sigma_\psi(A)(y, 0, \eta,
       D_t)\e^+ : H^s(\R_+) \to H^{s- \mu}(\R_+)			  
\end{equation}       
is a family of Fredholm operators (in general, neither 
surjective nor
injective) for $(y, \eta) \in T^*(\partial X) \setminus 0)$. Then,
elliptic conditions may exist both of trace and potential type in a
way that the associated boundary symbols
\[  \sigma_\partial(T)(y, \eta) : H^s(\R_+) \to J_{+, y}, \quad
    \sigma_\partial(K)(y, \eta) : J_{-, y} \to H^{s- \mu}(\R_+)   \]
for suitable vector bundles $J_\pm$ on $\partial X$ (for algebraic
reasons combined with a family of maps $\sigma_\partial(Q)(y, \eta)
:= \sigma_\psi(Q)(y, \eta) : J_{-, y} \to J_{+, y}$ for an operator $Q
\in L_{\cl}^\mu(\partial X; J_-, J_+)$) fill up the Fredholm family
\eqref{1021.Fred0806.eq} to a family of isomorphisms    
\begin{equation}
\label{1021.sigma(calA)0806.eq}
\begin{pmatrix}
\sigma_\partial(A) & \sigma_\partial(K)   \\
\sigma_\partial(T) & \sigma_\partial(Q)
\end{pmatrix}     \; (y, \eta): \;
   \Hsum{H^s(\R_+)}            {J_{-, y}}    \to
   \Hsum{H^{s- \mu}(\R_+)}     {J_{+, y}}
\end{equation}
for all $(y, \eta) \in T^*(\partial X) \setminus 0$.

Both trace and potential symbols may be 
necessary at the same time for
obtaining an isomorphism \eqref{1021.sigma(calA)0806.eq}. Locally 
the
operator families
$\sigma_\partial(T)$, $\sigma_\partial(K)$ have the structure of
homogeneous
principal components of trace and potential symbols as in
Example \ref{1013.G0206.exa} (of some type $\d \in \N$ in the
case of $\sigma_\partial(T)$). More generally, instead
of \eqref{1021.sigma(calA)0806.eq} we can consider isomorphisms
of the kind
\begin{equation}
\label{1021.sigmaG0906.eq}
\begin{pmatrix}
\sigma_\partial(A) + \sigma_\partial(G) & \sigma_\partial(K) \\
\sigma_\partial(T)     &     \sigma_\partial(Q)
\end{pmatrix} (y, \eta) :      
  \Hsum{H^s(\R_+)}                {J_{-, y}}     \to              
  \Hsum{H^{s- \mu}(\R_+)}         {J_{+, y}}
\end{equation}
with a Green symbol $\sigma_\partial(G)(y, \eta)$ of analogous
structure as in \eqref{1021.sigmapart(P)0806.eq} (it 
takes values in the space of compact operators $H^s(\R_+) \to
H^{s- \mu}(\R_+)$). 
Green symbols are generated in compositions of block matrices of the form 
\eqref{1021.sigma(calA)0806.eq} and also in inverses.
We now pass to an operator
\begin{equation}
\label{1021.calA0906.eq}
{\cal A} :=
 \begin{pmatrix}
 A+G  &  K   \\  T  &  Q
 \end{pmatrix} :
   \Hsum{H^s(\Int X)} 
        {H^s(\partial X, J_-)}    \to
   \Hsum{H^{s- \mu}(\Int X)}      
        {H^{s- \mu}(\partial X, J_+)},
\end{equation}
where $A$ is the original elliptic operator, and $G,T,K,Q$ 
are the
extra operators which constitute an elliptic boundary value
problem \eqref{1021.calA0906.eq} for $A$ with the principal
symbolic structure
\begin{equation}
\label{1021.sigma(calA)906.eq}
\sigma({\cal A}) = (\sigma_\psi({\cal A}), \sigma_\partial({\cal
                   A})),   
\end{equation}
for $\sigma_\psi({\cal A}) := \sigma_\psi(A)$ and
$\sigma_\partial({\cal A})$ given by \eqref{1021.sigmaG0906.eq}.
Note that
\begin{equation}
\label{1021.ho0906.eq}
\sigma_\partial({\cal A})(y, \lambda \eta) = \lambda^\mu
   \begin{pmatrix}
   \kappa_\lambda & 0  \\  0  &  1
   \end{pmatrix}
   \sigma_\partial({\cal A})(y, \eta)
   \begin{pmatrix}
   \kappa_\lambda & 0  \\  0  &  1
   \end{pmatrix}^{-1}
\end{equation}
for all $(y, \eta) \in T^*(\partial X) \setminus 0$, $\lambda 
\in \R_+$.   

Operators of the form
\begin{equation}
\label{1021.eq0906.eq}
{\cal A} = \begin{pmatrix}
           A+G  &  K  \\  T  &  Q
	   \end{pmatrix}
\end{equation}
constitute what is also called Boutet de Monvel's calculus (of
pseudo-differential boundary value problems with the
transmission property), cf. \cite{Bout1}. 
Operators of that kind also make sense on a not necessarily compact $C^\infty$ manifold with boundary.

We now  propose `answer number 1' to the question `what is a
boundary value problem' to an operator $A$, namely, such a 

\begin{quote}
 `{$2
   \times 2$}  block matrix \eqref{1021.eq0906.eq} with $A$ in the
   upper left corner, where the extra operators 
    $G,  T, K, Q$ 
   are an additional information from the boundary'    
   \end{quote}
%\end{align*}   
of a specific nature (roughly speaking, pseudo-differential operators on the
boundary with operator-valued symbols as in Example
\ref{1013.G0206.exa}).

Let $\got{B}^{\mu,d}(X)$ denote the space of all operator block
matrices of the form \eqref{1021.eq0906.eq} that are of the
structure as mentioned before, in particular, $A$ is of
order $\mu \in \Z$ and has the transmission property, and the
other operators are of order $\mu$ and type $\d \in \N$.

In this block matrix set-up the multiplicativity of the
principal symbols \eqref{1021.sigma(calA)906.eq} is again
restored; the only condition for a composition ${\cal A} {\cal
B}$ (say, for compact $X$, otherwise combined with a
localisation) is that rows and columns in the middle fit
together (more precisely, the bundles on the boundary), and we
then have
\begin{equation}
\label{1021.Pro0906.eq}
\sigma({\cal A} {\cal B}) = \sigma({\cal A}) \sigma({\cal B}),
\end{equation}
where the multiplication is componentwise, i.e.,
$\sigma_\psi({\cal A} {\cal B}) = \sigma_\psi({\cal A})
\sigma_\psi({\cal B})$, $\sigma_\partial({\cal A} {\cal B}) =
\sigma_\partial({\cal A}) \sigma_\partial({\cal B})$.

\begin{Remark}
\label{1021.1506re.re}
The definitions and results about operators 
\eqref{1021.eq0906.eq} including 
Definition {\st{\ref{1021.df0906.de}}}, and Theorems
{\st{\ref{1021.1006th.th}}}, {\st{\ref{1021.the1106.th}}} below
easily extend to operators between
distributional  sections of vector bundles $E, F \in \Vect(X)$ and
$J_\pm \in \Vect(Y)$. In this case instead of
\eqref{1021.calA0906.eq} we have the continuity
\begin{equation}
\label{1021.As1506.eq}
{\cal A} :
    \Hsum{H^s(\Int X,E)}        
         {H^s(\partial X, J_-)}       \to
    \Hsum{H^{s- \mu}(\Int X,F)}	 
	 {H^{s- \mu}(\partial X, J_+)}
\end{equation}
for all $s > \d - \frac{1}{2}$ when $\d \in \N$ denotes the type of
the involved Green and  trace operators.	 
\end{Remark}

Let us now enlarge Definition \ref{1021.de2805.de} as follows.

\begin{Definition}
\label{1021.df0906.de}
An operator
is called elliptic if both components of its principal symbol
\eqref{1021.sigma(calA)906.eq} are bijective, i.e., for the
principal interior symbol of ${\cal A}$ we have
$\sigma_\psi({\cal A}) \not= 0$ on $T^*X \setminus 0$, and the
principal boundary symbol $\sigma_\partial({\cal A})$ defines
isomorphisms \eqref{1021.sigmaG0906.eq} for all $(y, \eta) \in
T^*(\partial X) \setminus 0$ and any 
{\st{(}}sufficiently large{\st{)}} $s \in \R$. The isomorphism
\eqref{1021.sigmaG0906.eq} is also called the Shapiro-Lopatinskij
condition {\st{(}}for the elliptic operator $A${\st{)}}.
\end{Definition}

\begin{Theorem}
\label{1021.1006th.th}
Let $X$ be a compact $C^\infty$ manifold with boundary and
${\cal A}$ be an operator \eqref{1021.eq0906.eq} which
represents a boundary value problem for $A$ in the upper left
corner. Then the following properties are equivalent:
\begin{enumerate}
\item The operator ${\cal A}$ is elliptic in the sense of
      Definition {\st{\ref{1021.df0906.de}}};
\item ${\cal A}$ is Fredholm as an operator
      \eqref{1021.calA0906.eq} for some fixed 
      {\st{(}}sufficiently large{\st{)}} $s \in \R$.      
\end{enumerate}
\end{Theorem}      

\begin{Theorem}
\label{1021.the1106.th}
Let $X$ be a $C^\infty$ manifold with boundary and ${\cal A} \in
\got{B}^{\mu,\d}(X)$ an elliptic operator. Then there is a parametrix
${\cal P} \in \got{B}^{\mu, (\d- \mu)^+}(X)$ in the sense that the
remainders in the relations
\[  {\cal P} {\cal A} = {\cal I} - {\cal C}_\l, \quad
    {\cal A} {\cal P} = {\cal I} - {\cal C}_\r     \]
are operators ${\cal C}_\l \in \got{B}^{- \infty, \d_\l}(X), 
{\cal C}_\r
\in \got{B}^{- \infty, \d_\r}(X)$ where $\d_\l = 
\max(\mu,\d), 
\d_\r = (\d-
\mu)^+$, and ${\cal I}$ are corresponding identity operators. Here
$\nu^+ := \max(\nu, 0)$ for any $\nu \in \R$.    
\end{Theorem} 

Summing up the calculus of operators \eqref{1021.eq0906.eq} with its
symbolic structure solves the problem to find an operator algebra
that contains all elliptic boundary value problems
\eqref{1021.calA2805.eq} for differential operators together 
with their parametrices
\eqref{1021.calP0206.eq}. Block matrices appear, for
instance, in
compositions  when we form
\begin{equation}
\label{1021.eq1106.eq}
\binom{A}{T} (P \quad K)   = 
                           \begin{pmatrix}
                           AP & AK   \\   TP & TK
			   \end{pmatrix}
\end{equation} 
for different elliptic operators $\displaystyle \binom{A}{T}$ 
and $(P \quad K)$
(not necessarily being a parametrix of each other).

In the special case that $(P \quad K)$ is the parametrix of an
elliptic boundary value problem $\displaystyle \binom{A}{T}$, 
and if 
$\displaystyle \binom{A}{\widetilde{T}}$ is another elliptic 
boundary value
problem for the same operator $A$, then we have
\begin{equation}
\label{1021.Mult1506.eq}
\binom{A}{\widetilde{T}} (P \quad K) = 
                  \begin{pmatrix}
                      1 & 0 \\ \widetilde{T} P & \widetilde{T} K
		      \end{pmatrix}
\end{equation}
(modulo a compact operator in Sobolev spaces).

\begin{Remark}
\label{1021.red1506.re}
The operator $\widetilde{T} K$ is a classical elliptic
pseudo-differential operator on $\partial X$, called the reduction
of $\widetilde{T}$ to the boundary {\st{(}}by means of
$T${\st{)}}, and we have
\begin{equation}
\label{1021.Ag1506.eq}
\ind \binom{A}{T} - \ind \binom{A}{\widetilde{T}} = \ind
     \widetilde{T} K.
\end{equation}
The relation \eqref{1021.Ag1506.eq} is also called the
Agranovich-Dynin formula. It compares the indices of elliptic
boundary value problems for the same elliptic operator $A$ in
terms of an elliptic pseudo-differential operator on the
boundary. A result of that kind is also true for boundary value
problems of general $2 \times 2$ block matrix form, cf. \stcite{Remp2}.
\end{Remark}     

This is one of the occasions where pseudo-differential operators
are really useful to understand the nature of elliptic boundary
value problems for differential operators (apart from the 
aspect
of expressing parametrices). Elliptic pseudo-differential
operators on the boundary `parametrise' via the formula
\eqref{1021.Mult1506.eq} the set of all possible elliptic boundary
value problems for an elliptic operator $A$ on a compact
$C^\infty$ manifold with boundary.

This shows, in particular, that there are many different elliptic
boundary value problems for $A$ (which is also evident by
the above filling up procedure of $\sigma_\partial(A)$ to an
isomorphism). Of course, it is not so clear at the first glance
how many elliptic problems \eqref{1021.calA2805.eq} exist for an
elliptic differential operator $A$ with differential boundary
conditions of the kind \eqref{1021.Tj2805.eq} (up to the
pseudo-differential order reduction on the boundary that we
admitted for simplifying the formulation in the sense of
\eqref{1021.calAJ2805.eq}). 
An answer is given in Agmon, Douglis, and Nirenberg \cite{Agmo1}.
There are also elliptic differential operators $A$ that do not admit at all elliptic
boundary value problems in the sense \eqref{1021.eq0906.eq} (for
instance, Dirac operators in even dimensions and other
interesting geometric operators). 
Later on we will return to this aspect from the point of view of edge conditions.

We will discuss this problem in Section 5 in more detail. 
At least, the existence of operators of that kind
shows that regular boundaries are not harmless from such a
point of view.
It turns out that nevertheless there are other kinds of
elliptic boundary value problems rather than
Shapiro-Lopatinskij elliptic ones; in that framework we may
admit arbitrary elliptic operators $A$, cf. \cite{Schu37},
\cite{Schu50}.

Let us consider the case when a differential operator $A$
admits two different (Shapiro-Lopatinskij) elliptic problems
$\displaystyle \binom{A}{T_+}$ and $\displaystyle \binom{A}{T_-}$ 
with trace operators
$T_+$ and $T_-$. An example is the Laplace operator $\Delta$
on a (say, compact) $C^\infty$ manifold $X$ with boundary and
$T_-$ the Dirichlet, $T_+$ the Neumann condition.

An interesting category of boundary value problems, are mixed
problems, where the boundary $\partial X$ is subdivided into
two (say, $C^\infty$) submanifolds $Y_+, Y_-$ with common
boundary $Z$ (of codimension $1$ on $\partial X$) such that
$\partial X = Y_- \cup Y_+$, $Z = Y_- \cap Y_+$. Let us
slightly change notation and identify $T_\pm$ with the
restriction of the former $T_\pm$ to $\Int Y_\pm$. Then we
obtain an operator
\begin{equation}
\label{1021.calAT1606.eq}
{\cal A} := {}^\t(A \quad T_- \quad T_+)
\end{equation}
which represents a mixed problem
\begin{equation}
\label{1021.Nr1606.eq}
A u = f \quad \text{in}\quad \Int X, \quad
T_\mp u = g_\mp\quad \text{on}\quad \Int Y_\mp.
\end{equation}
The question is then which are the natural Sobolev spaces for
such problems and to what extent we can expect the Fredholm
property when $A$ is elliptic and $T_\mp$ satisfy the
Shapiro-Lopatinskij condition on $Y_\mp$ (up to $Z$ from the
respective sides).	    
\eqref{1021.calAT1606.eq} for the Laplace operator $A$ and
Dirichlet and Neumann conditions $T_\mp$ on $Y_\mp$
represents the so called Zaremba problem. 
Reducing $T_+$ to the boundary by means of $T_-$ gives rise to an operator $R$ on
$Y_+$ (which has of course an extension $\widetilde{R}$ to a
neighbourhood $\widetilde{Y}_+$ of $Y_+$ in $\partial X$) 
that has not the transmission property at
$Z$. This shows that the concept of boundary value problems
has to be generalised to the case without the transmission
property if one asks the solvability properties of mixed
problems \eqref{1021.Nr1606.eq}. 

The formulation of (pseudo-differential) boundary value problems
\eqref{1021.eq0906.eq} shows some specific features 
that should be carefully looked at.

\begin{Remark}
\label{1021.re1106.re}
The transmission property of symbols \eqref{1021.neu1106.eq} rules
out practically all symbols which are smooth up to boundary, except
for a thin set, defined by the condition \eqref{1021.neu1106.eq} for
all $j \in \N$. For instance, symbols which have $|\xi|^\mu$ as
their homogeneous principal part have the transmission property
only when $\mu \in 2 \Z$.
\end{Remark}

Observe that (up to a constant factor)
the absolute value  $|\eta|$ of the
covariable on the boundary is the homogeneous principal symbol 
  of the operator on $\partial X$ which follows from
the reduction of the Neumann problem for the Laplace operator to the
boundary by means of the potential belonging to the solution of the Dirichlet problem.
As such they fail to have the transmission
property at any hypersurface of codimension $1$ on the 
boundary.

\begin{Remark}
\label{1021.1106re.re}
Another remarkable point is that the operator convention
\eqref{1021.r+0906.eq} is not defined intrinsically on $\Int X$; it
employs the existence of a neighbouring manifold $\widetilde{X}$ and
the action of an operator $\widetilde{A}$ on $\widetilde{X}$,
combined with an extension operator $\e^+$ from $\Int X$ to the
other side and then the restriction $\r^+$ to $\Int X$.
Fortunately, despite of the jump of $\e^+u$ at $\partial X$ we have
the continuity of
\begin{equation}
\label{1021.+1106.eq}			        
\r^+ \widetilde{A} \e^+ : H^s(\Int X) \to H^{s- \mu}(\Int X)
\end{equation}
for $s > - \frac{1}{2}$ {\st{(}}when $X$ is compact, otherwise
between \textup{`$\comp/\loc$'} spaces{\st{)}}. In 
particular, it follows that
\begin{equation}
\label{1021.++1106.eq}
\r^+ \widetilde{A} \e^+ : C^\infty(X) \to C^\infty(X)
\end{equation}
is continuous.
Thus  the transmission property has the consequence that the
smoothness up to the boundary is preserved under the action.
\end{Remark}

We may ask to what extent a general pseudo-differential operator
$\widetilde{A} \in L_{\cl}^\mu(\widetilde{X})$ induces a controlled
mapping behaviour on $\Int X$ when we first realise $\r^+ A \e^+$
as a map \eqref{1021.cont0906.eq} and then try to extend it to
Sobolev spaces on $\Int X$ or to smooth functions up to the
boundary. The answer is disappointing, even in the simplest case on
the half-axis when we look at
\begin{equation}
\label{1021.op+1106.eq}
\op^+(a) := \r^+ \op(a) \e^+ : C_0^\infty(\R_+) \to
          C^\infty(\R_+)
\end{equation}	  
for a symbol $a(\tau) \in S_{\cl}^\mu(\R)$ with constant
coefficients. Taking into account that, at least for $\mu = 0$, the
operator \eqref{1021.op+1106.eq} induces a continuous map
\begin{equation}
\label{1021.L21106.eq}
\op^+(a) : L^2(\R_+) \to L^2(\R_+),
\end{equation}
there is no continuous extension as
\begin{equation}
\label{1021.1706neu.eq}  
\op^+(a) : H^s(\R_+) \to H^s(\R_+)     
\end{equation}
for arbitrary $s$ and hence no control of smoothness up to $0$. An
example where this smoothness  fails to hold is
\begin{equation}
\label{1021.1706tildeneu.eq}  
a(\tau) = \chi(\tau) \bigl( \theta^+(\tau) - 
    \theta^-(\tau) \bigr)  
\end{equation}
when $\theta^\pm(\tau)$ is the characteristic function of
$\R_\pm$ and $\chi(\tau)$ any excision function. 
The transmission property at $t = 0$ is violated in a spectacular way:
Instead of
$a_{(0)}(+1) = a_{(0)}(-1)$   
we have in this case
\[  \hspace*{3.5mm} a_{(0)}(+1) = - a_{(0)}(-1),     \]
which is a kind of `anti-transmission property'.

Let us set
\begin{equation}
\label{1021.Formelneu1606.eq}  
\Op_x(a) u(x) = \iint e^{i(x-x')\xi} a(x, \xi) u(x')dx' \dbar \xi.      
\end{equation}

Apart from the `brutal' operator convention with $\r^+$ and $\e^+$,
say, in the half-space
\begin{equation}
\label{1021.Kasten1106.eq}
\Op_y(\op^+(a)(y, \eta)) = \r^+ \Op_x(a) \e^+
\end{equation}
for symbols $a(x, \xi)= a(y,t,\eta, \tau) \in S_{\cl}^\mu(\R^{n-1}
\times \ol{\R}_+ \times \R_{\eta, \tau}^n)$ (where we omit
indicating an extension $\tilde{a}$ of $a$ to the opposite 
side, since the
choice does not affect the action on $\R^{n-1} \times \R_+
= \R_+^n$),
the question is which are the natural substitutes of the
Sobolev spaces $H^s(\R_+^n)$ which are the right choice for the case
with the transmission property.	 

This brings us back to the question of Section 1.3. As observed
before, symbols without the transmission property at the boundary
have played a role in mixed elliptic problems,
e.g., the Zaremba problem. In classical papers of Vishik and
Eskin \cite{Vivs2}, \cite{Vivs3} and the book of Eskin \cite{Eski2}
it was decided to realise \eqref{1021.Kasten1106.eq} as continuous
operators
\begin{equation}
\label{1021.bo1106.eq}
\r^+ \Op_x(a) : H_0^s(\ol{\R}_+^n) \to H^{s- \mu}(\R_+^n)
\end{equation}
(e.g., under the assumption that the symbols are independent of
$x$ for large $|x|$). Here
$H_0^s(\ol{\R}_+^n) = \{ u \in H^s(\R^n) : 
                         \supp u \subseteq \ol{\R}_+^n \}$    
and
$H^s(\R_+^n) = \{ u \big|_{\R_+^n} : u \in H^s(\R^n) \}$.
			 
There is a natural identification between $H_0^s(\ol{\R}_+^n)$ 
and
$H_0^s(\ol{\R}_+^n) \big|_{\R_+}$ for $s > - \frac{1}{2}$. 
Hence, for those $s$ we can identify $\r^+ \Op(a)$ with $\r^+ 
\Op(a)
\e^+$. However, the operator convention \eqref{1021.bo1106.eq} is
not symmetric with respect to the spaces in the preimage and 
the image;
this makes the composition of operators to a problem. However, for the half-axis case and for $s = \mu =
0$ the book \cite{Eski2} gave a completely different operator
convention rather than $\op^+(a)$, based on the Mellin transform on
$\R_+$. In the following section we say more about Mellin operator
conventions. This will show why there is no hope for a continuous
restriction of \eqref{1021.L21106.eq} to a continuous map 
between Sobolev spaces $H^s(\R_+)$ for arbitrary $s > 0$ or to a continuous map
\begin{equation}
\label{1021.neue11706.eq}  
\op^+(a) : {\cal S}(\ol{\R}_+) \to {\cal S}(\ol{\R}_+)    
\end{equation}
which preserves smoothness up to zero. This answers the question of
Section \ref{s.10.2} as follows:
\begin{equation}
\label{1021.h1106.eq}
\textup{`regular boundaries are not harmless'}
\end{equation}
in the context of boundary value problems, even if the boundary is 
a single point $\{ 0 \} = \partial \ol{\R}_+$.

Nevertheless, the way out is very beautiful, and we meet old
friends: Operators of the kind
\eqref{1021.L21106.eq} belong to the cone algebra on $\ol{\R}_+$,
cf. \cite{Schu31}, where $\ol{\R}_+$ is regarded as a manifold with
conical singularity $\{ 0 \}$.

What concerns the half-space, (or, more generally, a $C^\infty$
manifold with boundary) the answer is not less surprising. The
`right' Sovolev spaces are weighted edge  spaces
\begin{equation}
\label{eq.123}
  {\cal W}^{s, \gamma}(\R_+^n) := {\cal W}^s(\R^{n-1},
    {\cal K}^{s, \gamma; g}(\R_+))    
\end{equation}
for any $g \in \R$ (in the local description near the boundary). As the 
`answer number $2$' to the question of Section 2.1 we
offer:
\begin{equation}
\label{1021.be1106.eq}
\textup{`boundary value problems are edge problems'}
\end{equation}    
in the sense of a corresponding edge pseudo-differential calculus,
cf. Rempel and Schulze \cite{Remp1}, the monograph \cite{Schu31}, 
as well as Schulze and Seiler
\cite{Schu41}. The nature
of edge problems will be discussed in more detail in Section 3.1
below.

Also mixed elliptic boundary value problems of the type
\eqref{1021.calAT1606.eq} belong to the category of edge 
problems,
where the interface $Z$ on the boundary in the above description
plays the role of an edge. The same is true of crack problems with
smooth crack boundaries as mentioned at the beginning of Section
1.3.

The case of non-smooth interfaces or boundaries (in the sense of
`higher' edges and corners) requires more advanced tools, cf.
Section \ref{s.10.5} below.
\begin{Remark}
\label{r.2.17}
It can easily be proved that 
$H^s_\comp (\R_+^n) \subset \s{W}^{s,\gamma} (\R^n_+) \subset H^s_\loc (\R_+^n)$
for every $s,\gamma \in \R$, cf., analogously, the relation \eqref{eq.65.}.
Thus, if $X$ is a {\em(}say, compact{\em)} $C^\infty$ manifold with boundary {\em(}with a fixed collar
neighbourhood of $\partial X$,  locally identified  with $\ol{\R}_+^n \ni (y,t)${\em)} from \eqref{eq.123} we obtain global spaces on $X$ that we denote by $\s{W}^{s,\gamma}(X)$.
For simplicity, in the global definition we assume the coordinate diffeomorphisms to be independent of the normal variable $t$ for small $t$.
Then, given an asymptotic type $P = \{( m_j, p_j) \}_{ j \in \N}$ as in Section {\em 1.2 }with 
$\pi_\C P \subset \{ w \in \C: \re w < \frac{1}{2} - \gamma \}$, we can also define subspaces 
$\s{W}^{s,\gamma}_P (X)$ locally near $\partial X$ based on
$\s{W}^s (\R^{n-1}, \s{K}^{s,\gamma}_P (\R_+))$.
\end{Remark}

Let us briefly return to \eqref{1021.be1106.eq}.
What we suggest (and what is really the case) is that, when we interpret a manifold with $C^\infty$ boundary as a manifold with (regular) edge (where the boundary is the edge and $\ol{\R}_+$, the inner normal, the model cone of local wedges), boundary value problems are a special case of edge problems.                                                                        The edge calculus should contain all elements of the calculus of boundary value problems in generalised form, including edge conditions of trace and potential type, as analogues of boundary conditions.
Moreover, parametrices of elliptic edge problems should contain analogues of Green's function in elliptic boundary value problems.
Those appear in parametrices, even when we ignore non-vanishing 
edge / boundary data.
If we perform the edge calculus on a manifold with boundary, where the typical differential operators $A$ are edge-degenerate, i.e.,
$A = r^{-\mu} \sum_{j+|\alpha|\leq \mu}
       a_{j\alpha} (r,y) 
          \big(-r \frac{\partial}{\partial r} \big)^j
                    (rD_y)^\alpha$
in a coordinate neighbourhood $\cong \ol{\R}_+ \times \Omega$ of the boundary,
$\Omega \subseteq \R^q$ open, 
$a_{j\alpha}   \in C^\infty
                  (\ol{\R}_+ \times \Omega)$,
then there is the following chain of proper inclusions: 
\begin{eqnarray*}
&&
\ \ \
\big\{ 
\text{bvp's with the transmission property at the boundary} 
\big\} \\
&&
\subset 
 \big\{ 
\text{bvp's without (or with) the transmission property at the boundary} 
\big\}   \\
&&
\subset 
\big\{
 \text{edge problems} \big\};
\end{eqnarray*}
here `bvp's, standards for `boundary value problems'.

\subsection{Quantisation}
\label{s.10.2.2}

Quantisation in a pseudo-differential scenario means a rule to
pass from a symbol function to an operator.
This notation comes from quantum mechanics with its
relationship between Hamilton functions on phase spaces and
associated operators in Hilbert spaces.

\begin{Definition}
\label{1022.neu0607.de}
In the pseudo-differential terminology the map
\begin{equation}
\label{1022.sy1606.eq}
\Op: \textup{symbol $\to$ operator}
\end{equation}
is called an operator convention.
\end{Definition}

Rules of that kind can be organised in terms of the Fourier
transform
$F u(\xi) = \int_{\R^n} e^{- ix \xi} u(x)dx$.     
Given a symbol $a(x, \xi)$ on the `phase space' $\R^n \times \R^n
\ni (x, \xi)$ we obtain an associated operator by $\Op(a) =
F_{\xi \to x}^{-1} \bigl\{ a(x, \xi) F_{x \to \xi} \}$, cf. the
formula \eqref{1021.Formelneu1606.eq}.

If a symbol is involved in this form we also call $a(x, \xi)$ a
`left symbol'. More generally, we may admit `double symbols'
$a(x,x', \xi)$, and especially `right symbols' $a(x', \xi)$;
then we have
 \begin{equation}
\label{1022.Opx1606.eq}
\Op(a) u(x) := \iint e^{i(x-x')\xi} a(x,x', \xi) u(x')
               dx' \dbar \xi.
\end{equation}
Concerning $x,x'$ we do not insist on the full $\R^n$ but also
admit $x,x'$ to vary in an open subset $\Omega$. Then we obtain
a continuous map 
\[  \Op(a) : C_0^\infty(\Omega) \to C^\infty(\Omega),     \]
provided that $a(x,x', \xi) \in C^\infty(\Omega \times \Omega
\times \R^n)$ belongs to a reasonable symbol class.

Here we take H\"ormander's classes
$S_{(\cl)}^\mu(\Omega \times \Omega \times \R^n)$,
cf. Definition \ref{1013.Sy2005.de} (for the case $E =
\widetilde{E} = \C$ and trivial group actions).	       
The possibility to give $a(x, \xi)$ the meaning of a left or a
right symbol (where the resulting operators are different) shows
that the quantisation process is not canonical.

\begin{Remark}
\label{1022.re0607.re}
A map
\begin{equation}
\label{1022.symb0507.eq}
\textup{symb : operator $\to$ symbol}
\end{equation}
which is a right inverse of \eqref{1022.sy1606.eq}
{\st{(}}possibly up to negligible terms{\st{)}} may be
interpreted as an analogue of semi-classical asymptotics:
Objects of classical mechanics are recovered from their
quantised versions. In pseudo-differential terms we can
construct such a map
\begin{equation}
\label{1022.sF0607.eq}
\symb : L_{(\cl)}^\mu(\Omega)_{\proper} \to S_{(\cl)}^\mu(\Omega
        \times \R^n)
\end{equation}
on the space $L_{(\cl)}^\mu(\Omega)_{\proper}$ of properly
supported elements of $L_{(\cl)}^\mu(\Omega)$ by the rule
\begin{equation}
\label{1022.Axi0607.eq}
A \to e_{- \xi} A e_\xi =: a(x, \xi)
\end{equation}
for $e_\xi := e^{i x \xi}$. This follows from the Fourier
inversion formula $u(x) = \int e^{i x \xi} \hat{u}(\xi) \dbar
\xi$ by applying $A$ on both sides with respect to $x$, which
yields $A u (x) = \int e^{i x \xi}(e_{- \xi}(x) A e_\xi(.))
\hat{u}(\xi)\dbar \xi$.
\end{Remark}

A generalisation of \eqref{1022.Opx1606.eq} is the expression
\begin{equation}
\label{1022.Opvarphi1606.eq}
\Op(a; \varphi) u(x) := \iint e^{i \varphi(x,x', \xi)}
    a(x,x', \xi) u(x') dx' \dbar \xi,
\end{equation}
$a(x,x', \xi) \in S_{(\cl)}^\mu(\Omega \times \Omega \times \R^n)$.
Here $\varphi(x,x', \xi) \in C^\infty(\Omega \times \Omega \times
\R^n)$ is a real-valued (so called pseudo-differential phase) 
function of the form
\[  \varphi(x,x', \xi) = \sum_{j=1}^{n} \varphi_j(x,x')\xi_j  \]
with coefficients $\varphi_j(x,x') \in C^\infty(\Omega \times
\Omega)$, such that $\grad_{x,x', \xi} \varphi(x,x', \xi) \not=
0$ for $\xi \not= 0$ and $\grad_\xi \varphi(x,x', \xi) = 0
\Leftrightarrow x = x'$. In particular, $\varphi(x,x',
\xi) = (x-x')\xi$ is an admitted choice. Then, as is well known,
also \eqref{1022.Opvarphi1606.eq} represents a
pseudo-differential operator $\Op(a; \varphi) \in
L_{(\cl)}^\mu(\Omega)$. The relation 
\begin{equation}
\label{1022.neu++1706.eq}
a(x, \xi) \to \Op(a; \varphi) 
\end{equation}
may be interpreted as an operator convention. It is
known to induce an isomorphism
\[  S_{(\cl)}^\mu(\Omega \times \R^n)/S^{- \infty}(\Omega \times
      \R^n) \to L_{(\cl)}^\mu (\Omega)/
      L^{- \infty}(\Omega).     \]
As a consequence we  have the following result:    

\begin{Theorem}
\label{1022.the1606.th}
Let \, $\varphi(x,x', \xi)$ \, and \, $\tilde{\varphi}(\tilde{x},
\tilde{x}', \tilde{\xi})$\,  be pseudo-differential 
phase functions. Then there is a map
\begin{equation}
\label{1022.SO1606.eq}
S_{(\cl)}^\mu(\Omega \times \R^n) \to S_{(\cl)}^\mu
      (\Omega \times \R^n),
\end{equation}
$a(x, \xi) \to \tilde{a}(\tilde{x}, \tilde{\xi})$, such that
\begin{equation}
\label{1022.neuLo1706.eq}  
\Op(a; \varphi) = \Op({\tilde{a}}; \tilde{\varphi})  
       \bmod L^{- \infty}(\Omega).      
\end{equation}
The relation \eqref{1022.SO1606.eq} induces an isomorphism
\begin{equation}
\label{1022.neu+1706.eq}  
S_{(\cl)}^\mu(\Omega \times \R^n)/S^{- \infty}(\Omega \times
      \R^n) \to S_{(\cl)}^\mu(\Omega \times \R^n)/S^{-
      \infty}(\Omega \times \R^n).     
\end{equation}
\end{Theorem}          

The map \eqref{1022.neu+1706.eq} incorporates a change of the
operator convention \eqref{1022.neu++1706.eq} with the phase
function $\varphi$ to the one with the phase function
$\tilde{\varphi}$. The corresponding map
\[  a(x, \xi) \to \tilde{a}(\tilde{x}, \tilde{\xi})    \]
between (left) symbols is not canonical insofar in the
preimage we may add any $c(x, \xi) \in S^{- \infty}(\Omega \times \R^n)$ and
in the image any $\tilde{c}(\tilde{x}, \tilde{\xi}) \in S^{-
\infty}(\Omega \times \R^n)$ without violating
\eqref{1022.neuLo1706.eq}. In a more precise version of such
operator conventions we may ask whether there is more control of
smoothing operators (under suitable assumptions on 
the behaviour of the phase functions near the
boundary $\partial \Omega$).

The following discussion can be subsumed under the following
question:
Let $\Omega \subset \R^n$ be an open set, let $\varphi(x,x', 
\xi) \in C^\infty(\Omega \times \Omega
\times \R^n)$ be a pseudo-differential phase function, and let
$a(x, \xi) \in S^\mu(\Omega \times \R^n)$ be a symbol.
Do there exist `natural' scales of subspaces 
${\cal H}^s(\Omega)$, $\widetilde{\cal H}^s(\Omega)$ of
$H_{\loc}^s(\Omega)$ such that 
$\Op(a; \varphi) : C_0^\infty(\Omega) \to C^\infty(\Omega)$ 
extends to a continuous operator
$\Op(a; \varphi) : {\cal H}^s(\Omega) \to \widetilde{\cal
        H}^{s- \mu}(\Omega)$
for every $s \in \R$ (or, if necessary, for certain 
specific $s$)?	

To illustrate the point let us consider the operator
\begin{equation}
\label{1022.opneu0607.eq}  
\op(a)u(t) = \iint e^{i(t-t')\tau} a(t, \tau) u(t')
          dt' \dbar \tau,      
\end{equation}
$a(t, \tau) \in S_{\cl}^\mu(\R_+ \times \R)$, first for $u \in
C_0^\infty(\R_+)$. If $a$ belongs to 
$S_{\cl}^\mu(\ol{\R}_+ \times \R) =
S_{\cl}^\mu(\R \times \R) |_{\ol{\R}_+ \times \R}$ and has 
the
transmission property at $t = 0$, there is an
extension of $\op(a)$ as a continuous map	   
\[  \op^+(a) : H^s(\R_+) \to H^{s- \mu}(\R_+)     \]
for every $s > - \frac{1}{2}$. However, if we change the phase
function, i.e., replace $\varphi(t,t', \tau) = (t-t')\tau$ by
another pseudo-differential phase function $\tilde{\varphi}
(r,r',\varrho)$, the corresponding operator
\[  \op(a; \tilde{\varphi}) : u(r) \to
    \iint e^{i \tilde{\varphi}(r,r', \varrho)}
    a(r, \varrho) u(r') dr' \dbar \varrho     \]
is not necessarily extendible in that way.

Let us now consider the case that $a(t, \tau) \in
S_{\cl}^\mu(\ol{\R}_+ \times \R)$ has not the transmission
property at $0$. Assume for the moment $\mu = 0$ and $a$
independent of $t$. Recall that the operator $\op^+(a)$ is
continuous as a map \eqref{1021.L21106.eq} but (in general) not as
\eqref{1021.1706neu.eq} for all $s$ or as a continuous
operator \eqref{1021.neue11706.eq}.
Beautiful examples are the symbols
\[  a_\pm(\tau) = \chi(\tau) \theta_\pm(\tau),    \]
cf. \eqref{1021.1706tildeneu.eq}. Observe that the operators
\begin{equation}
\label{1022.remneu2206.eq}  
\op^+((1-\chi)\theta_\pm) : L^2(\R_+) \to L^2(\R_+)    
\end{equation}
have kernels in ${\cal S}(\ol{\R}_+ \times \ol{\R}_+) (= {\cal
S}(\R \times \R) \big|_{\ol{\R}_+ \times \ol{\R}_+} )$.
Thus the essential properties of $\op^+(a_\pm)$ are reflected
by
\begin{equation}
\label{1022.a1706.eq}  
\op^+(\theta_\pm) : L^2(\R_+) \to L^2(\R_+).
\end{equation}    
The following result may be found in Eskin's book
\cite{Eski2} (see also \cite{Schu31}).

\begin{Proposition}
\label{1022.Me2206.pr}
We have {\st{(}}as an equality of continuous operators
$L^2(\R_+) \to L^2(\R_+)${\st{)}}
\begin{equation}
\label{1022.nichts2206.eq}
\op^+(\theta_\pm) = \op_M(g_\pm)
\end{equation}
for the functions
$g_+(w) := (1-e^{-2 \pi i w})^{-1}, \quad
       g_-(w) := 1 - g_+(w) = (1-e^{2 \pi i w})^{-1}$.
\end{Proposition}

In other, words the pseudo-differential operator
$\op^+(\theta_\pm)$ on $\R_+$ based on the Fourier transform
(combined with the special 
precaution at $0$ in terms of $\e^+, \r^+$) is
equivalently expressed as a Mellin pseudo-differential
operator $\op_M(g_\pm)$ (cf. the formula
\eqref{1012.neu22206.eq}) with the symbol $g_\pm
\big|_{\Gamma_{\frac{1}{2}}}$. Moreover, we have
\begin{equation}
\label{1022.chitheta2206.eq}
\op^+(\chi\theta_\pm) = \op_M(g_\pm) + G,
\end{equation}
where $G$ is an operator with kernel in ${\cal S}(\ol{\R}_+
\times \ol{\R}_+)$, cf. the remainder term
\eqref{1022.remneu2206.eq}.

\begin{Remark}
\label{re2206.re}
We have
\[  g_\pm(w) \in {\cal M}_R^0      \]
for $R = \{(j,0) \}_{j \in \Z}$ {\st{(}}in the notation of
Section {\st{10.1.2}}{\st{)}}. More precisely, we have
\[  g_+(\beta + i \varrho) \to
    \begin{cases}
      0 & \text{for $\varrho \to + \infty$}, \\
      1 & \text{for $\varrho \to - \infty$}     
    \end{cases}                 \]
for all $\beta \in \R$, uniformly in compact
$\beta$-intervals, and the converse behaviour of $g_-(\beta +
i \varrho)$.    
\end{Remark}      

\begin{Corollary}
\label{1022.co2206.co}
The operators $\op^+(\theta_\pm), \op^+(\chi \theta_\pm) :
L^2(\R_+) \to L^2(\R_+)$ restrict to continuous maps
\begin{equation}
\label{1022.optheta2206.eq}
\op^+(\theta_\pm), \op^+(\chi \theta_\pm) : {\cal S}(\ol{\R}_+)
    \to {\cal S}_P^0(\R_+)
\end{equation}
for the asymptotic type $P = \{(j,1) \}_{j \in - \N}$. Note
that a function $f \in {\cal S}_P^0(\R_+)$ has an asymptotic
expansion
\[  f(t) \sim \sum_{j=0}^{\infty} \{ c_j t^j + d_j t^j
      \log t \} \quad \text{as $t \to 0$}     \]
with constants $c_j, d_j \in \R$. Thus the operators 
\eqref{1022.optheta2206.eq} cannot be extendible to continuous
maps $H^s (\R_+) \to H^s(\R_+)$ for all $s \in \R$.
\end{Corollary}

The relation \eqref{1022.chitheta2206.eq} gives us an idea of
how the operator $\op^+(a)$ for an arbitrary $a(\tau) \in
S_{\cl}^\mu(\R)$ can be expressed as a Mellin
pseudo-differential operator on $\R_\pm$, modulo a smoothing
operator of a controlled behaviour. Let us consider the case
$\mu \in \Z$ {\st{(}}the case $a(t, \tau) \in S_{\cl}^\mu(\ol{\R}_+
\times \R)$ for arbitrary $\mu \in \R$ is treated in
\stcite{Schu31}{\st{)}}. A classical symbol $a(\tau) \in
S_{\cl}^\mu(\R)$ has an asymptotic expansion 
\begin{equation}
\label{1022.neua0607.eq}  
a(\tau) \sim \sum_{j=0}^{\infty} \chi(\tau)
      a_{\mu-j}^\pm \theta_\pm(\tau) \tau^{\mu - j}\quad \textup{for
      $\tau \to \pm \infty$}     
\end{equation}
with constants $a_{\mu-j}^\pm \in \C$. 
That means, for every $k \in\N$ there is an $N = N(k)$ such that
$\op^+ \bigl(a(\tau) - \sum_{j=0}^{N} \chi(\tau) 
        a_{\mu-j}^\pm
        \theta_\pm(\tau) \tau^{\mu - j}\bigr)$        
has a kernel in $C^k(\ol{\R}_+ \times \ol{\R}_+)$. 
Thus the essential point is to reformulate the operators
$\op^+(\chi(\tau) \theta_\pm (\tau) \tau^l)$, $l \in \Z$, by means of
the Mellin transform. For
the case $l \in \N$ we can write
\begin{equation}
\label{1022.neuM0107.eq}  \op^+(\chi(\tau)\theta_\pm(\tau) \tau^l) =
      \op^+(\tau^l) \op^+(\chi(\tau) \theta_\pm(\tau)).    
\end{equation}

In order to express $\op^+(\tau^l)$ in Mellin terms we
observe that
$\op^+(\tau) = t^{-1} i \op_M(w)$   
on $C_0^\infty(\R_+)$, i.e.,
\[  
\op^+(\tau^l) = \prod_{j=0}^{l-1} 
        (t^{-1} i \op_M(w))
      = t^{-l} i^l \op_M\bigl( \prod_{j=0}^{l-1}(w+j)\bigr). \]

In the latter formula we employed the commutation rule
$\op_M(T^1f) = t \op_M(f)t^{-1}$, with the notation 
$(T^\beta f)(w) :=  f(w+\beta)$, for a holomorphic Mellin symbol $f(w)$, e.g., a
polynomial in $w$.      

Thus, setting $h_l(w) := i^l \prod_{j=0}^{l-1} (w+j)$ for $l \in
\N$ we obtain 
\[  \op^+(\chi(\tau)\theta_\pm(\tau)\tau^l) = t^{- l}
    \op_M(h_l g_\pm) + C_l     \]
for the smoothing operator $C_l = t^{-t} \op_M(h_l)G$. For the
case $- l \in \N$ we have
\begin{align*}
\op^+(\chi(\tau) \theta_\pm(\tau) \tau^l) & =
    \bigl(t^l \op_M(h_{- l}) \bigr)^{-1} \op^+
             (\chi(\tau)\theta_\pm(\tau))
         = \bigl(\op_M(h_{-l})\bigr)^{-1} t^{-l}
       	  \op^+(\chi(\tau) \theta_\pm(\tau))        \\
   & = t^{-l} \op_M(T^l h_{-l}^{-1}) 
              \op^+(\chi(\tau)\theta_\pm(\tau))
	 = t^{-l} \op_M \bigl( (T^l h_{-l})^{-1} g_\pm \bigr)
	   + C_l
\end{align*}
for the smoothing operator $C_l = t^{-l} \op_M \bigl( (T^l
h_{-l})^{-1} \bigr) G$.	        	     

Thus the formula \eqref{1022.neua0607.eq} gives us for every $k
\in \N$ the representation
\begin{equation}
\label{1022.rep0607.eq}
\op^+(a) = \op_M(m_k) + D_k,
\end{equation}
for $m_k(t,w) := \sum_{j=0}^{k} t^{- \mu + j} f_{\mu -j}(w)$,
\[  f_{\mu-j}(w) = \{ a_{\mu -j}^+ g_+(w) + a_{\mu -j}^- g_-(w)
    \} h_{\mu -j}(w)   \] 
for $j = 0, \ldots, \mu$, $h_l(w) = i^l \prod_{j=0}^{l-1} (w+j)$,
and
\[  f_{\mu -j}(w) = \{ a_{\mu - j}^+ g_+(w) + a_{\mu -j}^-
                  g_-(w) \} (T^{j- \mu} h_{\mu -j}^{-1})(w)   
		  \]
for $j > \mu$, and $D_k$ is an operator of a controlled
behaviour, explicitly given by the considerations before. Its
kernel belongs to $C^N(\R_+ \times \R_+)$ with $N = 
N(k) \to
\infty$ as $k \to \infty$. This concerns the case $\mu \in \Z$;
as mentioned before, analogous representations for $\mu
\in \R$ may be found in \cite{Schu31}.      	      	    

Of course, the formula \eqref{1022.rep0607.eq} is not a complete
reformulation of an operator from the Fourier to the Mellin
representation, although it is a good approximation, since we
can pass to an asymptotic sum $\sum_{j=0}^{\infty} t^{- \mu + j}
f_{\mu - j}(w)$.

However, as a corollary of Theorem \ref{1022.the1606.th} we
obtain Mellin representations immediately:

\begin{Proposition}
\label{1022.Me0607.pr}
For every $a(t,x, \tau, \xi) \in S_{(\cl)}^\mu (\R_+ \times
\Omega \times \R^{1+n})$ there is an $m(r,x, w, \xi)$ $\in
S_{(\cl)}^\mu(\R_+ \times \Omega \times \Gamma_{\frac{1}{2} -
\gamma} \times \R^n)$ such that
\begin{equation}
\label{1022.opt1708.eq}  
\Op_x(\op_t(a)) = \Op_x(\op_M^\gamma(m)) \bmod 
         L^{- \infty}(\R_+ \times \Omega).     
\end{equation}
\end{Proposition}

We want to illustrate Proposition {\st{\ref{1022.Me0607.pr}}} on
the half-axis $\R_+$ {\st{(}}the generalisation to $\R_+ \times
\Omega$ is trivial{\st{)}}. Let us admit double symbols on the
Fourier as well as on the Mellin side; if necessary
pseudo-differential generalities allow us to pass to
representations in terms of left symbols.

Consider the weighted Mellin pseudo-differential operator
\begin{align}
\label{1022.opM0607.eq}
\op_M^\gamma(f)u(r) & = \int \int_{0}^{\infty}
   \Bigl(\frac{r}{r'}\Bigr)^{-(\frac{1}{2} - \gamma + i
   \varrho)} f(r,r',\frac{1}{2} - \gamma + i \varrho) u(r')
   \frac{dr'}{r'} \dbar \varrho	   \\
 & = r^{- \frac{1}{2} + \gamma} \int \int_{0}^{\infty}
     e^{i \varrho(\log r' - \log r)} f(r,r', \frac{1}{2} -
     \gamma + i \varrho)(r')^{- \frac{1}{2} - \gamma}
     u(r') dr' \dbar \varrho      \nonumber
\end{align}
for an $f(r,r', w) \in S_{(\cl)}^\mu(\R_+ \times \R_+ \times
\Gamma_{\frac{1}{2} - \gamma})$. The operator 
\[  B : v \to \int \int_{0}^{\infty}
    e^{i\tilde{\varphi}(r,r',\varrho)}f(r,r', \frac{1}{2} -
    \gamma + i \varrho) u(r')dr' \dbar \varrho   \]
is an element of $L_{(\cl)}^\mu(\R_+)$, since
\begin{equation}
\label{1022.tildevarphi67.eq}
\tilde{\varphi}(r,r', \varrho) = \varrho(\log r' - \log r)
\end{equation}
is a pseudo-differential phase function.           
This implies $r^{- \frac{1}{2} + \gamma} B r^{- \frac{1}{2} -
\gamma} \in L_{(\cl)}^\mu(\R_+)$, and, according to Theorem
\ref{1022.the1606.th}, we find a representation of
$\op_M^\gamma(f)$ by a symbol
$a(t, \tau) \in S_{(\cl)}^\mu(\R_+ \times \R) \bmod L^{-
\infty}(\R_+)$, cf. the formula \eqref{1022.opneu0607.eq}.	 

\begin{Remark}
\label{1022.di0607.re}
Consider the diffeomorphism
\[  \chi: \R_+ \to \R, \quad \chi(r) := - \log r,    \]
and set $y := - \log r$, i.e., $r = e^{- y}$. Then the operator
push forward of $\op_M^\gamma(f)$ under $\chi$ has the form
\begin{equation}
\label{1022.chi20607.eq}
\iint e^{i(y-y') \varrho} \bigl\{ e^{(\frac{1}{2} -
      \gamma)(y-y')} f(e^{- y}, e^{- y'},\frac{1}{2} - \gamma +
      i \varrho) \bigr\} v(y') dy' \dbar \varrho.
\end{equation}
\end{Remark}      

Now, since $\chi_* : L_{(\cl)}^\mu(\R_+) \to L_{(\cl)}^\mu(\R)$
is an isomorphism, for every $a(t,t', \tau) \in
S_{(\cl)}^\mu(\R_+ \times \R_+ \times \R)$ we can form a
$b(y,y', \varrho) \in S_{(\cl)}^\mu(\R \times \R \times \R)$
such that
\[  \chi_* \op_t(a) = \op_y(b) \bmod L^{- \infty}(\R).    \]
Setting
\[  f(r,r', \frac{1}{2} - \gamma + i \varrho) :=
     \bigl(\frac{r}{r'}\bigr)^{- (\frac{1}{2} - \gamma)}
     b(- \log r, - \log r', \varrho)     \]
from \eqref{1022.chi20607.eq} it follows that $\chi_* \op_t(a) =
\chi_* \op_M^\gamma(f) \bmod L^{- \infty}(\R)$ and hence
\[  \op_M^\gamma(f) = \op_t(a) \bmod L^{- \infty}(\R_+).    \] 

A similar argument applies to symbols on $\R_+ \times \Omega \ni
(t,x)$ rather than on the half-axis $\R_+$. Then, if $m(r,x, w \xi) \in
S_{(\cl)}^\mu(\R_+ \times \Omega \times \Gamma_{\frac{1}{2} -
\gamma} \times \R^{n})$ denotes a left symbol associated with
$f(r,r', x, \frac{1}{2} - \gamma + i \varrho, \xi)$ we just
obtain Proposition \ref{1022.Me0607.pr}.    

\begin{Remark}
\label{1022.re0607..re}
Similarly as \eqref{1022.sF0607.eq} we can construct a map
\[  \symb_M : L_{(\cl)}^\mu(\R_+ \times \Omega)_{\proper} \to
    S_{(\cl)}^\mu(\R_+ \times \Omega \times \Gamma_{\frac{1}{2}
    - \gamma} \times \R^n)     \]
by using the inversion formula
\[  u(r,x) = \int_{\R^n} \int_{\Gamma_{\frac{1}{2} - \gamma}}
      r^{- w} e_\xi(x) (M_\gamma F u)(w, \xi) \dbar w 
      \dbar \xi    \]
for $\dbar w = (2 \pi i)^{-1} dw$ and applying $A \in
L_{(\cl)}^\mu(\R_+ \times \Omega)_{\proper}$ under the integral      
sign. This gives us
\[  \symb_M(A)(r,x,w,\xi) = r^w e_{- \xi}(x) A r^{-w} e_\xi(.) \in
       S_{(\cl)}^\mu(\R_+ \times \Omega \times
       \Gamma_{\frac{1}{2} - \gamma} \times \R^n).     \]
For $\gamma = \frac{1}{2}$ this is, of course, equivalent 
to the
formula \eqref{1022.Axi0607.eq}, cf. also Remark 
{\st{\ref{1022.di0607.re}}}.       
\end{Remark}    

Proposition \ref{1022.Me0607.pr} and Remarks
\ref{1022.di0607.re} reformulate operators from the Fourier 
to the
Mellin representation, modulo smoothing remainders. More
interesting are  reformulations with remainders of a
controlled behaviour near $r = 0$ as obtained in the formula
\eqref{1022.rep0607.eq}. Such results are known in many
special situations, cf. the monograph \cite{Schu31} or the
papers \cite{Schu41}, \cite{Dine1}. Precise reformulations have
been mentioned before in connection with edge-degenerate
operators \eqref{1011.A1905.eq} coming from `standard'
differentiaal operators $A \in \Diff^\mu(\R^n \times \Omega)$,
$\Omega \subseteq \R^q$ open, $q \geq 0$. For $q = 0$ we
obtain Fuchs type operators of the form
\eqref{1011.Fuchs0705.eq}. 
By virtue of $- r \partial_{r} = M^{-1} w M = \op_M(w)$ ($=
\op_M^\gamma(w)$ on functions with compact support in $r \in
\R_+$) we can write \eqref{1011.A1905.eq} in the form
\begin{equation}
\label{1022.Agamma1708.eq}
A = r^{- \mu} \Op_y (\op_M^\gamma(h))
\end{equation}
for every $\gamma \in \R$ for the $(y, \eta))$ depending
Mellin symbol 
\[  h(r,y,w,\eta) = \sum_{j+|\alpha|\leq \mu} a_{j
      \alpha}(r,y)w^j(r \eta)^\alpha,    \] 
$a_{j \alpha}(r,y) \in C^\infty(\R^q,
\Diff^{\mu-(j+|\alpha|)}(X))$; in this case $X$ is a 
sphere.

For $q = 0$ the action $\Op_y(.)$ is simply to be omitted,
i.e., we have 
\begin{equation}
\label{1022.1708eqr.eq}  
A = r^{- \mu} \op_M^\gamma(h)
\end{equation}
for $h(r,w) = \sum_{j=0}^{\mu} a_j(r)w^j$, $a_j(r) \in
C^\infty(\ol{\R}_+, \Diff^{\mu-j}(X))$. 

In general, if $A$ is not a differential operator but a
classical pseudo-differential; operator in $\R^m$ we can first
consider the push forward of $A \big|_{\R^m \setminus \{ 0
\}}$ under polar coordinates $\chi : \R^m \setminus \{ 0
\} \to \R_+ \times X, x \to (r, \phi)$ (for 
$X = S^{m-1}$) to a
pseudo-differential operator with (operator-valued) symbol of
the form
\[  r^{- \mu} p(r, \varrho)    \]
such that
\[  p(r, \varrho) = \tilde{p}(r, r \varrho)    \]
and $\tilde{p}(r, \tilde{\varrho}) \in C^\infty(\ol{\R}_+,
L_{(\cl)}^\mu(X; \R_{\tilde{\varrho}}))$, and then obtain
\[  \chi_*(A \big|_{\R^m \setminus \{ 0 \} }) = r^{- \mu}
      \op_r(p) \bmod L^{- \mu}(\R^m \setminus \{ 0 \}). \]
In a second step from $p(r, \varrho)$ we produce  a
Mellin symbol $h(r,w) \in C^\infty(\ol{\R}_+, L_{\cl}^\mu(X;
\C))$ such that
\[ \op_r(p) = \op_M^\gamma(h) \bmod L^{- \infty}(\R_+ \times X).
                                     \]
Here (for any $C^\infty$ manifold $X$)
\[  L_{\cl}^\mu(X; \C \times \R^q)   \]
denotes the space of all holomorphic $L_{\cl}^\mu(X;
\R^q)$-valued functions $h(w, \eta)$ such that
$h(\beta + i \varrho, \eta) \in L_{\cl}^\mu(X;\R_\varrho
\times \R_\eta^q)$ for every $\beta \in \R$, uniformly in
compact $\beta$-intervals (for $q = 0$ we write
$L_{\cl}^\mu(X; \C)$).		

\begin{Theorem}
\label{1022.Op1708.th}
\begin{enumerate}
\item Given an arbitrary $A \in L_{(\cl)}^\mu(\R^m \times
      \Omega)$, $\Omega \subseteq \R^q$ open; the push 
      forward of 
      $A \big|_{(\R^m \setminus \{ 0 \}) \times \Omega}$ under $\chi : (x,y) \to (r,\phi, y)$ has
      the form
   \begin{equation}
   \label{1022.chi*1708.eq}
   \chi_*(A \big|_{(\R^m \setminus \{ 0 \}) \times \Omega}) =
      r^{- \mu} \Op_y(\op_r(p)) \bmod L^{- \infty}(\R_+ \times
      X \times \Omega)
   \end{equation}
    {\st{(}}for $X = S^{m-1}${\st{)}} for a family
    \begin{equation}
    \label{1022.p1708.eq}
    p(r,y,\varrho, \eta) = \tilde{p}(r,y,r \varrho, r \eta),
    \end{equation}
    \begin{equation}
    \label{1022.p'1708.eq}
    \tilde{p}(r,y,\tilde{\varrho}, \tilde{\eta}) \in
    C^\infty(\ol{\R}_+ \times \Omega \
    L_{(\cl)}^\mu(X; \R_{\tilde{\varrho}} \times
    \R_{\tilde{\eta}}^q));
    \end{equation}               		     
\item for every operator function \eqref{1022.p1708.eq} with
      \eqref{1022.p'1708.eq} for a $C^\infty$ manifold $X$
      there exists an
   \[  \tilde{h}(r,y,w, \tilde{\eta}) \in C^\infty(\ol{\R}_+
         \times \Omega, L_{(\cl)}^\mu(X; \C \times
         \R_{\tilde{\eta}}^q))       \]
      such that for
      $h(r,y,w,\eta) := \tilde{h}(r,y,w, r \eta)$    
      we have
   \begin{equation}
   \label{1022.Me1708.eq}
   \Op_y(\op_r(p)) = \Op_y(\op_M^\gamma(h)) \bmod
                L^{- \infty}(\R_+ \times X \times \Omega)  
   \end{equation}
   for every $\gamma \in \R$.
\end{enumerate}
\end{Theorem}   		   	 

\begin{Remark}
\label{1022.re1708.re}
Note that, although in the relations \eqref{1022.chi*1708.eq}
or \eqref{1022.Me1708.eq} we may have smoothing remainders the
kernels of which are not specified near $(r,r') = 0$, the
choice of $\tilde{p}$ and $\tilde{h}$ is possible in such a
way that the dependence on $r$ is smooth up to $r = 0$. In
other words, from the relation \eqref{1022.chi*1708.eq}
for every $\gamma \in \R$
we see that the control of the operator convention is much
more precise than in the general set-up of Proposition
{\st{\ref{1022.Me0607.pr}}}.
\end{Remark}

\begin{Remark}
\label{1022.1708re.re}
Observe that there is also a variant of Proposition
{\st{\ref{1022.Me0607.pr}}} for symbols $a(t, x, \tau,
\xi) \in S_{(\cl)}^\mu(\ol{\R}_+ \times \Omega \times 
\R^{1+n})$ that
always admit a choice of $m(r,x,w, \xi) \in
S_{(\cl)}^\mu(\ol{\R}_+ \times \Sigma \times \C \times \R^n)$
such that \eqref{1022.opt1708.eq} holds. Here
$S_{(\cl)}^\mu(\ol{\R}_+ \times \Sigma \times \C \times \R^n)$
is the space of all holomorphic functions $m(r,x,w,\xi)$ in
$w \in \C$ with values in $S_{(\cl)}^\mu(\ol{\R}_+ \times \Sigma
\times \R^n)$ such that
\[  m(r,x, \beta + i \varrho, \xi) \in S_{(\cl)}^\mu(\ol{\R}_+
      \times \Sigma \times \R_{\varrho, \xi}^{1+n})     \]
for every $\beta \in \R$, uniformly in compact
$\beta$-intervals.
\end{Remark}

The relation \eqref{1022.Me1708.eq} is a generalisation of
\eqref{1022.Agamma1708.eq} to pseudo-differential operators.
Intuitively it tells us that a pseudo-differential operator
$A$ on $\R^n \times \Omega$ near a fictitious edge $\Omega$
(or on $\R^m$ near the fictitious conical singularity $0$)
feels like a (weighted) Mellin operator in model cone
direction transversal to the edge (or on the cone $\R^m
\setminus \{ 0 \} \cong X^\land$ for $q = 0$, $X = S^{m-1}$).
Another interpretation is that $A$ is edge-degenerate (or of
Fuchs type) with respect to every fictitious smooth edge (or
any fictitious conical singularity). We thus see that the
smooth pseudo-differential calculus is full of `singular
confessions': 
Smooth operators belongs to the
more distinguished world of singular (or degenerate)
operators, although  they are usually not recognised as
legitimate members of that society. After this presumption we
may conject that the ambitions are going much deeper.  In
fact, as we saw at the end of Section 2.4.2, the possibilities
of smooth differential operators to pretend to be singular are
only bounded by the dimension of the underlying space. Similar
observations are true of pseudo-differential operators with
respect to higher edges and corners.

Surprisingly enough, there are not only fictitious
difficulties connected with fictitious singularities, as
explained in \cite{Dine1} or \cite{Liu2}. Even in the case of
differential operators \eqref{1011.A1905.eq} we can ask
the properties of edge symbols
\begin{equation}
\label{1022.+1708.eq}
\sigma_\land(A)(y, \eta) : {\cal K}^{s, \gamma}(X^\land) \to
    {\cal K}^{s- \mu, \gamma - \mu}(X^\land)
\end{equation}
in connection with the families of subordinate conormal
symbols
\begin{equation}
\label{1022.++1708.eq}
\sigma_\c \sigma_\land(A)(y,w) = \sum_{j=0}^{\mu}
   a_{j0}(0,y)w^j : H^s(X) \to H^{s- \mu}(X).
\end{equation}
If $A \in \Diff^\mu(\R^m \times \Omega)$ is elliptic, it is
interesting to know for which weights $\gamma \in \R$ the
operators \eqref{1022.+1708.eq} are Fredholm for all $(y, 
\eta) \in T^*
\Omega \setminus 0$. Admissible weights in that sense are
determined by the condition that the weight line
$\Gamma_{\frac{n+1}{2} - \gamma}$ does not intersect the set
of points $w \in \C$ where \eqref{1022.++1708.eq} is not
bijective, for all $y \in \Omega$. If this is the case we may
hope to find vector bundles $J_\pm$ on the edge $\Omega$ and
a block matrix family of operators
\begin{equation}
\label{1022.+++1708.eq}
\begin{pmatrix}
\sigma_\land(A) & \sigma_\land(K)  \\
\sigma_\land(T) & \sigma_\land(Q)
\end{pmatrix}   (y, \eta) :
  \Hsum{{\cal K}^{s, \gamma}(X^\land)}
       {J_{-, y}}                        \to
  \Hsum{{\cal K}^{s- \mu, \gamma - \mu}(X^\land)}
       {J_{+,y}}
\end{equation}
which fills up \eqref{1022.+1708.eq} to a family of
isomorphisms. Let $A \in \Diff^\mu(M)$ for a closed compact
$C^\infty$ manifold $M$ (of dimension $m+q$) with an embedded
closed compact manifold $Y$ (of dimension $q$) as a
fictitious edge. Then, considering the former $A \in
\Diff^\mu(\R^m \times \Omega)$ as a local representive, the
existence of isomorphisms \eqref{1022.+++1708.eq} is a global
problem and equivalent to the condition that $\ind_{S^*Y}
\sigma_\land(A) \in K(S^*Y)$ is the pull back of an element of
$K(Y)$ (namely $[J_+] - [J_-]$) under the canonical projection
$S^*Y \to Y$. This is a topological obstruction for the
existence of additional edge conditions (of trace, potential,
etc., type) which complete $A$ to a Fredholm block matrix      
\begin{equation}
\label{1022.calA1708.eq}
{\cal A} =  \begin{pmatrix}
            A & K  \\  T & Q
	    \end{pmatrix}  :
	\Hsum{{\cal W}^{s, \gamma}(M \setminus Y)}
	     {H^s(Y, J_-)}        \to
	\Hsum{{\cal W}^{s- \mu, \gamma - \mu}(M \setminus Y)}                    
             {H^{s- \mu}(Y, J_+)},
\end{equation}
cf. the discussion in Section 10.5.1 below. 
Both the evaluation of the non-bijectivity points of
\eqref{1022.++1708.eq} and of $\ind_{S^*Y} \sigma_\land(A)$
may be serious problems that are far from being  trivial in
the case of fictitious edges.

The spaces ${\cal W}^{s, \gamma}(M \setminus Y)$ are global weighted
edge spaces on $M \setminus Y$, locally near $Y$ modelled
on ${\cal W}^s(\R^q, {\cal K}^{s, \gamma}(X^\land))$.	     

As is known for $s = \gamma$, $s- \mu > \frac{m}{2}$, 
$s - \frac{m}{2} \not\in \N$, cf. \cite{Dine1}, operators 
of the kind \eqref{1022.calA1708.eq} are equivalent 
reformulations of differential operators 
\[  A : H^s(M) \to H^{s- \mu}(M)      \]
by applying suitable isomorphisms 
\[ \Hsum{{\cal W}^{s,s}(M \setminus Y)}
         {H^s(Y, J(s))} \to H^s(M), \quad 
	  H^{s- \mu}(M) \to
   \Hsum{{\cal W}^{s- \mu, s- \mu}(M \setminus Y)}
        {H^{s- \mu}(Y, J(s- \mu))}      \]
for vector bundles  $J(s), J(s- \mu) \in \Vect(Y)$.
In particular, for $\codim Y = 1$ such a block matrix 
\eqref{1022.calA1708.eq}
corresponds to a reformulation of $A$ with respect to
the subdivision of $M$ by means of $Y$. It would be
interesting to achieve similar reformulations of $A$ in terms
of subdivisions with corners, e.g., triangulations of $M$.		       

Let us return to the relation \eqref{1022.Me1708.eq},
interpreted as a local result for a pseudo-differential
operator $A \in L_{(\cl)}^\mu(M)$ on a closed compact manifold
$M$ with an embedded fictitious edge $Y$ of dimension $q$. We
then obtain the following result:

\begin{Theorem}
\label{1022.theo1808.th}
For $A \in L_{(\cl)}^\mu(M)$ and every $\gamma \in \R$ there
exists an operator $C_\gamma \in L^{- \infty}(M \setminus Y)$
such that
$A_\gamma := A - C_\gamma$
has an extension to a continuous operator
\[  A_\gamma : {\cal W}^{s,\gamma}(M \setminus Y) \to
      {\cal W}^{s- \mu, \gamma - \mu}(M \setminus Y)    \]
for every $s \in \R$.
\end{Theorem}      

This is an immediate consequence of \eqref{1022.Me1708.eq}
together with the fact that the operators $\Op_y(\op_M^{\gamma
- \frac{n}{2}}(.))$ are continuous in the weighted spaces
${\cal W}^s(\R^q, {\cal K}^{s, \gamma}(X^\land))$ (in their 
`$\comp$' or `$\loc$' versions on open sets $\Omega$ with
respect to $y$).

\subsection{The conormal cage}
\label{s.10.2.3}

%10.2.3
Let $X$ be a compact manifold with boundary $\partial X$.
By the `conormal cage' we understand the set $S^* X \cup
N^*$, explained in Remark \ref{1021.con2206.re}, consisting
of the cosphere bundle $S^* X$ as the cage and the conormal
unit intervals over the boundary as the bars. 

Consider a pseudo-differential operator
\begin{equation}
\label{1023.Aneu1108.eq}  
A := \r^+ \widetilde{A} \e^+ : L^2(X) \to L^2(X)    
\end{equation}
for an $\widetilde{A} \in L_{\cl}^0(\widetilde{X})$, where
$\widetilde{X}$ is an neighbouring $C^\infty$ manifold of 
$X$ (for instance, $2 X$). If $A$ has the transmission
property at the boundary, the homogeneous principal symbol
$\sigma_\psi(A)$ of order zero has an extension to a
continuous function $\bs{\sigma}_\psi(A)$ on $S^* X \cup N^*$
which is automatically determined by the extension of
$\sigma_\psi(A) \big|_{S^* X |_{\partial X}}$ 
from the north and south poles by homogeneity
$0$ to $N^*$. This is just the explanation of the relation
\eqref{1021.bssigma2206.eq}. If $A$ has not the transmission
property this may be not the case, cf. Remark
\ref{1021.con2206.re}.

Let us ask what we have to ensure  about the symbolic structure
of the operator $A$ when we want to associate with
\eqref{1023.Aneu1108.eq} a Fredholm boundary value problem (with
extra conditions on $\partial X$). First we require the usual
ellipticity of $A$, i.e., $\sigma_\psi(A) \not= 0$ on $T^*X
\setminus 0$. In addition, after the experience of Section
10.2.1, we have to consider the principal boundary symbol
\begin{equation}
\label{1023.sigma1108.eq}
\sigma_\partial(A)(y,\eta) = \r^+
\sigma_\psi(A)(y,0,\eta, D_t)\e^+ : L^2(\R_+) \to L^2(\R_+)
\end{equation}
(in local coordinates $x = (y, t) \in \R_+^n$, with the
covariables $\xi = (\eta, \tau)$).

In order to fill up \eqref{1023.sigma1108.eq} to a family of
isomorphisms \eqref{1021.sigma(calA)0806.eq} (here for $s = 0$)
we need that \eqref{1023.sigma1108.eq} is a family of Fredholm
operators for $(y, \eta) \in T^* (\partial X) \setminus 0$.

\begin{Theorem}
\label{1023.theo1108.th}
For the Fredholm property of \eqref{1023.sigma1108.eq} for
all $(y, \eta) \in T^*(\partial X) \setminus 0$ it is
necessary and sufficient that $\sigma_\psi(A) \big|_{S^* X
|_{\partial X}} \not= 0$ and that
\[  \sigma_\c \sigma_\partial(A)(y,w) :=
       \sigma_\psi(A)(y,0,0, +1)g_+(w) + \sigma_\psi(A)
       (y,0,0,-1)g_-(w)          \]
does not vanish for all $w \in \Gamma_{\frac{1}{2}}$ and
$y  \in \partial X$.
\end{Theorem}

This result may be found in Eskin's book \cite{Eski2}, see 
also \cite{Schu31}.       

\begin{Remark}
\label{1023.Rc1108.re}
Observe that the set
\[  \{w \in \C : w = a_+ g_+ (\frac{1}{2} + i \varrho) + a_-
       g_-(\frac{1}{2} + i \varrho), \ \varrho \in \R \}     \]
is the straight connection of the points $a_\pm \in \C$ in the
complex plane. The numbers
\[  a_\pm(y) := \sigma_\psi(A)(y,0,0, \pm 1)     \]
are the values of $\sigma_\psi(A)$ on the north and the south
pole of $S^* X |_{\partial X}$. Given any $f(y,w) \in
C^\infty(\partial X, {\cal S}(\Gamma_{\frac{1}{2}}))$ the 
points
\begin{equation}
\label{1023.co1108.eq}
a_+(y)g_+(\frac{1}{2} + i \varrho) + a_-(y)g_-(\frac{1}{2} + i
   \varrho) + f(y, \frac{1}{2} + i \varrho)
\end{equation}
define another connection between $a_+(y)$ and $a_-(y)$ in the
complex plane.
Choosing any diffeomorphism $(-1,+1) \to \Gamma_{\frac{1}{2}}$,
$\tau \to \frac{1}{2} + i \varrho$, such that $\tau \to \pm 1$
corresponds to $\varrho \to \mp \infty$ the connection
\eqref{1023.co1108.eq} can be reformulated as
\begin{equation}
\label{1023.d1108.eq}
a_+(y) g_+(\frac{1}{2} + i \varrho(\tau)) + a_-(y)
     g_-(\frac{1}{2} + i \varrho(\tau)) + f(y, \frac{1}{2} + i
     \varrho(\tau)),   
\end{equation}
$\tau \in [-1,1]$, which represents together with the values of
$\sigma_\psi(A) \big|_{S^* X |_{\partial X}}$ a continuous
function on the conormal cage $S^*X \cup N^*$.
\end{Remark}       

The function $f(y,w)$ in the relation \eqref{1023.co1108.eq} can
be regarded as a Mellin symbol of a family of operators
\begin{equation}
\label{1023.m1108.eq}
m(y, \eta) := \omega(t[ \eta ])
         \op_M(f)(y)\tilde{\omega}(t[\eta])
\end{equation}
for any choice of cut-off functions $\omega(t),
\tilde{\omega}(t)$.
For the discussion here we may take Mellin symbols 
\[  f(y,w) \in C^\infty(\partial X, {\cal M}^{- \infty}),   \]
where ${\cal M}^{- \infty}$ is the union over $\varepsilon > 0$
of all spaces ${\cal M}_\varepsilon^{- \infty} := \{ f(w) \in
{\cal A}(\{ \frac{1}{2} - \varepsilon < \re w < \frac{1}{2} +
\varepsilon \}) : f(\beta + i \varrho) \in {\cal S}
(\R_\varrho)$ for
every $\beta \in (\frac{1}{2} - \varepsilon, \frac{1}{2} +
\varepsilon)$, uniformly in compact subintervals$\}$.	 

With any such Mellin symbol we can associate an operator
\[  M : L^2(X) \to L^2(X)   \]
which is locally on $\partial X$ defined by $\Op_y 
(\op_M(f))$ and then glued together by using a partition of unity
on $\partial X$. We then set
\[  \sigma_\partial(M)(y, \eta) := \omega(t|\eta|) \op_M(f)(y)
       \tilde{\omega}(t|\eta|)     \]
and
\begin{equation}
\label{1023.sigmaN1108.eq}
\sigma_\c( A + M)(y, w) := a_+(y)g_+(\frac{1}{2} +
   i \varrho) + a_-(y) g_-(\frac{1}{2} + i \varrho) +
   f(y, \frac{1}{2} + i \varrho)
\end{equation}   
for $w = \frac{1}{2} + i \varrho$.

\begin{Definition}
\label{1023.df1208.de}
The function \eqref{1023.sigmaN1108.eq} is called the
{\st{(}}principal{\st{)}} conormal symbol of the operator $A +
M$.
\end{Definition}

\begin{Remark}
\label{1023.rea1808.re}
The notation {\st{`}}conormal symbol{\st{'}} is motivated by the 
bijection
\begin{equation}
\label{1023.nu1208.eq}
\nu : N^* \to \partial X \times \Gamma_{\frac{1}{2}}, \quad
      (y, \tau) \to (y, \frac{1}{2} + i \varrho(\tau))
\end{equation}
which admits the interpretation of $\sigma_\c(A + M)(y,w)$ for
$w \in \Gamma_{\frac{1}{2}}$ as a function on the conormal unit
interval bundle $N^*$ of the boundary $\partial X$. 
\end{Remark}

What concerns the summand $A$ the notation is compatible with
the information of the preceding section. In fact, let us write
$A$ locally in the coordinates $(y, t) \in \Omega \times
\ol{\R}_+$ in the form
\[  A = \Op_y(\r^+ \op_t(a) (y, \eta)\e^+)     \]
for a symbol $a(y,t,\eta,\tau) \in S_{\cl}^\mu(\Omega
\times \ol{\R}_+ \times \R_{\eta, \tau}^n)$. Then we know that
\[  \op^+(a)(y, \eta) = \r^+ \op_t(a)(y, \eta) \e^+   \]
admits a Mellin representation near $t = 0$ with the principal
conormal symbol
\begin{equation}
\label{1023.co1208.eq}
\sigma_\c (\op^+(a))(y, w) =
 a_{(0)}(y,0,0,1) g_+(w) +
     a_{(0)}(y,0,0,-1)g_-(w),     
\end{equation}
where $a_{(0)}(y,t,\eta,\tau)$ is the homogeneous principal component of the symbol $a$. 
For the Mellin summand $M =
\Op_{(y)}(m)$ we employ such a notation anyway, namely,
\[  
\sigma_\c(M)(y,w) = \sigma_\c(m)(y,w) =     f(y,w) , 
\]
cf. also Section 10.1.2.       
The notation `conormal symbol' of an operator, originally introduced in \cite{Remp1}, is motivated by the relationship with the conormal bundle of a boundary and the inner normal, interpreted as a manifold with conical singularity.
The boundary symbol of \eqref{1023.Aneu1108.eq} generates a function \eqref{1023.d1108.eq} (for $f=0$) on the conormal interval $[-1, +1]$, and this function has just the meaning of the conormal symbol of the operator \eqref{1023.Aneu1108.eq} when it is regarded as an element of the cone calculus on the half-axis, see \cite{Schu31}, \cite{Schu41}.

\begin{Theorem}
\label{1023.F11108.th}
The conditions $\sigma_\psi(A)\big|_{S^*X |_{\partial X}} \not= 0$ and
$\sigma_\c( A + M)(y, w) \not= 0$ for all $y \in
\partial X$ and all $w \in \Gamma_{\frac{1}{2}}$ are necessary
and sufficient for the Fredholm property of the operators
\begin{equation}
\label{1023.F01108.eq}
\sigma_\partial(A + M)(y,\eta) : L^2(\R_+) \to
     L^2(\R_+)
\end{equation}
for all $(y, \eta) \in T^* (\partial X) \setminus 0$.
\end{Theorem}

This result is  an information from \cite{Eski2}, see also \cite{Schu31}.

By virtue of the homogeneity
$\sigma_\partial(A + M)(y, \lambda \eta) =
      \kappa_\lambda \sigma_\partial( A + M)(y, \eta)
      \kappa_\lambda^{-1}$
for all $\lambda \in \R_+$ the index of the Fredholm operators
\eqref{1023.F01108.eq} is determined by the operators for $(y,
\eta) \in S^*(\partial X)$, the unit cosphere bundle induced by
$T^*(\partial X)$. The space $S^*(\partial X)$ is compact, and
we have
\begin{equation}
\label{1023.**1108.eq}
\ind_{S^*(\partial X)} \sigma_\partial( A + M) \in
        K(S^*(\partial X))      
\end{equation}
(here $K(.$) denotes the $K$ group on the space in the
brackets; recall that $K(.)$ is the group of equivalence classes
of pairs $(J,G)$ of vector bundles on that space, where $(J,G)
\sim (\widetilde{J}, \widetilde{G}) \Leftrightarrow J \oplus
\widetilde{G} \oplus H \cong \widetilde{J} \oplus G \oplus H$ for
some vector bundle $H$ (we are talking about smooth complex
vector bundles when the underlying space is a $C^\infty$
manifold, otherwise about continuous complex vector bundles).

The element \eqref{1023.**1108.eq} is represented by the
families of kernels and cokernels of the operators
$\sigma_\partial(A + M)(y, \eta)$, $(y, \eta) \in S^*(\partial
X)$, when their dimensions do not depend on $(y, \eta)$,
otherwise by an easy algebraic construction which reduces the
general case to that of constant dimensions, see, e.g., \cite[Section 2.1.7]{Schu31}.

The canonical projection $\pi_1: S^*(\partial X) \to \partial X$
gives rise to a homomorphism $\pi_1^* : K(\partial X) \to
K(S^*(\partial X))$ induced by the bundle pull back, which is
compatible with the equivalence relation.

In order to pass from the operator
\[  A + M : L^2(X) \to L^2(X)     \]
to a block matrix Fredholm operator
\begin{equation}
\label{1023.calA1108.eq}
{\cal A} = \begin{pmatrix}
            A  + M & K  \\  T & Q
	   \end{pmatrix}   :
    \Hsum{L^2(X)}     {L^2(\partial X, J_-)}    \to	   
    \Hsum{L^2(X)}     {L^2(\partial X, J_+)}
\end{equation}	 
for suitable vector bundles $J_\pm$ on $\partial X$, where $T, K$
and $Q$ are of similar meaning as the corresponding operators in
\eqref{1021.eq0906.eq}, we have to require that
\begin{equation}
\label{1023.top1108.eq}
\ind_{S^*(\partial X)} \sigma_\partial( A + M) \in
     \pi_1^* K(\partial X)
\end{equation}
which is a pseudo-differential version for a topological
obstruction for the existence of elliptic boundary value 
problems of Atiyah and Bott \cite{Atiy5}. We will come back 
to the
nature of such obstructions in a more general context in Section 5.3 below. 
A special case is the following result:

\begin{Theorem}
\label{1023.1108th.th}
Let $ A + M$ satisfy the conditions of Theorem
{\st{\ref{1023.F11108.th}}}. Then there is an elliptic boundary
value problem of the form \eqref{1023.calA1108.eq} if and only
if the relation \eqref{1023.top1108.eq} holds.
\end{Theorem}

Of course, if we talk about the extra operators $K, Q$ in
\eqref{1023.calA1108.eq} we mean that they are of a similar
structure as those in \eqref{1021.eq0906.eq}. To be more precise,
the construction follows by filling up the Fredholm family
\eqref{1023.sigma1108.eq} to a family of isomorphisms
\[  \begin{pmatrix}
    \sigma_\partial(A + M) & \sigma_\partial(K) \\
    \sigma_\partial(T)   & \sigma_\partial(Q)
    \end{pmatrix} (y, \eta) :
  \Hsum{L^2(\R_+)}     {J_{-, y}}     \to
  \Hsum{L^2(\R_+)}     {J_{+, y}},             \]
where $J_{\pm,y}$ are the fibres of vector bundles $J_\pm$ over
$y \in \partial X$. Those bundles on $\partial X$ just represent
the element \eqref{1023.top1108.eq}, i.e., 
\[  \ind_{S^*(\partial X)} \sigma_\partial(A + M) =
      [J_+] - [J_-],     \]
where $[J_+] - [J_-]$ denotes the equivalence class of $(J_+,
J_-)$.      

\begin{Remark}
\label{r.2.38}
Let us consider, more generally, operators of the form $A+M+G$ for a so called Green operator $G \in \s{L} (L^2(X))$ which is defined by $G=G_0 +G_\infty$ where
$$
G_\infty : L^2 (X) \to
   \s{W}^{\infty,0}_P (X),
      G^*_\infty : L^2 (X) \to \s{W}^{\infty,0}_Q (X)
$$
are continuous for asymptotic types
$P = \{ (p_j, m_j) \}_{j\in \N}$ and $Q = \{ (q_j, n_j) \}_{j \in \N}$
as in Section  {\em 1.2}, 
$\pi_\C P, \pi_\C Q \subset \{ \re w < \frac{1}{2} \}$, cf. Remark {\em \ref{r.2.17}}, and $G_0$ is locally in coordinates 
$(y,t) \in \R^{n-1} \times \ol{\R}_+$ of the form $\Op (g)$ for an operator-valued symbol
$g(y,\eta) \in 
    S^0_\cl (\R^{n-1} \times \R^{n-1}; L^2 (\R_+), \s{S}^0_P (\R_+))$
such that the pointwise adjoint is a symbol
 $g^*(y,\eta) \in 
    S^0_\cl (\R^{n-1}\times \R^{n-1}; L^2 (\R_+), \s{S}^0_Q (\R_+))$.
Setting 
$\sigma_\partial (G) (y,\eta) = g_{(0)} (y,\eta)$,
$(y,\eta) \in T^* (\partial X) \setminus 0$,
{\em(}cf. Definition {\em \ref{1013.Sy2005.de}} and Remark {\em \ref{1013.re0206.re}}{\em)} then we obtain a family of compact operators
$$
\sigma_\partial (G) (y,\eta) :
    L^2 (\R_+) \to L^2 (\R_+).
$$
It  follows that 
$\ind_{S^*  (\partial X)} \sigma_\partial (A+M)= 
  \ind_{S^* (\partial X)} \sigma_\partial (A+M+G)$.
\end{Remark}

The operators $G$ of Remark \ref{r.2.38} play a similar role as the Green operators in boundary value problems \eqref{1021.eq0906.eq} with the transmission property.
In the latter case the Mellin operators $M$ are not necessary to generate an operator algebra.
In the case without the transmission property (here, for simplicity, in $L^2$ spaces and of order zero) 
boundary value problems have the form of matrices
\begin{equation}
\label{eq.172}
{\cal A} = \begin{pmatrix}
           A+M+G  &  K  \\  T  &  Q
	   \end{pmatrix}.
\end{equation}
 It also makes sense to consider operators between sections in bundles $E,F \in \Vect (X)$ also in the upper left cornes, i.e., to consider operators
\begin{equation}
\label{eq.173}
{\cal A} :
    \Hsum{L^2(X,E)} 
        {L^2(\partial X, J_-)}    \to
    \Hsum{L^2(X,F)} 
        {L^2(\partial X, J_+)}.    
\end{equation}
Let $\got{V}^0(X)$ denote the space of all such operators.
We then have $\got{B}^0(X) := \got{B}^{0,0} (X) \subset \got{V}^0(X)$, cf. the notation of Section 2.1.
Similarly as \eqref{1021.sigma(calA)906.eq} the principal symbolic hierarchy has two components, namely,
\begin{equation}
\label{eq.174}
\sigma (\s{A}) = (\sigma_\psi (\s{A}), \sigma_\partial (\s{A})),
\end{equation}
with the interior symbol
\begin{equation}
\label{eq.175}
\sigma_\psi (\s{A}) :=
    \sigma_\psi (\s{A}) : \pi^*_X E \to \pi^*_X F,
\end{equation}
$\pi_X : T^* X \setminus 0 \to X$, and the boundary symbol
\begin{equation}
\label{eq.176}
\sigma_\partial ({\cal A}) :
    \pi^*_{\partial X}
    \Hsump{E' \otimes L^2 (\R_+)} 
        { J_-} \to 
    \pi^*_{\partial X}
    \Hsump{F' \otimes L^2 (\R_+)}  
        {J_+}
\end{equation}
$\pi_{\partial X} : T^* (\partial X) \setminus 0 \to \partial X$, which are bundle morphisms,
$E' := E|_{\partial X}$, $F' := F|_{\partial X}$.

\begin{Remark}
\label{r.2.39}
\begin{enumerate}
\item
The operators \eqref{eq.172} form an algebra {\em(}algebraic operations are defined when the entries of the operators fit together{\em)}.
In particular, we have
$$
\sigma (\s{AB}) = \sigma (\s{A}) \sigma (\s{B})
$$
with componentwise multiplication;	   
\item
if $\sigma (\s{A}) =0$, then \eqref{eq.173} is compact.
\end{enumerate}
\end{Remark}
\begin{Definition}
\label{d.2.40}
An operator $\s{A}$ of the form  \eqref{eq.172} is called elliptic if both components of $\sigma (\s{A})$ are isomorphisms.
\end{Definition}
\begin{Theorem}
\label{t.2.41}
An operator \eqref{eq.172} is elliptic if and only if \eqref{eq.173} is a Fredholm operator.
\end{Theorem}
Given an $\s{A} \in \got{V}^0 (X)$ and bundles $H \in \Vect (X)$, $L \in \Vect (\partial X)$ we can pass to a stabilisation of $\s{A}$ by forming a larger block matrix
$$
\wt{\s{A}} :=
 \begin{pmatrix}
    A + M +G       & 0      & K   &   0       \\
          0        & \id_H  & 0   &   0       \\
          T        & 0      & Q   &   0       \\
          0        & 0      & 0   &   \id_L   \\
\end{pmatrix}  :
  \Hsum{L^2(X, E \oplus H)}     {L^2(\partial X, J_{-} \oplus L)}     \to
  \Hsum{L^2(X, F \oplus H)}     {L^2(\partial X, J_{+} \oplus L)} 
$$
which also belongs to $\got{V}^0(X)$.        
It is evident that the ellipticity of $\s{A}$ entails the ellipticity of $\wt{\s{A}}$.
If $\s{A,B} \in \got{V}^0 (X)$ are elliptic we say that $\s{A}$ is stable homotopic to $\s{B}$, if there are stabilisations $\wt{\s{A}}$ and $\wt{\s{B}}$ of $\s{A}$ and $\s{B}$ respectively, such that there is  a continuous map $\gamma : [0,1] \to \got{V}^0(X)$ such that $\gamma(t)$ is elliptic for every $t \in [0,1]$ and $\gamma (0) = \wt{\s{A}}$,   $\gamma (1) = \wt{\s{B}}$  
(here we tacitly use a natural locally convex topology of $\got{V}^0 (X)$).
In a similar manner we can define stable equivalence of pairs of symbols of elliptic operators.

Clearly the index of an elliptic $\s{A} \in \got{V}^0 (X)$ only depends on the stable equivalence class of its principal symbols $\sigma(\s{A})$.
The space $\got{V}^0(X)$ of boundary value problems on $X$ of order $0$ 
(as well the subspace $\got{B}^0(X)$) is an example of an operator algebra with a principal symbolic hierarchy, where several components participate in the ellipticity.
It is an interesting task to understand in which way the components contribute to the index and whether and how (analytically, i.e., in terms of symbols) the contribution from one component can be shifted to another one by applying a stable homotopy through elliptic symbols.
Questions of that kind are reasonable for every operator algebra with symbolic hierarchies.
In the present case of the algebra $\got{V}^0(X)$ the picture is particularly beautiful.
First, for the subalgebra $\got{B}^0(X)$ a stable homotopy classification of elliptic principal symbols was given by Boutet de Monvel \cite{Bout1}.
The nature of homotopies depends on whether or not we admit homotopies through elliptic symbols in 
$\got{V}^0(X)$, cf. Rempel and Schulze \cite{Remp1}.
We do not give the explicit answer here, but we want to make a few remarks.
If $\s{A} \in \got{V}^0(X)$ is $(\sigma_\psi, \sigma_\partial)$-elliptic, then the upper left corner of \eqref{eq.176} is a family of Fredholm operators.
Similarly as Theorem \ref{1023.F11108.th} that means that (in the bundle case)
$\sigma_\c (A+M) (y,w)$ is a family of isomorphisms parametrised by $w \in \Gamma_{\frac{1}{2}}$, or, alternatively, when we pass to the parametrisation as in Remark \ref{1023.rea1808.re},
$\sigma_\c (A+M)(y,\eta)$ connects the isomorphisms
\begin{equation}
\label{eq.178}
\sigma_\psi (A) (y,0,0,\pm 1) : E'_y \to F'_y
\end{equation}
for $y \in Y$ by a family of isomorphisms parametrised by $N^*$.
In other words, the ellipticity of $\s{A}$ gives  rise to an isomorphism between the pull backs of $E$ and $F$ to the conormal cage $S^* X \cup N^*$ with respect to the canonical projection
$$
\pi : S^* X \cup N^* \to X.
$$
In the case $\s{A} \in \got{B}^0(X)$ the isomorphisms \eqref{eq.178}
are the same for the `plus' and the `minus' sign, and by virtue of the homogeneity of order zero the above mentioned pull back, restricted to $N^*$, is nothing other than
\begin{equation}
\label{eq.179}
\sigma_\psi (A) (y,0,0,\tau) :E'_y \to F'_y
\end{equation}
for all $-1 \leq \tau \leq 1$.
In the case of a boundary value problem $\s{A} \in \got{V}^0$ we have to replace \eqref{eq.179} by a family of the form \eqref{1023.d1108.eq} with  $f$ coming from the (in general non-trivial) Mellin symbol $f(y,w)$ which can cause a non-trivial contribution to 
$\ind_{S^*Y} \sigma_\partial (\s{A})$, cf. also \cite[Section 2.1.9]{Schu31}.

\section{How interesting are conical singularities?}
\label{s.10.3}
\begin{minipage}{\textwidth}
\setlength{\baselineskip}{0cm}
\begin{scriptsize}
An example of a cone is what is given to children in
Germany on their first day at school, a large cornet
filled with sweets.
The tip of the cone (the `conical singularity') then appears not so interesting, essential
things in this connection should be of non-vanishing
volume, while the tip is an unwelcome end.

However, if we look at a piece of material with conical
singularities (e.g., glass or iron) and observe heat flow and
tension in the body, the physical effects near the
conical points can be very important (for instance,
destruct the material). Near the tips the solutions of 
corresponding partial differential equations may be singular 
in a specific way.

The analysis in a neighbourhood of a conical singularity is a 
first
necessary step for building up calculi on configurations
with higher (`polyhedral') singularities, when we
interpret wedges as Cartesian products of cones and
$C^\infty$ manifolds, or `higher' corners as cones with
base spaces of a prescribed singular geometry.
\end{scriptsize}
\end{minipage}
     
\subsection{The iterative construction of higher singularities}
\label{s.10.3.1}

%10.3.1

Intuitively, a manifold $B$ with conical singularities is
a topological space $B$ with a (finite) subset $B'$ of 
conical
points such that $B \setminus B'$ is a $C^\infty$ manifold,
and every $v \in B'$ has a neighbourhood $V$ in $B$ that is
modelled on a cone
\begin{equation}
\label{1031.neuX1210.eq}  
X^\Delta = (\ol{\R}_+ \times X)/(\{ 0 \} \times X)   
\end{equation}
with base $X$, where $X$ is a $C^\infty$ manifold. In order to
classify different possibilities of choosing `singular charts' 
\[  \chi : V \to X^\Delta     \]
on $B$ we only admit maps of a system of
singular charts such that for any other element
$\tilde{\chi} : \widetilde{V} \to X^\Delta$ of that system
the transition map
\[  \tilde{\chi}_{\reg} \circ \chi_{\reg}^{-1} : \R_+ 
          \times X \to \R_+ \times X   \]
(for $\chi_{\reg} := \chi \big|_{V \setminus \{ v \}}$, etc.)
is smooth up to $0$, i.e., the restriction of a diffeomorphism
$\R \times X \to \R \times X$ to $\R_+ \times X$.
In this way we distinguish a conical
singularity from an infinite variety of mutually
non-equivalent cuspidal singularities.

Let us assume that $B'$ only consists of a single
point $v$; many (not all) considerations for a finite set of
conical singularities are similar to the case of one conical
singularity. 	  

The impact of a conical singularity of a space $B$ can
easily be underestimated. At the first glance we might
think that the new effects (compared with the smooth case)
in connection with ellipticity 
and other structures around the Fredholm
property of a Fuchs type operator $A$ 
are of the same size as the singularities themselves.
However, as we already saw, there is suddenly a pair
$(\sigma_\psi(A), \sigma_\c(A))$ of principal symbols, with
the conormal symbol $\sigma_\c(A)$ as a new component, 
a family of elliptic operators on the base of the cone,
and, apart from all the other remarkable things in connection
with the
pseudo-differential nature of parametrices in the conical
case, the conormal
symbol has `hidden' spectral properties, i.e.,
non-bijectivity points in the complex plane $\C \ni w$
\[  \sigma_\c(A)(w) : H^s(X) \to H^{s- \mu}(X)    \]
(and also poles in the pseudo-differential case) that are
often not explicitly known (or extremly difficult to 
detect), even in the case of fictitious conical 
singularities.

Conical singularities are important to create higher order
`polyhedral' singularities. In fact, starting from a cone
$X^\Delta$ with a smooth base $X$ we can form Cartesian products
$X^\Delta \times \Omega$ with open sets $\Omega$ in an Euclidian
space $\R^q$. A manifold $W$ with smooth edge $Y$ is then
modelled on such wedges $X^\Delta \times \Omega$ near $Y$ (with
$\Omega$ corresponding to a chart on $Y$). Similarly as for
conical singularities we impose some condition on the nature of
transition maps between local wedges. More precisely, if
$\chi : V \to X^\Delta \times \Omega, \
    \tilde{\chi} : V \to X^\Delta \times \widetilde{\Omega}$   
are two singular charts on $W$ near a point $y \in Y$, and if we
set
\[  \chi_{\reg} := \chi \big|_{V \setminus Y} : V \setminus Y
     \to  \R_+ \times X \times \Omega, \quad
    \chi_{\reg} := \tilde{\chi} \big|_{V \setminus Y} : V
      \setminus Y \to \R_+ \times X \times \widetilde{\Omega}, \]   
then the transition map
\begin{equation}
\label{1031.neut1910.eq}
\tilde{\chi}_{\reg} \circ \chi_{\reg}^{-1} : \R_+ \times X
      \times \Omega \to \R_+ \times X \times \widetilde{\Omega}
\end{equation}      
is required to be the restriction of a diffeomorphism $\R \times
X \times \Omega \to \R \times X \times \widetilde{\Omega}$ to
$\R_+ \times X \times \Omega$. This allows us to invariantly
attach $\{ 0 \} \times X \times \Omega$ to the open stretched
wedge $\R_+ \times X \times \Omega$ which gives us $\ol{\R}_+
\times X \times \Omega$, the local description of the so called
stretched manifold $\W$ with edge, associated with $W$. The
stretched manifold $\W$ is a $C^\infty$ manifold with 
boundary, and
$\partial \W$ has the structure of an $X$-bundle over $Y$.

Manifolds with edges form a category $\got{M}_1$, with natural
morphisms, especially, isomorphisms. The manifolds with conical
singularities form a subcategory (with edges of dimension $0$).

From $\got{M}_1$ we can easily pass to the category $\got{M}_2$
of manifolds with singularities of order 2, locally near the
singular subsets modelled on
\[  \textup{cones $W^\Delta$ or wedges $W^\Delta \times
       \Omega$} \]     
for a manifold $W \in \got{M}_2$ and open $\Omega \subseteq
\R^{q_2}$. This concept has been carried out in a paper of 
Calvo, Martin, and Schulze \cite{Calv2}. In other words, by 
repeatedly
forming cones and wedges we can reach caterogies of manifolds
with singularities which contain many concrete stratified 
spaces that are interesting in applications.

\begin{Remark}
\label{1031.re1210.re}
The notion `manifold' in this connection is only used for
convenience. Although there are analogues of charts, here called
singular charts, the spaces are topological manifolds only in
exceptional cases, e.g., $X^\Delta$ is a topological manifold
when $X$ is a sphere but not when $X$ is a torus.     
\end{Remark}

Observe that the category $\got{M}_k$ of spaces $M$ of
singularity order $k \in \N$ (where $k = 0$ means the $C^\infty$
case) can also be generated as follows: A space $M$ belongs to
$\got{M}_k$ if there is a submanifold $Y \in \got{M}_0$ such that
$M \setminus Y \in \got{M}_{k-1}$, and every $y \in Y$ has a
neighbourhood $V$ modelled on a wedge $X_{(k-1)}^\Delta \times 
\Omega$
for a base $X_{(k-1)} \in \got{M}_{k-1}$, $\Omega \subseteq \R^q$
open, $q = \dim Y$, with similar requirements on the transition
maps as before, cf. \cite{Calv2}. For $\dim Y = 0$ we have a
corner situation, while $\dim Y > 0$ corresponds to a higher edge. 

Setting, for the  moment $Y^{(k)} := Y$ from $M \setminus Y^{(k)} \in
\got{M}_{k-1}$ we obtain in an  analogous manner a manifold
$Y^{(k-1)} \in \got{M}_0$ such that $(M \setminus Y^{(k)})
\setminus Y^{(k-1)} \in \got{M}_{k-2}$. By iterating this
procedure we obtain a sequence of disjoint $C^\infty$ manifold
$Y^{(l)}$, $l = 0, \ldots, k$, such that $M \setminus \bigl\{
\bigcup_{j=0}^{m} Y^{(k-j)} \bigr\} \in \got{M}_{k-(m+1)}$ for
every $0 \leq m < k$, and $Y^{(0)} := M \setminus 
\bigl\{ \bigcup_{j=0}^{k-1} Y^{(k-j)} \bigr\}$.

Then we have
$M = \bigcup_{l=0}^{k} Y^{(l)}$,      
and the spaces
\[  M^{(j)} := \bigcup_{l=j}^{k} Y^{(l)} \in \got{M}_{k-l}   \]
form a sequence
\begin{equation}
\label{1031.eq1210.eq}
M =: M^{(0)} \supset M^{(1)} \supset \ldots \supset M^{(k)}
\end{equation}
such that $Y^{(j)} = M^{(j)} \setminus M^{(j+1)}$, $j = 0, \ldots,
k-1$, and $Y^{(k)} = M^{(k)}$  are $C^\infty$ manifolds. 
Those may be interpreted as smooth edges of $M$ of different dimensions. 
In particular, $Y^{(0)}$ is the $C^\infty$ part of $M$ of highest dimension.
Incidentally we call $Y^{(0)}$ the main stratum of $M$ and set $\dim M := \dim Y^{(0)}$.
Moreover, we have 
$M^{(j)} \setminus  M^{(j+1)} \in \got{M}_j$, $j=0, \ldots , k-1$.
Locally near any $y \in Y^{(j)}$ the space $M$ is modelled on a wedge
\begin{equation}
\label{eq.new}
X^\triangle_{(j-1)} \times \Omega
\end{equation}
for an open $\Omega \subseteq \R^{\dim Y^{(j)}}$ and an element 
$X_{(j-1)} \in \got{M}_{j-1}$.
\begin{Example}
\label{1031.exa1210.ex}
\begin{enumerate}
\item 
      If $M$ is a $C^\infty$ manifold with boundary, we have $M
      \in \got{M}_1$ and $M^{(1)} = \partial M$, and $Y^{(0)} = M
      \setminus \partial M$.
\item 
      A manifold $M$ with conical singularity $\{ v \}$ belongs
      to $\got{M}_1$, and we have
      $M^{(1)} = \{ v \}$ and $Y^{(0)} = M \setminus \{ v \}$.
\item 
      Let $\Omega_j$, $j = 0,1,2$, be $C^\infty$ manifolds, and set
      $M := (\Omega_0^\Delta \times \Omega_1)^\Delta \times
          \Omega_2$.      
      Then we have $M \in \got{M}_2$ and
      $Y^{(0)} = \R_+ \times \R_+ \times \Omega_0 \times 
                    \Omega_1 \times \Omega_2, \;
      Y^{(1)} = \R_+ \times \Omega_1 \times \Omega_2, \;
                  Y^{(2)} = \Omega_2$.	
\item
     Another example of a manifold with singularities is a cube $M$ in $\R^3$ with its boundary $M^{(1)}$, the system $M^{(2)}$ of one-dimensional edges including the corners, and $M^{(3)}$ the set of corner points.
In this case we have 
$M \in \got{M}_3$, $M^{(1)} \in \got{M}_2$,	$M^{(2)} \in \got{M}_1$ and
$M^{(3)} \in \got{M}_0$.  	   
\end{enumerate}
\end{Example}

For the calculus of operators on an $M \in \got{M}_k$ it is
reasonable to have a look at the space of `adequate' differential
operators. For $M \in \got{M}_0$ we simply take
$\Diff^\mu(M)$, the space of differential operators of order
$\mu$ with smooth coefficients. For $M \in \got{M}_1$ we take
$\Diff_{\deg}^\mu(M)$, defined as the subspace of all $A \in
\Diff^\mu(M \setminus Y^{(1)})$ that have in a neighbourhood of
any $y \in Y^{(1)}$ in the splitting of variables $(r, x, y) \in
\R_+ \times X_{(0)} \times \Omega$ (with $X_{(0)} \in \got{M}_0$
being the base of the local model cone near $Y^{(1)}$ and
$\Omega \subseteq \R^{\dim Y^{(1)}}$ open) the form
\begin{equation}
\label{1031.A1210.eq}
r^{- \mu} \sum_{j+|\alpha|\leq \mu}a_{j \alpha}(r,y) 
      \bigl(- r
      \frac{\partial}{\partial r}\bigr)^j(r D_y)^\alpha
\end{equation}
with coefficients $a_{j \alpha}(r, y) \in C^\infty(\ol{\R}_+
\times \Omega, \Diff^{\mu-(j+|\alpha|)}(X_{(0)}))$. 
By induction we can define
\begin{equation}
\label{1031.Ak1210.eq}
A \in \Diff_{\deg}^\mu(M)
\end{equation}
for every $M \in \got{M}_k$ as follows. On $M \setminus Y^{(k)}
\in \got{M}_{k-1}$ we assume $A \big|_{M \setminus Y^{(k)}} \in
\Diff_{\deg}^\mu(M \setminus Y^{(k)})$ which is already defined,
and in the splitting of variables $(r,x,y) \in \R_+ \times
X_{(k-1)} \times \Omega$ near any point $y \in Y^{(k)}$, 
$\Omega \subseteq \R^{\dim Y^{(k)}}$ open, $X_{(k-1)} \in \frak{M}_{k-1}$, the operator $A$ is required to be of the form \eqref{1031.A1210.eq} with coefficients
\begin{equation}
\label{eq.185}
  a_{j \alpha}(r,y) \in C^\infty(\ol{\R}_+ \times \Omega, \;
    \Diff_{\deg}^{\mu-(j+|\alpha|)}(X_{(k-1)})).    
\end{equation}   
    
The definition of \eqref{1031.Ak1210.eq} gives rise to the
notion of a principal symbolic hierarchy
\begin{equation}
\label{1031.sigma1210.eq}  
\sigma(A) := (\sigma_j(A))_{j=0, \ldots, k},
\end{equation}
where $\sigma_0(A) = \sigma_\psi(A \bigl|_{M \setminus M'}\bigr)$
is the standard homogeneous principal symbol of $A \big|_{M
\setminus M'} \in \Diff^\mu(M \setminus M')$ (recall that $M
\setminus M' \in \got{M}_0$). More generally,
$(\sigma_j(A))_{j=0, \ldots, k-1}$ is the symbol of $A \big|_{M
\setminus Y^{(k)}} \in \Diff_{\deg}^\mu(M \setminus Y^{(k)})$ in
the sense of $M \setminus Y^{(k)} \in \got{M}_{k-1}$, while 
we set
\begin{equation}
\label{1031.neusigma2610.eq}  
\sigma_k(A)(y,\eta) := r^{- \mu} \sum_{j+|\alpha| \leq \mu}
        a_{j \alpha}(0,y) \bigl(- r \frac{\partial}{\partial r}
        \bigr)^j(r \eta)^\alpha,     
\end{equation}
$(y, \eta) \in T^* Y^{(k)} \setminus 0$, as a family of
operators between functions on the model cone
$X_{(k-1)}^\land$. The nature of those functions will be
explained in Section 5.1 in more detail.

Observe that differential operators \eqref{1031.Ak1210.eq}
can be
generated in connection with Riemannian metrics. Assume that
$X_{(k-1)} \in \got{M}_{k-1}$, and let $g_{(k-1)}$ be a
Riemannian metric on $X_{(k-1)} \setminus X_{(k-1)}' 
\in \got{M}_0$. Consider the Riemannian metric
\begin{equation}
\label{1031.Rie1310.eq}
dr_k^2 + r_k^2 g_{(k-1)} + dy_k^2
\end{equation}
on the stretched wedge $\R_+ \times (X_{(k-1)} \setminus
X_{(k-1)}') \times \Omega_k$, $\Omega_k \subseteq \R^{q_k}$
open, for $r_k \in \R_+$, $y_k := (y_{k,1}, \ldots, y_{k,q_k})
\in \Omega_k$. Then the Laplace-Beltrami operator associated
with  \eqref{1031.Rie1310.eq} has the form
\begin{equation}
\label{1031.Deltak1310.eq}
r_k^{-2} \Bigl( r_k^2 \frac{\partial^2}{\partial r_k^2} +
    \Delta_{g_{(k-1)}} + r_k^2 \Delta_{\Omega_k} \Bigr),
\end{equation}
where $\Delta_{g_{(k-1)}}$ is the Laplace-Beltrami operator on
$X_{(k-1)} \setminus X_{(k-1)}'$ associated with $g_{(k-1)}$
and $\Delta_{\Omega_k} = \sum_{j=1}^{q_k}
\displaystyle \frac{\partial^2}{\partial y_{k,j}^2}$ the Laplacian on
$\Omega_k$. 
Note that
$r_k^2 \frac{\partial^2}{\partial r_k^2} = \Bigl( r_k
      \frac{\partial}{\partial r_k} \Bigr)^2 + 
         r_k \frac{\partial}{\partial r_k}$.
Assume that $g_{(k-1)}$ is given as
\begin{equation}
\label{1031.Lk1310.eq}
g_{(k-1)} := dr_{k-1}^2 + r_{k-1}^2 g_{(k-2)} + dy_{k-1}^2
\end{equation}
when $X_{(k-1)} \in \got{M}_{k-1}$ is locally modelled near an
edge point on a wedge $X_{(k-2)}^\Delta \times \Omega_{k-1}$
for an $X_{(k-2)} \in \got{M}_{k-2}$, $\Omega_{k-1} \subseteq
\R^{q_{k-1}}$ open, $(r_{k-1}, x, y_{k-1}) \in \R_+ \times
X_{(k-2)} \times \Omega_{k-1}$, and $g_{(k-2)}$ a Riemannian
metric on $X_{(k-2)} \setminus X_{(k-2)}'$. Inserting the
Laplace-Beltrami operator to \eqref{1031.Lk1310.eq} (using
notation analogous to \eqref{1031.Deltak1310.eq}) into
\eqref{1031.Lk1310.eq} it follows that the Laplace-Beltrami
operator for the Riemannian metric
\[  dr_k^2 + r_k^2 \bigl\{ dr_{k-1}^2 + r_{k-1}^2 g_{(k-2)} +
       dy_{k-1}^2 \bigr\} + dy_k^2       \]
has the form       	
\[  r_k^{-2} \bigl\{ r_k^2 \frac{\partial^2}{\partial r_k^2} +
     r_{k-1}^{-2} \bigl\{ r_{k-1}^2 \frac{\partial^2}{\partial
     r_{k-1}^2} + \Delta_{g_{(k-2)}} + r_{(k-2)}^2
     \Delta_{\Omega_{k-1}} \bigr\} + r_{k-1}^2
     \Delta_{\Omega_k} \bigr\}.    \]
By iterating this process we finally arrive at an $X_0 \in
\got{M}_0$; if we prescribe a Riemannian metric $g_0$ on $X_0$
and insert one Laplacian into the other we obtain 
an element of $\Diff_{\deg}^2(M)$ on the
singular manifold
$M := X_{(k-1)}^\Delta \times \Omega_k \in \got{M}_k$.

\subsection{Operators with sleeping parameters}
\label{s.10.3.2}

%10.3.2

The (pseudo-differential) calculus of operators on a
manifold with conical singularities or edges, locally
modelled on
\begin{equation}
\label{1032.Delta1310.eq}
\textup{cones $X^\Delta$ or wedges $X^\Delta \times \Omega$},
\end{equation}
for a (say, closed and compact) $C^\infty$ manifold $X$,
gives rise to specific operator-valued amplitude functions,
taking values in operators on $X$ and $X^\Delta$,
respectively. For instance, the calculus on the (infinite
stretched) cone $X^\Delta = \R_+ \times X \ni 
(r,x)$ starts
from Fuchs type differential operators
\begin{equation}
\label{1032.neuA1410.eq}  
A = r^{- \mu} \sum_{j=0}^{\mu} a_j(r) \left(- r \partial_r \right)^j,      
\end{equation}
$a_j(r) \in C^\infty(\ol{\R}_+,
\Diff^{\mu-j}(X))$. The operator family
$f(r,w) := \sum_{j=0}^{\mu} a_j(r) w^j$
can be regarded as an element of $C^\infty(\ol{\R}_+,
L_{\cl}^\mu(X; \Gamma_\beta))$ for every $\beta \in \R$.
Then, if we fix $\beta := \frac{n+1}{2} - \gamma$ for
$n = \dim X$, we can interpret $A$ as an operator
\[  A = r^{- \mu} \op_M^{\gamma - \frac{n}{2}}(f) : {\cal
       K}^{s, \gamma}(X^\land) \to {\cal K}^{s- \mu, 
       \gamma - \mu}(X^\land)      \]
(under suitable assumptions on the $r$-dependence of the
coefficients $a_j(r)$ for $r \to \infty$, for instance,
independence of $r$ for large $r$).       

Alternatively, we can start from an operator family $p(r,
\varrho) := \tilde{p}(r, r \varrho)$ for any
\begin{equation}
\label{1032.tildep1310.eq}
\tilde{p}(r, \tilde{\varrho}) \in C^\infty(\ol{\R}_+,
    L_{\cl}^\mu(X; \R_{\tilde{\varrho}})).
\end{equation}
In the pseudo-differential case we apply suitable
quantisations which produce operators $C_\gamma \in L^{-
\infty}(X^\land)$ such that
\begin{equation}
\label{1032.N1310.eq}
A_\gamma := r^{- \mu} \op_r(p) - C_\gamma :
    {\cal K}^{s, \gamma}(X^\land) \to {\cal K}^{s- \mu,
    \gamma - \mu}(X^\land)
\end{equation}
is continuous. Such a quantisation can be obtained by
constructing a (non-canonical) map $\tilde{p}(r,
\tilde{\varrho}) \to f(r,w)$  for an 
$f(r,w) \in C^\infty (\ol{\R}_+, L^\mu_\cl (X; \Gamma_{\frac{n+1}{2}-\gamma}))$ such that
\begin{equation}
\label{1032.M1310.eq}
\op_r(p) = \op_M^{\gamma - \frac{n}{2}}(f) \bmod
    L^{- \mu}(X^\land).
\end{equation}  

In order to find $C_\gamma$ we choose cut-off functions
$\omega(r), \tilde{\omega}(r), \tilde{\tilde{\omega}}(r)$
such that $\tilde{\omega} = 1$ on $\supp \omega$, $\omega
\equiv 1$ on $\supp \tilde{\tilde{\omega}}$. Then, using
pseudo-locality, we obtain       
\[  r^{- \mu} \op_r(p) = \omega r^{- \mu} \op_r(p)
      \tilde{\omega} + (1- \omega) r^{- \mu} \op_r(p)(1-
      \tilde{\tilde{\omega}}) + C     \]
for some $C \in L^{- \infty}(X^\land)$. Now
\eqref{1032.M1310.eq} allows us to write
\begin{equation}
\label{1032.(q)1310.eq}
r^{- \mu} \op_r(p) = \omega r^{- \mu} \op_M^{\gamma -
    \frac{n}{2}}(f) \tilde{\omega} + (1- \omega) r^{- \mu}
    \op_r(p)(1- \tilde{\tilde{\omega}}) + C_\gamma
\end{equation}
for $C_\gamma := C + \omega r^{- \mu} \bigl\{ \op_r(p) -
\op_M^{\gamma - \frac{n}{2}}(f) \bigr\} \tilde{\omega} \in
L^{- \infty}(X^\land)$. This gives us the continuity of
\eqref{1032.N1310.eq}.
More precisely, $A_\gamma :
C_0^\infty(X^\land) \to C^\infty(X^\land)$ extends by
continuity to \eqref{1032.N1310.eq} ($C_0^\infty(X^\land)$
is dense in ${\cal K}^{s, \gamma}(X^\land)$ for every $s,
\gamma \in \R$). This is remarkable, since we have
\[  \omega {\cal K}^{s, \gamma}(X^\land) = \omega r^\gamma
      {\cal K}^{s, 0}(X^\land)      \] 
for every $\gamma \in \R$ and a cut-off function $\omega(r)$,
which shows that the argument functions may have a pole at
$r = 0$ of arbitrary order when $\gamma$ is negative enough
(cf., analogously, Theorem \stref{1022.theo1808.th}).     

The process of generating operators \eqref{1032.N1310.eq} in
terms of parameter-dependent families
\eqref{1032.tildep1310.eq} can be modified by starting 
from
an edge-degenerate family $p(r,y,\varrho, \eta) :=
\tilde{p}(r,y, r \varrho, r \eta)$ for
\begin{equation}
\label{1032.tildetildep1310.eq}
\tilde{p}(r,y, \tilde{\varrho}, \tilde{\eta}) \in
    C^\infty(\ol{\R}_+ \times \Omega, L_{\cl}^\mu(X;
    \R_{\tilde{\varrho}, \tilde{\eta}}^{1+q})),
\end{equation}
$\Omega \subseteq \R^q$ open.
We can interpret \eqref{1032.tildetildep1310.eq} also as a
family \eqref{1032.tildep1310.eq} with `sleeping parameters'
$(y, \tilde{\eta}) \in \Omega \times \R^q$, while
\eqref{1032.tildep1310.eq} itself consists of an operator in
$L_{\cl}^\mu(X)$ with sleeping parameters $(r,
\tilde{\varrho}) \in \ol{\R}_+$. These are waked up in the
process of cone quantisation $r^{- \mu} \op_r(p) \to
A_\gamma$. The remaining parameters $(y, \tilde{\eta}) \in
\Omega \times \R^q$ are waked up by means of a suitable edge
quantisation.

The latter step is organised by means of a reformulation 
\eqref{1032.(q)1310.eq} 
depending on the parameters $(y, \eta)$. According to 
Theorem \stref{1022.Op1708.th} we choose an operator
function
\[  f(r,y, w, \eta) := \tilde{f}(r,y,w, r \eta)   \]
for an $\tilde{f}(r,y,w,\tilde{\eta}) \in C^\infty(\ol{\R}_+
\times \Omega, L_{\cl}^\mu(X; \Gamma_{\frac{n+1}{2} -
\gamma} \times \R_{\tilde{\eta}}^q))$ such that
\[
\op_r(p)(y, \eta) = \op_M^{\gamma - \frac{n}{2}}(f)(y,
      \eta) \bmod C^\infty(\R_+ \times \Omega, L^{-
      \infty}(X^\land; \R_\eta^q)).    \]
Then we write $r^{- \mu} \op_r(p)(y, \eta)$ in the form
\[  r^{- \mu} \op(p)(y, \eta) = A_\gamma(y, \eta) +
     C_\gamma(y,\eta)    \]      
for
\begin{align}
\label{1032.AneuAgamma.eq}  
A_\gamma(y, \eta) := & \ \omega(r[\eta])r^{- \mu}
     \op_M^{\gamma - \frac{n}{2}} (f)(y, \eta)
                 \tilde{\omega}(r[\eta]) \nonumber \\
     & + (1- \omega(r[\eta])) r^{- \mu} \op_r(p)(y, \eta)(1-
     \tilde{\tilde{\omega}}(r[\eta])).    
\end{align}  
From the construction it follows that $C_\gamma(y, \eta) \in C^\infty(\Omega, L^{-
\infty}(X^\land; \R_\eta^q))$. Recall that $\eta \to [ \eta
]$ is a strictly positive $C^\infty$ function in $\R^q$
that is equal to $| \eta |$ for large $| \eta |$. 
Now
\[  A_\gamma(y, \eta) : {\cal K}^{s, \gamma}(X^\land) \to
      {\cal K}^{s- \mu, \gamma - \mu}(X^\land)   \]      
is again a family of continuous operators for all $s \in
\R$, provided that (what we tacitly assume) the operator
family \eqref{1032.tildetildep1310.eq} has a suitable
dependence on $r$ for $r \to \infty$ (e.g., independent of
$r$ for $r > \const$). In the edge quantisation (i.e.,
quantisation near $r = 0$) it is convenient instead of
$A_\gamma(y, \eta)$ to consider
the operator function
\[  a_\gamma(y, \eta) := \sigma(r) A_\gamma(y, \eta)
      \tilde{\sigma}(r)       \]
for some cut-off functions $\sigma, \tilde{\sigma}$ which is
completely sufficient, since far from $r = 0$ our operator
on a manifold with edge should belong to the standard
calculus of pseudo-differential operators (where $\sigma,
\tilde{\sigma}$ are localising functions in connection with
a partition of unity on the respective manifold). Summing
up it follows that
\[  \sigma r^{- \mu} \op_r(p)(y, \eta) \tilde{\sigma} =
       \sigma A_\gamma(y, \eta) \tilde{\sigma} \, \bmod
       C^\infty(\Omega, L^{- \infty}(X^\land; \R^q)).   \]

\begin{Theorem}
\label{1032.the1310.th}
We have
$\sigma A_\gamma(y, \eta) \tilde{\sigma} \in 
       S^\mu(\Omega \times \R^q; {\cal K}^{s, \gamma}(X^\land),
       {\cal K}^{s- \mu, \gamma - \mu}(X^\land))$
for every $s, \gamma \in \R$ {\st{(}}cf. Definition
\stref{1013.Sy2005.de}{\st{)}}.
\end{Theorem}            

The edge-quantisation itself associated with $r^{- \mu} p(r,y,\varrho,
\eta)$ now follows by applying $\Op_y$ which gives us continuous
operators
\[  \Op_y(\sigma A_\gamma \tilde{\sigma}) : {\cal W}_{\comp}^s(\Omega,
      {\cal K}^{s, \gamma}(X^\land)) \to {\cal W}_{\loc}^{s- \mu}
      (\Omega,
      {\cal K}^{s- \mu, \gamma - \mu}(X^\land))     \]
for all $s, \mu \in \R$.

In Theorem \ref{1032.the1310.th} and the subsequent application of the
Fourier operator convention along $\Omega \ni y$ we took operators
of $L_{\cl}^\mu(X)$ with sleeping parameters $(r,y, \varrho, \eta) \in
\R_+ \times \Omega \times \R^{1+q}$, combined with a specific  rule to
activate them. By a globalisation (with a partition of unity,  etc.)
we obtain operators on a manifold $X_1 \in \got{M}_1$ 
in the sense of Section 3.1. Again we can assume that our operators contain
sleeping parameters $(r_2, y_2, \varrho_2, \eta_2) 
\in \R_+ \times
\Omega_2 \times \R^{1+q_2}$ and apply a similar scheme for the next
quantisation.

It turns out that it is advisable for such a calculus on wedges
$X_1^\Delta \times \Omega_2$ of second generation  to slightly modify
the expression for $A_{\gamma_2}(y_2, \eta_2)$ (the analogue of
\eqref{1032.AneuAgamma.eq}) by an extra localising factor in the
second summand, on the diagonal with respect to the $r_2$-variable,
cf. \cite{Calv1} and Section 5.4, formula \eqref{eq.301}. 
This shows that in the iteration of this process one
has to be careful, because the infinite cone $X_1^\Delta$ has edges
with exit to infinity for $r_2 \to \infty$. The shape of quantisations
is worth to be analysed also for other reasons, see the paper
\cite{Gil2} for alternative edge quantisations and their role for a
transparent composition behaviour of edge symbols. In other words, the
idea of introducing sleeping parameters in iterated 
quantisations for higher
calculi should be combined with other technical inventions.     

\subsection{Smoothing operators who contribute to the index}
\label{s.10.3.3}

%10.3.3

Let $M$ be a smooth compact manifold and $L_{\cl}^\mu(M)$ the
space of classical pseudo-differential operators of order $\mu$
on $M$. Moreover, let $S^{(\mu)}(T^* M \setminus 0)$  denote
the set of all $a_{(\mu)}(x, \xi) \in C^\infty(T^* M \setminus
0)$ such that $a_{(\mu)}(x, \lambda \xi) = \lambda^\mu
a_{(\mu)}(x, \xi)$ for all $\lambda > 0, (x, \xi) \in T^* M
\setminus 0$. Then we have the principal symbolic map
\[  \sigma_\psi : L_{\cl}^\mu(M) \to S^{(\mu)}(T^* M \setminus
               0).     \]
Together with the canonical embedding $L_{\cl}^{\mu-1}(M) \to
L_{\cl}^\mu(M)$ we obtain an exact sequence
\[  0 \to L_{\cl}^{\mu-1}(M) \to L_{\cl}^\mu(M) \to
          S^{(\mu)}(T^* M \setminus 0) \to 0,     \]
in particular, $L_{\cl}^{\mu-1}(M) = \ker \sigma_\psi$.	Every $A
\in L_{\cl}^\mu(M)$ induces continuous operators
\begin{equation}
\label{1033.A1410.eq}
A : H^s(M) \to H^{s- \mu}(M),
\end{equation}
and \eqref{1033.A1410.eq} is compact for $A \in
L_{\cl}^{\mu-1}(M)$. In particular, $L^{- \infty}(M) \cong
C^\infty(M \times M)$ consists of compact operators. As  we know
the ellipticity of $A$ is equivalent to the Fredholm property of
\eqref{1033.A1410.eq}, and we have
\begin{equation}
\label{1033.(ind)1410.eq}
\ind A = \ind(A + C)
\end{equation}
for every $C \in L_{\cl}^{\mu-1}(M)$. Denoting by
$L_{\cl}^\mu(M)_{\elle}$ the set of all elliptic $A \in 
L_{\cl}^\mu(M)$ and $S^{(\mu)}(T^* M \setminus 0)_{\elle} :=
\sigma_\psi L_{\cl}^\mu(M)_{\elle}$,
this relation shows that the index
\[  \ind : L_{\cl}^\mu(M)_{\elle} \to \Z     \]
can be regarded as a map
\begin{equation}
\label{1033.inda1410.eq}
\ind : S^{(\mu)}(T^* M \setminus 0)_{\elle} \to \Z.
\end{equation}
 $S^{(\mu)}(T^* M \setminus 0)_{\elle} :=
\sigma_\psi L_{\cl}^\mu(M)_{\elle}$.   	       
As is known the index only depends on stable homotopy classes of
elliptic principal symbols (the above mentioned relations are
valid in analogous form for operators acting between Sobolev
spaces of sections of smooth complex vector bundles on $M$;
the direct sum of elliptic operators is again elliptic), and the
classical Atiyah-Singer index theorem just refers to these
facts.

The phenomena completely change if the underlying manifold is
not compact. A simple example is the case
\[  M := (0,1).    \]
Taking the identity operator $A = 1$ in $H^0(M) := H^0(\R)
|_{(0,1)}$ which belongs to $L_{\cl}^0(M)$, for every $k \in \N$
we find an operator $C_k \in L^{- \infty}(M) \cap {\cal
L}(H^0(M), H^0(M))$ such that $A + C_k : H^0(M) \to H^0(M)$ is a
Fredholm operator and
\[  \ind (A + C_k) = k.     \]
We can construct such $C_k$ in the form
\[  C_k = \omega \op_M(f_k) \tilde{\omega}     \]
for a suitable Mellin symbol $f_k(z) \in S^{-
\infty}(\Gamma_{\frac{1}{2}})$, with cut-off functions $\omega,
\tilde{\omega}$ vanishing in a neighbourhood of $1$. In this
case $k$ just coincides with the winding number of the curve 
\begin{equation}
\label{1033.Lneu2010.eq}  
L := \{ w \in \C : 1 + f_k(z), \ z \in \Gamma_{\frac{1}{2}} \}
\end{equation}
under the ellipticity assumption $0 \not\in L$.

This is a very special case of operators on a manifold with
conical singularities, here the unit interval with the end
points as conical singularities. 
In other words, in the
(pseudo-differential) calculus on such a manifold we find
smoothing operators that produce any other index when added to a
Fredholm operator. Clearly we can also destroy the Fredholm
property when the first summand $A$ is Fredholm, or may achieve
it when $A$ is not Fredholm before. In the present example
this is just determined by $0 \in L$ or $0 \not\in L$. Other 
examples are elliptic
operators on  more general manifolds $B$ with conical 
singularities. If we
take, for instance, $B = X^\Delta$, with a closed compact
$C^\infty$ manifold $X$, and start from an operator 
\eqref{1032.neuA1410.eq},
\begin{equation}
\label{1033.K1410.eq}
A : {\cal K}^{s, \gamma}(X^\land) \to {\cal K}^{s- \mu, \gamma -
        \mu}(X^\land),
\end{equation}
then \eqref{1033.K1410.eq} is Fredholm if and only if it is
elliptic with respect to the components of the principal
symbolic hierarchy
\begin{equation}
\label{1033.neusigma1910.eq}  
\sigma_\gamma(A) := (\sigma_\psi(A), 
      \sigma_\c(A) \big|_{\Gamma_{\frac{n+1}{2} - \gamma}},
          \sigma_{\st E}(A)).
\end{equation}	 	
Here $\sigma_\psi(A)$ is the principal interior	symbol with
ellipticity in the Fuchs type sense; moreover, $\sigma_\c(A)(z)$
is the principal conormal symbol with ellipticity in the sense
that
\[  \sigma_\c(A)(z) : H^s(X) \to H^{s- \mu}(X)     \]
is a family of bijections for all $z \in \Gamma_{\frac{n+1}{2}-
\gamma}$ and any $s \in \R$. Finally, $\sigma_{\st E}(A)$ is the principal
exit symbol. The meaning of $\sigma_{\st E}(A)$ is as 
follows.
Consider $A$ in any subset $\R_+ \times U
\ni (r,x)$ for $r \to \infty$, with $U$ being a coordinate 
neighbourhood on
$X$. We choose a chart $\chi : \R_+ \times U \to \Gamma$ to
a conical set $\Gamma \subset \R_{\tilde{x}}^{n+1} \setminus \{
0 \}$ in such a  way that
$\chi(r,x) = r \chi_1(x)$
for  a diffeomorphism $\chi_1 : U \to V$ to an open
subset $V \subset S^n$. Then, in Euclidean coordinates 
$\tilde{x} \in
\Gamma$ (induced by $\R^{n+1}$ and  related to $(r, \varphi)$
for $\varphi = \chi_1(x)$ via polar coordinates) the symbol of
$A$ takes the form 
\begin{equation}
\label{1033.(rotp)1410.eq}
p(\tilde{x}, \tilde{\xi}) = \sum_{|\alpha|\leq \mu}
      a_\alpha(\tilde{x}) \tilde{\xi}^\alpha,  
\end{equation}  
$a_\alpha \in C^\infty(\Gamma)$. Concerning
the precise behaviour
of that symbol with respect to $\tilde{x} \not= 0$ for
$|\tilde{x}| \to \infty$, in this discussion we are completely
free to make a convenient choice.

Assume, for simplicity, $\Gamma = \R^{n+1} \setminus \{ 0 \}$;
otherwise our considerations can easily be localised in
$\Gamma$. The condition is then
\begin{equation}
\label{eq.n208}  \chi(\tilde{x}) a_\alpha(\tilde{x}) \in
        S_{\cl}^0(\R_{\tilde{x}}^{n+1})      
\end{equation}
for any excision function $\chi$ in $\R^{n+1}$, i.e.,
$\chi(\tilde{x}) \equiv 0$ for $|\tilde{x}| < c_0$,
$\chi(\tilde{x}) \equiv 1$ for $| \tilde{x}| > c_1$ for certain
$0 < c_0 < c_1$.
Let $a_{\alpha,(0)}(\tilde{x})$ denote the homogeneous
principal symbol of $a_\alpha(\tilde{x})$ of order $0$ in
$\tilde{x} \not= 0$. Then we set
\[  \sigma_{\st E}(p)(\tilde{x}, \tilde{\xi}) :=
      (\sigma_\e(p)(\tilde{x}, \tilde{\xi}), 
      \sigma_{\psi, \e}(p)(\tilde{x}, \tilde{\xi}))       \]
for
\begin{align*}  
\sigma_\e(p)(\tilde{x}, \tilde{\xi}) &:=
      \sum_{|\alpha|\leq \mu}
      a_{\alpha,(0)}(\tilde{x})\tilde{\xi}^\alpha, \
      (\tilde{x}, \tilde{\xi}) \in (\R^{n+1} \setminus \{ 0
      \}) \times \R^{n+1},    \\
\sigma_{\psi,\e}(p)(\tilde{x}, \tilde{\xi}) & :=
      \sum_{|\alpha|= \mu}
      a_{\alpha,(0)}(\tilde{x})\tilde{\xi}^\alpha, \
      (\tilde{x}, \tilde{\xi}) \in (\R^{n+1} \setminus \{ 0
      \}) \times (\R^{n+1} \setminus \{ 0 \}).
\end{align*}
This construction has an invariant meaning, first, locally on
conical sets $\Gamma, \widetilde{\Gamma} \subset \R^{n+1}
\setminus \{ 0 \}$, under transition maps $\Gamma \to
\widetilde{\Gamma}$ that are homogeneous in the variable
$|\tilde{x}|$ of order $1$, and then globally on $\R_+ \times
S^n$. This gives us a pair of functions
\[ \sigma_\e(A)(r,x, \varrho, \xi) \in 
               C^\infty(T^*(\R_+ \times X)), \
\sigma_{\psi, \e}(A)(r,x, \varrho, \xi) \in
            C^\infty(T^*(\R_+ \times X) \setminus 0)    \]
with the homogeneity properties
$\sigma_\e(A)(\lambda r, x, \varrho, \xi) = \sigma_\e(A)(r,x,
\varrho, \xi)$, 
$\sigma_{\psi, \e}(A)(r,x, \varrho, \lambda \xi) =
      \lambda^\mu \sigma_{\psi, \e}(A)(r,x, \varrho, \xi)$,  
for $\lambda > 0$ (in particular, $\sigma_\e(A)$ does not
depend on $r$ in this case).

The pair
\[  \sigma_{\st E}(A) := (\sigma_\e(A), \sigma_{\psi,
                \e}(A))\]      
is called the principal exit symbol of $A$ (of order $(\mu;
0$). Now the ellipticity of $A$ with respect to
$\sigma_{\st E}(A)$ is defined as $\sigma_\e(A) \not= 0$ and
$\sigma_{\psi, \e}(A) \not= 0$. Together with the above
mentioned ellipticity conditions with respect to
$\sigma_\psi(A)$ and $\sigma_\c(A)$ we thus obtain the
ellipticity of $A$ with respect to $\sigma_\gamma(A)$, cf. the
formula \eqref{1033.neusigma1910.eq}.		

Let $F$ be a Fr\'echet space, and $\mu \in \R$. Then
$S_{(\cl)}^\mu(\R^m; F)$ denotes the space of all (classical
or non-classical) symbols $p(\eta)$ with values in $F$, i.e.,
if $(\pi_\iota)_{\iota \in \N}$ is a semi-norm system for the
Fr\'echet topology of $F$, the condition is
$\pi_\iota(D_\eta^\alpha p) \leq c_\alpha \langle \eta
       \rangle^{\mu-|\alpha|}$
for all $\alpha \in \N^m, \iota \in \N$, and, in the classical
case, $p(\eta) \sim  \chi(\eta) \sum_{j=0}^{\infty} p_{(\mu
-j)}(\eta)$ for homogeneous components $p_{(\mu -j)}(\eta) \in 
C^\infty(\R^m
\setminus \{ 0 \}, F)$, $p_{(\mu -j)}(\lambda \eta) =
\lambda^{\mu -j}p_{(\mu -j)}(\eta)$ for all $\lambda > 0$, $j
\in \N$.    

\begin{Remark}
\label{1033.A1910.re}
The condition on local symbols \eqref{1033.(rotp)1410.eq} of
an operator
\begin{equation}
\label{1033.A1910.eq}
A = r^{- \mu} \sum_{j=0}^{\mu} a_j(r) \bigl(- r
       \frac{\partial}{\partial r} \bigr)^j,      
\end{equation}
$a_j(r) \in C^\infty(\ol{\R}_+, \Diff^{\mu -j}(X))$, $j = 0,
\ldots, \mu$, can also be formulated as
\begin{equation}
\label{1033.co1910.eq}       
a_j(r) \in S_{\cl}^0(\R, \Diff^{\mu -j}(X)) \big|_{\ol{\R}_+},
   \ j = 0, \ldots, \mu.
\end{equation}
\end{Remark}   

According to Theorem \stref{1011.Elle2304.th} the ellipticity
of $A$ is equivalent to the Fredholm property of the map
\eqref{1033.inda1410.eq}.

The operator $A$ belongs to the cone algebra of
pseudo-differential operators on $X^\land$. The cone algebra 
is motivated by the problem to express parametrices of elliptic
differential operators. It contains operators of the kind
\begin{equation}
\label{eq.new70}
  C := r^{- \mu} \omega(r) \op_M^{\gamma - \frac{n}{2}}(f)
       \tilde{\omega}(r)     
\end{equation}
with meromorphic Mellin symbols $f(w) \in {\cal M}_R^{-
\infty}(X)$ and arbitrary cut-off functions $\omega,
\tilde{\omega}$, cf. Section 10.1.2. Clearly we have $C \in
L^{- \infty}(X^\land)$, and
\[  C : {\cal K}^{s, \gamma}(X^\land) \to {\cal K}^{\infty,
        \gamma - \mu}(X^\land)         \]
is continuous for every $s \in \R$.

Similarly as before in the special case $M = (0,1)$ or $M =
\R_+$ we have the following general theorem.

\begin{Theorem}
\label{1033.theo1910.th}
Let $A$ be an operator on $X^\land$ as in Remark
\stref{1033.A1910.re} which is elliptic with respect to
$\sigma_\gamma(A)$. Then for every $k \in \N$ there exists an
$f_k \in {\cal M}_R^{- \infty}(X)$ for some discrete
asymptotic type $R$ such that, when we set $C_k = r^{- \mu}
\omega(r) \op_M^{\gamma - \frac{n}{2}}(f_k)
\tilde{\omega}(r)$, we have
\[  \ind(A + C_k) = k    \]
as a Fredholm operator \eqref{1033.neusigma1910.eq}.
\end{Theorem}

The operator $A + C_k$ belongs to the cone algebra
on $X^\land$ (with discrete asymptotics), and the Fredholm
property in  general is equivalent to the ellipticity.
In the present case the principal conormal symbol
\[  \sigma_\c(A+C_k)(w) = \sigma_\c(A)(w) + f_k(w)    \]
is elliptic with respect to the weight $\gamma$, i.e., induces
a family of isomorphisms
\[  \sigma_\c(A+C_k)(w) : H^s(X) \to H^{s- \mu}(X)    \]
for $w \in \Gamma_{\frac{n+1}{2} - \gamma}$, $s \in \R$.
This is a generalisation of the above mentioned condition $0
\not\in L$ for the curve \eqref{1033.Lneu2010.eq}

\begin{Remark}
\label{r.new70}
The phenomenon that a calculus of operators contains non-compact operators that are smoothing on the main stratum is a hint that those smoothing operators contain a hidden extra principal symbolic structure.
In the case of operators in Theorem {\em\ref{1033.theo1910.th}} this is just the conormal symbolic structure which is non-vanishing on operators of the form \eqref{eq.new70}.
Other examples are the Green operators $G$ occurring in boundary value problems \eqref{1021.calA0906.eq} which are smoothing over $\Int X$, but their boundary symbol $\sigma_\partial (G)$ may be non-trivial.
\end{Remark}

\subsection{Are cylinders the true corners?}  %the beyond}
\label{s.10.3.4}

%22.07.04 - Datsei 10.3.4 jetzt geht's los
If we want to describe analytic phenomena on a non-compact
manifold, for instance, on $\R_+$, it appears advisable to do
that in intrinsic terms, not referring to the negative counterpart $\R_-$.

Another aspect is the requirement of invariance of the
calculus under diffeomorphisms, for instance,
\begin{equation}
\label{1034.chi1910.eq}
\chi : \R_+ \to \R, \ r \to - \log r.
\end{equation}
Diffeomorphisms may destroy the geometry of the underlying space. 
For instance, \eqref{1034.chi1910.eq} transforms the conical singularity $r=0$ to infinity, while $\R_-$ disappears in the beyond.
The main results, e.g., on  boundary value
problems in, say, a half-space $\R_+^n = \R^{n-1} \times \R_+
\ni (y, r)$, with $\R_+$ being the inner normal to the boundary
$\R^{n-1}$, can certainly be transformed to results in $\R^n$
by the substitution $(y,r) \to (y, \chi(r))$, but after such a
transformation we lose a
part of the feeling for some ingredients of such problems, for
instance,  for the operator of restriction  on Sobolev spaces
$u(y,r) \in H^s(\R_+^{n})$, $s > \frac{1}{2}$, $u(y,r) \to u(y,
0)$.
Moreover, if we define pseudo-differential actions in
$H^s(\R_+^n)$ by $\r^+ A \e^+$, where $A$ is a 
pseudo-differential operator in $\R^n$, and $\e^+$ the extension by zero from $\R^n_+$ 
to $\R^n$, $\r^+$ the restriction to the half-space, it is not very natural to transport
the boundary to infinity. 
Boundary value problems may be regarded as edge problems with all the aspects of interpreting
$\ol{\R}_+^n$ as a manifold with edge $r = 0$, $\ol{\R}_+$ as
the model cone of local  wedges, and $A$ as an edge-degenerate
operator, cf. \stcite{Schu31}, \stcite{Schu41}.

Also other information is better 
located `in the finite', for instance, on the precise behaviour
of (pseudo-)differential operators in $\R^n \ni x = (x_1,
\ldots, x_n)$ with respect to a fictitious conical singularity $x
= 0$, or a fictitious edge or corner, e.g., $(x_{q+1}, \ldots,
x_n) = 0$ for some $0 \leq q < n$. 
The various cone, edge or corner 
(pseudo-) differential operators in $\R^n$ with smooth symbols 
across the singularities belong to the more exclusive clubs of
Fuchs, edge, or corner operators near those singularities.

As we saw, the new interpretation has its price: the
quantisations produce a complex degenerate behaviour of the
resulting operators, see, for instance, the formulas
\eqref{1011.Adeg1905.eq}, \eqref{1011.Vect1905.eq}.

In Fuchs degenerate operators on $\R_+ \times X$ it is natural 
to employ the
Mellin transform instead of the Fourier transform in the
axial variable $r \in \R_+$. The corresponding Mellin symbols 
$r^{- \mu} f(r, w)$ for
$f \in C^\infty(\ol{\R}_+, L_{\cl}^\mu(X; \Gamma_{\frac{n+1}{2}
- \gamma}))$, $n = \dim X$, just produce operators in the cone
calculus
\[  r^{- \mu} \omega(r) \op_M^{\gamma - \frac{n}{2}}(f)
      \tilde{\omega}(r): {\cal K}^{s, \gamma}(X^\land) \to
      {\cal K}^{s- \mu, \gamma - \mu}(X^\land)    \]
and conormal symbols
\[  f(0, w) : H^s(X) \to H^{s- \mu}(X),    \]
$w \in \Gamma_{\frac{n+1}{2} - \gamma}$. More generally, in 
edge-degenerate
situations we have parameter-dependent Mellin symbols $r^{-
\mu} f(r, y, w, r \eta)$ for $f(r,y,w, \tilde{\eta}) \in
C^\infty(\ol{\R}_+ \times \Omega, L_{\cl}^\mu(X;
\Gamma_{\frac{n+1}{2} - \gamma} \times \R_{\tilde{\eta}}^q))$
and corresponding $y$-dependent conormal symbols
\[  f(0, y, w, 0) : H^s(X) \to H^{s- \mu}(X),     \]
$y \in \Omega$, $w \in \Gamma_{\frac{n+1}{2} - \gamma}$. In such a
connection it is also common to replace the stretched cone $\R_+
\times X$ by an infinite cylinder by applying the substitution
\eqref{1034.chi1910.eq}.      	

Considering the cylinder $X^\wedge$ instead of the `true' cone 
$X^{\Delta} = (\ol{\R}_+ \times X) \setminus (\{ 0 \} \times X)$ we already make a compromise insofar 
we give up the cone `as it is'.
On the cylinder we have to restore the information by declaring certain differential operators as the `natural' ones, namely, those which generalise the shape of the Laplace-Bertrami operators belonging to conical metrics, cf. \eqref{1032.tildep1310.eq}.
Transforming such an operator to the cylinder $\R \times X$ by the substitution $t =-\log r$ we obtain an operator of the form
\begin{equation}
\label{eq.208}
e^{-t\gamma} \sum^\mu_{j=0} b_j (t) \frac{\partial^j}{\partial t^j}
\end{equation}
with coefficients
$b_j (t) \in C^\infty (\R, \Diff^{\mu-j} (X))$ of a specific behaviour for $t \to \infty$.
The smoothness of $a_j (r)$ up to $r=0$ has an equivalent reformulation in terms of a corresponding property of $b_j (t)$ up to $t = \infty$, but, as noted before, this appears less intuitive.
Moreover, on a cylinder, regarded as the original configuration, we could find quite different operators more natural, for instance, when we identify 
$\R^{n+1}_{\wt{x}} \setminus \{ 0 \}$ with 
$\R_+ \times S^n \ni (t,x) $ via polar coordinates and transform standard differential operators
$\sum_{|\alpha|\leq \mu} a_\alpha (\wt{x}) D^\alpha_{\wt{x}}$ with coefficients
$a_\alpha (\wt{x}) \in S^\nu_\cl (\R^{n+1}_{\wt{x}})$
(for some $\nu \in \R$ and, say, $a_\alpha (\wt{x}) =0$ in a neighbourhood of $\wt{x}=0$; 
cf. also \eqref{eq.n208}) into the form
\begin{equation}
\label{eq.209}
\sum^\mu_{j=0} c_j (t) \frac{\partial^j}{\partial t^j}
\end{equation}
with certain resulting coefficients
$c_j (t) \in C^\infty (\R_+, \Diff^{\mu-j} (S^n))$
(the latter vanish for $t < \varepsilon$ for some $\varepsilon >0$ and are thus identified with functions on $\R \ni t$).
It is clear that the behaviour of \eqref{eq.209} for $t \to \infty$ is fairly different from that of \eqref{eq.208} (in the case $X = S^n$).
Moreover, considering a differential operator on an infinite cylinder $\R \times X \ni (t,x)$ in general, we can assume any other behaviour for $t \to \pm \infty$.
If the crucial point are the qualitative properties of solutions $u$ of the equation $Au=f$, the answer depends on those assumptions, and different classes of operators may have nothing to do with each other.
In other words, the consideration of a  `geometric' object alone (e.g., a cylinder, or a 
\{cone\}$\setminus$\{vertex\}, or a \{compact smooth manifold with boundary\}$\setminus$\{boundary\}, 
or another non-compact manifold which is diffeomorphic to that) implies nothing on the analysis there, unless 
we do not make a specific choice of the operators.
Many non-equivalent cases may be of interest, but it happens that the terminology is  confusing.
The first step of finding out who 
is speaking about what may be to be aware of the ambiguity of notation.
Not only the scenarios around conical singularities, or cylindrical ends,  
or  boundary value problems   came to a colourful terminology, also the higher floors 
of singular contemplation produced an impressive diversity of different things under the same headlines, cooked with corner manifolds, analysis on polyhedra, multi-cylinders, etc.

Genuine geometric corners with their non-complete metrics, induced by smooth geometries of ambient spaces (e.g., cubes in $\R^3$, with all the physical phenomena, such as heat diffusion in bodies like that, or deformation and tension in models of elasticity) also live somewhere in the singular labyrinth.
Although they have a very complex character, they are not the hidden beasts but the beauties, waiting for their hero.

\section{Is \textup{`degenerate'} bad?}     
\label{s.10.4}  
\begin{minipage}{\textwidth}
\setlength{\baselineskip}{0cm}
\begin{scriptsize}    
In partial differential equations a deviation from normality is often called `degenerate'.  Similarly as in
the pseudo-differential terminology this has not a negative undertone, at least what concerns importance and
relevance to understand phenomena in natural sciences.  The notation `pseudo' in connection with ellipticity
is motivated in a similar manner as `negative' in the context of numbers.  A negative number is not
necessarily bad, it makes the life with computations a little easier, but it is also a contribution to the
symmetry of the corresponding mathematical structure.

Necessity and beauty form a unity also in problems in partial
differential equations on manifolds with singularities, and, in fact, degenerate operators  satisfy such an expectation.
\end{scriptsize}
\end{minipage}
     
\subsection{Operators on stretched spaces}   
\label{s.10.4.1}

%10.4.1

By singularities we understand what is described in Section 3.1.
The main idea was to  identify a neighbourhood of lower 
dimensional
strata by wedges $X^\Delta \times \Omega$ for a manifold $X$ of
smaller singularity order.
 A more complete characterisation is to
say that every $M \in \got{M}_k$ has a subspace $Y \in \got{M}_0$ (equal to $Y^{(k)} = M^{(k)}$ in the meaning of \eqref{1031.eq1210.eq}) 
such that a neighbourhood $U$ of $Y$ in $M$ is isomorphic to an
$X^\Delta$-bundle over $Y$ for some $X \in \frak{M}_{k-1}$ (details may be found in \stcite{Calv2}
or \stcite{Calv3}).
In this connection it is natural to employ the so-called stretched manifolds. 
From the local description of $M \in \got{M}_k$  near
$Y$ by wedges $X^\Delta \times \Omega$, $\Omega \subseteq
\R^{\dim Y}$ open, we have (as a consequence of the precise
definitions) also a local characterisation of $M \setminus Y$ `near
$Y$' by open stretched wedges $X^\land \times \Omega$, 
$X^\land
= \R_+ \times X$ and  a cocycle of transition maps
\begin{equation}
\label{1041.g2010.eq}
\R_+ \times X \times \Omega \to \R_+ \times X \times
    \widetilde{\Omega}
\end{equation}
which are isomorphisms in the category $\got{M}_{k-1}$ and
represent an $\R_+ \times X$-bundle $\L_{\reg}$ over $Y$. By assumption
the maps \eqref{1041.g2010.eq} are restrictions of
$\got{M}_{k-1}$-isomorphisms
\begin{equation}
\label{1041.tildeg2010.eq}
\R \times X \times \Omega \to \R \times X \times \Omega
\end{equation}
to $\R_+ \times X \times \Omega$. Moreover,
\eqref{1041.tildeg2010.eq} also restrict to
$\got{M}_{k-1}$-isomorphisms
\begin{equation}
\label{1041.b2010.eq}
\{ 0 \} \times X \times \Omega \to \{ 0 \} \times X \times
          \widetilde{\Omega}    
\end{equation}	  
which form a cocycle of an $X$-bundle over $Y$ that we call $\L_{\sing}$. 
By invariantly attaching 
$\{ 0 \} \times X \times \Omega$ to $\R_+ \times X \times \Omega$ 
we obtain $\ol{\R}_+ \times X \times \Omega$. 
Then \eqref{1041.g2010.eq} and
\eqref{1041.b2010.eq} together give us a cocycle of maps
\[  \ol{\R}_+ \times X \times \Omega \to \ol{\R}_+ \times X \times
           \widetilde{\Omega}    \]
which represents an $\ol{\R}_+ \times X$-bundle $\L$ over 
$Y$. Let us form the disjoint union
\[  \L = \L_{\sing} \cup \L_{\reg}.     \]
The bijection
\begin{equation}
\label{1041.U2010.eq}
U \setminus Y \cong \L_{\reg}
\end{equation}
allows us to complete $U \setminus Y$ to a stretched
neighbourhood $\U$ by forming the disjoint union
\[  \U = \M_{\sing} \cup (U \setminus Y)     \]
for an $X$-bundle $\M_{\sing}$ over $Y$ which is
$\got{M}_{k-1}$-isomorphic to $\L_{\sing}$,
\begin{equation}
\label{1041.L2010.eq}
\M_{\sing} \cong \L_{\sing},
\end{equation}
such that there is a bijection
\[   \U \cong \L      \]
that restricts to \eqref{1041.U2010.eq} and \eqref{1041.L2010.eq}.
Since $U \setminus Y \subset M$ we obtain at the same time
\[  \M := \U \cup (M \setminus U),     \]
the so-called stretched manifold associated with $M$. We then
set 
\[  \M_{\reg} := \M \setminus \M_{\sing}.      \]	   
From the construction it follows a continuous map
$\pi : \M \to M$
which restricts to the bundle projection
$\pi \big|_{\M_{\sing}} : \M_{\sing} \to Y$
and  to an $\got{M}_{k-1}$-isomorphism
$\pi \big|_{\M_{\reg}} : \M_{\reg} \to M \setminus Y$.
An example is the space $M := X^\Delta \times \Omega$ for 
$Y := \Omega\subseteq \R^q$ open, $X \in \got{M}_{k-1}$. 
In this case we have
\[  \M = \ol{\R}_+ \times X \times \Omega, \
    \M_{\sing} = \{ 0 \} \times X \times \Omega, \
    \M_{\reg} = \R_+ \times X \times \Omega.     
\]

\begin{Remark}
\label{r.4.1}
Given an $M \in \got{M}_k$ with the stretched manifold $\M$ there is the double $2\M \in \got{M}_{k-1}$ obtained by gluing together two copies of $\M$ along $\M_\sing$.
The construction of $2\M$ can be explained in local terms as
$2 (\ol{\R}_+ \times X \times \Omega) = \R \times X \times \Omega$.
\end{Remark} 
As explained in Section 3.1 a natural way of choosing  differential operators on a space $M \in \frak{M}_k$ is to locally identify a neighbourhood of a point $z \in Z := M^{(k)}$ with a wedge $X^{\Delta} \times \Xi$, to pass to the open stretched wedge $X^\wedge \times \Xi \ni (t,x,z)$ and to generate operators
\begin{equation}
\label{eq.216}
A=
  t^{- \mu} \sum_{j+|\alpha|\leq \mu}a_{j \alpha}(t,z) 
      \big(- t
      \frac{\partial}{\partial t}\big)^j(t D_z)^\alpha
\end{equation}
with coefficients
$a_{j\alpha} \in
      C^\infty (\ol{\R}_+ \times \Xi, \Diff^{\mu-(j+|\alpha|)}_{\deg} (X))$,
 taking values in an already constructed class of operators on 
$X \in \frak{M}_{k-1}$.

By definition, $X$ contains an edge $Y := X^{(k-1)} \in \frak{M}_0$ such that a neighbourhood $V$ of $Y$ in $X$ is isomorphic to a $B^\Delta$-bundle over $Y$ for some $B \in \frak{M}_{k-2}$.
Again we can fix a neighbourhood of a point $y \in Y$ modelled on a wedge $B^\Delta \times \Xi$ for some open $\Xi \subseteq \R^p$, pass to the associated open stretched wedge 
$B^\wedge \times \Omega \ni (r,b,y)$ and write the coefficients \eqref{eq.185} in the form
$$
a_{j\alpha} (t,z) =
  r^{-(\mu-(j+|\alpha|))}
     \sum_{k+|\beta| \leq \mu -(j +|\alpha|)}
         c_{j \alpha;k\beta} (r,t,y,z) (-r \frac{\partial}{\partial r})^k (rD_y)^\beta
$$
with coefficients 
$c_{j \alpha; k \beta} (r,t,y,z)$ in
\begin{equation}
\label{eq.217}
C^\infty \big(\ol{\R}_+ \times \Omega,
   \Diff^{\mu-(j+|\alpha|) -(k +|\beta|)}_{\deg} (B)\big).
\end{equation}

By inserting $a_{j\alpha} (t,z)$ into \eqref{eq.216} we obtain
\begin{eqnarray}
A &=& r^{-\mu} t^{-\mu} 
     \sum_{j+|\alpha| \leq \mu} 
                       r^{j+|\alpha|}
     \sum_{k+|\beta| \leq \mu -(j+|\alpha|)}
        c_{j \alpha; k\beta} (r,t,y,z) 
           \big( -r \frac{\partial}{\partial r} \big)^k 
              (rD_y)^\beta
           \big( -t \frac{\partial}{\partial t} \big)^j
              (tD_z)^\alpha                    \nonumber    \\
  &=& r^{-\mu} t^{-\mu} 
     \sum_{j+|\alpha|+k+|\beta| \leq \mu} 
             d_{j \alpha; k\beta} (r,t,y,z) 
           \big( -r \frac{\partial}{\partial r} \big)^k 
              (r D_y)^\beta
           \big( -rt \frac{\partial}{\partial t} \big)^j
              (rtD_z)^\alpha                                   \label{eq.deg}
\end{eqnarray}
with coefficients $d_{j \alpha; k \beta} (r,t,y,z) \in$ \eqref{eq.217}.
This process can be iterated, and gives rise to a class of degenerate operators on a corresponding `higher' stretched wedge, more precisely, on its interior which is of the form
$({\R_+})^k \times \Sigma \times \prod^k_{l=1} \Omega_l \ni 
      (r_1, \ldots , r_k, x, y_1, \ldots , y_k)$
for open sets 
$\Sigma \subseteq \R^n$, $\Omega_l \subseteq \R^{q_l}$ for
some dimensions $n, q_l$, $l=1, \ldots , k$.
In the case $k=2$ we have (with the corresponding modified notation) $B \in \frak{M}_0$, locally identified with 
$\Sigma$, while $\Omega$ and $\Xi$ correspond to $\Omega_2$ and $\Omega_1$, respectively, and the variables
$(r_1, r_2, y_1, y_2)$ to $(r,t,y,z)$.

At the end  of the chain of substitutions the operator $A$ takes the form
\begin{equation}
\label{eq.219}
A= r^{-\mu}_1  r^{-\mu}_2 \cdot \ldots \cdot r^{-\mu}_k  \wt{A} (R,V,Y),
\end{equation}
where $\wt{A}$ is a polynomial of order $\mu$ in the vector fields
\begin{equation}
\label{eq.220}
R_1= r_1 \partial_{r_1}, 
R_2= r_1 r_2  \partial_{r_2}, \ldots ,
R_k = r_1 r_2 \cdot \ldots \cdot r_k \partial_{r_k},
\end{equation}

\begin{equation}
\label{eq.221}
V_j = \partial_{x_j}, 
j=1, \ldots , n,
\ \ \mbox{where} \ \ \
x=(x_1, \ldots , x_n) \in  \Sigma,
\end{equation}

\begin{equation}
\label{eq.222}
Y_1 =(r_1 \partial_{y_{1i}})_{i=1,\ldots,q_1},
Y_2 = (r_1 r_2 \partial_{y_{2i}})_{i=1, \ldots , q_2}, \cdots,
Y_k =(r_1 r_2 \cdot \ldots \cdot r_k \partial_{y_{ki}})_{i=1, \ldots , q_k},
\end{equation}
with coefficients in 
$C^\infty ((\ol{\R}_+)^k \times \Sigma \times \prod^k_{l=1} \Omega_l)$,
$R:=(R_1, \ldots , R_k)$,   
$V=(V_1, \ldots , V_n)$,   
$Y=(Y_1, \ldots , Y_k)$.
Similar operators have been discussed before in Section 1.1, cf. Remark \ref{1011.Delta21905.re}.

The operators \eqref{eq.219} are degenerate in the sense that the coefficients at the derivatives in 
$r_l \in \R_+$ or $y_{l} \in \Omega_l$ tend to zero when $r_j \to 0$ for $1 \leq j \leq l$.
Clearly they are much more `singular' at the face $(r_1, \ldots , r_l) =0$ than those in Section 1.1, because the latter ones were obtained by repeatedly introducing polar coordinates into `smooth' operators given in an ambient space.
However, this special case shows that the class of operators of the kind \eqref{eq.219} is far of being rare, since it already  contains the operators with smooth coefficients.
In any case the operators \eqref{eq.219} have a nice shape, and they are waiting to be accepted as the new beauties of a future singular world.

Moreover, as we saw, special such operators of this category (and their pseudo-differential analogues) are a useful frame to understand the calculus of elliptic boundary value problems (especially, without the transmission property at the boundary), and these operators are accompanied by a tail of other (operator-valued) symbols which encode in this case the ellipticity of boundary conditions.
We return in Section 5.2 below once again to the aspect of symbolic hierarchies.
Let us note in this connection that, since we have to be aware of the conormal symbolic structure, it is better to consider the operators in the form
$$
A= r^{-\mu}_1 \cdot \ldots \cdot r^{-\mu}_k \wt{A} (R,Y)
$$
for a polynomial $\wt{A}$ in the vector fields \eqref{eq.220} and \eqref{eq.222} and coefficients in
$$
C^\infty ((\ol{\R}_+)^k \times
    \prod^k_{l=1} \Omega_l, \Diff^{\mu-\nu} (X))
$$
for an $X \in \frak{M}_0$, where the former $\Sigma$ plays the role of local coordinates on $X$, and $\nu$ is the number of vector fields of the type \eqref{eq.220}, \eqref{eq.222}, composed with the corresponding coefficient.

Once we have chosen the `stretched' space
$$
(\ol{\R}_+)^k \times X \times \Omega, \ \Omega \in \frak{M}_0,
$$
 as an object of interest we can interpret this again as a manifold with corner and repeat the game of passing to associated stretched spaces; this gives rise to  an infinite sequence.
This is particularly funny when we start from $\ol{\R}_+ \times \ol{\R}_+$.
The stretched space has two corners, and each stretching doubles up the number of corner points.

There are many possible choices of degenerate operators on such configurations, e.g., based on the vector fields
\begin{equation}
\label{eq.223}
r^{\lambda_l}_l \partial_{r_l}, \ \ \ l=1, \ldots , k,
\end{equation}
for certain $\lambda_l \in \R$, together with other vector fields on $X$ and $\Omega$, and also the weight factors in front of the operators can be modified.
We do not discuss such possibilities here, but we want to stress that usually the properties of degenerate operators drastically change when we change the nature of degeneracy.
In particular, when we replace the components of \eqref{eq.220} by \eqref{eq.223} (say, for the case $\lambda_l=1$, $l = 1, \ldots , k$), the resulting operators have 
a quite different behaviour than the former ones, except for $k=1$.
In other words there are many singular futures.

To return to the question in the headline of Section 4, our answer is `no'.
Although corner geometries give rise to ugly technicalities if the structure ideas remain unclear, the calculus on a singular manifold may dissolve the difficulties.

\subsection{What is \textup{`smoothness'} on a singular manifold?}
\label{s.10.4.2}

%10.4.2 10.4.2

Smoothness of a function on a manifold $M$ with singularities
$M' \subset M$, cf. the notation of Section 3.1, should mean
smoothness on the $C^\infty$ manifold $M \setminus M'$,
together with some controlled behaviour close to $M'$. For
instance, if $M = [0,1]$ is the unit interval on the real axis,
we might talk about $C^\infty$ up to the end points $\{0\}$ and
$\{1\}$. More generally, if $M$ consists of a one-dimensional
net with a system $M'$ of knots, i.e., intersection
points of finitely many intervals, (for instance, $M$ may 
be the
boundary of a triangle in the plane, or the system of
one-dimensional edges of a cube in $\R^3$, including corners) we
could ask $C^\infty$ on the intervals up to the end points and
continuity across $M'$.

The `right' notion of smoothness depends on the expectations on
the role of that property. For the analysis of (elliptic)
operators on $M$ the above mentioned notion is not convenient.

Smoothness should survive when we ask the regularity of
solutions to an elliptic equation
\begin{equation}
\label{1042.+2110.eq}
Au = f,
\end{equation}
$A \in \Diff_{\deg}^\mu(M)$, for a smooth right hand side $f$.
To illustrate a typical phenomenon we want to formulate the
following slight modification of Theorem \stref{1012.th2505.th}
which refers to the case $M = \R_+ \cup \{ 0 \} \cup \{ +
\infty \}$, with $M' = \{ 0 \} \cup \{ + \infty\}$ being 
regarded as conical singularities. Consider an operator $A$
given by \eqref{1012.Eq0705.eq}.

\begin{Theorem}
\label{1042.reg2110.th}
Let $A$ be elliptic with respect to $(\sigma_\psi(A)$,
$\sigma_\c(A) \big|_{\Gamma_{\frac{1}{2}} - \gamma})$ for some
weight $\gamma \in \R$, i.e., $a_\mu \not= 0$ and
$\sigma_\c(A)(w) \not= 0$ on $\Gamma_{\frac{1}{2} - \gamma}$.
Then for every $f \in L_{Q_\gamma^0, Q_\gamma^\infty}^{2,
\gamma}(\R_+)
\cap C^\infty(\R_+)$ the equation \eqref{1042.+2110.eq} has a
unique solution
\[  u \in L_{P_\gamma^0, P_\gamma^\infty}^{2, \gamma}(\R_+)
       \cap C^\infty(\R_+)     \]
for every pair $(Q_\gamma^0, Q_\gamma^\infty)$ of discrete
asymptotic types with some resulting $(P_\gamma^0,
P_\gamma^\infty)$.
\end{Theorem}       

In other words, elliptic regularity in the frame of smooth
functions has three aspects:
\begin{align}
\label{1042.s2110.eq}
&\textup{standard smoothness on $M \setminus M'$}, \\
\label{1042.w2110.eq}
&\textup{weighted properties close to $M'$},  \\
\label{1942.as2110.eq}
&\textup{asymptotic properties close to $M'$}.
\end{align}

The individual weighted and asymptotic properties are
determined by the calculus of elliptic operators that we
choose on $M$. 
There are several choices, as we shall see 
by the following Theorem \stref{1042.the2110.th} and
Remark \stref{1042.re2110.re}.

\begin{Theorem}
\label{1042.the2110.th}
Let $A$ be an operator on $X^\land$ as in Section 
{\st{3.3}} that
is elliptic with respect to the principal symbol
\eqref{1033.neusigma1910.eq}. Then
\[  A u = f \in {\cal K}_Q^{s- \mu, \gamma - \mu}(X^\land)
      \]
and $u \in {\cal K}^{- \infty, \gamma}(X^\land)$ implies
\[  u \in {\cal K}_P^{s, \gamma}(X^\land)     \]
for every discrete asymptotic type $Q$ with some
resulting $P$, for every $s \in \R$ 
{\st{(}}in particular, this
also holds for $s = \infty${\st{)}}. Moreover,
\[  A u = f \in {\cal S}_Q^{\gamma - \mu}(X^\land)    \]
and $u \in {\cal K}^{- \infty, \gamma}(X^\land)$ entail
\[  u \in {\cal S}_P^\gamma(X^\land).     \]
\end{Theorem}      

\begin{Remark}
\label{1042.re2110.re}
There is also a theorem on elliptic  
regularity for operators \eqref{1033.A1910.eq} when the 
behaviour of coefficients for $r \to \infty$ is analogous
to that for $r \to 0$, namely,
$a_j(r^{\pm 1}) \in C^\infty(\ol{\R}_+, \Diff^{\mu -
        j}(X))$.     
In that case $r = 0$ and $r = \infty$ are treated as
conical singularities, cf. also Remark {\em\ref{r.1.22}}. 

If $A$ is elliptic in the
sense $\sigma_\psi(A) \not= 0$, and $\sigma_\psi(A)$ as
well as $\sigma_\psi(I^{- n} A I^n)$ are elliptic up to $r
= 0$ {\st{(}}in the Fuchs type sense{\st{)}} and
$\sigma_\c(A)(w) \big|_{\Gamma_{\frac{n+1}{2} - \gamma}}:
H^s(X) \to H^{s- \mu}(X)$ is a family of isomorphisms,
then
$A u = f \in {\cal H}_{Q_\gamma^0, Q_\gamma^\infty}^{s-
      \mu, \gamma - \mu}(X^\land)$
and $u \in {\cal H}^{- \infty, \gamma}(X^\land)$ imply
$u \in {\cal H}_{P_\gamma^0, P_\gamma^\infty}^{s,
        \gamma}(X^\land )$     
for every pair of discrete asymptotic types $(Q_\gamma^0,
Q_\gamma^\infty)$ with some resulting 
$(P_\gamma^0, P_\gamma^\infty)$.
Especially, for $s=\infty$ we see which kind of smoothness survives under the process of solving an elliptic equation.
\end{Remark}	
Once we arrived at the point to call a function smooth on a manifold $M$ with singularities $M'$ when $u$ is smooth on $M \setminus M'$ and of a similar qualitative behaviour near $M'$ as a solution of an elliptic equation (belonging to the calculus adapted  to $M$) we have a candidate  of a definition also for manifolds with edges and corners.
In Section 4.5 below we give an impression on the general asymptotic behaviour of a solution near a smooth edge.
The variety of possible `asymptotic configurations' in this case is overwhelming, and it is left to the individual feeling of the reader to see in this behaviour the opened door to an asymptotic hell or to a spectral paradise.

The functional analytic description of corner asymptotics for the singularity order $k \geq 2$ is another non-trivial part of the story.
For instance, if an edge has conical singularities (which corresponds to the case $k=2$) we have to expect asymptotics in different axial directions $(r_1, r_2) \in \R_+ \times \R_+$ near $r_1=0$ and $r_2=0$, and the description of the interaction of both contributions near the corner point $r_1 = r_2 =0$ requires corresponding inventions in terms of weighted distributions with asymptotics (especially, when the Sobolev smoothness $s$ is finite).
The asymptotics of solutions of elliptic equations in such corner situations has been investigated from different point of view in \cite{Schu30}, \cite{Schu10}, \cite{Schu27}, \cite{Krai3}. 

\subsection{Schwartz kernels and Green operators}
\label{s.10.4.3}

%10.4.3 hallo

The notation `Green operators' in cone and edge calculi is
derived from Green's function of boundary value problems. In
the most classical context we have Green's function of the
Dirichlet problem
\begin{equation}
\label{1043.D2210.eq}
\Delta u = f \;\textup{in $\Omega$,\quad $T u = g$ on $\partial
   \Omega$}, 
\end{equation}   
$T u := u |_{\partial \Omega}$, in a bounded smooth domain in
$\R^n$. For convenience we assume for the moment $f \in
C^\infty(\ol{\Omega})$, $g \in C^\infty(\partial \Omega)$. The
problem \eqref{1043.D2210.eq} has a unique solution $u \in
C^\infty(\ol{\Omega})$ of the form
\[  u = P f + K g.     \]
Here
$P : C^\infty(\ol{\Omega}) \to C^\infty(\ol{\Omega})$   
just represents Green's function of the Dirichlet problem, 
and
$K : C^\infty(\partial \Omega) \to C^\infty(\ol{\Omega})$
is a potential operator.

The operator $P$ is a parametrix of $\Delta$ in $\Omega$.
Every fundamental solution $E$ of $\Delta$ is a parametrix,
too. Thus the operator 
\[  P - E =: G    \]
has a kernel in $C^\infty(\Omega \times \Omega)$. The operator
$G$ is a Green operator in the sense of our notation. The
operator
\[  (P \quad K) = \binom{\Delta}{T}^{-1}   \]
belongs to the pseudo-differential calculus of boundary value
problems with the transmission property at the boundary.

From that calculus we know some very remarkable relations. Near
the boundary in local coordinates $x = (y, t) \in \Omega
\times \ol{\R}_+$, $\Omega \subseteq \R^q$ open, $q = n-1$,
the operator $G$ has the form
\[  G = \Op_y(g) + C    \]
for a symbol $g(y,\eta)$ such that
\begin{equation}      
\label{1043.g*2210.eq}
g(y, \eta), g^*(y,\eta) \in S_{\cl}^\mu(\Omega \times \R^q; \;
     L^2(\R_+), {\cal S}(\ol{\R}_+))
\end{equation}
for $\mu = -2$, cf. Example \stref{1013.G0206.exa}. 
The operator $C$ is   smoothing in the
calculus of boundary value problems (in this case with a
kernel in $C^\infty(\Omega \times \ol{\R}_+ \times \Omega
\times \ol{\R}_+)$). The structure of Green symbols $g(y,
\eta)$ is closely related to the nature of elliptic 
regularity of the homogeneous principal boundary symbol
\begin{equation}
\label{1043.R2210.eq}
\sigma_\partial(\Delta)(\eta) = - |\eta|^2 +
    \frac{\partial^2}{\partial t^2} : H^s(\R_+) \to
    H^{s-2}(\R_+)
\end{equation}
for $\eta \not= 0$, $s > \frac{3}{2}$, which is an operator
with the transmission property at $t = 0$, elliptic  as usual
in the finite (up to $t = 0$) and exit elliptic for $t \to
\infty$. Thus
\[  \sigma_\partial(\Delta)(\eta) u(t) = f(t) \in {\cal
        S}(\ol{\R}_+)          \]
implies $u(t) \in {\cal S}(\ol{\R}_+)$. In particular, we have
smoothness at $t = 0$ and the Schwartz property for $t \to
\infty$.	

As noted in Example \stref{1013.G0206.exa} the Green symbols
act as operators
\[  g(y, \eta)u(t) = \int_{0}^{\infty} f_G(t[\eta], t'[\eta];
            y, \eta)u(t')dt'     \]
for a function
$f_G(t,t; y, \eta) \in {\cal S}(\ol{\R}_+ \times
              \ol{\R}_+)$	    
for every fixed $y, \eta$. In addition, the property
\eqref{1043.g*2210.eq} reflects remarkable rescaling
properties, hidden in Green operators, here encoded by the
twisted homogeneity of the components of the corresponding classical symbol.

The question is now whether this behaviour is an accident, or a
typical phenomenon with a more general background. The answer
should be contained in the pseudo-differen\-tial algebras on
manifolds with singularities. Although many details on 
the higher singular
algebras are projects for the future, the expectation is as
follows. If $M \in \got{M}_k$ is a manifold with 
singularities of order $k$,
and $Y \in \got{M}_0$ such that $M \setminus Y \in
\got{M}_{k-1}$, with $Y$ being a corresponding higher edge,
then $Y$ has a neighbourhood $U$ in $M$ which is
$\got{M}_k$-isomorphic to an $X_{k-1}^\Delta$-bundle over $Y$
for an $X_{k-1} \in \got{M}_{k-1}$. This gives rise to an
axial variable $r_k \in \R_+$ of the cone $X_{k-1}^\Delta$,
or, if convenient, of the open stretched cone $X_{k-1}^\land =
\R_+ \times X_{k-1}$. Then, locally on $Y$, we can construct 
Green symbols $g(y, \eta)$ that are classical in the
covariables and take values in operators on weighted cone
Sobolev spaces where 
\begin{align}
\label{1043.k2210.eq}
g(y, \eta) & : {\cal K}^{s, \gamma}(X_{k-1}^\land) \to
           {\cal S}_P^{\gamma - \mu}(X_{k-1}^\land), \\
\label{1043.k*2210.eq}	   	      
g^*(y, \eta) & : {\cal K}^{s,-\gamma + \mu}(X_{k-1}^\land) \to
           {\cal S}_Q^{- \gamma}(X_{k-1}^\land).
\end{align}	   
Here $\gamma = (\gamma_1, \ldots, \gamma_k) \in \R^k$ is a
tuple of weights, where 
$\gamma - \mu := (\gamma_1 - \mu, \ldots, \gamma_k - \mu)$, and
\begin{equation}
\label{1043.calS2210.eq}   
{\cal S}_P^\varrho(X^\land)    
\end{equation}
are analogues of the spaces \eqref{eq.n56} with  `higher' asymptotic types $P$ that encode a
specific asymptotic behaviour for $r_k \to 0$. 
In Section 6.3 below we shall deepen the insight on the nature of
higher ${\cal K}^{s, \gamma}$- and ${\cal S}_P^\gamma$-spaces.

The mappings \eqref{1043.k2210.eq} are a generalisation of
\begin{equation}
\label{1043.f2210.eq}
g(y, \eta) : L^2(\R_+) \to {\cal S}(\ol{\R}_+),
\end{equation}
cf. \eqref{1043.g*2210.eq}. In \eqref{1043.f2210.eq} the
asymptotic  type $P$ means nothing other than smoothness up 
to
$t = 0$ (Taylor asymptotics). The Schwartz property at
infinity is typical in Green symbols. It comes from the role
of Green edge symbols to adjust operator families in the  
full calculus of homogeneous edge symbols by smoothing 
elements (in
particular, in the elliptic case in connection with kernels
and cokernels), taking into account that the calculus  
treats $t
\to \infty$ as an exit to infinity (for $\eta \not= 0$). 
In such
calculi the remainders near $\infty$ have Schwartz kernels.

Another interesting aspect on kernels is their 
behaviour near $t = 0$.

In the preceding section we tried to give an impression on the enormous variety of different asymptotic phenomena which may occur in smooth functions on a manifold with singularities.
In \eqref{1043.calS2210.eq} this is summarised under the notation `$P$'; it encodes not only asymptotic information at the tip of the corner with base $X$ but on all the edges of different dimension, generated by the singularities of $X$.
In particular, with such infinite edges of $X^\wedge$ also the asymptotic information is travelling to $\infty$, i.e., to the conical exit of $X^\wedge$ for $t \to \infty$.

Smoothness and asymptotics are not only an aspect of Green symbols but also of the global smoothing operators on a manifold $M$ with singularities $M'$ which are usually regarded as the simplest objects in an operator algebra on $M$.
They are defined, for instance, by their property to map weighted distributions on $M \setminus M'$ to smooth functions (and the same for the formal adjoints).
However, as we saw in Section 4.2, the notion of smoothness of a function on $M$ is itself a special invention and an input to the a priori philosophy of how the operators on $M$ in general (also those with non-vanishing symbols) have to look like.
Smoothness in that sense has to be compatible with pseudo-locality of operators which gives rise to smoothing operators by cutting out distributional kernels off the diagonal.
Their characterisation in terms of (say, tensor products of) smooth functions is an important aspect of the full calculus on $M$, and so we need to know what is smoothness on $M$.
As we see this is a substantial aspect.

\subsection{Pseudo-differential aspects, solvability of equations} 
\label{s.10.4.4}
%10.4.4

Pseudo-differential operators on a $C^\infty$ manifold $M$ can be motivated by parametrices of elliptic differential operators.
More precisely, there is a hull operation which extends the algebra $\bigcup_{\mu\in\N} \Diff^\mu (M)$
to a corresponding structure that is closed under forming parametrices of elliptic elements.
This process is natural for the same reason as the construction of multiplicative inverses of non-vanishing integers in the elementary calculus.
If $M$ is a manifold with singularities in the sense that there is a subset $M' \subset M$ of singular points such that $M \setminus M'$ is $C^\infty$, the hull operation makes sense both for 
$\Diff (M \setminus M'):= \bigcup_{\mu\in \N} \Diff^\mu (M \setminus M')$ (as before) and for suitable subalgebras of $\Diff (M \setminus M')$.
While for $\Diff^\mu (M \setminus M')$ the ellipticity is still expressed by $\sigma_\psi (A)$ (the homogeneous principal symbol of $A$ of order $\mu$), in subalgebras we may have additional principal symbolic information as sketched in Section 1.1.
The latter aspect is just one of the specific novelties with singularities.

Let us have look at some special cases.   

{\bf(A.1) Boundary value problems.} If $M$ is a $C^\infty$ manifold with boundary, the
task to complete classical differential boundary value
problems (e.g., Dirichlet or Neumann problems for Laplace
operators) gives rise to Boutet de Monvel's calculus of
pseudo-differential operators with the transmission property
at the boundary. The operators ${\cal A}$ are $2 \times 2$
block matrices, and the principal symbolic hierarchy consists
of pairs
\begin{equation}
\label{1044.sigmapartial2610.eq}
\sigma({\cal A}) = (\sigma_\psi({\cal A}),
        \sigma_\partial({\cal A})),
\end{equation}	
with the interior symbol $\sigma_\psi$ and boundary symbol
$\sigma_\partial$. 

{\bf(A.2) Operators on manifold with conical singularities.} 
Another case is  a manifold $M$ with conical
singularities. As the typical differential operators $A$ we
take the class $\Diff_{\deg}^\mu(M)$ of operators that are 
of Fuchs type near the conical
singularities (in stretched coordinates, and
including the weight factors $r^{- \mu}$ for $\mu = \ord A$).
The principal symbols consist of pairs
\begin{equation}
\label{1044.sigmac2610.eq}
\sigma(A) = (\sigma_\psi(A), \sigma_{\c, \gamma}(A))
\end{equation}
with the  (Fuchs type) interior symbol $\sigma_\psi$ and the
conormal symbol $\sigma_{\c, \gamma}$ (referring to the weight
line $\Gamma_{\frac{n+1}{2} - \gamma}$ as described before; 
$n$ is equal to the dimension of the base of the local cone, and
$\gamma \in \R$ is a weight). 
The associated pseudo-differential
calculus is called (in our terminology) the cone algebra,
equipped with the principal symbolic hierarchy
\eqref{1044.sigmac2610.eq}.

The stretched manifold $\M$ associated with a manifold $M$ 
with
conical singularities is a $C^\infty$ manifold with boundary 
(recall that the stretched coordinates $(r,x)$ just refer to a collar
neighbourhood of $\partial \M$ with $r$ being the normal 
variable). Nevertheless the cone calculus has a completely
different structure than the calculus of (A.1) of boundary 
value
problems with the transmission property. This shows that when
$\M$ means a stretched manifold to a manifold with conical
singularities the notation $C^\infty$ `manifold with boundary'
does not imply a canonical choice of a calculus (although
there are certain relations  between the calculi of (A.1) and
(A.2)). 
The cone algebra solves the problem of expressing
parametrices of elliptic differential operators $A \in
\Diff_{\deg}^\mu(M)$, and it is closed under parametrix
construction for elliptic elements, also in the
pseudo-differential case.

\begin{Remark}
\label{1044.cone2610.re}
On a manifold $M$ with conical singularities there are many
variants of `cone algebras':
\begin{enumerate}
\item The weight factor $r^{- \mu}$ can be replaced by any
      other factor $r^{- \beta}$, $\beta \in \R$, without an
      essential change of the calculus;
\item the ideals of smoothing operators depend on the choice
      of asymptotics near the tip of the cone, with finite or
      infinite asymptotic expansions and discrete or
      continuous asymptotics; this affects the nature of
      smoothing Mellin operators {\st{(}}with lower order conormal
      symbols{\st{)}} and of Green operators;
\item there is a cone algebra on the infinite cone $M =
      X^\Delta$ with an extra control at the conical exit 
      to infinity $r \to \infty$. In that case we have a 
      principal symbolic
      hierarchy with three components $\sigma(A) =
      (\sigma_\psi(A), \sigma_{\c, \gamma}(A),
      \sigma_{\st E}(A))$;
\item in cone algebras which are of interest in applications the
      base $X \cong \partial \M$ of the cone may have a
      $C^\infty$ boundary; we then have a cone calculus of
      boundary value problems in the sense of {\st (A.1)},
      i.e., $2 \times 2$ block matrices ${\cal A}$, with
      principal symbolic hierarchies
      \[  \sigma({\cal A}) = (\sigma_\psi({\cal A}),
              \sigma_\partial({\cal A}), \sigma_{\c, \gamma}
	      ({\cal A})),     \]          
      or, in the case $X^\Delta$ for $\partial X \not=
      \emptyset$, with exit calculus at $\infty$,
      \[  \sigma({\cal A}) = (\sigma_\psi({\cal A})),
            \sigma_\partial({\cal A}), \sigma_{\c,
	    \gamma}({\cal A}), \sigma_{\st E}({\cal A}),
	    \sigma_{E'}({\cal A}))      \]
      {\st{(}}with exit symbols $\sigma_{\st E}$ and
      $\sigma_{E'}$ from the interior and the boundary, 
      respectively{\st{)}}.	    
\end{enumerate}
\end{Remark}

{\bf (A.3) Operators on manifolds with edges.} 
Let $M$ be a manifold with smooth edge $Y$. As the
typical differential operators we take $\Diff_{\deg}^\mu(M)$
as  explained in Section 3.1. In this case the weight 
factor $r^{- \mu}$ in front of the operator (in stretched
coordinates) is essential for our edge
algebra. Similarly as (A.1), the edge algebra consists of $2
\times 2$ block matrices ${\cal A}$ with extra edge conditions
of trace and potential type. Instead of the principal boundary
symbol in \eqref{1044.sigmapartial2610.eq} (which is a $2 
\times 2$ block matrix family on $\R_+$, the inner normal to
the boundary) we now have a principal edge symbol 
$\sigma_{\land,
\gamma}({\cal A})$ which takes values in the cone algebra on
the infinite model cone $X^\Delta$ of local wedges, as
described in Remark \stref{1044.cone2610.re}. The weight
$\gamma \in \R$ is inherited from the cone algebra;
$\sigma_{\land, \gamma}({\cal A})$ as a $2 \times 2$ block
matrix family of operators ${\cal K}^{s, \gamma}(X^\land) \to
{\cal K}^{s- \mu, \gamma - \mu}(X^\land)$, parametrised by
$ T^* Y \setminus 0$. The principal symbolic hierarchy in the
edge algebra has again two components
\[  \sigma({\cal A}) = (\sigma_\psi({\cal A}), \sigma_{\land,
          \gamma}({\cal A})),      \]
with the (edge-degenerate) interior symbol $\sigma_\psi$ and
the principal edge symbol $\sigma_{\land, \gamma}$.	  
The edge algebra solves the problem of expressing parametrices
for elliptic elements with an operator $A \in
\Diff_{\deg}^\mu(M)$ in the upper left corner, and it is
closed under constructing parametrices of elliptic elements
also in the pseudo-differential case.

\begin{Remark}
\label{1044.edge2610.re}
On a manifold $M$ with edge $Y$ there are many variants of
`edge algebras', similarly as Remark \stref{1044.cone2610.re}      
for the case $\dim Y = 0$.
\begin{enumerate}
\item The edge algebra very much depends on the choice of the
      ideal of smoothing operators on the level of edge
      symbols, cf. Remark \stref{1044.cone2610.re} {\st{(ii)}}.
\item It is desirable to have an  edge algebra on the infinite cone
      $M^\Delta$ with a corresponding exit symbolic 
      structure, also in the variants of boundary
      value problems, i.e., a combination of {\st{(A.2)}} and
      {\st{(A.3)}}, when we have $\partial X \not= 0$ for the
      base $X$ of local model cones. 
      We then  have to expect corresponding larger principal symbolic
      hierarchies.
\item The edge algebra in the `closed case' {\st{(}}i.e.,
      $\partial X = \emptyset${\st{)}} is a generalisation of          
      the algebra of boundary value problems in the sense
      of {\st{(A.1)}}; the edge plays the role of the boundary
      and the local model cone of the inner normal. The 
      operators in the upper left corner have not 
      necessarily the
      transmission property at the boundary {\st{(}}they may
      even be edge degenerate{\st{)}}.
\end{enumerate}      
\end{Remark}

\begin{Remark}
\label{r.4.7}
The manifolds $M$ of Remarks {\em\ref{1044.cone2610.re}} and 
{\em\ref{1044.edge2610.re}}
{\em(}in the case without boundary{\em)} belong to $\got{M}_1$.
For $M \in \got{M}_2$ we also talk about the calculus of second generation.
The papers
{\em\cite{Schu29},
\cite{Schu25},
\cite{Schu27},
\cite{Haru11},
\cite{Haru12},
\cite{Mani2},
\cite{Calv3}} belong to this program.
\end{Remark}

The precise calculus of higher corner algebras, i.e., for $M
\in \got{M}_k$ for $k \geq 3$ is a program of future research,
although there are partial partial results, cf. \cite{Schu25}, \cite{Calv2}, and  Section  5 below.

The nature of a parametrix of an elliptic operator $A$ characterises to some extent the solvability of the equation
\begin{equation}
\label{eq.247}
Au=f.
\end{equation}
To illustrate that we consider the simplest case of 
$A \in \Diff^\mu (M)$ for a closed compact $C^\infty$ manifold $M$ and the scale 
$H^s (M)$, $s \in \R$, of standard Sobolev spaces on $M$.
The way to derive elliptic regularity of solutions $u$ is as follows.
By virtue of the ellipticity of $A$ there is a parametrix
$P \in L^{-\mu}_\cl (M)$, and 
$P: H^r(M) \to H^{r+\mu} (M)$ is continuous for every $r \in \R$.
Then \eqref{eq.247} gives us
$$
PA u = (1-G) u = Pf
$$
for an operator $G \in L^{-\infty} (M)$.
Because of $G H^{-\infty} (M) \subset H^\infty (M) \subset H^s (M)$ for every $s \in \R$ it follows that 
$u \in H^{-\infty} (M)$, $f \in H^r (M)$ implies 
$u =Pf +Gu \in H^{r+\mu} (M)$.
Observe that the existence of a parametrix also gives rise to so-called a priori estimates for the solutions.
That means, for every $r \in \R$ we have
\begin{equation}
\label{eq.248}
\| u \|_{H^s(M)} \leq
     c \big( \| f \|_{H^{s-\mu} (M)} +
             \| u \|_{H^r (M)} 
       \big)
\end{equation}
when $u \in H^{-\infty} (M)$ is a solution of 
$Au =f \in H^{s-\mu} (M)$, for a constant $c = c(r,s)>0$.
In fact, we have
$$
\| u \|_{H^s(M)}  \leq
   \| Pf \|_{H^s (M)} +
   \| Gu \|_{H^s (M)},
$$
and the right hand side can be estimated by \eqref{eq.248}, since
$P : H^{s-\mu} (M) \to H^s (M)$ and $G : H^r (M) \to H^s (M)$
are continuous for all $s,r \in \R$.

Similar conclusions make sense for elliptic operators on a manifold with singularities.
In other words, to characterise the solvability of the equation \eqref{eq.247} it is helpful to have the following structures.   

{\bf (S.1)   Operator algebras, symbols.}  
Construct an algebra of operators with a principal symbolic structure that defines operators modulo lower order terms.                 

{\bf (S.2)  Ellipticity, parametrices.} 
Define ellipticity in terms of the principal symbols (and, if necessary, kinds of Shapiro-Lopatinskij or global projection data)
and construct parametrices within the algebra.               

{\bf (S.3)  Smoothing operators.} 
Establish an ideal of smoothing operators to characterise the left over terms. 

{\bf (S.4)  Scales of spaces.} 
Introduce natural scales of distribution spaces such that the elements of the algebra induce continuous operators.

These aspects together with other features, such as asymptotic summation, formal Neumann series constructions, operator conventions (quantisations) and recovering of symbols from the operators, or kernel characterisations, belong to the desirable elements of calculi, also on manifolds with higher singularities.
As we saw this can be a very complex program.
However, the effort is justified.
The characterisation of the operators in the algebra reflects the internal structure of parametrices, while the functional analytic features of adequate scales of distributions describe in advance the nature of elliptic regularity.

In addition the algebra aspects appear at once in connection with single operators.
In order to treat any fixed operator of interest on a singular manifold of higher order we are faced with operator-valued symbolic components which are operator functions with values in the algebras of lower singularity order that may range over these algebras in (nearly) full generality.
Also from that point of view we need the constructions within a calculus with the features 
(S.1) - (S.4).

\subsection{Discrete, branching, and continuous asymptotics} 
\label{s.4.5}

%10.4.5

An interesting problem in partial differential equations near geometric singularities is the asymptotic behaviour of solutions close to the singularities.
For instance, in classical elliptic boundary value problems (with smooth boundary) the smoothness of the right hand sides and boundary data entails the smoothness of solutions up to the boundary (of course, there are also other features of elliptic regularity, e.g., in Sobolev spaces).

The latter property can already be observed on the level of operators on the half-axis 
$\R_+ \ni r$ for an elliptic operator $A$ of the form
\begin{equation}
\label{eq.249}
A = 
    \sum^{\mu}_{j=0}
        c_j \frac{d^j}{dr^j}
\end{equation}
with (say, constant) coefficients $c_j$.
We can rephrase $A$ as
\begin{equation}
\label{eq.250}
A = r^{-\mu}
    \sum^{\mu}_{j=0}
        a_j \big( -r \frac{d}{dr} \big)^j
\end{equation}
with other coefficients $a_j \in \C$.        
Assume that $c_0 \not= 0$ which is equivalent to $a_0 \not= 0$.
The asymptotics of solutions of an equation
\begin{equation}
\label{eq.251}
Au =f
\end{equation}
for $r \to 0$ when $f$ is smooth up to $r=0$ can be obtained in a similar manner as in Section 1.2.
In this case the resulting asymptotic type of $u$ is again of the form
$P = \big\{ (-j,0) \big\}_{j\in \N}$,
i.e., represents Taylor asymptotics.
Observe that the weight factor $r^{-\mu}$ in \eqref{eq.250} does not really affect the consideration; in Section 1.2 we could have considered the case with weight factors as well (as we saw the weight factors are often quite natural).

The transformation from \eqref{eq.249} to \eqref{eq.250} can be identified with a map
$(c_j)_{0 \leq j \leq \mu} \to 
    (a_j)_{0 \leq j \leq \mu}$,
$\C^{\mu+1} \to \C^{\mu+1}$, which is not surjective for $\mu >0$.
From Section 1.2 we know that when the coefficients $a_j$ in \eqref{eq.250} are arbitrary, Taylor asymptotics of solutions $u$ of \eqref{eq.251} is an exceptional case.
In fact, even for right hand sides $f$ that are smooth up to $r=0$ we obtain solutions $u$ with asymptotics of other types $P$, determined by the poles of the inverse of the conormal symbol
$\big( \sum^{\mu}_{j=0} a_j w^j \big)^{-1}$.
The asymptotic behaviour of solutions becomes much more complex when the equation \eqref{eq.251} is given on a (stretched) cone 
$X^\wedge = \R_+ \times X$ with non-trivial base, say, for a closed compact $C^\infty$ manifold $X$.
Then the resulting asymptotic types may be infinite, and it is interesting to enrich the information by finite-dimensional spaces
$L_j \subset C^\infty (X)$, i.e., to consider sequences
\begin{equation}
\label{eq.252}
P = \big\{ (p_j, m_j, L_j) \big\}_{j \in \N},
\end{equation}
$\pi_\C P = \{ p_j \}_{j \in \N} \subset 
     \{ w \in \C : \re w < \frac{n+1}{2} - j \}$
for $n = \dim X$ and some weight $j \in \R$,
$\re p_j \to -\infty$ as $j \to \infty$.
Recall that an 
$u(r,x) \in \s{K}^{s,\gamma} (X^\wedge)$ has asymptotics for $r \to 0$ of type $P$, if for every  $\beta >0$ there is an $N = N(\beta)$ such that
\begin{equation}
\label{eq.253}
u(r,x) - \omega(r)
  \sum^{N}_{j=0} \sum^{m_j}_{k=0}  c_{jk} (x) r^{-p_j} \log^k r \in 
      \s{K}^{s,\beta} (X^{\wedge})     
\end{equation} 
with coefficients 
$c_{jk} \in L_j$, $0 \leq k \leq m_j$.
This condition just defines the space $\s{K}^{s,\gamma}_P (X^\wedge)$.
Observe that when we set $\Theta := (\vartheta,0]$ for some finite 
$\vartheta < 0$ and
\begin{equation}
\label{eq.254}
P_\Theta :=
   \big\{ (p,m,L) \in P : \ \ 
    \re p > \frac{n+1}{2} -\gamma + \vartheta \big\},
\end{equation}  
 
\begin{equation}
\label{eq.255}
\s{K}^{s,\gamma}_\Theta (X^\wedge) :=
   \limpr_{k \in \N} 
      \s{K}^{s, \gamma-\vartheta-k^{-1}} (X^\wedge),
\end{equation}   
\begin{equation}
\label{eq.256}
\s{E}_{P_\Theta} (X^\wedge) 
:=
   \Big\{
      \sum^{N}_{j=0} \sum^{m_j}_{k=0}  c_{jk} (x)  
                                       \omega (r) r^{-p_j} \log^k r  :  
c_ {jk} \in L_j , \ \ 
 0 \leq k \leq m_j, \ \ 0 \leq j \leq N  \Big\},   
 \end{equation}
 the direct sum
\begin{equation}
\label{eq.257}
\s{K}^{s,\gamma}_{P_\Theta} (X^\wedge) := 
\s{K}^{s,\gamma}_{\Theta} (X^\wedge) + 
\s{E}_{P_\Theta} (X^\wedge) 
\end{equation}
is a Fr\'echet space, and we have 
$\s{K}^{s,\gamma}_{P} (X^\wedge) =
    \limpr\limits_{\vartheta \to -\infty} 
      \s{K}^{s, \gamma}_{P_\Theta} (X^\wedge)$.
  
Operators of the form \eqref{eq.250} (in general with $r$-dependent coefficients $a_j$) occur as the (principal) edge symbols of edge-degenerate operators
\begin{equation}
\label{eq.258}
{\bsA} :=
  r^{-\mu} \sum_{j+|\beta| \leq \mu}
      b_{j\beta} (r,y) 
                 \Big( -r \frac{\partial}{\partial r} \Big)^j
                 (r D_y)^\beta
\end{equation}
on a (stretched) wedge $X^\wedge \times \Omega$, $\Omega \subseteq \R^q$ open,
with coefficients 
$b_{j\beta} (r,y) \in 
 C^\infty \big( \ol{\R}_+ \times \Omega$, $\Diff^{\mu-(j+|\beta|)} (X)
          \big)$.          
The principal edge symbol of \eqref{eq.258} is defined as 
\begin{equation}
\label{eq.259}
\sigma_\wedge (\bsA) (y,\eta) :=
   r^{-\mu} \sum_{j+|\beta| \leq \mu}
      b_{j\beta} (0,y) 
                 \Big( -r \frac{\partial}{\partial r} \Big)^j
                 (r \eta)^\beta,
\end{equation}                 
$(y,\eta) \in T^* \Omega \setminus 0$, and represents a family of continuous operators
\begin{equation}
\label{eq.260}
\sigma_\wedge (\bsA) (y,\eta) :
   \s{K}^{s,\gamma} (X^\wedge) \to
      \s{K}^{s-\mu, \gamma -\mu} (X^\wedge)
\end{equation}
for every $s,\gamma \in \R$.
It turns out that the asymptotics of solutions $\bsu$ of an elliptic equation
\begin{equation}
\label{eq.261}\bsA \bsu = \bsf
\end{equation}
on $X^\wedge \times \Omega$ for $r \to 0$ is determined by  \eqref{eq.260}, more precisely, by the inverse of the conormal symbol of \eqref{eq.260}, namely, 
$\sigma_{\c} \sigma_\wedge (\bsA)^{-1} (y, w)$, where
$$
\sigma_{\c} \sigma_\wedge (\bsA) (y, w) =
    \sum^\mu_{j=0} b_{j0} (0,y) w^j        
$$
which is a family of continuous operators
$$
\sigma_{\c} \sigma_\wedge (\bsA) (y, w) :
   H^s (X) \to H^{s-\mu} (X)
$$
for every $s \in \R$, smooth in $y \in \Omega$, holomorphic in $w \in \C$.
However, the question is:
`what means asymptotics in the edge case?'   
The answer is far from being straightforward, and, what concerns the choice of spaces that contain the solutions, we have a similar problem as above in connection with the `right approach' to Sobolev spaces, discussed in Section 1.3.
Here, in connection with asymptotics, this problem appears in refined form, because the choice of similar terms of asymptotics requires a confirmation of the formulation of the spaces.
In order to illustrate some of the asymptotic phenomena, for convenience, we consider the case 
$\Omega = \R^q$ and assume the coefficients $b_{j \beta} (r,y)$ to be independent of $y$ when $|y| >C$ for some $C>0$ and independent of $r$ for $r >R$ for some $R>0$.
In addition we assume that for some $\gamma \in \R$ the operators \eqref{eq.260} define isomorphisms for all $s \in \R$ (in general, we can only expect Fredholm operators; for the ellipticity those are to be filled up to $2 \times 2$ block matrices of isomorphisms by extra entries of trace, potential, etc., type with respect to the edge $\R^q$).
For the operator $\bsA$ we assume $\sigma_\psi$-ellipticity in the sense that the homogeneous principal symbol
$$
\sigma_\psi (\bsA) (r,x,y, \rho, \xi, \eta)
$$
does not vanish for $(\rho, \xi, \eta) \not= 0$ and that
$r^\mu \sigma_\psi (\bsA) (r,x,y,r^{-1} \rho, \xi, r^{-1} \eta) \not= 0$
for $(\rho, \xi,\eta) \not= 0$,
up to $r=0$.

The family of operators 
$$
\bsa (y,\eta) :=
  r^{-\mu} \sum_{j+|\beta| \leq \mu}
      b_{j \beta} (r,y) \Big(-r \frac{\partial}{\partial r} \Big)^j (r \eta)^\beta :
          \s{K}^{s,\gamma} (X^\wedge) \to
            \s{K}^{s-\mu,\gamma-\mu} (X^\wedge)
$$
can be interpreted as an element 
$\bsa(y,\eta)
      \in S^\mu \big(\R^q \times \R^q; 
             \s{K}^{s,\gamma} (X^\wedge),
              \s{K}^{s-\mu,\gamma-\mu} (X^\wedge)\big)$ 
for every $s \in \R$, cf. Definition \ref{1013.1305de.de} and Remark \ref{1011.1305re.re}, which gives us a continuous operator
$$
\bsA = \Op_y (\bsa) :
   \s{W}^s (\R^q, \s{K}^{s,\gamma} (X^\wedge)) \to
      \s{W}^{s-\mu} (\R^q, \s{K}^{s-\mu,\gamma-\mu} (X^\wedge)) 
$$
for every $s \in \R$.
Now the pseudo-differential calculus of edge-degenerate operators allows us to construct a symbol
\begin{equation}
\label{eq.262}
\bsp (y,\eta) \in
    S^{-\mu} \big(\R^q \times \R^q; \s{K}^{s-\mu,\gamma-\mu} (X^\wedge),
                                 \s{K}^{s,\gamma} (X^\wedge)\big) 
\end{equation}
such that the operator
$\bsP := \Op_y (\bsp) : \s{W}^{s-\mu} (\R^q, \s{K}^{s-\mu,\gamma-\mu} (X^\wedge))  \to
                       \s{W}^s (\R^q, \s{K}^{s,\gamma} (X^\wedge))$
is a parametrix of $\bsA$ in the sense that there is an $\varepsilon>0$ such that
\begin{equation}
\label{eq.263}
\bsP \bsA - \bsI :
  \s{W}^s (\R^q, \s{K}^{s,\gamma} (X^\wedge)) \to
                                 \s{W}^{\infty} (\R^q, \s{K}^{\infty,\gamma+\varepsilon}(X^\wedge)) 
\end{equation}  
is continuous for all $s$, and, similarly,
$\bsP \bsA - \bsI$.
The relation \eqref{eq.263} gives us elliptic regularity of solutions
$\bsu \in   
   \s{W}^{-\infty} (\R^q, \s{K}^{-\infty,\gamma}(X^\wedge))$ to \eqref{eq.261} for
\begin{equation}
\label{eq.264}
\bsf \in
  \s{W}^{s-\mu} (\R^q, \s{K}^{s-\mu, \gamma-\mu}(X^\wedge)), 
\end{equation}
namely,
\begin{equation}
\label{eq.265}
\bsu \in
  \s{W}^{s} (\R^q, \s{K}^{s, \gamma}(X^\wedge)). 
\end{equation}  
The precise nature of $\bsP$ is a subtle story, contained in the analysis of the edge algebra, cf. \cite{Schu2}, or \cite{Schu20}.
The question is, do we have an analogue of elliptic regularity with edge asymptotics.
Here, an $\bsu \in  \s{W}^{s} (\R^q, \s{K}^{s, \gamma} (X^\wedge))$ is said to have discrete edge asymptotics of type $P$ if $\bsu \in  \s{W}^{s} (\R^q, \s{K}^{s, \gamma}_P (X^\wedge))$.
In the following theorem we assume that the coefficients $b_{j \beta}$ are independent of $y$ everywhere.
\begin{Theorem}
\label{t.4.3}
Let $\bsA$ satisfy the above conditions, and let the coefficients $b_{j\beta}$ be independent of $y \in \R^q$.
Then for every discrete asymptotic type 
$Q = \{ (q_j, n_j, M_j) \}_{j \in \N}$
{\em(}with $\pi_\C Q \subset 
   \{ w \in \C : \re w < \frac{n+1}{2} - (\gamma -\mu) \}${\em)}
   there exists a $P$ as in {\em(10.4.32)} such that 
   $\bsu \in  \s{W}^{-\infty} (\R^q, \s{K}^{-\infty, \gamma} (X^\wedge))$ and
 $\bsA \bsu = \bsf \in 
    \s{W}^{s-\mu} (\R^q, \s{K}^{s-\mu, \gamma-\mu}_Q (X^\wedge))$ implies
$\bsu \in \s{W}^s (\R^q, \s{K}^{s,\gamma}_P (X^\wedge))$.
\end{Theorem}

This result may be found in \cite{Schu20}.
It is based on the fact that there is a parametrix $\bsP = \Op_y (\bsp)$ for an amplitude function \eqref{eq.262} that restricts to elements
\begin{equation}
\label{eq.266}
\bsp (y,\eta) \in
  S^{-\mu} (\R^q \times \R^q;
   \s{K}^{s-\mu, \gamma-\mu}_S (X^\wedge), \s{K}^{s,\gamma}_B (X^\wedge))
\end{equation}
(in this special case independent of $y$)
for every discrete  asymptotic type $S$ with some resulting $B$ and such that the remainder
$\bsG := \bsP \bsA - \bsI$ defines continuous operators
$$
\bsG : \s{W}^s (\R^q, \s{K}^{s,\gamma} (X^\wedge)) \to
        \s{W}^\infty (\R^q, \s{K}^{\infty, \gamma}_{\wt{B}} (X^\wedge))
$$
for some discrete asymptotic type $\wt{B}$.

The rule in Theorem \ref{t.4.3} to find  $P$ in terms of $Q$ is very close to that discussed before in Section 1.2.
The essential observation is that there is a Mellin asymptotic type $R$ 
(see the formula \eqref{eq.51}) such that
$$
\sigma_\c \sigma_\wedge (\bsA)^{-1} (w) \in 
      \s{M}^{-\mu}_R (X),
$$
and $R$ in this case is independent of $y$.
Unfortunately, this conclusion does not work in general, when the coefficients depend on $y$.
Although we also have
$$
\sigma_\c \sigma_\wedge (\bsA)^{-1} (y, w) \in 
     \s{M}^{-\mu}_{R(y)} (X)
$$
for every fixed $y$, the asymptotic type $R(y)$ may depend on $y$, and we cannot expect any property like \eqref{eq.266}.
The $y$-dependence of $R$ means that all components of \eqref{eq.51} depend on $y$; in particular, the numbers $n_j$ which  encode the multiplicities of poles, may jump with varying $y$, and there are no smooth `paths' of poles $r_j (y)$, $y \in \R^q$, in the complex plane, but, in general, irregular clouds of points
$\{ \pi_\C R(y) : y \in \R^q \}$.
Then, even if we can detect some $y$-dependent families of discrete asymptotic types $Q{(y)}$, 
$P (y)$, with the hope to discover a rule as in Theorem \ref{t.4.3} also in the general case, the first question is, what are the spaces 
$\s{W}^{s,\gamma} (\R^q, \s{K}^{s, \gamma}_{P(y)} (X^\wedge))$?
An answer for the case $\dim X=0$ is given in \cite {Schu34}, \cite{Schu36}.
The point is to encode somehow the expected variable discrete and branching patterns of poles (that appear after Mellin transforming a function with such asymptotics).
We do not discuss here all the details up to the final conclusions; this would go beyond the scope of this exposition.
More information may be found in \cite{Schu2}, or \cite{Schu20}, see also \cite{Kapa10}.
We want to give an idea of how discrete and  branching asymptotics are organised in such a way that the concept admits edge spaces together with continuity results for pseudo-differential operators in such a framework.

The key word in this connection is `continuous asymptotics'.
The notion is based on analytic functionals in the complex plane.
We do not recall here too much material on this topic.
Let us content ourselves with some notation.
If $\s{A} (U)$, $U \subseteq \C$ open, is the space of all holomorphic functions in $U$, endowed with the Fr\'echet topology of uniform convergence on compact subsets, we have the space $\s{A}' (U)$ of all linear continuous functionals
$$
\zeta : \s{A} (U) \to \C,
$$
the so-called analytic functionals in $U$.
For every open $U \subseteq V$ we have a restriction operator $\s{A}' (U) \to \s{A}' (V)$.
Given $\zeta \in \s{A}' (\C)$, an open set $U \subseteq \C$ is called a carrier of $\zeta$, if there is an element $\zeta_U \in \s{A}' (U)$ which is the restriction of $\zeta$ to $U$.
A compact subset $K \subset \C$ is said to be a carrier of $\zeta \in \s{A}' (\C)$, if every open $U \supset K$ is a carrier of $\zeta$ in the former sense.

By $\s{A}' (K)$ we denote the subspace of all $\zeta \in \s{A}' (\C)$ carried by the compact set $K$.
It is known that  $\s{A}'(K)$ is a nuclear Fr\'echet space in a natural system of semi-norms.

It also makes sense to talk about analytic functionals with values in a, say, Fr\'echet space $E$, i.e., we have the spaces $\s{A}' (K,E) = \s{A}' (K) \hat{\otimes}_\pi E$ of $E$-valued analytic functionals, carried by $K$.
We may take, for instance, $E = C^\infty (X)$.
\begin{Example}
\label{ex.4.10}
Let $K \subset \C$ be a compact set, and let $C$ be a smooth compact curve in $\C \setminus K$ surrounding the set $K$ counter-clockwise.
In addition we assume that there is a diffeomorphism $\kappa : S^1 \to C$ such that, when we identify any $w \in K$ with the origin in $\C$, the corresponding  winding number of $\kappa$ is equal to $1;$ this is required for every $w \in K$.
It can be proved that for every $\varepsilon >0$ there exists a curve $C$ of this kind such that $\dist (K,C) < \varepsilon$, cf. {\em\cite[Theorem 13.5]{Rudi1}}.
Let $f \in \s{A} (\C \setminus K)$, and form
\begin{equation}
\label{eq.267}
\langle \zeta, h \rangle :=
   \frac{1}{2\pi i}
      \int_C f(w) h(w) dw
\end{equation}
for $h \in \s{A} (\C)$.
Then we have $\zeta \in \s{A}' (K)$.
More generally, considering an $f(y,w) \in C^\infty (\Omega, \s{A} (\C \setminus K, E))$
for an open set  $\Omega \subseteq \R^q$ and a Fr\'echet space $E$, by
\begin{equation}
\label{eq.268}
\langle \zeta (y), h \rangle :=
   \frac{1}{2\pi i}
      \int_C f(y,w) h(w) dw
\end{equation}
we obtain an element $\zeta (y) \in C^\infty (\Omega, \s{A}' (K,E))$.
Clearly \eqref{eq.268} is independent of the choice of $C$.
\end{Example}

In \eqref{eq.267} we can take, for instance,
$f(w) = M_{r \to w} (w(r) r^{-p} \log^k r) (w)$ for any $p \in \C$, $k \in \N$, with $M$ being the weighted Mellin transform {\em(}with any weight $\gamma \in \R$ such that $\re p < \frac{1}{2} - \gamma${\em)} and a cut-off function $\omega (r)$.
Then \eqref{eq.267} takes the form
\begin{equation}
\label{eq.269}
\langle \zeta, h \rangle =
  (-1)^k \frac{d^k}{dw^k} h(w) |_{w =p},
\end{equation}
$h \in \s{A} (\C)$.
This corresponds to the $k$-th derivative of the Dirac measure at the point $p$, and we have $\zeta \in \s{A}' (\{ p \})$.
Inserting $h(w) := r^{-w}$ in \eqref{eq.269}   it follows that 
$\langle \zeta, r^{-w} \rangle =
           r^{-p} \log^k r$.
More generally, we have the following proposition.
 
\begin{Proposition}
\label{p.4.11}
Let 
$K := \{ p_0 \} \cup \{ p_1 \} \cup \ldots \cup \{ p_N \}$
for $p_j \in \C$, $j=0, \ldots, N$, and let $f \in \s{A} (\C \setminus K, E)$
be a meromorphic function with poles at $p_j$ of multiplicity $m_j+1$, and let 
$(-1)^k k! c_{jk} \in E$ be the Laurent coefficients at $(z - p_j)^{-(k+1)}$, $0 \leq k \leq m_j$.
Then the formula \eqref{eq.267} represents an element of $\s{A}' (K,E)$
which is of the form
\begin{equation}
\label{eq.270}
\langle \zeta, h \rangle =
   \sum_{j=0}^N \sum_{k=0}^{m_j}  (-1)^k
        c_{jk} \frac{d^k}{dw^k} h(w) \big|_{w=p_j},
\end{equation}
and we have
$\langle \zeta, r^{-w} \rangle =
   \sum_{j=0}^N \sum_{k=0}^{m_j}  
        c_{jk} r^{-p_j} \log^k r$.
\end{Proposition}      
An analytic functional of the form \eqref{eq.270} will be  called discrete (and of finite order).
In particular, if \eqref{eq.262} is a discrete asymptotic type, the relation \eqref{eq.253} can be interpreted as follows.
There is a sequence $\zeta_j \in \s{A}' (\{p_j \}, L_j)$ of discrete analytic functionals such that for every $\beta > 0$ there is an $N = N (\beta) \in \N$ such that 
\begin{equation}
\label{eq.271}
u(r,x) - \omega (r)
     \sum^N_{j=0}
        \langle \zeta_j, r^{-w} \rangle \in
             \s{K}^{s,\beta} (X^\wedge).
\end{equation}
This definition of the space $\s{K}^{s,\gamma}_P (X^\wedge)$ admits a generalisation as follows.
We replace $\{ p_j \}$ by arbitrary compact sets
$K_j \subset \{ w \in \C:
      \re w < \frac{n+1}{2} - \gamma \}$, $j \in \N$, such that 
$\sup \{ \re w : w \in K_j \} \to - \infty$ as $j \to \infty$.
Then an element $u(r,x) \in \s{K}^{s,\gamma} (X^\wedge)$ is said to have continuous asymptotics for $r \to 0$, if there is a sequence
$\zeta_j \in \s{A}' (K_j, C^\infty (X))$ such that the relation \eqref{eq.271} holds for every $\beta >0$ with some $N = N(\beta)$.

The notion of continuous asymptotics has been introduced in Rempel and Schulze \cite{Remp6} and then investigated in detail in \cite{Schu28}, \cite{Schu30}, \cite{Schu2}, \cite{Schu16}, \cite{Schu34}, \cite{Schu36}, see also \cite{Schu20}, or \cite[Section 2.3.5]{Kapa10}.
The original purpose was to find a way to express variable discrete asymptotics.
We do not develop here the full story but only give the main idea.
Intuitively, a family $u(r,x,y) \in  C^\infty ( \Omega, \s{K}^{s,\gamma} (X^\wedge))$ should have asymptotics of that kind, if there is a family
\begin{equation}
\label{eq.272}
P(y) =
   \{ (p_j (y), m_j (y), L_j (y) ) \}_{j\in \N}
\end{equation}
of discrete asymptotic types such that for every compact subset $M \subset \Omega$ and every 
$\beta >0$ there is an $N = N(\beta)$ with the property
\begin{equation}
\label{eq.273}
u(r,x,y) - \omega (r)
   \sum_{j=0}^N \sum_{k=0}^{m_j}
       c_{jk} (x,y) r^{-p_j (y)} \log^k r \in
           \s{K}^{s,\beta} (X^\wedge)
\end{equation}
with coefficients $c_{jk} (x,y) \in L_j (y)$, $0 \leq k \leq m_j$, for every fixed $y \in M$.

The nature of the family \eqref{eq.272} which appears in `realistic' pointwise discrete and branching asymptotic types belonging to solutions $\bsu$ of \eqref{eq.261} in the general case can be described as follows.
For every open set $U \subset \Omega$ such that $\ol{U} \subset \Omega$, $\ol{U}$ compact, there exists a sequence of compact sets $K_j \subset \C$, $j \in \N$, with the above-mentioned properties and a sequence $\zeta_j \in C^\infty (U, \s{A}' (K_{j}, C^\infty (X)))$, $j \in \N$, such that 
$\zeta_j (y) \in \s{A}' (\{ p_j (y) \}, L_j (y))$ is discrete for every fixed $y \in U$, and for every $\beta >0$ there is an $N = N (\beta) \in \N$ such that 
$$
u(r,x,y) - \omega (r) 
  \sum_{j=0}^N 
    \langle \zeta_j (y), r^{-w} \rangle \in
       \s{K}^{s,\beta} (X^\wedge)
$$
for every $y \in U$.
Observe that this notion really admits branchings of the exponents in \eqref{eq.273} and jumping $m_j (y)$ and $c_{jk} (x,y)$ with varying $y$.
Examples of such $\zeta_j(y)$ (say, in the scalar case) are functionals of the form \eqref{eq.268} for
$$
f(y,w) :=
  \frac{c(y) -z}{(a(y) -z) (b(y) -z)}
$$
with coefficients $a,b,c \in C^\infty (\Omega)$ taking values in $K_j$.

The characterisation of elements of $\s{W}^s_\loc (\Omega, \s{K}^{s,\gamma} (X^\wedge))$, $s \in \R$, with branching discrete asymptotics is also interesting.
The details are an excellent excersise for the reader.
In order to have an impression what is going on we want to consider once again the case of constant (in $y$) discrete asymptotics.
Let us give a notion of singular functions of the edge asymptotics of elements in 
$\s{W}^s (\R^q, \s{K}^{s,\gamma} (X^\wedge))$ of type \eqref{eq.252}.
To this end we fix any $-\infty < \vartheta <0$, form $P_\Theta$ by \eqref{eq.254} for $\Theta = (\vartheta, 0]$, and consider the decomposition \eqref{eq.257} which gives rise to a decomposition
\begin{equation}
\label{eq.274}
\s{W}^s (\R^q, \s{K}^{s,\gamma}_{P_\Theta} (X^\wedge)) =
           \s{W}^s (\R^q, \s{K}^{s,\gamma}_{\Theta} (X^\wedge)) +
              \s{V}^s (\R^q, \s{E}_{P_\Theta} (X^\wedge)),
\end{equation}
see Remark \ref{1013.re2005.re}, which is valid in analogous form also for the Fr\'echet space 
$E := \s{K}^{s,\gamma}_{P_\Theta} (X^\wedge)$ with the subspaces
$$
L := \s{K}^{s,\gamma}_\Theta (X^\wedge),
M := \s{E}_{P_\Theta} (X^\wedge).
$$
Thus every
$$
u(r,x,y) \in
 \s{W}^{s} (\R^q, \s{K}^{s,\gamma}_P  (X^\wedge) \subset
    \s{W}^{s} (\R^q, \s{K}^{s,\gamma}_{P_\Theta}  (X^\wedge))
$$
can be written as
$$
u(r,x,y) = u_\flat (r,x,y) + u_\sing (r,x,y)
$$
for an element
$$
u_\flat (r,x,y) \in 
   \s{W}^{s} (\R^q, \s{K}^{s,\gamma}_\Theta  (X^\wedge))
$$
of edge flatness $\Theta$ (relative to the weight $\gamma$) and a 
$$
u_\sing (r,x,y) \in
   \s{V}^{s} (\R^q, \s{E}_P  (X^\wedge))    \equiv
     F^{-1} \kappa_{\langle \eta \rangle}
       F (H^s (\R^q, \s{E}_{P_\Theta} (X^\wedge))),
$$
with the Fourier transform $F = F_{y \to \eta}$.
The space $\s{E}_{P_\Theta} (X^\wedge)$ is of finite dimension, cf. the formula \eqref{eq.258}.
The space $\s{V}^s (\R^q, \s{E}_{P_\Theta} (X^\wedge))$ consists of all linear combinations of functions
\begin{equation}
\label{eq.275}
F^{-1}
 \Big\{ \langle \eta \rangle^{\frac{n+1}{2}}
   \hat{v} (\eta) c_{jk} (x) \omega (r \langle \eta \rangle)
         (r \langle \eta \rangle)^{-p_j} \log^k 
         (r \langle \eta \rangle ) \Big\}
\end{equation}
for arbitrary $v \in H^s (\R^q)$, $\hat{v}(\eta) = ( F_{y \to \eta} v) (\eta)$.
In other words, \eqref{eq.275} describes the shape of the singular functions of the edge asymptotics of constant (in $y$) discrete type $P$.
In particular, we see (say, for the case $k=0$) how the Sobolev smoothness in $y \in \R^q$ of the coefficients of the asymptotics depends on $\re p_j$.
Note that decompositions of the kind \eqref{eq.274} have a nice analogue in classical Sobolev spaces $H^s (\R^{d+q})$ relative to a hypersurface $\R^q$, cf. the paper \cite{Dine1}.
The singular functions \eqref{eq.275} can also be written as
\begin{equation}
\label{eq.276}
F^{-1}
 \Big\{ \langle \eta \rangle^{\frac{n+1}{2}}
   \hat{v} (\eta) 
        \langle \zeta, (r \langle \eta \rangle)^{-w}  \rangle \Big\}
\end{equation}        
for suitable discrete $\zeta \in \s{A}' (\{ p_j \}, L_j)$ and $v \in H^s (\R^q)$.
The generalisation to continuous asymptotics is based on singular functions of the form
$$
F^{-1}
 \Big\{ \langle \eta \rangle^{\frac{n+1}{2}}
        \langle \zeta(\eta), (r \langle \eta \rangle)^{-w}  \rangle \Big\}
$$
for $\zeta (\eta) \in 
   \s{A}' (K, C^\infty (X) \hat{\otimes}_\pi \hat{H}^s (\R^q_\eta))$, $K \subset \C$ compact, where
$\hat{H}^s (\R^q_\eta) = F_{y \to \eta} H^s (\R^q_y)$.
Edge asymptotics in the $y$-wise discrete case on a wedge $X^\wedge \times \Omega$ can be modelled on
\begin{equation}
\label{eq.277}
F^{-1}_{y \to \eta}
 \Big\{ \langle \eta \rangle^{\frac{n+1}{2}} 
        \langle \zeta (y,\eta), (r \langle \eta \rangle)^{-w}  \rangle \Big         \}
\end{equation}                
for functions
\begin{equation}
\label{eq.278}
\zeta (y,\eta) \in
   C^\infty 
    (U, \s{A}' (K, C^\infty (X) 
       \hat{\otimes}_\pi
         \hat{H}^s (\R^q_\eta))),
\end{equation}     
$U \subset \Omega$ open, $\ol{U} \subset \Omega$, $\ol{U}$ compact, $K = K(U) \subset \C$ compact, where
$\zeta (y,\eta)$ is as in \eqref{eq.276} for every fixed $y \in U$.

Now the singular functions with variable discrete and branching asymptotics are formulated as 
\eqref{eq.277}     where \eqref{eq.278} are pointwise discrete and of finite order, i.e., pointwise of the form \eqref{eq.270}, with coefficients
$c_{jk} (x,y,\eta) \in L_j (y) \otimes F (H^{s-\frac{n+1}{2}} (\R^q))$,
$p_j =p_j (y)$, $m_j = m_j (y)$, cf. the expression \eqref{eq.272}.
Edge asymptotics in such a framework is a rich program, partly for future research.
Elliptic regularity of solutions to elliptic edge problems with continuous asymptotics is carried out in different contexts, see, e.g.,
\cite{Schu2}, \cite{Schu20}, or \cite{Kapa10}.
Variable discrete asymptotics for boundary value problems have been studied in \cite{Schu34}, \cite{Schu36} and by Bennish \cite{Benn1}.
Here we only want to mention that such a program requires the preparation of Mellin and Green symbols which also reflect such asymptotics, similarly as in the discrete case  (with Mellin symbol spaces consisting of meromorphic operator functions).

\section{Higher generations of calculi}   
\label{s.10.5}  
\begin{minipage}{\textwidth}
\setlength{\baselineskip}{0cm}
\begin{scriptsize}    
Manifolds with singularities of order $k$ form a category $\frak{M}_k$ ($\frak{M}_0$ is the category of $C^\infty$ manifolds, $\frak{M}_1$ the one of manifolds with  conical singularities or smooth edges, etc.).
The elements of $\frak{M}_{k+1}$ can be defined in terms of $\frak{M}_k$ by an iterative process.
Every $M \in \frak{M}_k$ supports an algebra of natural differential operators, with principal symbolic hierarchies and notions of ellipticity.
It is an interesting task to construct associated algebras, as outlined in Section 4.4.
The answers that are already given for $\frak{M}_1$ and $\frak{M}_2$, see, for instance, \cite{Schu32}, \cite{Schu29}, or \cite{Schu27}, show that the structures on the level $k+1$ require the parameter-dependent calculus from the level $k$, together with elements of the index theory and many other features that are also of interest on their own right.
The analysis on manifolds with singularities is not a simple induction from $k$ to $k+1$, although some general observations seem to be clear in `abstract terms'.
\end{scriptsize}
\end{minipage}

\subsection{Higher generations of weighted corner spaces}
\label{s.10.5.1}

%10.5.1

One of the main issues of the analysis on manifolds $M$ with
higher singularities is the character of weighted Sobolev
spaces on such manifolds. 
According to the general
principle of successively generating cones and wedges and
then to globalising the distributions on $M$ we mainly have to
explain the space
\begin{equation}
\label{1051.calK2610.eq}
{\cal K}^{s, \gamma}(X^\land), \quad (t,x) \in X^\land,
\end{equation}
for a (compact) manifold $X \in \got{M}_k$, $s \in \R$, for
a weight tuple $\gamma \in \R^{k}$, and what is the
weighted wedge space
$$
  {\cal W}^{s, \gamma}(X^\land \times \R^q), \quad
     (t,x,y) \in X^\land \times \R^q.      
$$
Before we give an impression on how these spaces are organised, we want to recall that every $M \in \got{M}_k$ is connected with a chain of subspaces \eqref{1031.eq1210.eq}, $M^{(j)} \in \got{M}_{k-j}$,
$j=0, \ldots, k$, where $M^{(0)} =M$.
Let us assume for the moment that $M^{(j)}$ is compact for every $j$ (otherwise, when we talk about weighted corner spaces, we will also have variants with subscript `(comp)' and `(loc)').
Moreover, for simplicity, we first consider `scalar' spaces; the case of distributional sections of vector bundles will be a modification.

For $k=0$ we take the standard Sobolev spaces $H^s (M)$, $s \in \R$.
If $M \in \got{M}_1$ has conical singularities, we have our weighted cone spaces $\s{H}^{s,\gamma} (\M)$ for $s,\gamma \in \R$, on the corresponding stretched manifold $\M$.
For $M \in \got{M}_1$ with smooth edge we can take the weighted edge spaces $\s{W}^{s,\gamma} (\M)$ for $s,\gamma \in \R$, on the stretched manifold $\M$ associated with $M$.
Those are subspaces of $H^s_\loc (\Int \M)$, modelled on 
\begin{equation}
\label{eq.f0}
\s{W}^s (\R^q, \s{K}^{s,\gamma} (X^\wedge)),
\end{equation}
locally in a neighbourhood of $\partial \M = \M_\sing$.
The invariance of these spaces refers to an atlas on $\M$, where the transition maps near $r=0$ are independent of $r$.
Recall from Section 1.3 that we can also form the spaces
\begin{equation}
\label{eq.g.}
\s{W}^s (\R^q, \s{K}^{s,\gamma;g} (X^\wedge))
\end{equation}
for every $s,\gamma,g \in \R$, based on the group action \eqref{eq.new11}.
Let $\s{W}^{s,\gamma;g} (\M)$ denote the corresponding global spaces on $\M$ (they make sense for similar reasons as before with an atlas as for $g=0$).
\begin{Remark}
\label{r.inv.}
The spaces $\s{W}^{s,\gamma;s-\gamma} (\M)$ are invariantly defined for particularly natural charts on $\M$, namely, those mentioned at the beginning of Section {\em 4.1}, here for the case $k=1$ {\em(}cf. also {\em\cite{Tark4})}.
For simplicity, in the following discussion we return to the case $g=0$ and ignore this extra information.
\end{Remark}

For the higher calculi it seems better to modify some notation and to refer to the singular manifolds $M$ themselves rather than their stretched versions, although the distributions are always given on 
$M \setminus M^{(1)} \in \got{M}_0$.
So we replace notation as follows:
\begin{equation}
\label{eq.an}
\s{H}^{s,\gamma} (\M) \to  \s{H}^{s,\gamma} (M), \ \
\s{W}^{s,\gamma} (\M) \to  \s{H}^{s,\gamma} (M), \ \
\s{K}^{s,\gamma} (X^\wedge) \to  \s{K}^{s,\gamma} (X^\Delta).
\end{equation}
We only preserve the $\s{W}^{s,\gamma}$-notation in wedges $X^\Delta \times \R^q$, in order to keep in mind the edge-definition in the sense of Definition \ref{1013.1305de.de}.
In other words, we set
\begin{equation}
\label{eq.ed}
\s{W}^{s,\gamma} (X^\Delta \times \R^q) :=
     \s{W}^s (\R^q, \s{K}^{s,\gamma} (X^\Delta))
\end{equation}
which is equal to \eqref{eq.f0}.

For $M \in \got{M}_k$, $k \geq 1$, the weights will have the meaning of tuples
\begin{equation}
\label{eq.wei}
\gamma = (\gamma_1, \ldots , \gamma_k) \in \R^k.
\end{equation}
Here $\gamma_k$ is the `most singular' weight.
For the subspaces $M^{(j)} \in \got{M}_{k-j}$, $j=0, \ldots$, $k-1$, we take the subtuples
\begin{equation}
\label{eq.sub}
\gamma^{(j)} :=
   (\gamma_{j+1}, \ldots , \gamma_k) \in \R^{k-j}.
\end{equation}
The weighted corner space on $M^{(j)}$ of smoothness $s \in \R$ and weight $\gamma^{(j)}$ will be denoted by
$$
\s{H}^{s,\gamma^{(j)}} (M^{(j)}), \ \ \
    j=0, \ldots , k-1.
$$
Occasionally, in order to unify the picture, we also admit the case $\gamma^{(k)}$ as the empty weight tuple  and set in this case
$$
\s{H}^{s, \gamma^{(k)}} (M^{(k)}) := H^s (M^{(k)})
$$
(recall that $M^{(k)} \subset M$ is a $C^\infty$ manifold, cf. Section 3.1).
Knowing the meaning of the spaces 
\begin{equation}
\label{eq.men}
\s{H}^{s,\gamma} (M) \ \ \ \text{for} \ \ \ M \in \got{M}_k
\end{equation}
and of
\begin{equation}
\label{eq.it}
\s{H}^{s,\gamma} (X^\Delta), \s{K}^{s,\gamma} (X^\Delta) \ \ \
\text{for} \ \ \  X \in \got{M}_{k-1},
\end{equation} 
$\gamma =(\gamma_1, \ldots , \gamma_k)$, for a given $k \geq 1$, the question  is how to pass to the corresponding spaces for $k+1$.
An answer is given in \cite{Calv2}, and we briefly describe the result.
We keep in mind the group of isomorphisms 
$\{ \kappa_\lambda \}_{\lambda \in \R_+}$,
on $\s{K}^{s,\gamma} (X^\Delta)$, defined by
\begin{equation}
\label{eq.kappa.}
(\kappa_\lambda u) (r,x) =
   \lambda^{\frac{1+\dim X}{2}} 
     u(\lambda r, x), \ \lambda \in \R_+,
\end{equation}
which allows us to form the spaces \eqref{eq.ed}.
The space \eqref{eq.men} is locally near $Y = M^{(k)}$ modelled on spaces \eqref{eq.ed} when $\dim Y>0$ and on $\s{H}^{s,\gamma}(X^\Delta)$ for $\dim Y =0$ when (for simplicity) $Y$ consists of a single corner point.
The definition of $\s{H}^{s,\gamma} (X^\Delta)$ is as follows:
$$
\s{H}^{s,\gamma} (X^\Delta) :=
   S^{-1}_{\gamma_{k}-\frac{\dim X}{2}} (\s{H}^{s,\gamma'} (\R \times X))
$$
when we write 
$\gamma = (\gamma', \gamma_k)$ for 
$\gamma' := (\gamma_1, \ldots , \gamma_{k-1})$,
and we employ the induction assumption that the cylindrical space 
$\s{H}^{s,\gamma'} (\R \times X)$ is already known.
Here
$$
(S_\beta u) (\bsr, .) :=
     e^{-(\frac{1}{2}-\beta) \bsr} u(e^{-\bsr},.)
$$
for $u(r,.) \in \s{H}^{s,\gamma} (X^\Delta)$, $(r, .)\in \R_+ \times X$, 
$(\bsr, .) \in \R \times X$, $r:= e^{-\bsr}$.
A similar description of  $\s{H}^{s,\gamma} (M)$ holds locally near 
$Y^{(i)} = M^{(j)} \setminus M^{(j+1)}$ for every $j=0, \ldots , k-1$.
Then the space
\begin{equation}
\label{eq.Hj}
\s{H}^{s,\gamma} (M)
\end{equation}
itself may be obtained by gluing together the local pieces by a construction in terms of singular charts and a partition of unity on $M$.
In order to define the spaces
$$
\s{H}^{s,(\gamma, \theta)} (M^\Delta)
$$
with $\theta \in \R$ being the weight belonging to the new axial variable $t \in \R_+$ we need again cylindrical spaces
\begin{equation}
\label{eq.ft}
\s{H}^{s,\gamma} (\R_{\bst} \times M)
\end{equation}
which are locally near any $y \in Y^{(j)}$, $j=0, \ldots, k$, modelled on
$$
\s{W}^s \big(\R_{\bst} \times \R^{\dim Y^{(j)}},
               \s{K}^{s, \gamma^{(j)}} (X^\Delta_{(j-1)}) \big)
\ \ \text{for} \ \ 
X_{(j-1)} \in \got{M}_{j-1}
$$
which refers to the representation of a neighbourhood of $y$ in $M$ as 
\eqref{eq.new}.
Since $Y^{(0)} \in \got{M}_0$, we have the standard cylindrical  Sobolev spaces contributing to \eqref{eq.ft} over 
$\R_{\bst} \times Y^{(0)}$.
By virtue of $M = \bigcup^k_{j=0} Y^{(j)}$ the space $\R \times M$ can be covered by cylindrical neighbourhoods. 
Gluing  together the local spaces, using  singular charts along $M$ and a partition of unity, yields the space \eqref{eq.ft}.
We then set
$$
\s{H}^{s,(\gamma, \theta)} (M^\Delta) :=
    S_{\theta-\frac{\dim M}{2}}
       (\s{H}^{s, \gamma} (\R \times M)),
$$
where $\theta$ plays the role of the new weight $\gamma_{k+1}$, belonging to $t \in \R_+$.
Moreover, by forming arbitrary locally finite sums of elements of $\s{H}^{s,\gamma} (\R \times M)$ with compact support in $\bst \in \R$ we obtain a space that we  denote by
\begin{equation}
\label{eq.locH}
\s{H}^{s,\gamma}_{\loc (\bst)} (\R_{\bst} \times M).
\end{equation}
Now let us form the spaces
$$
\s{W}^s \big( 
          \R^{1+\dim Y^{(j)}}_{t,\wt{y}},
            \s{K}^{s,\gamma^{(j)}}_{\wt{r},x} (X^\Delta_{(j-1)}) \big),
$$
$j=1, \ldots ,k$, and $H^s (\R_t \times \R^{\dim Y^{(0)}})$ for $j=0$.
Moreover, we employ the spaces \eqref{eq.locH} with $t$ instead of $\bst$.
Let $\omega (t)$ be any cut-off function on the half-axis and interpet $1-\omega (t)$ in the following notation as a function in $t \in \R$ that vanishes for $t \leq 0$.
By $\s{W}^{s,\gamma}_\cone (\R_+ \times M)$ we denote the set of all 
$u \in \s{H}^{s,\gamma}_{\loc (t)} (\R_t \times M) |_{\R_{+,t} \times M}$
such that $(1-\omega (t)) u(t, .)$ (for any cut-off function $\omega$)
expressed in local coordinates on $M$ in the wedge $\R \times \R_+ \times X_{(j-1)} \times \R^{\dim Y^{(j)}}_y \ni (t,r,x,y)$, cf. the formula 
\eqref{eq.new}, has the form
$$
v (t,tr,x,ty)
$$
for some $v(t,\wt{r},x, \wt{y}) \in \s{W}^s (\R_{t,\wt{y}}, 
                        \s{K}^{s,\gamma^{(j)}}_{\wt{r},x} (X^\Delta_{(j-1)}))$
for all $j=1, \ldots , k$, and 
$v(t, \wt{y}) \in 
     H^s (\R_t \times \R^{\dim Y^{(0)}})$ for $j=0$.
The invariance of our spaces under (adequate) coordinate transformations is not completely trivial.
We do not deepen this aspect here.
Let us only mention that we have to specify the charts with the 
cocycle of transition maps, and we do not necessarily admit arbitrary isomorphisms of respective local wedges (as manifolds of the corresponding singularity order).

Coming back to Remark \ref{r.inv.}, we could also employ modified definitions of higher wedge spaces, based on the spaces
     $\s{K}^{s,\gamma;g} (X^\Delta) :=
           \langle r \rangle^{-g} \s{K}^{s,\gamma} (X^\Delta)$
with group actions
$$
(\kappa^g_\lambda u) (r,x) =
   \lambda^{g+\frac{1+\dim X}{2}}
        u(\lambda r,x), \ \ \ 
            \lambda \in \R_+,
$$
for $u(r,x) \in \s{K}^{s,\gamma;g} (X^\Delta)$
instead of \eqref{eq.kappa.}.
However, this has a chain of consequences 
which we do not discuss in more detail here.
Now we set
$$
\s{K}^{s,(\gamma,\theta)} (M^\Delta) :=
   \big\{ \omega u + (1-\omega) v :
       u \in \s{H}^{s,(\gamma,\theta)} (M^\Delta),
       v \in \s{W}^{s,\gamma}_\cone (\R_+ \times M) \big\}
$$
for some cut-off function $\omega (t)$.

Summing up we have constructed spaces of the kind \eqref{eq.Hj}, namely,
\begin{equation}
\label{eq.M}
\s{H}^{s,(\gamma,\theta)} (M^\Delta),
   \s{K}^{s,(\gamma,\theta)} (M^\Delta)\ \ \text{for} \ \
      M \in \got{M}_k.
\end{equation}
For arbitrary $N \in \got{M}_{k+1}$ we obtain the spaces
$$
\s{H}^{s,(\gamma,\theta)} (N)
$$
by gluing together local spaces over wedges, similarly as above \eqref{eq.Hj} for the case $M \in \got{M}_k$.
All these spaces are Hilbert spaces with adequate scalar products.
On $\s{K}^{s,(\gamma,\theta)} (M^\Delta)$ we have a strongly continuous group of isomorphisms, similarly as \eqref{eq.kappa.}.
This allows us to form the spaces 
$\s{W}^s (\R^q, \s{K}^{s,(\gamma,\theta)} (M^\Delta))$, and we thus have again the raw material for the next generation of weighted corner spaces.

Let us finally note that the constructions also make sense for non-compact $M \in \got{M}_k$; we assume, for instance, that $M$ is a countable  union of compact sets and that such $M$ are embedded in a compact 
$\wt{M} \in \got{M}_k$.
Then we can talk about $\s{H}^{s,\gamma}_{(\comp)} (M)$, defined to be the set of all elements of $\s{H}^{s,\gamma} (\wt{M})$, supported by a compact subset $K \subseteq M$ and about $\s{H}^{s,\gamma}_{(\loc)} (M)$ to be the set of all locally finite sums of elements in $\s{H}^{s,\gamma}_{(\comp)} (M)$.
The notation `(comp)' and `(loc)' in parentheses is motivated by the fact that, although the distributions are given on $M \setminus M^{(1)}$, the support refers to $M$ (e.g., if $M =\ol{\R}_+, M^{(1)} = \{ 0 \}$, we talk about compact support in $\ol{\R}_+$).    
     
\subsection{Additional edge conditions in higher corner algebras}
\label{s.10.5.2}

As we saw in boundary value problems a basic idea to
complete an elliptic operator $A$ to a Fredholm operator
between Sobolev spaces is to formulate additional boundary
conditions. This can be done on the level of symbols by
filling up the boundary symbol \eqref{1021.Fred0806.eq} to a
family \eqref{1021.sigma(calA)0806.eq} of isomorphisms. In
general it is necessary to admit vector bundles $J_\pm$ on
the boundary, even if the operator $A$ itself is scalar. 
We also can start from operators acting between distributional
sections of vector bundles $E$ and $F$, and we then have
Fredholm operators as in Remark \stref{1021.1506re.re}.

In a similar manner we proceed for a manifold $M$ with  edge
$Y=M^{(1)}$. 
If $A \in \Diff_{\deg}^\mu(\M; E, F)$ is an edge-degenerate operator 
(between weighted edge space of (distributional) sections of 
vector bundles $E, F$) ellipticity requires filling up the 
homogeneous principal edge symbol
\begin{equation}
\label{1052.sigmaland2610.eq}
\sigma_\land(A)(y, \eta) : {\cal K}^{s, \gamma}(X^\Delta,
   E_y) \to {\cal K}^{s- \mu, \gamma - \mu}(X^\Delta, F_y))
\end{equation}   
to a $2 \times 2$ block matrix family of isomorphisms
\begin{equation}
\label{1052.Asigma2610.eq}
\sigma_\land({\cal A})(y, \eta) :
   {\cal K}^{s, \gamma}(X^\Delta, E_y)  \oplus  J_{-, y}      \to
   {\cal K}^{s- \mu, \gamma- \mu}(X^\Delta, F_y) \oplus   J_{+,y}
\end{equation}
for suitable $J_\pm \in \Vect(Y)$, $(y, \eta) \in T^* Y
\setminus 0$.		
	
Here, in abuse of notation, $E_y, F_y \in \Vect(X^\land)$
denote bundles that  are obtained as follows. 
First consider the $X$-bundle $\M_{\sing}$ over $Y$ and the
associated $X^\land$-bundle $\M_{\sing}^\land$, with the
canonical projection $p : \M_{\sing}^\land \to \M_{\sing}$,
induced by $X^\land \to X$. For every $E \in \Vect(\M)$ we
obtain a bundle $p^*(E \big|_{\M_{\sing}})$; then the
restriction of the latter bundle to the fibre of
$\M_{\sing}^\land$ over $y \in Y$ is  denoted again 
by $E_y$.	

The construction of \eqref{1052.Asigma2610.eq} for a
$\sigma_\psi$-elliptic element $A$ of the edge algebra 
is also meaningful for pseudo-differential
operators. 
The weight $\gamma \in \R$ is kept fixed; in general there
are many admissible weights for which
\eqref{1052.sigmaland2610.eq} is a Fredholm family, but the
bundles $J_\pm$ depend on $\gamma$.

In any case we obtain $2 \times 2$ block matrices
\begin{equation}
\label{1052.calA2610.eq}
{\cal A} :
  \Hsum{{\cal H}^{s, \gamma}(M, E)}
       {H^s(Y, J_-)}       \to
  \Hsum{{\cal H}^{s- \mu, \gamma - \mu}(M, F)}
       {H^{s- \mu}(Y, J_+)}
\end{equation}
which are Fredholm operators as soon as ${\cal A}$ is
elliptic with respect to $\sigma_\psi({\cal A})$ and
$\sigma_\land({\cal A})$; the latter condition is just the
bijectivity of \eqref{1052.Asigma2610.eq} for all $(y, \eta)
\in T^* Y \setminus 0$ for some $s \in \R$.

Let $\got{A}^\mu(M, \bs{g})$ denote the space of all
\eqref{1052.calA2610.eq} belonging to the weight data
$\bs{g} := (\gamma, \gamma - \mu)$ (the bundles are assumed
to be known for every concrete ${\cal A}$, otherwise we
write $\got{A}^\mu(M, \bs{g}; \bs{v})$ for $\bs{v} := (E, F;
J_-, J_+)$).            

On a manifold $M$ in the category $\got{M}_k$ with the
sequence of subspaces \eqref{1031.eq1210.eq}, $M^{(j)} \in
\got{M}_{k-j}$, $0 \leq j \leq k$, the picture (in simplified
form) is a follows.
We have weighted Sobolev spaces
\[  
{\cal H}_{(\comp)}^{s, \gamma^{(j)}}(M^{(j)}, E^{(j)}) \;
      \textup{and ${\cal H}_{(\loc)}^{s,
      \gamma^{(j)}}(M^{(j)}, E^{(j)})$}      
\]
with weights 
$\gamma^{(j)} = (\gamma_{j+1}, \ldots, \gamma_{k}) \in \R^{k-j}$ for $0 \leq j \leq k-1$ and 
$H^s_{(\comp/\loc)} (M^{(k)}, E^{(k)})$ 
for vector bundles $E^{(j)}$. 
(For $M^{(j)}$ compact we omit subscripts
`$(\comp)$' and `$(\loc)$'.)

The higher corner operator space $\got{A}^\mu(M, \bs{g};
\bs{v})$ of operators of order $\mu$ on $M$ then consists of
$(k+1) \times (k+1)$ block matrices
\begin{equation}
\label{1052.calAE2710.eq}
{\cal A} : \bigoplus_{j=0}^{k} {\cal H}_{(\comp)}^{s,
     \gamma^{(j)}} (M^{(j)}, E^{(j)})\to
     \bigoplus_{l=0}^{k} {\cal H}_{(\loc)}^{s-\mu,
     \gamma^{(l)}- \mu}(M^{(l)}, E^{(l)}),
\end{equation}
with weight data 
$\bs{g} := (\gamma^{(j)}, \gamma^{(j)} - \mu)_{j = 0, \ldots, k-1}$, $\gamma^{(j)} - \mu := (\gamma_{j+1} -\mu, \ldots, \gamma_{k}- \mu)$, 
and tuples 
$\bs{v} := (E^{(j)}, F^{(j)})_{j = 0, \ldots, k}$, 
$E^{(j)},F^{(j)} \in \Vect(M^{(j)})$. 
The fibre dimensions of the involved bundles may be zero. 
In that case the corresponding spaces are omitted.
Writing
$$
 {\cal A} = ({\cal A}_{ij})_{i,j = 0, \ldots, k}   
$$
we have
$$
  ({\cal A}_{ij})_{i,j = l, \ldots, k} \in
      \got{A}^\mu(M^{(l)}, \bs{g}^{(l)}; \bs{v}^{(l)})    
$$   
for every $0 \leq l \leq k$, 
with weight and bundle data $\bs{g}^{(l)}$ and $\bs{v}^{(l)}$,
respectively, that follow from $\bs{g}$ and $\bs{v}$ by
omitting corresponding components. 
Set $\ulc \s{A} := (\s{A}_{ij})_{i,j=0, \ldots ,k-1}$.
The principal symbolic hierarchy 
\begin{equation}
\label{1052.sigma2710.eq}
\sigma({\cal A}) := (\sigma_j({\cal A}))_{0 \leq j \leq k}
\end{equation}
is defined inductively, where 
$(\sigma_j({\cal A}))_{0 \leq j \leq k-1}$ is the
principal symbol of $\ulc {\cal A} \big|_{M \setminus M^{(k)}}$ with
$M \setminus M^{(k)} \in \got{M}_{k-1}$, such that the symbols
up to the order $k-1$ are known, while
\begin{equation}
\label{1052.sigmaland2710.eq}
\sigma_k({\cal A})(y, \eta) : 
  \Hsum{\oplus_{j=0}^{k-1} 
         {\cal K}^{s, \gamma^{(j)}} 
          ( (X_{k-1}^{(j)})^\Delta, E_y^{(j)} )}
       {E_y^{(k)}}     \to
  \Hsum{\oplus_{l=0}^{k-1} {\cal K}^{s- \mu,
                \gamma^{(l)}-\mu} ( (X_{k-1}^{(l)})^\Delta, F_y^{(l)} )}
       {F_y^{(k)}}
\end{equation}
for $(y, \eta) \in T^* M^{(k)} \setminus 0$ is the highest
principal symbol of order $k$, cf.
\eqref{1031.neusigma2610.eq} for the case when ${\cal A}$ 
consists of an upper left corner which is a differential
operator $A$.
Here $X_{k-1} \in \got{M}_{k-1}$ (by assumption,
compact) is the fibre of the $X_{k-1}$-bundle $\M_{\sing}$ over
$M^{(k)} \ni y, \; M^{(k)} \in \got{M}_0$, and
\[  X_{k-1}  = X_{k-1}^{(0)} \supset X_{k-1}^{(1)} \supset 
      \ldots \supset X_{k-1}^{(k-1)}		   \]   
is the chain of subspaces, analogously as
\eqref{1031.eq1210.eq}. 
In this discussion we tacitly	
assumed $\dim M^{(k)} > 0$. In the case $\dim M^{(k)} = 0$ 
the space $M^{(k)}$ consists of
corner points. Then \eqref{1052.sigmaland2710.eq} is to be
replaced by an analogue of the former conormal symbols, 
namely,	      
\begin{equation}
\label{1052.con2710.eq}
\sigma_\c({\cal A})(w) :
   \Hsum{\oplus_{j=0}^{k-2} {\cal H}^{s, \gamma^{\prime (j)}}
           (X_{k-1}^{(j)}, E^{(j)} )}
        {E^{(k-1)}}      \to
   \Hsum{\oplus_{l=0}^{k-2} {\cal H}^{s- \mu, \gamma^{\prime (l)}-
           \mu} ( X_{k-1}^{(l)}, E^{(l)} )}
	{F^{(k-1)}},   	
\end{equation}
$w \in \Gamma_{\frac{\dim X_{k-1}+ 1}{2}- \gamma_k}$ where 
$\gamma' =(\gamma_1, \ldots , \gamma_{k-1})$.
Clearly, \eqref{1052.con2710.eq}
depends on the discrete corner points $y \in M^{(k)}$; for
simplicity, we assume that $M^{(k)}$ consists of a single
point (subscripts `$y$' are then omitted). 

\begin{Definition}
\label{1052.def2710.de}
An operator ${\cal A} \in \got{A}^\mu(M, \bs{g}; \bs{v})$ for
$M \in \got{M}_k$ is called elliptic of order $\mu$, if 
${\cal
A}' := {\cal A} \big|_{M \setminus M^{(k)}}$ is elliptic as an
element of $\got{A}^\mu(M \setminus M^{(k)}, \bs{g};
\bs{v}')$ for $\bs{g}' := (\gamma^{(j)}, \gamma^{(j)} - 
\mu)_{j
= 1, \ldots, k-1}$, $\bs{v}' = (E^{(j)}, F^{(j)})_{j = 1,
\ldots, k-1}$, and if \eqref{1052.sigmaland2710.eq} for $\dim
M^{(k)} > 0$ 
is a family of isomorphisms for
all $(y, \eta) \in T^* M^{(k)} \setminus 0$
{\st{(}}or \eqref{1052.con2710.eq} for $\dim M^{(k)} = 0$, for
all $w \in \Gamma_{\frac{\dim X_{k-1} 
+ 1}{2} - \gamma_k}${\st{)}}. 
\end{Definition}

\begin{Theorem}
\label{1052.the2710.th}
Let ${\cal A} \in \got{A}^\mu(M, \bs{g}; \bs{v})$ be elliptic
and $M \in \got{M}_k$ compact. Then \eqref{1052.calAE2710.eq}
is a Fredholm operator for every $s \in \R$, and ${\cal A}$ 
has a parametrix ${\cal P} \in \got{A}^{-\mu} (M, \bsg^{-1}; \bsv^{-1})$.
\end{Theorem}

\subsection{A hierarchy of topological obstructions}
\label{s.10.5.3}

%10.5.3
 
Looking at the constructions of Section 2.1 for an elliptic
operator $A$ on $X$ in connection with the process of filling up
the Fredholm family \eqref{1021.Fred0806.eq} to a family of
isomorphisms \eqref{1021.sigma(calA)0806.eq} we did not
emphasise that the existence of the vector bundles $J_\pm \in
\Vect(\partial X)$ is by no means automatic. To illustrate that
we first recall the homogeneity
\[  \sigma_\partial(A)(y, \lambda \eta) = \lambda^\mu
      \kappa_\lambda \sigma_\partial(A)(y, \eta)
      \kappa_\lambda^{-1}     \]
for all $\lambda \in \R_+, (y, \eta) \in T^*(\partial X)
\setminus 0$ which shows that we may consider
\eqref{1021.Fred0806.eq} for $(y, \eta) \in S^*(\partial X)$
(the unit cosphere bundle induced by $T^*(\partial X)$) which
is a compact topological space (when $X$ compact).
It suffices to construct \eqref{1021.sigma(calA)0806.eq}
first for $(y, \eta) \in S^*(\partial X)$ and then to extend
it by homogeneity $\mu$  to arbitrary $(y, \eta) \in
T^*(\partial X) \setminus 0$, setting
\begin{equation}
\label{1053.hom2810.eq}
\sigma_\partial({\cal A})(y, \eta) := | \eta |^\mu
    \begin{pmatrix}
    \kappa_{| \eta |} & 0  \\  0 & 1
    \end{pmatrix} \sigma_\partial({\cal A})(y, \frac{\eta}{|
    \eta |})
    \begin{pmatrix}
    \kappa_{| \eta |} & 0  \\  0 & 1
    \end{pmatrix}^{-1},
\end{equation}
cf. the relation \eqref{1021.ho0906.eq}. From the fact that
\eqref{1021.Fred0806.eq} is a family of Fredholm operators,
parametrised by the compact space $S^*(\partial X)$ we have
a $K$-theoretic index element
\begin{equation}
\label{1053.ind2810.eq}
\ind_{S^*(\partial X)} \sigma_\partial(A) \in K(S^*(\partial X)).
\end{equation}
Recall that the $K$-group $K(.)$ (for a compact topological
space in the parentheses) is a group of equivalence classes
$[G_+] - [G_-]$ of pairs $(G_-, G_+)$ of vector bundles $G_-,
G_+ \in \Vect(.)$.         

If \eqref{1021.sigma(calA)0806.eq} is a family of isomorphisms,
the index element \eqref{1053.ind2810.eq} is  equal to $[J_+]
- [J_-]$ which means
\begin{equation}
\label{1053.pull2810.eq}
\ind_{S^*(\partial X) } \sigma_\partial(A) \in \pi_1^* K(\partial X)
\end{equation}
when $\pi_1 : S^*(\partial X) \to \partial X$ denotes the
canonical projection. In general, we only have
\eqref{1053.ind2810.eq}, i.e., \eqref{1053.pull2810.eq} is a
topological obstruction for the existence of an elliptic
operator \eqref{1021.eq0906.eq} with $A$ in the upper left
corner. The condition has been studied in \stcite{Atiy5} for
elliptic differential operators $A$ and in \stcite{Bout1} for
pseudo-differential operators $A$ with the transmission
property at $\partial X$. Dirac operators (in even dimensions)
and other interesting geometric operators belong the cases
where this obstruction does not vanish.

The following discussion is partly  hypothetical, it formulates
expectations that are not completely worked out in detail, 
except
for the obvious things such as the following observation.

Set
\begin{equation}
\label{1053.K2810.eq}
\bs{\cal K}^{s, \gamma}(X^\Delta, E_y) := \oplus_{j=0}^{k-1}
    {\cal K}^{s, \gamma^{(j)}}((X_{k-1}^{(j)})^\Delta, 
      E_y^{(j)}),
\end{equation}
\begin{equation}
\label{1053.L2810.eq}
\bs{\cal K}^{s- \mu, \gamma - \mu}(X^\Delta, F_y) :=
    \oplus_{l=0}^{k-1} {\cal K}^{s- \mu, \gamma^{(l)}-\mu}
((X_{k-1}^{(l)})^\Delta, F_y^{(l)}).
\end{equation}
Moreover, let $\ulc {\cal A} := ({\cal A}_{ij})_{i,j = 0,
\ldots, k-1}$.
If \eqref{1052.sigmaland2710.eq} is a family of isomorphisms,
the $k \times k$ upper left corner
\begin{equation}
\label{1053.Fre2810.eq}
\sigma_\land (\ulc {\cal A}) (y, \eta) : \bs{\cal K}^{s,
   \gamma}(X^\Delta, E_y) \to \bs{\cal K}^{s- \mu, \gamma -
   \mu}(X^\Delta, F_y)
\end{equation}
is a family of Fredholm operators for all $(y, \eta) \in T^*
M^{(k)} \setminus 0$. In addition, using the natural group 
actions $\{ \bs{\kappa}_\lambda \}_{\lambda \in \R_+}$ and
$\{ \bs{\tilde{\kappa}}_\lambda \}_{\lambda \in \R_+}$  
on the spaces
\eqref{1053.K2810.eq} and \eqref{1053.L2810.eq}, respectively,
we have the homogeneity
\[ \sigma_\land(\ulc {\cal A})(y, \lambda \eta) = \lambda^\mu
         \bs{\tilde{\kappa}}_\lambda \sigma_\land(\ulc 
	 {\cal A})(y, \eta) \bs{\kappa}_\lambda^{-1}    \]      
for all $\lambda \in \R_+, \ (y, \eta) \in T^* M^{(k)}
\setminus 0$. This allows us to interpret
\eqref{1053.Fre2810.eq} as a Fredholm family on $S^* M^{(k)}$,
and we have
\[  \ind_{S^*M^{(k)}} \sigma_\land(\ulc {\cal A}) \in K(S^*
         M^{(k)}).     \]	 
It is again a necessary and sufficient condition for the
existence of a block matrix family \eqref{1052.sigmaland2710.eq}
of isomorphisms with vector bundles $E^{(k)}, \ F^{(k)} \in
\Vect(M^{(k)})$ that
\begin{equation}
\label{1053.top2810.eq}
\ind_{S^*M^{(k)}} \sigma_\land(\ulc {\cal A}) \in \pi_1^*
     K(M^{(k)}),	 
\end{equation}
$\pi_1 : S^*M^{(k)} \to M^{(k)}$.     

If \eqref{1053.top2810.eq} holds we find the additional
entries for $\sigma_\land({\cal A})(y, \eta)$ in the $k$ th
row and column, first for $(y,
\eta) \in S^* M^{(k)}$ and then for all $(y, \eta) \in T^*
M^{(k)} \setminus 0$ by an extension by homogeneity, similarly
as \eqref{1053.hom2810.eq}.

The condition
\eqref{1053.top2810.eq} is a topological obstruction for the
existence of an elliptic element ${\cal A} \in \got{A}^\mu(M,
\bs{g}; \bs{v})$ for a given  operator of the form $\ulc {\cal A} \in \got{A}^\mu(M,
\bs{g}; \bs{w})$, $\bs{w} := (E^{(j)}, F^{(j)})_{j=0, \ldots,
k-1}$, that is elliptic with respect to  $(\sigma_j(.))_{j= 0, \ldots, k-1}$. 

If \eqref{1053.top2810.eq} is violated, it should be possible to modify the procedure of filling up \eqref{1053.Fre2810.eq} to a family of isomorphisms \eqref{1052.sigmaland2710.eq} by completing $\luc \s{A}$ to a Fredholm operator $\s{A}$, by using global projection data in analogy of the constructions of \cite{Schu42} for the case $k=1$ (see also \cite{Schu37}, \cite{Schu44} for the case boundary value problems).
The extra entries of $\s{A}$ (compared with $\ulc \s{A}$) then refer to subspaces of the standard Sobolev spaces on $M^{(k)}$ which are the image under a pseudo-differential projection.
In the opposite case, i.e., when \eqref{1053.top2810.eq} holds we obtain ellipticity of $\s{A}$ in the sense of Definition \ref{1052.def2710.de} which is an analogue of the Shapiro-Lopatinskij ellipticity, known from boundary value problems.
In that case it is interesting to talk about different possibilities of filling up the operator $\s{A}_{11}$ to a Fredholm operator $\s{A}$ in the above mentioned way.
Let $\s{B}$ be another operator containing $\s{A}_{11}$ in the upper left corner, and let $\s{B}$ be also elliptic in the sense of Definition \ref{1052.def2710.de}.
There is then a reduction of the conditions 
$(\s{B}_{ij})_{i,j=0,\ldots ,k, (i,j) \not= (1,1)}$ to the subspace $M^{(1)} \in \got{M}_{k-1}$ by means of 
$\s{A} =(\s{A}_{ij})_{i,j=0, \ldots ,k}$.
The algebraic process is similar to that in \cite[Section 3.2.1.3]{Remp2} for the case of boundary value problems.
In other words, there exists an elliptic element  $\s{R} \in \got{A}^0 (M^{(1)})$
(for brevity, weight and bundle data are omitted in the notation) such that
$$
\ind \s{B} - \ind \s{A} = \ind \s{R}.
$$
 This relation is an analogue of the Agranovich-Dynin formula, cf. Remark \ref{1021.red1506.re}, and
\cite[Section 3.2.1.3]{Remp2}.

The latter observation can also be interpreted as follows.
The elliptic operators on $M^{(1)}$ parametrise the elliptic operators in $\got{A}^\mu (M)$, apart from the ellipticity condition  for $\s{A}_{11}$ itself (which means, e.g., for $k=1$, that $\s{A}_{11}$ is elliptic of Fuchs type or in the edge-degenerate sense). 

\subsection{The building of singular algebras}
\label{s.10.5.4}

%10.5.4
If $M \in \got{M}_k$ is given, we assume to have constructed an algebra of operators $\got{A} (M) := \bigcup_\mu \got{A}^\mu (M)$ for 
$\got{A}^\mu (M) := \bigcup_{\bsg,\bsv} \got{A}^\mu (M,\bsg;\bsv)$,
cf. the notation in Section 5.2, with a principal symbolic structure 
$(\sigma_j (\s{A}))_{0 \leq j \leq k}$.
For $M \in \got{M}_0$ we may take, for instance, the algebra of classical pseudo-differential operators on $M$.
The program of the iterative calculus on $\got{M}_{k+1}, \got{M}_{k+2}, \ldots$, is to organise a natural scenario to pass from $\got{A} (M)$ to corresponding higher  generations of calculi.
Spaces in $\got{M}_{k+1}$ can be obtained from $M \in \got{M}_k$ by pasting together local cones $M^\Delta$ or wedges $M^\Delta \times \Omega$, $\Omega \subseteq \R^q$ open.
Analytically, the main steps (apart from invariance aspects) consist of understanding the correspondence between $\got{A} (M)$ and the next higher algebras
$$
\got{A} (M^\Delta) \ \ \ \text{and} \ \ \ \got{A} (M^\Delta \times \Omega).
$$
The way which is suggested here will be called conification and edgification of the calculus on $M$.
The experience from the cone and edge algebras of first generation leads to the following ingredients.        

{\bf (C.1) Parameter-dependent calculus}.
Establish $\got{A} (M; \R^l)$, a parameter-depen\-dent version
 of $\got{A}(M)$ with parameters $\lambda = (\lambda_1, \ldots , \lambda_l) \in \R^l$ of dimension $l \geq 1$.
Here $\zeta := (\lambda_2, \ldots , \lambda_l) \in \R^{l-1}$ may be treated as sleeping parameters in the sense of Section 3.2.
In the process of the iterative construction it becomes clear how the parameters are successively activated, cf. the points (C.2) - (C.4) below.
In this context we assume that $\got{A} (M; \R^l)$ is constructed for every $M \in \got{M}_k$; thus since $M^\wedge = \R_+ \times M$ also belongs to $\got{M}_k$ (with $\R_+$ being regarded as a $C^\infty$ manifold) we also have $\got{A} (M^\wedge)$ and $\got{A} (M^\wedge; \R^l)$.
If $\got{A}^{-\infty} (M)$ denotes the space of smoothing elements in the algebra $\got{A} (M)$ (defined by their mapping properties in weighted corner spaces), we set 
$$
\got{A}^{-\infty} (M; \R^l) := \s{S} (\R^l, \got{A}^{-\infty} (M))
$$
for every $M \in \got{M}_k$.  

{\bf (C.2) Holomorphic Mellin symbols and kernel cut-off}.
Generate an analogue of $\got{A} (M; \R^l)$, namely, $\got{A} (M; \C \times \R^{l-1})$ of holomorphic families in the complex parameter $v \in \C$ by applying a kernel cut-off procedure to elements of $\got{A} (M; \R^l)$ with respect to $\lambda_1$.
Here $\s{A} (v, \zeta) \in \got{A} (M; \C \times \R^{l-1})$ is holomorphic in $v = \beta + i\tau \in \C$ with values in $\got{A} (M; \R^{l-1}_\zeta)$ such that
$$
\s{A} (\beta + i \tau, \zeta) \in
    \got{A} (M; \R^l_{\tau,\zeta})
 $$
for every $\beta \in \R$, uniformly in finite $\beta$-intervals.
The holomorphy of operator families can be defined in terms of holomorphic families of the underlying local symbols (the notion directly follows in terms of the spaces of symbols, plus holomorphic families of smoothing operators which is also an easy notion, taking into account their mapping properties between global weighted spaces, or subspaces with asymptotics).
In a similar sense we can form the spaces
$$
C^\infty \big(\ol{\R}_+ \times \Xi, \got{A} (M; \R^l) \big)   \ \ \ \text{and} \ \ \
C^\infty \big(\ol{\R}_+ \times \Xi, \got{A} (M; \C \times \R^{l+1}) \big),
$$
respectively.         

{\bf (C.3) Mellin quantisation}. Given a 
$$
\wt{p} (t,z,\wt{\tau}, \wt{\zeta}) \in
   C^\infty 
         \big( \ol{\R}_+ \times \Xi,
               \got{A} (M; \R^l_{\wt{\tau}, \wt{\zeta}}) \big)
$$
we find an
$$
\wt{h} (t,z,v, \wt{\zeta}) \in
   C^\infty 
       \big( \ol{\R}_+ \times \Xi,
           \got{A} (M; \C \times \R^{l-1}_{\wt{\zeta}}) \big)
$$
such that for
$p(t,z,\tau,\zeta) := \wt{p} (t,z,t\tau, t\zeta)$,
$h(t,z,v,\zeta) := \wt{h} (t,z,v,t \zeta)$
we have 
$$
\op^\gamma_M (h) (z,\zeta) = \op_t (p) (t,\zeta)  \ \
\mod C^\infty (\Omega, \got{A}^{-\infty} (M^\wedge; \R^{l-1}))
$$ 
for every $\gamma \in \R$.
The correspondence $p \to h$ may be achieved by a combination of a transformation from the Fourier phase function $(t-t') \tau$ to the Mellin phase function $(\log t' - \log t) \tau$ with a kernel cut-off construction.   

{\bf (C.4) Edge quantisation}.
We start from a family
$\wt{p} (t,z,\wt{\tau}, \wt{\zeta}) \in
   C^\infty 
         \big( \ol{\R}_+ \times \Xi,
               \got{A} (M;$ $\R^{1+q}_{\wt{\tau}, \wt{\zeta}}) \big)$
               for $\Xi \subseteq \R^q$ open and obtain
$$
p(t,z,\tau, \zeta) := \wt{p} (t,z,t\tau, t\zeta), \ \ \ 
p_0(t,z,\tau, \zeta) := \wt{p} (0,z,t\tau, t\zeta), 
$$
$$
h(t,z,v, \zeta) := \wt{h} (t,z,v, t\zeta), \ \ \ 
h_0(t,z,v, \zeta) := \wt{h} (0,z,v, t\zeta)
$$
by Mellin quantisation.
Moreover, we fix cut-off functions 
$\omega$, $\wt{\omega}$, $\wt{\wt{\omega}}$ such that 
$\wt{\omega} \equiv 1$ on supp $\omega$, $\omega \equiv 1$ on supp $\wt{\wt{\omega}}$, and cut-off  functions
$\sigma, \wt{\sigma}$.
We set
\begin{equation}
\label{eq.294}
a_M (z, \zeta) := t^{-\mu} \omega 
       (t[\zeta]) \op^{\Theta-\frac{n}{2}}_M
           (h) (z, \zeta) \wt{\omega} (t' [\zeta])
\end{equation}
for a $\theta \in \R$ and $n=\dim M$,
\begin{equation}
\label{eq.295}
a_\psi (z, \zeta) := t^{-\mu} (1-\omega  (t[\zeta]))
                     \omega_0(t[\zeta],t'[\zeta]) \op_t (p) (z,\zeta) 
                       (1- \wt{\wt{\omega}}  (t'[\zeta]))
\end{equation}
($t' \in \R$ is the variable under the Mellin transform),
$\omega_0 (t,t') := \psi \Big(
                           \frac{(t-t')^2}{1+ (t-t')^2} \Big)$
for every $\psi \in C^\infty_0 (\ol{\R}_+)$
such that $\psi(t) =1$ for $t < \frac{1}{2}$, $\psi (t) =0$ for $t> \frac{2}{3}$,
cf. \cite[Lemma 2.10]{Calv1}, and form the operator-valued amplitude function
\begin{equation}
\label{eq.296}
a(z, \zeta) := \sigma \{ a_M (z,\zeta) + a_\psi (z,\zeta) \} \wt{\sigma}
\end{equation}
which belongs to
$S^\mu \big(\Xi \times \R^d; \bssK^{s,(\gamma,\theta)} (M^\Delta),
                         \s{K}^{s-\mu,(\gamma-\mu,\theta-\mu)} (M^\Delta)\big)$
for every $s \in \R$.
As above, $\theta \in \R$ plays the role of the additional weight 
$\gamma_{k+1}$.                

{\bf (C.5) Mellin plus Green symbols.}
Compositions of symbols of the kind \eqref{eq.296} and computations in connection with ellipticity and parametrices generate a further class of symbols, namely, Mellin plus Green symbols
\begin{equation}
\label{eq.297}
m(z, \zeta) + g(z,\zeta) \in
      S^\mu_\cl \big(
          \Xi \times \R^d; 
              \bssK^{s,(\gamma,\theta)} (M^\Delta),
              \bssS^{(\gamma-\mu,\theta-\mu)} (M^\Delta)\big).
\end{equation}
Here 
$$
\bssS^{(\gamma,\theta)} (M^\Delta ) :=
     \lim\limits_{\substack{\longleftarrow \\ N \in \N}}
         \langle t \rangle^{-N} 
              \bssK^{N,(\gamma,\theta)} (M^\Delta).
$$
There are many possible variants for \eqref{eq.297}, for instance, symbols referring to the weight line 
$\Gamma_{\frac{n+1}{2}-\theta}$ itself, $n=\dim M$, or to an $\varepsilon$-strip around this for a small $\varepsilon >0$, or to a larger strip in the complex $v$-plane in which asymptotic phenomena are encoded in terms of meromorphy.
Let as content ourselves here with the $\varepsilon$-strip.
In that case we choose a function $f (z,v)$ which is $C^\infty$ in $z \in \Xi$ and holomorphic in 
$\{ v \in \C : \frac{n+1}{2} - \theta - \varepsilon 
                 < \re v< \frac{n+1}{2} - \theta + \varepsilon \}$,
taking values in $\got{A}^{-\infty} (M)$, such that
$f(z,v) \in
      C^\infty (\Xi, \got{A}^{-\infty} (M; \Gamma_\beta))$ 
for  every
$\frac{n+1}{2} -\theta -\varepsilon < \beta < \frac{n+1}{2} - \theta + \varepsilon$,
uniformly in compact $\beta$-intervals.
Then we set 
\begin{equation}
\label{eq.298}
m(z,\zeta) := t^{-\mu} 
    \omega (t[\zeta])
       \op^{\theta-\frac{n}{2}}_M (f) (z)
          \wt{\omega} (t' [\zeta])
\end{equation}
for an arbitrary choice of cut-off functions $\omega, \wt{\omega}$.
A Green symbol $g(z,\zeta)$ is defined by 
\begin{equation}
\label{eq.299}
g(z, \zeta)  \in
      S^\mu_\cl \big(
          \Xi \times \R^d; 
              \bssK^{s,(\gamma,\theta)} (M^\Delta),
              \bssS^{(\gamma-\mu+\delta,\theta-\mu+\delta)} (M^\Delta)\big).
\end{equation}
for some $\delta >0$, for all $s$, together with a similar condition on the $(z,\zeta)$-wise formal adjoints.
Varying $\omega, \wt{\omega}$ in \eqref{eq.298} we only obtain a Green remainder.
The symbols \eqref{eq.299} take values in compact operators 
$ \bssK^{s,(\gamma,\theta)} (M^\Delta) \to
              \bssK^{s,(\gamma-\mu,\theta-\mu)} (M^\Delta)$;
the operator-valued symbols \eqref{eq.298} have not such a property.
We have
$$
\Op_z (m), \ \Op_z (g) \in \got{A}^{-\infty} (M^\wedge).
$$
Recall that $M^\wedge = \R_+ \times M$ belongs to $\got{M}_k$; therefore, the smoothing operators on $M^\wedge$ are already known by induction.
Nevertheless, $\Op_z(m), \Op_z (g)$ take part as non-smoothing contributions in the algebra $\got{A} (M^\Delta)$, cf. (C.7) below.          

{\bf (C.6) Global smoothing operators}.
Formulate the spaces 
$\got{A}^{-\infty} (N) \ni C$ for arbitrary $N \in \got{M}_{k+1}$ by requiring the mapping properties
\begin{equation}
\label{eq.300}
C : \s{H}^{s,(\gamma, \theta)}_{(\comp)} (N) \to
    \s{H}^{\infty,(\gamma-\mu+\delta, \theta-\mu+\delta)}_{(\loc)} (N),
\end{equation}
$s \in \R$, for some $\delta = \delta(C) >0$ and, analogously, for the formal adjoints $C^*$.
Here we fix weight data $((\gamma, \theta), (\gamma-\mu, \theta-\mu))$ for arbitrary weights and orders $\mu$ (the spaces in the relation \eqref{eq.300} are an abbreviation for the direct sums occurring before, with $k+1$ instead of $k$).                     

{\bf (C.7) Global corner operators of ${\bf (k+1)}$-th generation}.
An operator $A \in \got{A}^\mu (N)$ for $N \in \got{M}_{k+1}$, associated with weight data 
$((\gamma, \theta), (\gamma-\mu, \theta-\mu))$ is defined as follows:

We first choose cut-off functions $\sigma, \wt{\sigma}, \wt{\wt{\sigma}}$ on $N$ that are equal to $1$ in a small neighbourhood of $Z := N^{(k+1)}$ and vanish outside another such neighbourhood, such that $\wt{\sigma} =1$ on $\supp \sigma$, $\sigma =1$ on $\supp \wt{\wt{\sigma}}$.
Then $\got{A}^\mu (N)$ consists of all
$$
A = A_\sing + A_\reg +C
$$
such that
\begin{itemize}
\item [(i)]
$C \in \got{A}^{-\infty} (N)$;
\item[(ii)]
$A_\reg := (1 -\sigma) A_{\Int} (1 - \wt{\wt{\sigma}})$ for
$A_{\Int} := A |_{N \setminus Z} \in \got{A}^\mu (N \setminus Z)$,
also associated with the weight data $(\gamma, \gamma -\mu)$;
\item[(iii)]
$A_\sing$ (modulo pull backs to the manifold) is a locally finite sum of operators of the form
$$
\varphi \{ \Op_z (a+m+g) \} \psi
$$
referring to the local description of $N$ near $Z$ by wedges $M^\Delta \times \Xi$, $\Xi \subseteq \R^q$ open $(d= \dim Z)$, for arbitrary symbols $a,m,g$ as in (C.4), (C.5) and functions 
$\varphi, \psi \in C^\infty_0 (\Xi)$, $\varphi$ belonging to a partition of unity on $Z$ and 
$\psi \equiv 1$ on $\supp \varphi$.
\end{itemize}

{\bf (C.8) The principal symbolic hierarchy}.
For $A \in \got{A}^\mu (N)$ we set 
$$
\sigma (A) := (\sigma_{\Int} (A), \sigma_{\wedge_{k+1}}(A)),
$$
where
$\sigma_{\Int} (A) = \sigma (A_{\Int})$ 
is the symbol which is known from the step before, since 
$N \setminus Z \in \got{M}_k$, and 
\begin{eqnarray}
\label{eq.301}
\sigma_{\wedge_{k+1}} (A) (z, \zeta)   
&:=&
t^{-\mu} \{ \omega (t |\zeta|) 
   \op^{\theta - \frac{n}{2}}_M (h_0) (z,\zeta)
        \wt{\omega} (t |\zeta|)                    \nonumber \\
& + &
(1-\omega (t |\zeta|)) 
  \omega_0 (t |\zeta|, t' |\zeta|)
     \op_t (p_0) (z, \zeta) 
       (1-\wt{\wt{\omega}} (t' |\zeta|)) \}    \nonumber  \\
& +&
\sigma_{\wedge_{k+1}} (m+g) (z, \zeta),
\end{eqnarray}
where $\sigma_{\wedge_{k+1}} (m+g) (z,\zeta)$ is the homogeneous principal part of $m+g$ in the sense of (operator-valued) classical symbols of order $\mu$.
The edge symbol \eqref{eq.301} is interpreted as a family of operators
\begin{equation}
\label{eq.new98}
\sigma_{\wedge_{k+1}} (A) (z, \zeta) :
   \bssK^{s,(\gamma,\theta)} (M^\Delta) \to
    \bssK^{s-\mu,(\gamma-\mu,\theta-\mu)} (M^\Delta),
\end{equation}
$(z,\zeta) \in T^* Z \setminus 0$, and we have
$$
\sigma_{\wedge_{k+1}} (A) (z, \lambda \zeta) =
   \lambda^\mu \kappa_\lambda \sigma_{\wedge_{k+1}} (A)) (z, \zeta) \kappa^{-1}_\lambda
$$
for all $\lambda \in \R_+$.

\section{Historical background and future program}   
\label{s.10.6}  

\begin{minipage}{\textwidth}
\setlength{\baselineskip}{0cm}
\begin{scriptsize}    
The analysis on manifolds with singularities has a long
history. Motivations and models from the applied sciences go
back to the 19 th century. 
There are deep connections with 
pure mathematics, e.g., complex analysis, geometry, and topology. 
Numerous authors have contributed to the field.
We outline here a few aspects of the development and sketch
some challenges and open problems.
\end{scriptsize}
\end{minipage}

\subsection{Achievements of the past development} 
\label{s.10.6.1}

%10.6.1

The analysis on manifolds with singularities is inspired by ideas and achievements from classical areas of mathematics, such as singular integral operators, Toeplitz operators, elliptic boundary value problems, Sobolev problems, from the applied sciences with edge and corner geometries, crack problems, numerical computations, pseudo-differential calculus, asymptotic analysis and Mellin operators with meromorphic symbols, parameter-dependent ellipticity, spectral theory, ellipticity on non-compact manifolds, expecially, with conical exits to infinity, Dirac operators and other geometric operators, Hodge theory, index theory, spectral theory, functional calculus, and many other areas.

Elliptic boundary value problems (e.g., Dirichlet or Neumann for
the Laplace operator) in a smooth bounded domain in $\R^n$ are
often studied directly, not necessarily in the framework of a
voluminous calculus. However, it may be instructive to consider the
class of all elliptic boundary value problems for elliptic
differential operators at the same time.

The history of elliptic boundary value problems is well known;
there are many stages and numerous applications. In the present
exposition we will not give a complete list of merits and
achievements of the general development.

We mainly focus on ideas that played a role for the iterative
calculus of edge and corner problems. A classical reference is 
the
work of Lopatinskij \stcite{Lopa1} who introduced a general concept of
ellipticity of boundary conditions for 
an elliptic differential operator. 
We are
talking here about Shapiro-Lopatinskij conditions. The operators 
representing boundary conditions are also called trace operators. 

An algebraic characterisation of elliptic differential trace
operators may be found in Agmon, Douglis, and Nirenberg
\stcite{Agmo1}, the complementing condition. Let
us also mention the works of Schechter \stcite{Sche1},
Solonnikov \stcite{Solo1}, \stcite{Solo2}, and the
monograph of Lions and Magenes \stcite{Lion1}. Moreover,
Solonnikov \stcite{Solo3} studied parabolic problems in
such a framework. The Sixtees of the past century were
also a period of intensive development of the
pseudo-differential calculus, cf. Kohn and Nirenberg
\stcite{Kohn1}, H\"ormander \stcite{Horm7}, \stcite{Horm2}.
Ideas and sources of this theory (especially, of singular
integral operators) are, in fact, much older.

Wiener-Hopf operators became an important model for
different kinds of operator algebras with symbolic
structures, ellipticity, and Fredholm property. In higher
dimensions they played an essential role in the theory of
Vishik and Eskin on pseudo-differential boundary value
problems without (or with) the  transmission property at
the boundary, cf. Vishik and Eskin \stcite{Vivs2},
\stcite{Vivs3} and Eskin's monograph \stcite{Eski2}.
An algebra of pseudo-differential operators with the 
transmission property at the boundary was established by
Boutet de Monvel \stcite{Bout1}. 
This algebra is closed under constructing parametrices of
elliptic elements.

Similarly as in the work of Vishik and Eskin, the operators
in Boutet de Monvel's algebra have a $2 \times 2$ block
matrix structure with additional trace and potential
entries. Moreover, there appear extra Green operators in
the upper left corners which are indispensable in
compositions. Apart from the standard ellipticity of the
upper left corner there is a notion of ellipticity of the remaining
entries which is an analogue of the Shapiro-Lopatinskij
condition, a bijectivity condition for a second
(operator-valued) symbolic component.

It turned out very early that the ellipticity of the upper
left corner does not guarantee the existence of a
Shapiro-Lopatinskij elliptic $2 \times 2$ block matrix
operator, cf. Atiyah and Bott \stcite{Atiy5}. 
Despite of the general index theory, cf. Atiyah and Singer
\stcite{Atiy6} and the subsequent development, which is 
also an important source of the analysis on manifolds with
singularities, it remained unclear for a long time how to
complete boundary value problems for such operators to a
calculus which is closed under parametrix construction for
elliptic elements (answers are given in \stcite{Schu37},
\stcite{Schu41}, \stcite{Schu50}).
The case of differential boundary value problems of that type was widely studied by numerous authors, see Seeley \cite{Seel2}, Booss-Bavnbek and Wojciechowski \cite{Boos1}, or the author's joint paper with Nazaikinskij, Sternin, and Shatalov \cite{Naza1}, see also \cite{Naza7}, jointly with Nazaikinsij, Savin, and Sternin.

Interpreting a (smooth) manifold with boundary as a manifold
with edge (with the boundary as edge and the inner normal
as the model cone of local wedges) boundary value problems
have much in common with edge problems. This is
particularly typical for the theory of Vishik and Eskin
where the operators on the half-axis are Wiener-Hopf and
Mellin operators that belong (in the language of
\stcite{Schu2}, \stcite{Schu31}) to the cone algebra on the
half-axis, cf. Eskin's monograph        %[Eski..,\S 15].....
\textup{\cite[\S 15]{Eski2}}. 
Let us also mention in this connection the work of Cordes and Herman \cite{Cord3} and 
Gohberg and Krupnik \stcite{Gohb2}, \stcite{Gohb5}.
(The calculus of Vishik and Eskin was completed to an
algebra in \stcite{Remp1}.)
There is another category of problems with `edges', the so-called Sobolev problems, where elliptic conditions are
posed on submanifolds of codimension $\geq 1$, embedded in
a given manifold. This type of problems has been
systematically studied by Sternin \stcite{Ster2},
\stcite{Ster5}, including conditions of trace and potential
type. In this case the embedded manifolds can also be
interpreted as edges (cf. the recent papers
\stcite{Naza13}, \stcite{Dine1} and \stcite{Liu3}).

Boundary value problems for differential operators in domains with
conical singularities in weighted Sobolev spaces have been studied
by Kondrat'ev \stcite{Kond1} and by many other authors. 
The Fredholm property in \cite{Kond1} was obtained
under the condition of Fuchs type ellipticity together with the
ellipticity of the principal conormal symbol with respect to a
chosen weight. At the same time the asymptotics of solutions at
the tip of the cone was characterised in terms of the
non-bijectivity points of the principal conormal symbol which
gives rise to meromorphic operator functions, operating in Sobolev
spaces on the base of the local cones. (The notation `conormal
symbol' was introduced in \stcite{Remp1} in a situation of
boundary value problems, where `conormal' comes from the conormal bundle 
of a domain that corresponds to a cone, see also Section 2.3; other authors
speak about operator pencils or indicial families. Our notation 
is
motivated by their role of a principal symbolic component
in a hierarchy.)

Such conormal symbols fit into the frame of parameter-dependent
operators and parameter-dependent ellipticity on a manifold. This
is an aspect of independent importance. Agmon \cite{Agmo4}
interpreted a
spectral parameter as an additional covariable; a similar concept
was applied
by Agranovich and Vishik \stcite{Agra1} to parabolic problems, and
it played an important role in Seeley's work \stcite{Seel1} on
complex powers of an elliptic operator. 
Later on, parameter-dependent boundary value problems in the
technique of Boutet de Monvel's calculus were investigated by Grubb
\stcite{Grub1} with a more general dependence on parameters.

Parameters in the singular analysis appear in a very simple way.
If $A$ is a (say, differential) operator on a singular
configuration $M$ and if we analyse $A$ in a neighbourhood of a
(smooth) stratum $Y$ then we can freeze variables on $Y$ and
consider the cotangent variables $\eta$ to $Y$ (in the symbol of
$A$) as parameters. We then obtain an operator function $a(y,
\eta)$ on a cone $X^\Delta$ transversal to $Y$. In this connection
it is natural to accept $X^\Delta$ as an infinite cone and to
interpret $a(y, \eta)$ as an operator-valued symbol of $A$.
This is just the idea of boundary symbols on a manifold
$M$ with smooth boundary; the transversal cone in this case is
$\ol{\R}_+$.
%\textup{\cite[Proposition 2.1.136]{Schu31}} - Probe

In Shapiro-Lopatinskij ellipticity there is an automatic
control of operators for $r \to \infty$ when $\eta \not= 0$.
Similarly, also for $\dim X > 0$, it is interesting to observe
the behaviour of operators near the conical exit of
$X^\Delta$ to infinity.

The simplest model of such a manifold is the Euclidean space
$\R^n$ which
corresponds to $(S^{n-1})^\Delta$ in polar coordinates.
Ellipticity up to infinity in the case of differential
operators was studied by Nirenberg and Walker \cite{Nire1}.
The pseudo-differential calculus of such operators was independently
developed by Shubin \cite{Shub1}, Parenti \cite{Pare1} and  by
Cordes \cite{Cord2}. It is essential here that the manifolds
at infinity have as specific structure, i.e., there is a
`metric' background which leads to standard Sobolev spaces
up to infinity. For the singular analysis near $r = 0$ it
is also important to study operators on finite cylinders $\R
\times X$ between `cylindrical' Sobolev spaces. Although
infinite cones and infinite cylinders geometrically are
nearly the same, the ellipticities are quite different.
Ellipticity referring to the cylindrical metric was
investigated by Sternin \cite{Ster3}. The corresponding
results are close to the ones for weighted Sobolev spaces
near conical singularities.

Classical operator calculi with symbolic structures usually
contain the equivalence between ellipticity and Fredholm
property in the chosen Sobolev spaces on the given configuration.
This is a starting point of many beautiful connections to
index theories. 
Although this is an interesting side of the history, it 
goes beyond the scope of this exposition  which is focused more
on `analytic' aspects. 
Geometric and topological relations
are discussed in detail in a new monograph of Nazaikinskij,
Savin, Schulze and Sternin \cite{Naza14}. 

\subsection{Conification and edgification}  
\label{s.10.6.2}

   %10.6.2

By `iterative calculus' we understand a program to successively generating operator structures on manifolds with higher singularities, such that ellipticity of the operators, parametrices, and index theory make sense.
Let us first recall that a manifold $M \in \frak{M}_{k+1}$, $k\in \N$, can be generated by repeatedly forming cones
$X^{\Delta} =(\ol{\R}_+ \times X) /(\{0\} \times X)$ and wedges
$X^\Delta \times \Omega$, starting from elements $X \in \frak{M}_k$ and open $\Omega \subseteq \R^q$ (local edges), combined with pasting constructions to reach the `global' object $M$.
In the case of a $C^\infty$ manifold $X$, i.e., $X \in \frak{M}_0$, we obtain in this way manifolds with conical singularities and edges, i.e., objects in $\frak{M}_1$; a next step gives us corner manifolds in $\frak{M}_2$, i.e., of second generation, and so on.

Now the program of the iterative calculus is as follows.
Given a (pseudo-differential) operator algebra on $X \in \frak{M}_k$, apply a `conification' to generate a so-called cone algebra on $X^\Delta$, then an `edgification' to obtain a corresponding edge algebra on $X^\Delta \times \Omega$, and then past together the obtained local cone and edge algebras to the next higher algebra on $M \in \frak{M}_{k+1}$.
The question is now how to organise such conifications and edgifications.
Answers  of different generality may be found in the papers and monographs
\cite{Schu32}, \cite{Schu29}, \cite{Schu11}, \cite{Schu2}, \cite{Schu20}, \cite{Schu27}, 
as well as in the author's joint works with 
Rempel \cite{Remp1}, \cite{Remp2}, \cite{Remp3},
Egorov \cite{Egor1},
Kapanadze \cite{Kapa10}, or
Nazaikinskij, Savin, and Sternin \cite{Naza14}.
In order to make the conification and edgification idea transparent we try to give an impression of how the first cone, edge, and corner algebras were originally found.
(The following discussion has some intersection with the previous section).

First, in the context of the early achievements of the calculus of pseudo-differential operators, see Kohn and Nirenberg \cite{Kohn1},  H\"ormander \cite{Horm7}, \cite{Horm2}, and of the index theory, see Atiyah and Singer \cite{Atiy6}, it became standard to establish operator algebras with a principal symbolic structure, closed under the construction of parametrices of elliptic elements, and containing a minimal class of `desirable' elements, such as differential operators (see also the discussion in Section 4.4 before).
However, already for boundary value problems on a $C^\infty$ manifold with boundary this concept leads to `unexpected' difficulties.
Vishik and Eskin \cite{Vivs2}, \cite{Vivs3} established a very general calculus of pseudo-differential boundary value problems, but the `calculus answer' was not so smooth as in the boundaryless case; compositions and parametrices were not given within the calculus.
A `smooth' calculus of boundary value problems in that desirable sense was obtained later on by Boutet de Monvel \cite{Bout1}, however  under two severe restrictions.
The symbols are required to have the transmission property at the boundary (these symbols form a thin set in the space of all pseudo-differential symbols which are smooth up to the boundary).
Moreover elliptic operators (such as Dirac operators in even dimensions or other important geometric operators) are excluded (for topological reasons) from the notion of Shapiro-Lopatinskij ellipticity of boundary conditions, see also Atiyah and Bott \cite{Atiy5}, and the discussion in Section 5.3.
In any case, both Vishik, Eskin and Boutet de Monvel stressed the role of a second principal symbolic component, namely, the boundary symbol which encodes the Shapiro-Lopatinskij ellipticity of the boundary conditions and refers to the entries of a $2 \times 2$ block matrix with trace and potential operators.
The latter kind of operators (together with Green operators) was added as a contribution of the boundary.
(Note that an operator algebra for boundary value problems without any topological restriction (such as for geometric operators mentioned before) was given
in \cite{Schu37}, see also \cite{Schu50}.)
An algebra of boundary value problems that admits all smooth symbols (also those without the transmission property at the boundary, as in Vishik and Eskin's work), closed under parametrix construction of Shapiro-Lopatinskij-elliptic elements, was constructed by Rempel and Schulze \cite{Remp1}.
However, the structure of lower order terms was not yet analysed in \cite{Remp1}; this came later in the frame of the edge calculus.
A crucial role for \cite{Remp1} played a specific algebra on the half-axis from Eskin's book \cite{Eski2}, namely, a pseudo-differential algebra of operators of order zero on $\R_+$, without any condition of transmission property at $0$, formulated by means of the Mellin transform.
Lower order terms in this algebra in Eskin's formulation are Hilbert Schmidt operators in $L^2 (\R_+)$.
From the point of view of conical singularities this half-axis-algebra can be interpreted as a substructure of the `cone algebra', see also \cite{Schu31}, while the operators in \cite{Remp1} could be seen as edge operators with the boundary being interpreted as an edge and $\ol{\R}_+$, the inner normal, as the model cone of local wedges.
In that sense \cite{Remp1} gave a first example of an edgification of a cone algebra which is, roughly speaking, a pseudo-differential calculus along the edge with amplitude functions taking values in the cone algebra on the model cone, here $\ol{\R}_+$.
Of course, also Boutet de Monvel's algebra can be interpreted as an edgification of its boundary symbolic calculus, though the operators in this case form a narrower subalgebra of the cone algebra on $\R_+$.

In order to really recognise the algebras on the  half-axis in connection with conical singularities,
another input was necessary, namely, the analysis of operators of Fuchs, type, which are of independent interest on manifolds with conical singularities in general.
It was the work of Kondratyev \cite{Kond1} which motivated the author together with Rempel to try to carry out the hull operation, discussed in Section 4.4, i.e., to complete the Fuchs type differential operators to a corresponding pseudo-differential algebra.
This was first done in \cite{Remp7} for the case of a closed manifold with conical singularities, then in \cite{Remp3} for the case of boundary value problems on a manifold with conical singularities and boundary, see also \cite{Remp11}.
Another orientation (from the point of methods) have the works of Plamenevskij \cite{Plam1}, 
\cite{Plam2}, Derviz \cite{Derv1} and Komech \cite{Kome1} (the latter is close to technique of Vishik and Eskin.
Independently, also Melrose and Mendoza \cite{Melr1} constructed a pseudo-differential calculus for Fuchs type symbols, see also Melrose \cite{Melr4}.

The cone calculus of \cite{Remp7} refers to weighted Mellin Sobolev spaces and subspaces with discrete asymptotics, using suitable classes of meromorphic Mellin symbols with values in pseudo-differential operators on the base $X$ of the cone.
The cone calculus for $\dim X=0$ was formulated for the purposes of boundary value problems; compared with Eskin's algebra on $\R_+$ the cone algebra of \cite{Remp7} is not restricted to operators of order zero and to principal conormal symbols of order zero and to Hilbert Schmidt operators as the ideal of smoothing operators.
It contains operators of any order with coefficients that are smooth up to $0$, modulo a possible weight factor, and also lower order conormal symbols;   the smoothing operators are Green operators in the sense that they map $\s{K}^{s,\gamma} (\R_+)$ to  spaces of the kind $\s{S}^{\gamma-\mu}_P (\R_+)$ for some discrete asymptotic type $P$ and, analogously, the adjoints.
The details of this calculus were elaborated in  \cite{Remp9}, see also \cite{Remp3}, or \cite{Schu31}, \cite{Schu2}.
We stress these features here because the choice of the cone algebra for $\dim X=0$ is crucial for the nature of the `conification' of the pseudo-differential calculus on an arbitrary base $X$.
One step is to fix a choice of an algebra of pseudo-differential operators (with `sleeping parameters') on $X$, say 
$L^\mu_\cl (X; \R)$ in the case of smooth compact (the space of all classical parameter-dependent pseudo-differential operators of order $\mu$ on $X$), and then to organise the calculus with Mellin symbols 
$h(r,w) \in C^\infty (\ol{\R}_+, L^\mu_\cl (X; \Gamma_{\frac{n+1}{2}-\gamma}))$  
(for $n=\dim X$), along the lines of the cone algebra on $\R_+$.
The full structure is, of course, rich in details, for instance, we can take holomorphic (in $w \in \C$) non-smoothing symbols (reached by kernel cut-off constructions), meromorphic smoothing symbols, and, moreover, for $r \to \infty$ impose extra assumptions when we intend to edgify the obtained cone calculus.
Summing up, `conification'
 means to pass from a prescribed pseudo-differential algebra on a base $X$ (first smooth and compact) to a cone algebra on $X^\Delta$ by taking the former cone algebra on $\R_+$, but now with symbols taking values in the given algebra on $X$.
 
 This cone algebra near the tip of the cone (in the variant of a base $X$ with smooth boundary) is just what completes Kondratyev's theory \cite{Kond1} to an algebra with the above mentioned properties.
 At the same time, during  this period of the development  there was another main motivation for the refinement of Eskin's algebra to the cone algebra on $\R_+$, namely, the aim to generalise the boundary symbolic calculus of Boutet de Monvel's algebra to a boundary symbolic calculus of a future algebra of boundary value problems for symbols which have not necessarily the transmission property at the boundary.
 That algebra of boundary value problems itself was intended to be obtained as a corresponding edgification.
 This program finally created the calculus of boundary value problems without the transmission property as a substructure of a corresponding edge algebra, cf. \cite{Schu2}, \cite{Schu31}, \cite{Schu20}, 
\cite{Schu41}.
At that time also the structures of the edge algebra in general were invented, in which $\ol{\R}_+$, the model cone of the case of boundary value problems, was replaced by an arbitrary cone $X^\Delta$ with a compact manifold $X$ without (and with smooth) boundary, cf. \cite{Remp3}, \cite{Remp11}, \cite{Schu32}, \cite{Schu11}.
Moreover, the methods have been developed under the  aspect of the general idea of generating operator algebras in terms of the successive procedure of `conification' and `edgification' of already achieved structures.

Elements of this approach are sketched in Section 5.4.
As noted before, the `final' structures and many interesting details are a program of future research, cf. also Section 6.4 below.
But also the development up to the present state of the calculus contained some surprising elements.
One of them was the invention of abstract edge Sobolev spaces
$\s{W}^s (\R^q, H)_\kappa$ with Hilbert spaces $H$, endowed with the action of a strongly continuous group of isomorphisms
$\kappa = \{ \kappa_\lambda \}_{\lambda \in \R_+}$,
cf. \cite{Schu32}.
From the impression on anisotropic reformulations of standard Sobolev spaces, on the role of fictitious conical points and edges, and from the experience in boundary value problems it seemed quite canonical to take for $H$ a weighted space on an infinite cone with conical exit to infinity, with the `unspecific' weight $0$ at infinity.
A few years ago ($\cong$ 2001) I. Witt (who was at that time in Potsdam) suggested to admit also spaces $H$ with another weight at infinity with an adjusted variant of the group action $\kappa$.
This idea has been used by Airapetyan and Witt in \cite{Aira2}.
Later on Tarkhanov realised such an idea in \cite{Tark4}, see also \cite{Schu53} for the case of boundary value problems.
It turns out, see, for instance, \cite[Section 7.1.2]{Haru13}, that there is a continuum of different  edge spaces which all localise outside the edge to standard Sobolev spaces and admit an edge pseudo-differential calculus for the same classes of typical edge-degenerate differential operators.
Thus, the problem of `edge-quantising' edge-degenerate (pseudo-differential) symbols and of carrying out  a hull operation as discussed in Section 4.4 has many solutions.

\subsection{Similarities and differences between ellipticity and parabolicity}    %Hierarchies of operator algebras}   
\label{s.10.6.3}

%10.6.3
In this exposition we mainly focused on the concept of ellipticity.
Of course, also other types of equations are of interest on a manifold with singularities, for instance, parabolic or hyperbolic ones.
Many problems in this connection occur in models of physics.

We want to discuss here a few aspects on parabolic operators.
The simplest example is the heat operator
\begin{equation}
\label{eq.302}
A := \partial_t - \Delta
\end{equation}
with the Laplacian $\Delta$ on a Riemanian manifold $X$, $n=\dim X$, with $t \in \R$ being the time variable.
In local coordinates $x \in \R^n$ the operator \eqref{eq.302} has an anisotropic homogeneous principal symbol
$$
\sigma_\psi (A) (\tau, \xi) := i \tau + |\xi|^2
$$
which is anisotropic homogeneous of order 2, i.e., satisfies the relation
$$
\sigma_\psi (A) (\lambda^2 \tau, \lambda \xi) =
   \lambda^2 \sigma_\psi (A) (\tau, \xi)
$$
for every $\lambda \in \R_+$.
The operator \eqref{eq.302} is anisotropic elliptic of order 2 in the sense of the property
\begin{equation}
\label{eq.303}
\sigma_\psi (A) (\tau,\xi) \not= 0 \ \ \ \text{for all} \ \ \
(\tau, \xi) \in \R^{1+n} \setminus \{0\}.
\end{equation}
Parabolicity means that $\sigma_\psi (A)$ has an extension 
$\sigma_\psi (A) (\zeta,\xi) := i \zeta + |\xi|^2$ 
into the lower complex $\zeta$ half-plane $\C_-$ with respect to the time covariable, such that
$$
\sigma_\psi (A) (\zeta, \xi) \not= 0 \ \ \ \text{for all} \ \ \
(\zeta, \xi) \in (\ol{\C}_- \times \R^n) \setminus \{ 0\}.
$$
With 
$\langle \tau, \xi \rangle_{(l)} :=
   (1+|\tau|^2 + |\xi|^{2l} )^{1/2l}$, $l \in \N \setminus \{ 0 \}$,
we can form anisotropic Sobolev spaces 
$H^{s,(l)} (\R \times \R^n)$ with the norm
\begin{equation}
\label{eq.304}
\| u(t,x) \|_{H^{s,(l)} (\R \times \R^n)} =
   \Big\{ \int  \langle \tau, \xi \rangle^{2s}_{(l)}
          | \hat{u} (\tau, \xi) |^2 d\tau d \xi \Big\}^{\frac{1}{2}},
\end{equation}
$s \in \R$, or, more generally,
$H^{s,(l)} (\R \times X)$ on the infinite cylinder $\R \times X$.
Let us set 
$H^{s,(l)}_0 (\R \times X) :=
    \big\{ u \in H^{s,(l)} (\R \times X) :
     u |_{\R_- \times X} =0 \big\}$
and 
$$
H^{s,(l)}_0 ((0,T) \times X) :=
    \big\{ u |_{(-\infty,T) \times X} :
      u \in H^{s,(l)}_0 (\R \times X) \big \}
$$
for every $T >0$.
The operator \eqref{eq.302} defines continuous maps
\begin{equation}
\label{eq.305}
 A : H^{s,(2)}_0 ((0,T) \times X) \to
     H^{s-2, (2)}_0 ((0,T) \times X)
\end{equation}
for all $s \in \R$, and it is a reasonable problem to ask the solvability of the equation
$$
Au =f
$$
in this scale of spaces, more precisely, to find a solution 
$u(t,x) \in H^{s,(2)}_0 ((0,T) \times X)$ for every 
$f(t,x) \in H^{s-2,(2)}_0 ((0,T) \times X)$ and to construct a parametrix (or the inverse) of the operator \eqref{eq.305} within a corresponding anisotropic calculus of pseudo-differential operators on the cylinder.
An answer was given by Piriou \cite{Piri1} in the framework of a Volterra pseudo-differential calculus, not only for the anisotropy $l=2$, but for arbitrary even $l$.
Corresponding differential operators may have the form
$$
A := (\partial_t - D)^m,
$$
$m \in \N$, for an elliptic differential operator $D$ on $X$ of order $l$, for instance,
$$
D = (-1)^{1+\frac{l}{2}} \Delta^{\frac{l}{2}}.
$$
Any such operator induces continuous maps
$$
A : H^{s,(l)}_0 ((0,T) \times X) \to
     H^{s-m, (l)}_0 ((0,T) \times X)
$$
for all $s \in \R$.
The solvability problem is as before, cf. \cite{Piri1}.
More generally, this also concerns operators that are locally on $X$ of the form
\begin{equation}
\label{eq.306}
A = \sum_{|\alpha| \leq m}
     a_\alpha (t,x) D^{\alpha_0}_t D^{\alpha'}_x
\end{equation}
for $\alpha := (\alpha_0, \alpha') \in \N^{1+n}$,
$|\alpha |_l := l \alpha_0 + |\alpha' |$, $a_\alpha (t,x) \in C^\infty (\R \times X)$.
The anisotropic homogeneous principal symbol of $A$ is defined by
\begin{equation}
\label{eq.307}
\sigma_\psi (A) (t,x,\tau,\xi) :=
   \sum_{|\alpha|_l =m} 
      a_\alpha (t,x) \tau^{\alpha_0} \xi^{\alpha'}.
\end{equation}
It satisfies the identity
\begin{equation}
\label{eq.308}
\sigma_\psi (A) (t,x, \lambda^l \tau, \lambda \xi) =
    \lambda^m \sigma_\psi (A) (t,x,\tau, \xi)
\end{equation}
for all $\lambda \in \R_+$.
Parabolicity of \eqref{eq.306} means that the extension of 
$\sigma_\psi (A)$ to $(\zeta, \xi) \in (\ol{\C}_- \times \R^n) \setminus \{ 0 \}$
satisfies the condition
\begin{equation}
\label{eq.309}
\sigma_\psi (A) (t,x,\zeta,\xi) \not= 0 \ \ \ \text{for all} \ \ \
(t,x,\zeta,\xi) \in \R \times \R^n \times
    (\ol{\C}_- \times \R^n) \setminus \{ 0 \} ).
\end{equation}
Observe that then
$$
p(t,x,\tau, \xi) :=
   \sigma_\psi (A)^{-1} (t,x, \tau - i \varepsilon, \xi)
$$
for any fixed $\varepsilon >0$ belongs to 
$C^\infty (\R \times \R^n \times \R \times \R^n)$,
and extends to a function 
$p(t,x, \zeta, \xi)$ in $C^\infty (\R \times \R^n \times (\ol{\C}_- \times \R^n))$
which is holomorphic in $\zeta \in \C_-$ and satisfies the anisotropic symbolic estimates
\begin{equation}
\label{eq.310}
| D^\alpha_{t,x} D^\beta_{\zeta, \xi} p(t,x,\zeta,\xi)| \leq
    c  \langle \zeta, \xi \rangle^{\mu -|\beta|_l}_{(l)}
\end{equation}
for every $\alpha, \beta \in \N^{1+n}$, $(t,x) \in K_0 \times K'$, $K_0 \subset \subset \R$,
$K' \subset \subset \R^n$ compact,
$(\zeta, \xi) \in \ol{\C}_- \times \R^n$, $\mu := -m$, with constants $c =c (\alpha, \beta, K_0, K') > 0$.

Note that analogous considerations make sense for a (in the simplest case) smooth compact manifold $X$ with boundary.
Then, together with the operator $A$ we consider (first differential) boundary conditions on $(0,T) \times X$,
represented by a continuous operator
$$
T : H^{s,(l)}_0 ((0,T) \times X) \to
    \oplus^{N}_{j=1}
         H^{s-m_j -\frac{1}{2}, (l)}_0 ((0,T) \times \partial X)
$$
of the form $T = {}^\t (T_1, \ldots , T_N)$ for
$$
T_j u (t,y) :=
   B_j u(t,x) |_{(0,T) \times \partial X},
$$
$(t,y) \in (0,T) \times \partial X$, with differential operators
$$
B_j := \sum_{|\beta|_l \leq m_j}
     b_{j,\beta} (t,x) D^{\beta_0}_t D^{\beta'}_x,
$$
$b_{j,\beta} \in C^\infty (\R \times X)$.
Locally near $\partial X$, in a splitting  $x =(y, x_n) \in \partial X \times [0,1)$ in tangential and normal variables near the boundary, and covariables $(\eta, \xi_n)$, we have the boundary symbols
$$
\sigma_\partial (A) (t,y, \tau, \eta) := 
      \sigma_\psi (A) (t,y,0, \tau, \eta, D_{x_n}) :
         H^s (\R_+) \to H^{s-m} (\R_+),
$$
$$
\sigma_\partial (T_j) (t,y, \tau, \eta)  :
      H^s (\R_+) \to \C,
$$
defined by
\begin{eqnarray*}
\sigma_\partial (A) (t,y, \tau, \eta) u
&:=&
\sigma_\psi (A) (t,y,0, \tau, \eta, D_{x_n}) u,   \\   
\sigma_\partial (T_j) (t,y, \tau, \eta) u
&:=&
\sigma_\psi (B_j) (t,y,0, \tau, \eta, D_{x_n}) u|_{x_n=0},
\ \ \  j=1, \ldots, N,
\end{eqnarray*} 
for $(\tau, \eta) \not= 0$, 
$s > \max \{ m - \frac{1}{2}, m_1 + \frac{1}{2}, \ldots, m_N + \frac{1}{2} \}$.
Setting
$\s{A} := {}^\t (A \ \ T)$ for 
$T := {}^t (T_1, \ldots , T_N)$ 
we thus obtain the principal boundary symbol
\begin{equation}
\label{eq.312}
\sigma_\partial (\s{A}) (t,y,\tau,\eta) :=
  \left( \begin{array} {c}
            \sigma_\partial  (A) \\
            \sigma_\partial (T) 
         \end{array} \right) 
 (t,y,\tau, \eta) : H^s (\R_+) \to 
 \Hsum{H^{s-m} (\R_+)}{\C^N}.
\end{equation} 
Observe that for $(\kappa_\lambda u) (x_n) := 
    \lambda^{1/2} u (\lambda x_n)$, $\lambda \in \R_+$,
    we have anisotropic homogeneity, namely,
$$
\sigma_\partial (A) (t,y, \lambda^l \tau, \lambda \eta) =
     \lambda^m \kappa_\lambda 
        \sigma_\partial (A) (t,y, \tau, \eta) \kappa^{-1}_\lambda,   
$$
$$
\sigma_\partial (T_j) (t,y, \lambda^l \tau, \lambda \eta) =
     \lambda^{m_j +\frac{1}{2}} 
        \sigma_\partial (T_j) (t,y, \tau, \eta) \kappa^{-1}_\lambda,   
$$
$\lambda \in \R_+$, $j=1, \ldots, N$.
\begin{Definition}
\label{d.6.1}
The boundary value problem
{\em
\begin{equation}
\label{eq.313}
Au = f \ \ \text{in} \ \ (0,T) \times X, \ \
Au = g \ \ \text{in} \ \ (0,T) \times \partial X
\end{equation}
}
is called parabolic, if $A$ is parabolic in the sense of \eqref{eq.309}, and if the boundary symbol has an extension to an invertible family of operators
\begin{equation}
\label{eq.314}
\sigma_\partial (A) (t,y, \zeta, \eta) : H^s (\R_+) \to 
   \Hsum{H^{s-m} (\R_+)}{\C^N}
\end{equation}
in $(\zeta, \eta) \in (\ol{\C}_- \times \R^{n-1} ) \setminus \{ 0 \}$,
holomorphic in $\zeta \in \C_-$, 
$s > \max \{ m-\frac{1}{2}, m_1 + \frac{1}{2}, \ldots , m_N + \frac{1}{2} \}$.
\end{Definition} 
Similarly as in the elliptic theory, the number $N$ is determined by the parabolic operator $A$.
\begin{Theorem}
\label{t.6.2}
\begin{itemize}
\item[{\em(i)}]
Let $X$ be a compact $C^\infty$ manifold with boundary.
A parabolic boundary value problem $\s{A} = {}^\t (A \  \ T)$ induces isomorphisms
$$
\s{A} : H^{s,(l)}_0 ((0,T) \times X) \to
    \Hsum{H^{s-m,(l)}_0 ((0,T) \times X)}
         { \oplus^{N}_{j=1}  H^{s-m_j-\frac{1}{2}, (l)}_0  ((0,T) \times \partial X)}
$$
for all $s > \max \{ m-\frac{1}{2}, m_1 + \frac{1}{2}, \ldots , m_N + \frac{1}{2} \}$ and
$0 < T < \infty$.
The inverse operator belong to an anisotropic analogue of Boutet de Monvel's calculus on the cylinder and is parabolic within that framework.
\item[{\em(ii)}]
If $X$ is a closed compact $C^\infty$ manifold, then a parabolic operator $A$, {\em(}locally{\em)} of the form \eqref{eq.306}, induces isomorphisms
$$
A : H^{s,(l)}_0 ((0,T) \times X) \to
   H^{s-m,(l)}_0 ((0,T) \times X)
$$
for all $s \in \R$ and $0 < T < \infty$.
The inverse operator belongs to an anisotropic analogue of the calculus of classical pseudo-differential operators and is parabolic in this class.
\end{itemize}
\end{Theorem}
  
  A reference for Theorem \ref{t.6.2}  is Agranovich and Vishik \cite{Agra1} and Krainer \cite{Krai5} (in a slight modification for finite cylinders).
  The assertion (ii) may be found in the paper \cite{Piri1} of Piriou.
  It is also interesting to consider parabolicity on the infinite  half-cylinder $\R_+ \times X$ with special attention for $t \to \infty$ and to establish invertibility of the corresponding operators in weighted analogues of the above mentioned spaces with exponential weight up to $t=\infty$.
  Corresponding results for the case of closed compact $X$ are given in Krainer and Schulze \cite{Krai4}, see also Krainer \cite{Krai9}, \cite{Krai10}, and for the case of compact $X$ with $C^\infty$ boundary, in the framework of (pseudo-differential) boundary value problems, in Krainer \cite{Krai5}.
  
  In this approach the idea is to interpret the infinite time-space cylinder for $t \to \infty$ as a transformed anisotropic cone, obtained for $t >c$ for some $c>0$ by the substitution $t =-\log r$, 
$r \in (0, e^{-c})$, cf. also the discussion in Section 3.4, especially, the form of the operators \eqref{eq.208} which comes from operators of Fuchs type (in the parabolic case from anisotropic ones). 
Remember that in Fuchs type operators we imposed smoothness of the coefficients for $r \to 0$ (up to a possible weight factor).
That means, for the transformed operator in $t$ we impose a corresponding behaviour of the coefficients for $t\to \infty$.

The above mentioned results on infinite cylinders just express the inverses of such parabolic operators in the framework of an anisotropic analogue of the cone algebra, referring to a conical singularity at infinity, more precisely, within an anisotropic version of the cone algebra with a control of the Volterra property up to infinity.
Clearly at infinity an analogue of the principal conormal symbol is required to be bijective in Sobolev spaces on the cross section $X$.
This causes a discrete set of forbidden (exponential) weights at infinity, similarly as in the cone calculus at the tip of the cone (for the corresponding exponents in power weights).

Similarly as in the `usual' cone algebra it is also interesting to observe asymptotics of solutions, here interpreted as long-time asymptotics, coming from the meromorphic structure of the inverse of the principal conormal symbol,
see Krainer \cite{Krai9}, \cite{Krai10}, \cite{Krai5}. 

\begin{Remark}
\label{r.6.3}
In parabolic problems it is also common to pose {\em(}non-trivial{\em)} initial conditions at the bottom of the cylinder.
In the case of boundary value problems {\em(}see, Agranovich and Vishik {\em\cite{Agra1})} one usually assumes that the initial values are compatible with the values on the boundary of the cylinder.
We do not discuss the details here but return below to such problems from a more general point of view.
\end{Remark}
Parabolicity in the framework of algebras of anisotropic pseudo-differential operators and the computation of long-time asymptotics is also interesting in connection with special configurations $X$ with singularities.
Looking at simple models of heat flow in media with singularities of that kind we immediately see the relevance of such a generalisation.

For instance, if $X$ has conical singularities, the additional time variable generates an edge.
Then, considering long-time asymptotics for $t \to \infty$ we are faced with a corner problem in the category $\got{M}_2$, where $t$ plays the role of a corner axis variable.
The same is true when $X$ is a manifold with smooth edges.
Long-time asymptotics for the latter case have been studied by Krainer and Schulze in \cite{Krai3}.
Earlier, parabolicity in an anisotropic analogue of the edge algebra in a finite time interval $(0,T)$, i.e., when $X$ is a manifold with edge, was investigated in \cite{Buch1}.
Let us also mention that parabolic boundary value problems in the pseudo-differential set-up of Vishik and Eskin's  technique have been investigated in \cite{Can1}, \cite{Vivs5}.

Also for parabolic operators in cylinders with singularities (in the spatial variables) it is natural to pose additional data of trace and potential type along the lower-dimensional strata of the configuration, satisfying a parabolic analogue of the Shapiro-Lopatinskij condition.
That means, that the symbols admit holomorphic extensions into the lower complex half-plane of the time covariable, required to be invertible there.
If we have posed such conditions, then it is clear that the inverses of the corresponding operators again belong to the Volterra calculus of such operators.
However, in contrast to the analogous task in the elliptic theory, the explicit construction of extra Shapiro-Lopatinskij parabolic conditions seems to be not so easy, although it should be always possible.
Some results in this direction for specific parabolic problems may be found in 
\cite{Mika1}.
Also initial-edge conditions (in analogy of initial-boundary conditions) with non-trivial initial data on the bottoms of cylinders associated with lower-dimensional strata of the spatial configurations belong to the natural tasks in parabolicity on singular manifolds, both under the condition of compatibility between initial and edge data as well as of non-compatibility (cf. Remark \ref{r.6.3}).
As far as we know there is nothing done yet in this direction, and it is certainly interesting to know more on the nature of solvability of such problems.
Note that initial-boundary value problems with non-compatible data, even in simplest cases of the heat operator with Dirichlet or Neumann data (or even mixed data of Zaremba type) on the boundary of the cylinders, together with initial conditions on the bottom of the cylinder have a simple physical meaning.
In the non-compatible case those represent kinds of mixed problems, combined with corner singularities when the boundary of the cylinder is not smooth, or if the boundary data are mixed (e.g., of Zaremba type).

\subsection{Open problems and new challenges}   %Complexes}
\label{s.10.6.4}

%10.6.4
In the singular analysis (similarly as in other areas of mathematics)
it is difficult to give reasonable criteria on what is an `open problem'.
The solution may depend on the  person who finds something open or not.
It also happens that crucial notions  in this field (such as `ellipticity' or `corner manifolds') occur in quite different meanings.
Being aware of this uncertain background we want to discuss a few aspects of the singular analysis that contain challenges for the future research.
First of all the known elements of the elliptic (and also the parabolic) theory (including boundary value problems) on smooth configurations are of interest also in the singular case.
This concerns, in particular, the points (S.1) - (S.4) of Section 4.4 which can be specified for singular manifolds by the discussion in Section 5.
A number of new challenges can be summarised under the following key words.   

{\bf (F.1) Operator algebras.}
Given manifolds $M \in \got{M}_k$, $k \in \N$, $k \geq 2$, study the natural analogues of the (known for $k=0,1$) pseudo-differential algebras, including the principal symbolic hierarchies and additional data (of trace and potential type) on the lower-dimensional strata, and complete necessary elements of the conification and edgification process.

As we pointed out in different considerations before, the higher pseudo-differential algebras on stratified spaces are more general than everything what is usually contained in theories of boundary value problems (when we consider a boundary as a realisation of a smooth edge), including the case of symbols without the transmission property at the boundary.
Even if we ignore for the moment the aspect of existence or non-existence of Shapiro-Lopatinskij edge conditions (and assume, for instance, the case that the existence is guaranteed) there is a large variety of `technical' elements of a calculus to be established in the future in such a way that the theory on a space $M \in \got{M}_{k+1}$ is really a simple iteration of steps up to $\got{M}_k$.
There is the system of quantisations which contains anisotropic reformulations of isotropic (though degenerate in stretched variables) symbolic information in terms of various operator-valued symbols with twisted homogeneity, combined with the `right' choices of weighted distribution spaces.
In our exposition we discussed spaces based on `$L^2$-norms'.
In applications to non-linear operators it is often necessary to treat the `$L^p$-case' for $p \not= 2$.
This is one of the problems with is essentially open.

Another interesting aspect is the problem of variable and branching asymptotics that should thoroughly be investigated, cf. Section 4.5 for smooth edges and \cite{Schu34}, \cite{Schu36} for the case of boundary value problems without the transmission property.

There is also the question of `embedding' the calculus for $\got{M}_k$ as a subcalculus for $\got{M}_{k+1}$, for instance, by `artificially' seeing an $M \in \got{M}_k$ as an element of $\got{M}_{k+1}$.
This appears in connection with the following problem.
Take an elliptic operator $A$ on a closed compact $C^\infty$ manifold, fix a triangulation, and rephrase $A$ as an elliptic corner operator $\s{A}$ on the arising manifold with edges and corners, where $\ind A = \ind \s{A}$, then pass to a more refined triangulation and formulate $\s{A}$ again as a corresponding corner operator $\got{A}$ such that $\ind \s{A} = \ind \got{A}$, and so on.
The investigations in the author's joint papers with Dines \cite{Dine1}, Dines and Liu \cite{Dine3} can be seen as a contribution to this aspect.

Let us also point out that here we always speak about regular singularities. 
The various cuspidal cases may be of quite different character, and also here the main structures on operator algebras from the point of view of asymptotics in distribution spaces, possible extra edge conditions, adequate quantisations, construction of parametrices within the calculus, remain to be achieved.   

{\bf (F.2) Higher corner spaces.}
Complete and deepen the investigation of the higher generations of weighted Sobolev spaces that fit to the operator algebras of (F.1).

The choice of weighted edge spaces on a manifold with smooth edge that we discussed in Section 1.3 is not entirely canonical.
We saw that there is (at least) one continuously parametrised family of such (mutually non-equivalent) spaces which all admit the edge calculus, although there is a candidate which seems to be the most `natural' one.
Also on manifolds with higher corners we have many choices and one possible preferable one, which is for integer smoothness directly connected with degenerate vector fields on the respective stretched manifold that generate the space of typical differential operators.
In the higher corner cases it seems some work to be done to completely organise the variety of anisotropic reformulations in connection with higher $\s{K}^{s,\gamma} (X^\wedge)$-spaces on respective model cones, equipped with several necessary and useful characterisations in terms of degenerate families of pseudo-differential (corner-) operators on $X$ which take into account also the presence  of the conical exit of $X^\wedge$ to infinity.
It will also be useful to single out subspaces with asymptotic information and to establish analogues of the kernel characterisations of Green operators on $X^\wedge$.
Also the $L^p$-analogues for $p \not= 2$ should be investigated, especially, from the point of view of anisotropic corner representations and of the continuity of operators in the algebras between such spaces.          

{\bf (F.3) Ellipticity under extra conditions on lower-dimensional strata.}
Study ellipticity and parametrices as well as the Fredholm property, not only from the point of view of Shapiro-Lopatinskij ellipticity of conditions on the lower-dimensional strata but of conditions, partly (or mainly) to be invented for operators which violate the topological criterion for the existence of Shapiro-Lopatinskij elliptic data.

Having organised pseudo-differential algebras $\got{A} (N)$ on $N \in \got{M}_{k+1}$ in the spirit of (F.1) we have operators $A$ together with their symbolic hierarchies $\sigma (A)$ as in (C.8), Section 5.4.
The idea of (Shapiro-Lopatinskij) ellipticity with respect to $\sigma (A)$ is that $A|_{N \setminus Z}$ is required to be elliptic in $\got{M}_k$, i.e., with respect to $\sigma_{\Int} (A)$, and that, in addition, \eqref{eq.new98} is a family of isomorphisms for all $(z, \zeta) \in T^* Z \setminus 0$.
However, from the edge algebra for smooth edges (or from boundary value problems) we know that the `interior ellipticity', i.e., the one with respect to $\sigma_{\Int} (.)$, does not guarantee the bijectivity of \eqref{eq.new98}, but only the Fredholm property.
This is just the occasion to come back to the possible topological obstruction for the existence of Shapiro-Lopatinskij elliptic conditions on $Z$, see Section 5.3.
The role of such conditions is to fill up \eqref{eq.new98} to a $2 \times 2$ block matrix family of isomorphisms 
$\sigma_{\wedge_{k+1}} (\s{A}) (z,\zeta)$ by extra entries of trace and potential type.
On the level of operators they belong to an elliptic element $\s{A}$ in $\got{A} (N)$ which itself is Fredholm as soon as $N$ is compact (and otherwise has a parametrix within the calculus).
If the topological obstruction does not vanish, then, when the interior ellipticity of $A$ refers to Shapiro-Lopatinskij elliptic data on the lower-dimensional edges of $N \setminus Z$, it should be possible to perform again the machinery of global projection conditions on $Z$ along the lines of \cite{Schu42} (which treats the case of smooth edges and is a generalisation of \cite{Schu37} and \cite{Schu44}).
Vanishing or non-vanishing of the obstruction with respect to $Z$ might depend on the choice of the extra edge conditions on the steps for $N \setminus Z$ before.
It is completely open whether that happens and how it is to be controlled.
Another interesting point is, whether the idea of global projection (or Shapiro-Lopatinskij elliptic) conditions on $Z$ is also possible, if in the steps before, i.e., within $\sigma_{\Int} (.)$ on the edges of $N \setminus Z$ there are already involved global projection conditions on a lower level of singularity.
At this moment we have to confess that in the principal symbolic hierarchy $\sigma (A)$ which was defined in (C.8) Section 5.4 we tacitly assumed the symbolic components of $\sigma_{\Int} (A)$ to consist of contributions of Shapiro-Lopatinskij type on all the lower-dimensional  singular strata of $N \setminus Z$.

Moreover, an inspection of the methods of \cite{Schu37}, \cite{Schu44}, \cite{Schu42} shows that in the global projection case the symbolic data which define the ellipticity and then the Fredholm index of the resulting operator can be enriched by the choice of the respective pairs of global projections, i.e., `simply' considered as a symbolic contribution, too.
The open question is whether this is really fruitful and whether then, having done that to generalise $\sigma_{\Int} (.)$ on $N \setminus Z$, the construction for $Z$ can be continued again with two possible outcomes, vanishing or non-vanishing  of another topological obstruction.
Let us note at this point that, in order to carry out details of this kind, we have to refer all the times to background information on ellipticity in algebras on $M^\Delta$ for $M \in \got{M}_k$, including effects from the conical exits to infinity with the corresponding symbolic structures, similarly as is done for the edge calculus of second generation in \cite{Calv1}, \cite{Calv3}.
If the symbolic machinery in such a sense could be successfully established, there remain other beautiful tasks in connection with operator conventions, moreover, with relative indices under changing weights on different levels of singularity, and with the investigation of the system of ideals in the full algebras that are determined by vanishing of some components of $\sigma(A)$.
Parametrix constructions always belong to the main issues; because of the complexity of the involved structures, this should be done in a careful manner, and the work for the next singularity order is waiting.
Nedless to say that for all components of the symbolic structures one should show the necessity of ellipticity for the Fredholm property of the associated operator.
May be, this is straightforward (the necessity of ellipticity in the framework of global projection conditions can be found in \cite{Schu50}; the idea of how to do it in this case goes back to a private communication with Savin, Sternin and Nazaikinskij during their work in Potsdam 2004; it was used again in \cite{Schu42} in the edge case).                            

{\bf (F.4) Index theory.}
Establish index theories in the algebras of (F.1), both for Shapiro-Lopatinskij ellipticity and other types of ellipticity of the extra data in the sense of (F.3).

Ellipticity of an operator on a compact configuration (or a compactified one, where a certain specific behaviour near the non-compact exits is encoded by the nature of amplitude functions, leading to specific extra principal symbolic objects at the exits) is expected to guarantee the Fredholm property in natural distribution spaces.
If the calculus is well organised, both properties (i.e., ellipticity and Fredholm property) are equivalent.
This may be a starting point of index theories on singular manifolds.
After the eminent influence of the classical index theory to modern mathematics it is generally accepted that index theories should be created also for singular manifolds.
Speculations about that could fill several books; so we can only make a few remarks here.
Index theories can have many faces, and predictions on what is most fruitful very much depends on individual priorities.
As soon as we find some operator algebras (or single operators) to be of sufficient interest we can ask to what extent the index can be expressed purely in terms of symbols (or other data contributed by the notion of ellipticity, e.g., global projection conditions).
In the Shapiro-Lopatinskij set-up this aspect is quite natural, and, as a general property of the operator theories, the index only depends on the stable homotopy class of the symbols (through elliptic ones).
In ellipticities with symbolic hierarchies we have here a first essential problem.
The symbols have operator-valued components which can be interpreted as semi-classical objects, i.e., as operator families with amplitude functions, where a quantisation is applied with respect to a part of the covariables, while other covariables remained as parameters, see the Sections 2.2 and 3.2.
The ellipticity of the corresponding component  (i.e., the invertibility for all remaining variables and covariables, say, in the cotangent bundle minus zero section of the corresponding stratum) is a kind of parameter-dependent ellipticity of operators on an infinite singular cone.
There is a subordinate principal symbolic hierarchy with ellipticity in the algebra on the corresponding infinite cone, again with operator-valued components, again with subordinate symbols belonging to corresponding algebras where those symbols take their values, and so on.
Thus every symbolic component of the original operator induces tails of subordinate symbols who all participate in a well organised way in the structure of the operators, especially, in homotopies through elliptic elements.
It is of quite practical importance for the basic understanding of the algebras on manifolds with singularities (and a reason to discuss the index problems here anyway) to know things about the index (better the kernels and cokernels) of operator families in algebras on infinite model cones, since this just affects the number of additional conditions of trace and potential type on the corresponding strata.

What concerns homotopies through elliptic symbolic tuples it is interesting to understand to what extent different components may exchange `index information' along the path that determines the homotopy.
In optimal cases the stable homotopy classes can be represented by very specific ones with particularly simple or `standard' components.
In this connection one may ask, whether in some  cases (apart from the smooth compact case) it suffices to mainly look at kinds of Dirac operators as the elliptic operators on the main statum.
In the singular situation, including the `simple' conical case, we think that it is not true that an arbitrary elliptic operator on the main stratum, when all the extra conditions participate in the ellipticity (apart from the effect that global projection conditions influence the situation anyway) can be stably homotopied through elliptic operators (in the singular operator algebra) to an operator of Dirac type in the upper left corner.
In other words, if one identifies index theory in a 
singular case with doing things for Dirac operators, it is necessary to explain why this has something to do with the index of an arbitrary elliptic operator on the singular manifold in consideration.
This, of course, also needs to fix an operator algebra context in which all elliptic operators have their right to exist and to fix the meaning of homotopies of symbols which should lead to homotopies of Fredholm operators.

Another aspect of the index theory  (at least in the classical context of the work of Atiyah and Singer \cite{Atiy7}, \cite{Atiy6}) is to study external products $A \boxtimes B$ of elliptic operators $A$ and $B$ over different manifolds $M$ and $N$, respectively.
The product then lives on the Cartesian product $M \times N$, it should be elliptic, and we should have 
$\ind A \boxtimes B = \ind A \ind B$, see also Rodino \cite{Rodi1}.
For singular $M$ and $N$ the Cartesian product $M \times N$ has higher singularities, and a reasonable formulation of the multiplicativity of the index requires the corresponding calculus of operators for the resulting order of singularity.
At least such a question may motivate to seriously promote the calculus of operators on manifolds of arbitrary singularity order.
The problem of multiplicativity itself is sufficiently complex and far from being understood.
From the experience with the classical context it is also clear that we should study elliptic complexes on singular manifolds, Hodge theory, K\"unneth formulas, and other things, known in analogous form from the smooth compact case.
Also this is a wide field, and only partial results are known, see, for instance, \cite{Shaw1}, \cite{Pill1}, \cite{Schu45}, \cite{Mant2}, \cite{Schu6}, \cite{Gries1}.

Let us finally consider the problem of expressing the index in terms of symbols.
An interesting aspect in this connection are so-called analytic index formulas which may consist of expressions that directly compute the index by the symbol.
Here by `symbol' we understand the principal symbol which is, for instance, for conical singularities, the pair of interior principal symbol and the conormal symbol on a given weight line.
Even in that case the problem of deriving analytic index formulas (in analogy of Fedosov's analytic index formulas in the smooth case, see \cite{Fedo6}) is open, with the exception of some particular cases, while analytic index expressions in which lower order terms also participate are apparently easier to organise.

By this remark we stop the index discussion here.
As noted before, geometric or topological aspects of ellipticity on a singular manifold are not the main issue of this exposition; for that we refer   to the new monograph \cite{Naza14}, together with the bibliography there.

\subsection{Concluding remarks}   %Complexes}
\label{s.10.6.5}

%10.6.5
The structures that we discussed here can be motivated by a quite classical question, namely, what has to happen in a (pseudo-differential) scenario on a manifold with a non-complete geometry (for instance, a polyhedron embedded in an Euclidean space) such that parametrices of elliptic operators belong to the calculus.
More precisely, a starting point may be boundary value problems, say, in a cube in $\R^3$, with piecewise smooth Shapiro-Lopatinskij elliptic data on the boundary, for instance, Dirichlet on some of the faces, Neumann on the others.
Examples are also mixed elliptic problems on a $C^\infty$ manifold $X$ with boundary $Y$, with elliptic conditions on different parts $Y_{\pm}$ of the boundary, where $Y_\pm$ subdivide the boundary, i.e.,
$Y=Y_- \cup Y_+$, and $Z := Y_- \cap Y_+$ is of condimension 1.
In the simplest case $Z$ is $C^\infty$, in other cases $Z$ may have singularities, for instance, conical points or edges.

The answer consists of the corner pseudo-differential calculus of boundary value problems with the transmission property at the smooth part of the boundary.
The general technique and many details may be found in the author's joint monographs with  Kapanadze \cite{Kapa10} or Harutjunjan \cite{Haru13} which are based on 
\cite{Schu29}, see also 
\cite{Schu27}, 
\cite{Schu51}, or 
\cite{Calv1},
\cite{Calv2},
\cite{Calv3}, and the papers 
\cite{Dine2}, 
\cite{Haru6},
\cite{Haru11}.

At present there is an increasing stream of investigations in the literature on pseudo-differential theories which claim a relationship with analysis on polyhedral or corner manifolds.
As noted in Section 3.4 a problem is the terminology.
In many cases the investigations are focused on operators on  non-compact manifolds with complete Riemannian metrics and not to configurations in classical boundary value problems, for instance, manifolds with smooth boundary and operators with the transmission property at the boundary (cf. the calculus outlined in Section 2.1).
Parametrix constructions for the above mentioned boundary value problems require careful work with the trace and potential data occurring on the several faces of the configuration, cf. Section 5.2.

One can discover many `unexpected' relations between the analytic machineries on complete or incomplete Riemannian manifolds.
One example is the connection between pseudo-differential operators on the half-axis, with standard symbols, smooth up to $0$, based on the Fourier transform, and Mellin pseudo-differential operators, i.e., operators of Fuchs type, with Mellin symbols that are smooth up to zero, cf. \cite{Eski2},  \cite{Schu31}, \cite{Schu41} and the discussion in Section 2.2 around Mellin quantisation.
Another example is the possible embedding of  elliptic boundary problems with global projection data (`APS' and generalisations on a manifold with smooth boundary) into the pseudo-differential framework, cf. \cite{Schu37}, \cite{Schu50}.
Also the discussion of Section 1.1 on fictitious singularities which makes sense in the pseudo-differential context as well, shows that `usual' pseudo-differential operators which are smooth (across a fictitious conical singularity) may suddenly discover their affection to Fuchs type operators or other societies of corner operators, cf. the general class of Theorem \ref{t.E} below.
One key word is the blow up of singularities which gives rise to degenerate symbols which can be taken as a starting point for Mellin quantisations, cf. Theorem \ref{1022.Op1708.th}. 

In Section 4.1 we saw that there are many kinds of differential operators subsumed under the category `degenerate' with a completely different behaviour.
If they are the result of a blowing up process of singularities, applied to originally given differential operators $D$ on a singular configuration $M$ (minus $M'$, the set of singularities; see the considerations of Section 1.1), then also $M$ itself remains a source of interesting questions.
Also for the above mentioned boundary value problems in polyhedral domains it is helpful to carry out blow ups and to basically deal with the resulting edge- or corner-degenerate operators; which are as in the formulas \eqref{eq.219}, \eqref{eq.220}, \eqref{eq.221}.
Although the calculus mainly refers to such objects, we do not ignore what we want to achieve for the corner singularity itself.
The general scheme of constructing parametrices can be described in terms of a  continuation of the axiomatic approach of Section 5.4, see also \cite{Schu25}.
Many elements on what we understand by ellipticity (here in the sense of the Shapiro-Lopatinskij ellipticity of edge conditions or the ellipticity with respect to the conormal symbols) are described in Section 5.2.
Let us now  consider operators in the upper left corners, i.e., operators on the main stratum.
Those are known in advance, i.e., before we add any extra edge conditions.
Looking at a `higher' stretched corner of the form
$$
K := (\ol{\R}_+)^k \times \Sigma \times \Pi^k_{l=1} \Omega_l
$$
for open sets $\Sigma \subseteq \R^n$, $\Omega_l \subseteq \R^{q_l}$,
$l = 1, \ldots , k$, we first have the space $L^\mu_\cl (\Int K)$ of standard classical pseudo-differential operators on $\Int K$.
As such they have left symbols
$$
a(r,x,y,\rho, \xi, \eta) \in
    S^\mu_\cl (\Int K \times \R^{k+n+q})
$$
for $q := \sum^k_{l=1} q_l$,
$r = (r_1, \ldots , r_k) \in (\R_+)^k$,
$x \in \Sigma$, $y \in \Omega := \Pi^k_{l=1} \Omega_l$,
 with the covariables 
$\rho = ( \rho_1, \ldots , \rho_k) \in \R^k$, $\xi \in \R^n$,
$\eta = (\eta_1, \ldots , \eta_k) \in \R^q$,
$\eta_j \in \R^{q_j}$,
$j = 1, \ldots , k$.
Every $A \in L^\mu_\cl (\Int K)$ has the form 
\begin{equation}
\label{eq.K}
A = \Op (a) \ \ \ \mod \ \ \  L^{-\infty} (\Int K)
\end{equation}
for such a symbol $a$.
Now a first task to treating corner pseudo-differential operators which are related to parametrices of differential operators of the form \eqref{eq.219} with the vector fields \eqref{eq.220}, \eqref{eq.221} is to be aware that $\bigcup_\mu L^\mu_\cl (\Int K)$ contains lots of interesting subalgebras.
In the present case it is adequate take operators with left symbols of the form
\begin{equation}
\label{eq.332}
a (r,x,y,\rho, \xi, \eta) :=
   r^{-\mu} \wt{p} (r,x,y, \wt{\rho}, \xi, \wt{\eta})
\end{equation}
where $r^{-\mu} := r^{-\mu}_1 \cdot \ldots \cdot r^{-\mu}_k$ and
$$
\wt{\rho} :=
(r_1 \rho_1, r_1 r_2 \rho_2 , \ldots , r_1 r_2 \ldots r_k \rho_k),  \ \
\wt{\eta} :=
(r_1 \eta_1, r_1 r_2 \eta_2 , \ldots , r_1 r_2 \ldots r_k \eta_k).
$$
Let $L^\mu_\cl (\Int K)_{\corner}$ denote the subset of all $A \in L^\mu_\cl (\Int K)$ of the form \eqref{eq.K} with symbols \eqref{eq.332} for arbitrary
$$
\wt{p} (r,x,y, \wt{\rho}, \xi, \wt{\eta}) \in
  S^\mu_\cl 
           \big ((\ol{\R}_+)^k \times \Sigma 
                              \times \Omega 
                              \times \R^{k+n+q}_{\wt{\rho}, \xi, \wt{\eta}} \big).
$$    
As we know every $A \in L^\mu_\cl (\Int K)$ can be represented by a properly supported operator $A_0$ modulo an element $C \in L^{-\infty} (\Int K)$.
In particular, this is the case for $A \in L^\mu_\cl (\Int K)_{\corner}$.

An element $A \in L^\mu_\cl (\Int K)_{\corner}$ is called $\sigma_\psi$-elliptic if 
$$
\wt{p}_{(\mu)} (r,x,y, \wt{\rho}, \xi, \wt{\eta}) \not= 0 
\ \ \ \text{for all} \ \ \ 
(r,x,y) \in K, (\wt{\rho}, \xi, \wt{\eta} ) \not= 0,
$$
where $\wt{p}_{(\mu)}$ is the homogeneous principal symbol of $\wt{p}$ in $(\wt{\rho}, \xi, \wt{\eta}) \not= 0$ of order $\mu$.

\begin{Theorem}
\label{t.E}
\begin{itemize}
\item[{\em(i)}]
Let $A \in L^\mu_\cl (\Int K)_{\corner}$, $B \in L^\nu_\cl (\Int K)_{\corner}$, and
$A$ or $B$ properly supported.
Then we have $AB \in L^{\mu+\nu}_\cl (\Int K)_{\corner}$.
\item[{\em(ii)}]
Let $A \in L^\mu_\cl (\Int K)_{\corner}$ be $\sigma_\psi$-elliptic.
Then there is a {\em(}properly supported{\em)} para\-metrix $P_0 \in L^{-\mu}_\cl (\Int K)_{\corner}$
in the sense
$$
I - P_0 A, \ \
I - AP_0 \in L^{-\infty} (\Int K).
$$
\end{itemize}
\end{Theorem}
Moreover, $L^\mu_\cl (\Int K)_{\corner}$ is closed under the operation of formal adjoints.
The proof is elementary and essentially based on the fact that the spaces of involved symbols are closed under asymptotic summation (modulo symbols of order $-\infty$), especially, Leibniz multiplication and Leibniz inversion under the condition of  $\sigma_\psi$-ellipticity.

There are many other variants of Theorem \ref{t.E}, for instance, for operators with symbols with other weight factors instead of $r^{-\mu}$ (e.g., without weight factors).

Another aspect of the parametrix construction is to quantise the obtained Leibniz inverted symbols in such a way that there arise continuous operators in higher weighted corner spaces.
This was outlined in Section 5.4.
The nature of those spaces gives a hint about the adequate smoothing operators in the final corner pseudo-differential algebra.
They can be defined through their mapping properties (and their formal adjoints), namely, to continuously map weighted spaces of any smoothness $s$ to other weighted spaces of  smoothness $s = \infty$.
The latter aspect is a contribution to the discussion in Section 4.2.

Having a parametrix $P$ of $A$ in the corner algebra of the type of an upper left corner (the notation $P$ instead of $P_0$ indicates the chosen quantisation in order to reach an operator in the corner calculus), we can try to add extra elliptic conditions, according to Section 5.2, and to obtain a block matrix operator $\s{P}$ with $P$ in the upper left corner.
Then, if $A$ is the elliptic operator in the given boundary value problem 
$\s{A} = {}^\t (A \ \  T)$ (say, the Laplacian in a cube $M$ with the Dirichlet/Neumann conditions on the faces of the boundary $M'$, indicated by $T$) then the operator $\s{P}$ with can be employed to reduce $\s{A}$ to the boundary $M'$.
The result is an elliptic operator on $M'$ which can be treated on the level of operators on the corner manifold $M'$ without boundary.
The resulting operator $\s{R}$ on the boundary, in general being again an elliptic block matrix operator with an upper left corner $R$, can be interpreted as a transmission problem for the elliptic pseudo-differential operator $R$ on $M'$ with a jumping behaviour across the interfaces $Z = M''$  of $M'$ (in the case of a cube $M$ the interfaces $M''$ consist of the system of one-dimensional edges plus the corner points).
To treating $R$ is now a beautiful task in the framework of boundary value problems for pseudo-differential operators without the transmission property at the smooth part of the boundary, where the boundary itself may have corner points $M'''$.
Although the method to carry out all this is clear in principle, many details, refinements and more explicit information should be worked out in future.
By that we mean, in particular, computing the admissible weights in the weighted corner spaces, the number of extra interface conditions (of trace and potential type) depending on the weights, and the explicit (corner-) asymptotics of solutions.
In this context there are lots of other things worth to be developed, for instance, the calculus of operators on manifolds with conical exits to infinity, modelled on a cylinder with cross section that has itself singularities.
Other useful details to be completed and deepened are Green formulas of several kind, the kernel cut-off and corner quantisation, or potentials of densities on a manifolds with corners, embedded in an ambient smooth manifold, with respect to a fundamental solution of an elliptic operator.

\addcontentsline{toc}{section}{References} 
    
\bibliographystyle{plain}
\bibliography{master}

\end{document}